\documentclass[leqno,draft]{amsart}
\usepackage{amssymb,enumerate, latexsym, mathrsfs, mdwlist}

\swapnumbers
\theoremstyle{definition}
\newtheorem{nul}{}[section]
\newtheorem{dfn}[nul]{Definition}
\newtheorem{axm}[nul]{Axiom}
\newtheorem{cnstr}[nul]{Construction}
\newtheorem{ntn}[nul]{Notation}
\newtheorem{exm}[nul]{Example}

\newtheorem{rec}[nul]{Recollection}

\newtheorem*{dfn*}{Definition}
\newtheorem*{axm*}{Axiom}
\newtheorem*{ntn*}{Notation}
\newtheorem*{exm*}{Example}
\newtheorem*{exr*}{Exercise}
\newtheorem*{int*}{Intuition}
\newtheorem*{qst*}{Question}

\theoremstyle{plain}

\newtheorem{thm}[nul]{Theorem}
\newtheorem{prp}[nul]{Proposition}
\newtheorem{lem}[nul]{Lemma}

\newtheorem{cor}{Corollary}[nul]
\newtheorem*{thm*}{Theorem}
\newtheorem*{mainthm*}{Universal Additivity Theorem}
\newtheorem*{structurethm*}{Structure Theorem}
\newtheorem*{prp*}{Proposition}
\newtheorem*{cor*}{Corollary}
\newtheorem*{lem*}{Lemma}
\newtheorem*{cnj*}{Conjecture}

\numberwithin{equation}{nul}

\DeclareMathOperator{\cof}{cof}

\DeclareMathOperator{\colim}{colim}
\DeclareMathOperator{\Colim}{Colim}
\DeclareMathOperator{\creff}{cr}
\DeclareMathOperator{\End}{End}
\DeclareMathOperator{\ev}{ev}

\DeclareMathOperator{\Ext}{Ext}
\DeclareMathOperator{\Exc}{Exc}

\DeclareMathOperator{\Fun}{Fun}
\DeclareMathOperator{\Ho}{Ho}

\DeclareMathOperator{\id}{id}
\DeclareMathOperator{\Ind}{Ind}

\DeclareMathOperator{\Map}{Map}
\DeclareMathOperator{\Mor}{Mor}
\DeclareMathOperator{\Obj}{Obj}
\DeclareMathOperator{\pr}{pr}
\DeclareMathOperator{\Spec}{Spec}
\DeclareMathOperator{\Thy}{Thy}

\newcommand{\BB}{\mathbf{B}}
\newcommand{\CC}{\mathbf{C}}
\newcommand{\DD}{\mathbf{D}}

\newcommand{\FF}{\mathbf{F}}

\newcommand{\KK}{\mathbf{K}}

\newcommand{\SSS}{\mathbf{S}}

\newcommand{\UU}{\mathbf{U}}
\newcommand{\VV}{\mathbf{V}}

\newcommand{\XX}{\mathbf{X}}
\newcommand{\YY}{\mathbf{Y}}
\newcommand{\ZZ}{\mathbf{Z}}

\newcommand{\Add}{\mathrm{Add}}
\newcommand{\Alg}{\mathbf{Alg}}
\newcommand{\Cat}{\mathbf{Cat}}
\newcommand{\coWald}{\mathbf{coWald}}

\newcommand{\Grp}{\mathbf{Grp}}

\newcommand{\Kan}{\mathbf{Kan}}
\newcommand{\lWald}{\ell\mathbf{Wald}_{\infty}}
\newcommand{\Mod}{\mathbf{Mod}}
\newcommand{\Pair}{\mathbf{Pair}}
\newcommand{\Perf}{\mathbf{Perf}}

\renewcommand{\Pr}{\mathbf{Pr}}

\newcommand{\QCoh}{\mathbf{QCoh}}
\newcommand{\Rex}{\mathbf{Rex}}
\newcommand{\Sch}{\mathbf{Sch}}

\newcommand{\Set}{\mathbf{Set}}
\newcommand{\Sp}{\mathbf{Sp}}
\newcommand{\Stk}{\mathbf{Stk}}

\newcommand{\VaddWald}{\mathrm{D}_{\mathrm{fiss}}(\mathbf{Wald}_{\infty})}
\newcommand{\VWald}{\mathrm{D}(\mathbf{Wald}_{\infty})}
\newcommand{\Wald}{\mathbf{Wald}_{\infty}}

\newcommand{\add}{\mathrm{add}}
\newcommand{\cart}{\mathrm{cart}}
\newcommand{\cocart}{\mathrm{cocart}}

\newcommand{\fiss}{\mathrm{fiss}}

\newcommand{\op}{\mathrm{op}}

\newcommand{\coloneq}{\mathrel{\mathop:}=}

\makeatletter 
\def\revddots{\mathinner{\mkern1mu\raise\p@ 
\vbox{\kern7\p@\hbox{.}}\mkern2mu 
\raise4\p@\hbox{.}\mkern2mu\raise7\p@\hbox{.}\mkern1mu}} 
\makeatother

\usepackage{tikz}
\usetikzlibrary{matrix,arrows,decorations}

\pgfdeclaredecoration{single line}{initial}{\state{initial}[width=\pgfdecoratedpathlength-1sp]{\pgfpathmoveto{\pgfpointorigin}}\state{final}{\pgfpathlineto{\pgfpointorigin}}} 

\newcommand{\fromto}[2]{{#1}\ \tikz[baseline]\draw[>=stealth,->](0,0.5ex)--(0.5,0.5ex);\ {#2}}
\newcommand{\into}[2]{{#1}\ \tikz[baseline]\draw[>=stealth,right hook->](0,0.5ex)--(0.5,0.5ex);\ {#2}}

\newcommand{\cofto}[2]{{#1}\ \tikz[baseline]\draw[>=stealth,>->](0,0.5ex)--(0.5,0.5ex);\ {#2}}
\newcommand{\fibto}[2]{{#1}\ \tikz[baseline]\draw[>=stealth,->>](0,0.5ex)--(0.5,0.5ex);\ {#2}}
\newcommand{\equivto}[2]{{#1}\ \tikz[baseline]\draw[>=stealth,->,font=\scriptsize,inner sep=0.5pt](0,0.5ex)--node[above]{$\sim$}(0.5,0.5ex);\ {#2}}
\newcommand{\goesto}[2]{{#1}\ \tikz[baseline]\draw[>=stealth,|->](0,0.5ex)--(0.5,0.5ex);\ {#2}}
\newcommand{\adjunct}[4]{{#1}\colon{#2}\ \begin{tikzpicture}[baseline] \draw[>=stealth,->] (0,1ex) -- (0.75,1ex); \draw[>=stealth,->] (0.75,0.25ex) -- (0,0.25ex); \end{tikzpicture}\ {#3}\colon{#4}}

\let\oldtocsection=\tocsection

\let\oldtocsubsection=\tocsubsection

\renewcommand{\tocsection}[2]{\hspace{0em}\oldtocsection{#1}{#2}}
\renewcommand{\tocsubsection}[2]{\hspace{1.5em}\oldtocsubsection{#1}{#2}}

\title{On the algebraic $K$-theory of higher categories}

\author{Clark Barwick}
\address{Massachusetts Institute of Technology, Department of Mathematics, Building 2, 77 Massachusetts Avenue, Cambridge, MA 02139-4307, USA}
\email{clarkbar@gmail.com}

\dedicatory{In memoriam Daniel Quillen, 1940--2011, with profound admiration.}

\begin{document}

\begin{abstract} We prove that Waldhausen $K$-theory, when extended to a very general class of quasicategories, can be described as a Goodwillie differential. In particular, $K$-theory spaces admit canonical (connective) deloopings, and the $K$-theory functor enjoys a simple universal property. Using this, we give new, higher categorical proofs of the Approximation, Additivity, and Fibration Theorems of Waldhausen in this context. As applications of this technology, we study the algebraic $K$-theory of associative rings in a wide range of homotopical contexts and of spectral Deligne--Mumford stacks. 
\end{abstract}

\maketitle

\setcounter{tocdepth}{2}
\tableofcontents


\setcounter{section}{-1}
\section{Introduction} We characterize algebraic $K$-theory as a \emph{universal homology theory}, which takes suitable higher categories as input and produces either spaces or spectra as output. What makes $K$-theory a \emph{homology theory} is that it satisfies an excision axiom. This excision axiom is tantamount to what Waldhausen called \emph{additivity}, so that an excisive functor is precisely one that splits short exact sequences. What makes this homology theory \emph{universal} is this: if we write $\iota$ for the functor that carries a higher category to its moduli space of objects, then algebraic $K$-theory is initial among homology theories $F$ that receive a natural transformation $\fromto{\iota}{F}$. In the lingo of Tom Goodwillie's calculus of functors \cite{MR1076523,MR1162445,MR2026544}, $K$ is the \emph{linearization} of $\iota$. Algebraic $K$-theory is thus the analog of stable homotopy theory in this new class of categorified homology theories. From this we obtain an explicit universal property that completely characterizes algebraic $K$-theory and permits us to give new, conceptual proofs of the fundamental theorems of Waldhausen $K$-theory.

To get a feeling for this universal property, let's first contemplate $K_0$. For any ordinary category $C$ with a zero object and a reasonable notion of ``short exact sequence'' (e.g., an exact category in the sense of Quillen, or a category with cofibrations in the sense of Waldhausen, or a triangulated category in the sense of Verdier), the abelian group $K_0(C)$ is the universal target for Euler characteristics. That is, for any abelian group $A$, the set $\mathrm{Hom}(K_0(C),A)$ is in natural bijection with the set of maps $\phi\colon\fromto{\Obj C}{A}$ such that $\phi(X)=\phi(X')+\phi(X'')$ whenever there is a short exact sequence
\begin{equation*}
X'\ \tikz[baseline]\draw[>=stealth,>->](0,0.5ex)--(0.5,0.5ex);\ X\ \tikz[baseline]\draw[>=stealth,->>](0,0.5ex)--(0.5,0.5ex);\ X''.
\end{equation*}

We can reinterpret this as a universal property on the entire functor $K_0$, which we'll regard as valued in the category of \emph{sets}. Just to fix ideas, let's assume that we are working with the algebraic $K$-theory of categories with cofibrations in the sense of Waldhausen. If $E(C)$ is the category of short exact sequences in a category with cofibrations $C$, then $E(C)$ is also a category with cofibrations. Moreover, for any $C$,
\begin{enumerate}[(1)]
\item the functors
\begin{equation*}
\goesto{[X'\ \tikz[baseline]\draw[>=stealth,>->](0,0.5ex)--(0.5,0.5ex);\ X\ \tikz[baseline]\draw[>=stealth,->>](0,0.5ex)--(0.5,0.5ex);\ X'']}{X'}\textrm{\quad and\quad}\goesto{[X'\ \tikz[baseline]\draw[>=stealth,>->](0,0.5ex)--(0.5,0.5ex);\ X\ \tikz[baseline]\draw[>=stealth,->>](0,0.5ex)--(0.5,0.5ex);\ X'']}{X''}
\end{equation*}
together induce a bijection $\equivto{K_0(E(C))}{K_0(C)\times K_0(C)}$.
\suspend{enumerate}
The functor $\goesto{[X'\ \tikz[baseline]\draw[>=stealth,>->](0,0.5ex)--(0.5,0.5ex);\ X\ \tikz[baseline]\draw[>=stealth,->>](0,0.5ex)--(0.5,0.5ex);\ X'']}{X}$ now gives a commutative monoid structure $K_0(C)\times K_0(C)\cong K_0(E(C))\ \tikz[baseline]\draw[>=stealth,->](0,0.5ex)--(0.5,0.5ex);\ K_0(C)$. With this structure, $K_0$ is an abelian group. We can express this sentiment diagrammatically by saying that
\resume{enumerate}[{[(1)]}]
\item the functors
\begin{equation*}
\goesto{[X'\ \tikz[baseline]\draw[>=stealth,>->](0,0.5ex)--(0.5,0.5ex);\ X\ \tikz[baseline]\draw[>=stealth,->>](0,0.5ex)--(0.5,0.5ex);\ X'']}{X'}\textrm{\quad and\quad}\goesto{[X'\ \tikz[baseline]\draw[>=stealth,>->](0,0.5ex)--(0.5,0.5ex);\ X\ \tikz[baseline]\draw[>=stealth,->>](0,0.5ex)--(0.5,0.5ex);\ X'']}{X}
\end{equation*}
also induce a bijection $\equivto{K_0(E(C))}{K_0(C)\times K_0(C)}$.
\end{enumerate}
Now our universal characterization of $K_0$ simply says that we have a natural transformation $\fromto{\Obj}{K_0}$ that is initial with properties (1) and (2).

For the $K$-theory \emph{spaces} (whose homotopy groups will be the higher $K$-theory groups), we can aim for a homotopical variant of this universal property. We replace the word ``bijection'' in (1) and (2) with the words ``weak equivalence;'' a functor satisfying these properties is called an \emph{additive functor}. Instead of a map from the set of objects of the category with cofibrations $C$, we have a map from the \emph{moduli space} of objects --- this is the classifying space $N\iota C$ of the groupoid $\iota C\subset C$ consisting of all isomorphisms in $C$. An easy case of our main theorem states that algebraic $K$-theory is initial in the homotopy category of (suitably finitary) additive functors $F$ equipped with a natural transformation $\fromto{N\iota}{F}$.

Now let's enlarge the scope of this story enough to bring in examples such as Waldhausen's algebraic $K$-theory of spaces by introducing homotopy theory in the source of our $K$-theory functor. We use $\infty$-categories that contain a zero object and suitable cofiber sequences, and we call these \emph{Waldhausen $\infty$-categories}. Our homotopical variants of (1) and (2) still make sense, so we still may speak of \emph{additive functors} from Waldhausen $\infty$-categories to spaces. Moreover, any $\infty$-category has a moduli space of objects, which is given by the maximal $\infty$-groupoid contained therein; this defines a functor $\iota$ from Waldhausen $\infty$-categories to spaces. Our main theorem (\S \ref{sect:univpropKthy}) is thus the natural extension of the characterization of $K_0$ as the universal target for Euler characteristics:

\begin{mainthm*}[\S \ref{sect:univpropKthy}] Algebraic $K$-theory is homotopy-initial among \textup{(}suitably finitary\textup{)} additive functors $F$ equipped with a natural transformation $\fromto{\iota}{F}$.
\end{mainthm*}

It is well-known that algebraic $K$-theory is hair-raisingly difficult to compute, and that various theories that are easier to compute, such as forms of $\mathrm{THH}$ and $\mathrm{TC}$, are prime targets for ``trace maps'' \cite{MR1474979}. The Universal Additivity Theorem actually classifies all such trace maps: for any additive functor $H$, the space of natural transformations $\fromto{K}{H}$ is equivalent to the space of natural transformations $\fromto{\iota}{H}$. But since $\iota$ is actually represented by the ordinary category $\Gamma^{\op}$ of pointed finite sets, it follows from the Yoneda lemma that the space of natural transformations $\fromto{K}{H}$ is equivalent to the space $H(\Gamma^{\op})$. In particular, by Barratt--Priddy--Quillen, we compute the space of ``global operations'' on algebraic $K$-theory:
\begin{equation*}
\End(K)\simeq QS^0.
\end{equation*}

\bigskip

The proof of the Universal Additivity Theorem uses a new way of conceptualizing functors such as algebraic $K$-theory. Namely, we regard algebraic $K$-theory as a homology theory on Waldhausen $\infty$-categories, and we regard additivity as an excision axiom. But this isn't just some slack-jawed analogy: we'll actually pass to a homotopy theory on which functors that are $1$-excisive in the sense of Goodwillie (i.e., functors that carry homotopy pushouts to homotopy pullbacks) correspond to additive functors as described above. (And making sense of this homotopy theory forces us to pass to the $\infty$-categorical context.)

The idea here is to regard the homotopy theory $\Wald$ of Waldhausen $\infty$-categories as formally analogous to the ordinary category $V(k)$ of vector spaces over a field $k$. The left derived functor of a right exact functor out of $V(k)$ is defined on the \emph{derived category} $\mathrm{D}_{\geq 0}(k)$ of chain complexes whose homology vanishes in negative degrees. Objects of $\mathrm{D}_{\geq 0}(k)$ can be regarded as formal geometric realizations of simplicial vector spaces. Correspondingly, we define a \emph{derived $\infty$-category} $\VWald$ of $\Wald$, whose objects can be regarded as formal geometric realizations of simplicial Waldhausen $\infty$-categories. This entitles us to speak of the \emph{left derived functor} of a functor defined on $\Wald$. Then we suitably localize $\VWald$ in order to form a universal homotopy theory $\VaddWald$ in which exact sequences split; we call this the \emph{fissile derived $\infty$-category}. Our Structure Theorem (Th. \ref{thm:additiveequiv}) uncovers the following relationshp between excision on $\VaddWald$ and additivity:

\begin{structurethm*}[Th. \ref{thm:additiveequiv}] A \textup{(}suitably finitary\textup{)} functor from Waldhausen $\infty$-categories to spaces is additive in the sense above if and only if its left derived functor factors through an excisive functor on the fissile derived $\infty$-category $\VaddWald$.
\end{structurethm*}

This Structure Theorem is not some dreary abstract formalism: the technology of Goodwillie's calculus of functors tells us that the way to compute the universal excisive approximation to a functor $F$ is to form the colimit of $\Omega^nF\Sigma^n$ as $\fromto{n}{\infty}$. This means that as soon as we've worked out how to compute the suspension $\Sigma$ in $\VaddWald$, we'll have an explicit description of the \emph{additivization} of any functor $\phi$ from $\Wald$ to spaces, which is the universal approximation to $\phi$ with an additive functor. And when we apply this additivization to the functor $\iota$, we'll obtain a formula for the very thing we're claiming is algebraic $K$-theory: the initial object in the homotopy category of additive functors $F$ equipped with a natural transformation $\fromto{\iota}{F}$.

So what \emph{is} $\Sigma$? Here's the answer: it's given by the formal geometric realization of Waldhausen's $S_{\bullet}$ construction (suitably adapted for $\infty$-categories). So the universal homology theory with a map from $\iota$ is given by the formula
\begin{equation*}
\goesto{\mathscr{C}}{\colim_{n}\Omega^n|\iota S_{\ast}^n(\mathscr{C})|}.
\end{equation*}
This is exactly Waldhausen's formula for algebraic $K$-theory, so our Main Theorem is an easy consequence of our Structure Theorem and our computation of $\Sigma$.

Bringing algebraic $K$-theory under the umbrella of Goodwillie's calculus of functors has a range of exciting consequences, which we are only able to touch upon in this first paper. In particular, three key foundational results of Waldhausen's algebraic $K$-theory --- the Additivity Theorem \cite[Th. 1.4.2]{MR86m:18011} (our version: Cor. \ref{for:preaddsareeasy}), the Approximation Theorem \cite[Th. 1.6.7]{MR86m:18011} (our version: Pr. \ref{prp:approx}), the Fibration Theorem \cite[Th. 1.6.4]{MR86m:18011} (our version: Pr. \ref{thm:fibration}), and the Cofinality Theorem \cite[Th. 2.1]{MR990574} (our version: Th. \ref{thm:cofinality}) --- are straightforward consequences of general facts about the calculus of functors combined with some observations about the homotopy theory of $\Wald$.

To get a glimpse of various bits of our framework at work, we offer two examples that exploit certain features of the algebraic $K$-theory functor of which we are fond. First (\S \ref{sect:example1}), we apply our foundational work to the study of the connective $K$-theory of $E_1$-algebras in suitable ground $\infty$-categories. We define a notion of a \emph{perfect} left module over an $E_1$-algebra (Df. \ref{ntn:perfectmodule}). In the special case of an $E_1$ ring spectrum $\Lambda$, for any set $S$ of homogenous elements of $\pi_{\ast}\Lambda$ that satisfies a left Ore condition, we obtain a fiber sequence of connective spectra
\begin{equation*}
\KK(\mathbf{Nil}_{(\Lambda,S)}^{\omega})\ \tikz[baseline]\draw[>=stealth,->](0,0.5ex)--(0.5,0.5ex);\ \KK(\Lambda)\ \tikz[baseline]\draw[>=stealth,->](0,0.5ex)--(0.5,0.5ex);\ \KK(\Lambda[S^{-1}]),
\end{equation*}
in which the first term is the $K$-theory of the $\infty$-category of $S$-nilpotent perfect $\Lambda$-modules (Pr. \ref{prp:protolocforE1rings}). (Note that we only work with connective $K$-theory, so this is only a fiber sequence in the homotopy theory of connective spectra; in particular, the last map need not be surjective on $\pi_0$.) Such a result --- at least in special cases --- is surely well-known among experts; see for example \cite[Pr. 1.4 and Pr. 1.5]{BM}.

Finally (\S \ref{sect:example2}), we introduce $K$-theory in derived algebraic geometry. In particular, we define the $K$-theory of quasicompact nonconnective spectral Deligne--Mumford stacks (Df. \ref{dfn:algKthyDMstacks}). We prove a result analogous to what Thomason called the ``proto-localization'' theorem \cite[Th. 5.1]{MR92f:19001}; this is a fiber sequence of connective spectra
\begin{equation*}
\KK(\mathscr{X}\setminus\mathscr{U})\ \tikz[baseline]\draw[>=stealth,->,font=\scriptsize](0,0.5ex)--(0.75,0.5ex);\ \KK(\mathscr{X})\ \tikz[baseline]\draw[>=stealth,->,font=\scriptsize](0,0.5ex)--(0.75,0.5ex);\ \KK(\mathscr{U})
\end{equation*}
corresponding to a quasicompact open immersion $j\colon\fromto{\mathscr{U}}{\mathscr{X}}$ of quasicompact, quasiseparated spectral algebraic spaces. Here $\KK(\mathscr{X}\setminus\mathscr{U})$ is the $K$-theory of the $\infty$-category perfect modules $\mathscr{M}$ on $\mathscr{X}$ such that $j^{\star}\mathscr{M}\simeq0$ (Pr. \ref{prp:protolocforDMstacks}). Our proof is new in the details even in the setting originally contemplated by Thomason (though of course the general thrust is the same).

\subsection*{Relation to other work} Our universal characterization of algebraic $K$-theory has probably been known --- perhaps in a more restrictive setting and certainly in a different language --- to a variety of experts for many years. In fact, the universal property stated here has endured a lengthy gestation: the first version of this characterization emerged during a question-and-answer session between the author and John Rognes after a talk given by the author at the University of Oslo in 2006.

The idea that algebraic $K$-theory could be characterized via a universal property goes all the way back to the beginnings of the subject, when Grothendieck defined what we today call $K_0$ of an abelian or triangulated category just as we describe above \cite{MR0116023,MR50:7133}. The idea that algebraic $K$-theory might be expressible as a linearization was directly inspired by the ICM talk of Tom Goodwillie \cite{MR1159249} and the remarkable flurry of research into the relationship between algebraic $K$-theory and the calculus of functors --- though of course the setting for our Goodwillie derivative is more primitive than the one studied by Goodwillie et al.

But long before that, of course, came the foundational work of Waldhausen \cite{MR86m:18011}. Since it is known today that relative categories comprise a model for the homotopy theory of $\infty$-categories \cite{BarKan1011}, the work of Waldhausen can be said to represent the first study of the algebraic $K$-theory of higher categories. Furthermore, the idea that the defining property of this algebraic $K$-theory is additivity is strongly suggested by Waldhausen, and this point is driven home in the work of Randy McCarthy \cite{MR1217072} and Ross Staffeldt \cite{MR990574}, both of whom recognized long ago that the Additivity Theorem is the ur-theorem of algebraic $K$-theory.

In a parallel development, Amnon Neeman has advanced the \emph{algebraic $K$-theory of triangulated categories}  \cite{MR1491990,MR1604910,MR1656552,MR1724625,MR1793672,MR1798824,MR1798828,MR1828612,MR2181838} as a way of extracting $K$-theoretic data directly from the triangulated homotopy category of a stable homotopy theory. The idea is that the algebraic $K$-theory of a ring or scheme should by approximation depend (in some sense) only on a derived category of perfect modules; however, this form of $K$-theory has known limitations: an example of Marco Schlichting \cite{MR1930883} shows that Waldhausen $K$-theory can distinguish stable $\infty$-categories with equivalent triangulated homotopy categories. These limitations are overcome by passing to the \emph{derived $\infty$-category}.

More recently, Bertrand To\"en and Gabriele Vezzosi showed \cite{MR2061207} that the Waldhausen $K$-theory of many of the best-known examples of Waldhausen categories is in fact an invariant of the simplicial localization; thus To\"en and Vezzosi are more explicit in identifying higher categories as a natural domain for $K$-theory. In fact, in the final section of \cite{MR2061207}, the authors suggest a strategy for constructing the $K$-theory of a Segal category by means of an ``$S_{\scriptstyle{\bullet}}$ construction up to coherent homotopy.'' The desired properties of their construction are reflected precisely in our construction $\mathscr{S}$. These insights were explored more deeply in the work of Blumberg and Mandell \cite{MR2764905}; they give an explicit description of Waldhausen's $S_{\scriptstyle{\bullet}}$ construction in terms of the mapping spaces of the simplicial localization, and they extend Waldhausen's approximation theorem to show that in many cases, equivalences of homotopy categories alone are enough to deduce equivalences of $K$-theory spectra.

Even more recent work of Andrew Blumberg, David Gepner, and Gon{\c c}alo Tabuada \cite{BGT} has built upon brilliant work of the last of these authors in the context of DG categories \cite{MR2451292} to produce another universal characterization of the algebraic $K$-theory of stable $\infty$-categories. One of their main results may be summarized by saying that the algebraic $K$-theory of stable $\infty$-categories is a universal additive invariant. They do not deal with general Waldhausen $\infty$-categories, but they also study nonconnective deloopings of $K$-theory, with which we do not contend here.

Finally, we recall that Waldhausen's formalism for algebraic $K$-theory has of course been applied in the context of associative $S$-algebras by Elmendorf, Kriz, Mandell, and May \cite{MR97h:55006}, and in the context of schemes and algebraic stacks by Thomason and Trobaugh \cite{MR92f:19001}, To\"en \cite{MR2000h:14010}, Joshua \cite{MR1942183,MR1949305,MR1967388}, and others. The applications of the last two sections of this paper are extensions of their work.

\subsection*{A word on higher categories} When we speak of $\infty$-categories in this paper, we mean $\infty$-categories whose $k$-morphisms for $k\geq 2$ are weakly invertible. We will use the \emph{quasicategory} model of this sort of $\infty$-categories. Quasicategories were invented in the 1970s by Boardman and Vogt \cite{MR0420609}, who called them \emph{weak Kan complexes}, and they were studied extensively by Joyal \cite{Joyal, Joyal08} and Lurie \cite{HTT}. We emphasize that quasicategories are but one of an array of equivalent models of $\infty$-categories (including simplicial categories \cite{MR81h:55018,MR81h:55019,MR81m:55018,MR2276611}, Segal categories \cite{math.AG/9807049,MR2341955,simpson:book}, and complete Segal spaces \cite{MR1804411,MR2439415}), and there is no doubt that the results here could be satisfactorily proved in any one of these models. Indeed, there is a canonical equivalence between any two of these homotopy theories \cite{JT,MR2321038,Bergner} (or any other homotopy theory that satisfies the axioms of \cite{Toen} or of \cite{BSP}), through which one can surely translate the main theorems here into theorems in the language of any other model. To underscore this fact, we will frequently use the generic term \emph{$\infty$-category} in lieu of the more specialized term \emph{quasicategory}.

That said, we wish to emphasize that we employ many of the technical details of the particular theory of quasicategories as presented in \cite{HTT} in a critical way in this paper. In particular, beginning in \S \ref{sect:Waldfib}, the theory of fibrations, developed by Joyal and presented in Chapter 2 of \cite{HTT}, is instrumental to our work here, as it provides a convenient way to finesse the homotopy-coherence issues that would otherwise plague this paper and its author. Indeed, it is the convenience and relative simplicity of this theory that compelled us to work with this model.


\subsection*{Acknowledgments} There are a lot of people to thank. Without the foundational work of Andr\'e Joyal and Jacob Lurie on quasicategories, the results here would not admit such simple statements or such straightforward proofs. I thank Jacob also for generously answering a number of questions during the course of the work represented here.

My conversations with Andrew Blumberg over the past few years have been consistently enlightening, and I suspect that a number of the results here amount to elaborations of insights he had long ago. I have also benefitted from conversations with Dan Kan.

In the spring of 2012, I gave a course at MIT on the subject of this paper. During that time, several sharp-eyed students spotted errors, including especially Rune Haugseng, Luis Alexandre Pereira, and Guozhen Wang. I owe them my thanks for their scrupulousness.

John Rognes has declined to be listed as a coauthor, but his influence on this work has been tremendous. He was present at the conception of the main result, and this paper is teeming with insights I inherited from him.

I thank Peter Scholze for noticing an error that led to the inclusion of the Cofinality Theorem \ref{thm:cofinality}.

Advice from the first referee and from Mike Hopkins has led to great improvements in the exposition of this paper. The second referee expended huge effort to provide me with a huge, detailed list of little errors and omissions, and I thank him or her most heartily.

On a more personal note, I thank Alexandra Sear for her unfailing patience and support during this paper's ridiculously protracted writing process.


\part{Pairs and Waldhausen $\infty$-categories} In this part, we introduce the basic input for additive functors, including the form of $K$-theory we study. We begin with the notion of a \emph{pair} of $\infty$-categories, which is nothing more than an $\infty$-category with a subcategory of \emph{ingressive morphisms} that contains the equivalences. Among the pairs of $\infty$-categories, we will then isolate the \emph{Waldhausen $\infty$-categories} as the input for algebraic $K$-theory; these are pairs that contain a zero object in which the ingressive morphisms are stable under pushout. This is the $\infty$-categorical analogue of Waldhausen's notion of categories with cofibrations.

We will also need to speak of families of Waldhausen $\infty$-categories, which are called \emph{Waldhausen (co)cartesian fibrations}, and which classify functors valued in the $\infty$-category $\Wald$ of Waldhausen $\infty$-categories. We study limits and colimits in $\Wald$, and we construct the $\infty$-category of \emph{virtual Waldhausen $\infty$-categories}, whose homotopy theory serves as the basis for all the work we do in this paper.


\section{Pairs of $\infty$-categories} The basic input for Waldhausen's algebraic $K$-theory \cite{MR86m:18011} is a category equipped with a subcategory of weak equivalences and a subcategory of cofibrations. These data are then required to satisfy sundry axioms, which give what today is often called a \emph{Waldhausen category}.

A category with a subcategory of weak equivalences (or, in the parlance of \cite{BarKan1011}, a \emph{relative category}) is one way of exhibiting a homotopy theory. A \emph{quasicategory} is another. It is known \cite[Cor. 6.11]{BarKan1011} that these two models of a homotopy theory contain essentially the same information. Consequently, if one wishes to employ the flexible techniques of  quasicategory theory, one may attempt to replace the category with weak equivalences in Waldhausen's definition with a single quasicategory.

But what then is to be done with the cofibrations? In Waldhausen's framework, the specification of a subcategory of cofibrations actually serves two distinct functions.
\begin{enumerate}[(1)]
\item First, Waldhausen's Gluing Axiom \cite[\S 1.2, Weq. 2]{MR86m:18011} ensures that pushouts along these cofibrations are compatible with weak equivalences. For example, pushouts in the category of simplicial sets along inclusions are compatible with weak equivalences in this sense; consequently, in the Waldhausen category of finite spaces, the cofibrations are monomorphisms.
\item Second, the cofibrations permit one to restrict attention to the particular class of cofiber sequences one wishes to split in $K$-theory. For example, an exact category is regarded as a Waldhausen category by declaring the cofibrations to be the admissible monomorphisms; consequently, the admissible exact sequences are the only exact sequences that algebraic $K$-theory splits.
\end{enumerate}

In a quasicategory, the first function becomes vacuous, as the only sensible notion of \emph{pushout} in a quasicategory must preserve equivalences. Thus only the second function for a class of cofibrations in a quasicategory will be relevant. This means, in particular, that we needn't make any distinction between a cofibration in a quasicategory and a morphism that is \emph{equivalent} to a cofibration. In other words, a suitable class of cofibrations in a quasicategory $C$ will be uniquely specified by a subcategory of the homotopy category $hC$. We will thus define a \emph{pair} of $\infty$-categories as an $\infty$-category along with a subcategory of the homotopy category. (We call these \emph{ingressive morphisms}, in order to distinguish it from the more rigid notion of cofibration.) Among these pairs, we will isolate the Waldhausen $\infty$-categories in the next section.

In this section, we introduce the homotopy theory of pairs as a stepping stone on the way to defining the critically important homotopy theory of Waldhausen $\infty$-categories. As many constructions in the theory of Waldhausen $\infty$-categories begin with a construction at the level of pairs of $\infty$-categories, it is convenient to establish robust language and notation for these objects. To this end, we begin with a brief discussion of some set-theoretic considerations and a reminder on constructions of $\infty$-categories from simplicial categories and relative categories. We apply these to the construction of an $\infty$-category of $\infty$-categories and --- following a short reminder on the notion of a \emph{subcategory} of an $\infty$-category --- an $\infty$-category $\Pair_{\infty}$ of pairs of $\infty$-categories. Finally, we relate this $\infty$-category of pairs to an $\infty$-category of functors between $\infty$-categories; this permits us to exhibit $\Pair_{\infty}$ as a \emph{relative nerve}.


\subsection*{Set theoretic considerations} In order to circumvent the set-theoretic difficulties arising from the consideration of these $\infty$-categories of $\infty$-categories and the like, we must employ some artifice. Hence to the usual Zermelo--Frankel axioms \textsc{zfc} of set theory (including the Axiom of Choice) we add the following \emph{Universe Axiom} of Grothendieck and Verdier \cite[Exp I, \S 0]{MR50:7130}. The resulting set theory, called \textsc{zfcu}, will be employed in this paper.

\begin{axm}[Universe] Any set is an element of a universe.
\end{axm}

\begin{nul} This axiom is independent of the others of \textsc{zfc}, since any universe $\UU$ is itself a model of Zermelo--Frankel set theory. Equivalently, we assume that for any cardinal $\tau$, there exists a strongly inaccessible cardinal $\kappa$ with $\tau<\kappa$; for any strongly inaccessible cardinal $\kappa$, the set $\VV_{\kappa}$ of sets whose rank is strictly smaller than $\kappa$ is a universe \cite{MR0244035}.
\end{nul}

\begin{ntn} In addition, we fix, once and for all, three uncountable, strongly inaccessible cardinals $\kappa_0<\kappa_1<\kappa_2$ and the corresponding universes $\VV_{\kappa_0}\in\VV_{\kappa_1}\in\VV_{\kappa_2}$. Now a set, simplicial set, category, etc., will be said to be \textbf{\emph{small}} if it is contained in the universe $\VV_{\kappa_0}$; it will be said to be \textbf{\emph{large}} if it is contained in the universe $\VV_{\kappa_1}$; and it will be said to be \textbf{\emph{huge}} if it is contained in the universe $\VV_{\kappa_2}$. We will say that a set, simplicial set, category, etc., is \textbf{\emph{essentially small}} if it is equivalent (in the appropriate sense) to a small one.
\end{ntn}


\subsection*{Simplicial nerves and relative nerves} There are essentially two ways in which $\infty$-categories will arise in the sequel. The first of these is as \emph{simplicial categories}. We follow the model of \cite[Df. 3.0.0.1]{HTT} for the notation of simplicial nerves.
\begin{ntn}\label{ntn:superscriptsscats} A \textbf{\emph{simplicial category}} --- that is, a category enriched in the category of simplicial sets --- will frequently be denoted with a superscript $(-)^{\Delta}$.

Suppose $\CC^{\Delta}$ a simplicial category. Then we write $(\CC^{\Delta})_{0}$ for the ordinary category given by taking the $0$-simplices of the $\Mor$ spaces. That is, $(\CC^{\Delta})_{0}$ is the category whose objects are the objects of $\CC$, and whose morphisms are given by
\begin{equation*}
(\CC^{\Delta})_{0}(x,y)\coloneq\CC^{\Delta}(x,y)_0.
\end{equation*}
If the $\Mor$ spaces of $\CC^{\Delta}$ are all fibrant, then we will often write
\begin{equation*}
\CC\textrm{\quad for the simplicial nerve\quad}N(\CC^{\Delta})
\end{equation*}
\cite[Df. 1.1.5.5]{HTT}, which is an $\infty$-category \cite[Pr. 1.1.5.10]{HTT}.
\end{ntn}

It will also be convenient to have a model of various $\infty$-categories as \emph{relative categories} \cite{BarKan1011}. To make this precise, we recall the following.
\begin{dfn}\label{dfn:relnerve} A \textbf{\emph{relative category}} is an ordinary category $C$ along with a subcategory $wC$ that contains the identity maps of $C$. The maps of $wC$ will be called \textbf{\emph{weak equivalences}}. A \textbf{\emph{relative functor}} $\fromto{(C,wC)}{(D,wD)}$ is a functor $\fromto{C}{D}$ that carries $wC$ to $wD$.

Suppose $(C,wC)$ a relative category. A \textbf{\emph{relative nerve}} of $(C, wC)$ consists of an $\infty$-category $A$ equipped and a functor $p\colon\fromto{NC}{A}$ that satisfies the following universal property. For any $\infty$-category $B$, the induced functor
\begin{equation*}
\fromto{\Fun(A,B)}{\Fun(NC,B)}
\end{equation*}
is fully faithful, and its essential image is the full subcategory spanned by those functors $\fromto{NC}{B}$ that carry the edges of $wC$ to equivalences in $B$. We will say that the functor $p$ \textbf{\emph{exhibits $A$ as a relative nerve of $(C,wC)$}}.
\end{dfn}

Since relative nerves are defined via a universal property, they are unique up to a contractible choice. Conversely, note that the property of being a relative nerve is invariant under equivalences of $\infty$-categories. That is, if $(C,wC)$ is a relative category, then for any commutative diagram
\begin{equation*}
\begin{tikzpicture} 
\matrix(m)[matrix of math nodes, 
row sep=4ex, column sep=4ex, 
text height=1.5ex, text depth=0.25ex] 
{&NC&\\ 
A'&&A\\}; 
\path[>=stealth,->,font=\scriptsize] 
(m-1-2) edge node[above]{$p'$} (m-2-1) 
edge node[above]{$p$} (m-2-3) 
(m-2-1) edge node[below]{$\sim$} (m-2-3); 
\end{tikzpicture}
\end{equation*}
in which $\equivto{A'}{A}$ is an equivalence of $\infty$-categories, the functor $p'$ exhibits $A'$ as a relative nerve of $(C,wC)$ if and only if $p$ exhibits $A$ as a relative nerve of $(C,wC)$.

\begin{rec}\label{nul:constructrelnerve} There are several functorial constructions of a relative nerve of a relative category $(C,wC)$, all of which are (necessarily) equivalent.
\begin{enumerate}[(\ref{nul:constructrelnerve}.1)]
\item One may form the hammock localization $\mathrm{L}^{\!\mathrm{H}}(C,wC)$ \cite{MR81h:55019}; then a relative nerve can be constructed as the simplicial nerve of the natural functor $\fromto{C}{R(\mathrm{L}^{\!\mathrm{H}}(C,wC))}$, where $R$ denotes any fibrant replacement for the Bergner model structure \cite{MR2276611}.
\item One may mark the edges of $NC$ that correspond to weak equivalences in $C$ to obtain a marked simplicial set \cite[\S 3.1]{HTT}; then one may use the cartesian model structure on marked simplicial sets (over $\Delta^0$) to find a marked anodyne morphism
\begin{equation*}
\fromto{(NC,NwC)}{(N(C,wC),\iota N(C,wC))},
\end{equation*}
wherein $N(C,wC)$ is an $\infty$-category. This map then exhibits the $\infty$-category $N(C,wC)$ as a relative nerve of $(C,wC)$.
\item A relative nerve can be constructed as a fibrant model of the homotopy pushout in the Joyal model structure \cite[\S 2.2.5]{HTT} on simplicial sets of the map
\begin{equation*}
\fromto{\coprod_{\phi\in wC}\Delta^1}{\coprod_{\phi\in wC}\Delta^0}
\end{equation*}
along the map $\fromto{\coprod_{\phi\in wC}\Delta^1}{NC}$.
\end{enumerate}
\end{rec}


\subsection*{The $\infty$-category of $\infty$-categories} The homotopy theory of $\infty$-categories is encoded first as a simplicial category, and then, by application of the simplicial nerve \cite[Df. 1.1.5.5]{HTT}, as an $\infty$-category. This is a pattern that we will follow to describe the homotopy theory of pairs of $\infty$-categories below in Nt. \ref{ntn:pair}.

To begin, recall that an ordinary category $C$ contains a largest subgroupoid, which consists of all objects of $C$ and all isomorphisms between them. The higher categorical analogue of this follows.
\begin{ntn}\label{ntn:interior} For any $\infty$-category $A$, there exists a simplicial subset $\iota A\subset A$, which is the largest Kan simplicial subset of $A$ \cite[1.2.5.3]{HTT}. We shall call this space the \textbf{\emph{interior $\infty$-groupoid of $A$}}. The assignment $\goesto{A}{\iota A}$ defines a right adjoint $\iota$ to the inclusion functor $u$ from Kan simplicial sets to $\infty$-categories.
\end{ntn}

We may think of $\iota A$ as the \emph{moduli space of objects} of $A$, to which we alluded in the introduction.

\begin{ntn}\label{ntn:catofcat} The large simplicial category $\Kan^{\Delta}$ is the category of small Kan simplicial sets, with the usual notion of mapping space. The large simplicial category $\Cat^{\Delta}_{\infty}$ is defined in the following manner \cite[Df. 3.0.0.1]{HTT}. The objects of $\Cat^{\Delta}_{\infty}$ are small $\infty$-categories, and for any two $\infty$-categories $A$ and $B$, the morphism space
\begin{equation*}
\Cat^{\Delta}_{\infty}(A,B)\coloneq\iota\Fun(A,B)
\end{equation*}
is the interior $\infty$-groupoid of the $\infty$-category $\Fun(A,B)$.

Similarly, we may define the huge simplicial category $\Kan(\kappa_1)^{\Delta}$ of large simplicial sets and the huge simplicial category $\Cat_{\infty}(\kappa_1)^{\Delta}$ of large $\infty$-categories.
\end{ntn}

\begin{rec}\label{nul:Catrelnerve} Denote by
\begin{equation*}
w(\Kan^{\Delta})_0\subset(\Kan^{\Delta})_0
\end{equation*}
the subcategory of the ordinary category of Kan simplicial sets (Nt. \ref{ntn:superscriptsscats}) consisting of weak equivalences of simplicial sets. Then, since $(\Kan^{\Delta},w(\Kan^{\Delta})_0)$ is part of a simplicial model structure, it follows that $\Kan$ is a relative nerve of $((\Kan^{\Delta})_0,w(\Kan^{\Delta})_0)$. Similarly, if one denotes by
\begin{equation*}
w(\Cat_{\infty}^{\Delta})_{0}\subset(\Cat_{\infty}^{\Delta})_{0}
\end{equation*}
the subcategory of categorical equivalences of $\infty$-categories, then $\Cat_{\infty}$ is a relative nerve (Df. \ref{dfn:relnerve}) of $(\Cat^{\Delta}_{\infty})_{0},w(\Cat^{\Delta}_{\infty})_{0})$. This follows from \cite[Pr. 3.1.3.5, Pr. 3.1.3.7, Cor. 3.1.4.4]{HTT}.

Since the functors $u$ and $\iota$ (Nt. \ref{ntn:interior}) each preserve weak equivalences, they give rise to an adjunction of $\infty$-categories \cite[Df. 5.2.2.1, Cor. 5.2.4.5]{HTT}
\begin{equation*}
\adjunct{u}{\Kan}{\Cat_{\infty}}{\iota}.
\end{equation*}
\end{rec}


\subsection*{Subcategories of $\infty$-categories} The notion of a \emph{subcategory} of an $\infty$-category is designed to be completely homotopy-invariant. Consequently, given an $\infty$-category $A$ and a simplicial subset $A'\subset A$, we can only call $A'$ a subcategory of $A$ if the following condition holds: any two equivalent morphisms of $A$ both lie in $A'$ just in case either of them does. That is, $A'\subset A$ is completely specified by a subcategory $(hA)'\subset hA$ of the homotopy category $hA$ of $A$.

\begin{rec}\label{rec:subcats} Recall \cite[\S 1.2.11]{HTT} that a \textbf{\emph{subcategory}} of an $\infty$-category $A$ is a simplicial subset $A'\subset A$ such that for some subcategory $(hA)'$ of the homotopy category $hA$, the square
\begin{equation*}
\begin{tikzpicture} 
\matrix(m)[matrix of math nodes, 
row sep=4ex, column sep=4ex, 
text height=1.5ex, text depth=0.25ex] 
{A'&A\\ 
N(hA)'&N(hA)\\}; 
\path[>=stealth,->,font=\scriptsize] 
(m-1-1) edge[right hook->] (m-1-2) 
edge (m-2-1) 
(m-1-2) edge (m-2-2) 
(m-2-1) edge[right hook->] (m-2-2); 
\end{tikzpicture}
\end{equation*}
is a pullback diagram of simplicial sets. In particular, note that a subcategory of an $\infty$-category is uniquely specified by specifying a subcategory of its homotopy category. Note also that any inclusion $\into{A'}{A}$ of a subcategory is an inner fibration \cite[Df. 2.0.0.3, Pr. 2.3.1.5]{HTT}.

We will say that $A'\subset A$ is a \textbf{\emph{full subcategory}} if $(hA)'\subset hA$ is a full subcategory. In this case, $A'$ is uniquely determined by the set $A'_0$ of vertices of $A'$, and we say that $A'$ is \textbf{\emph{spanned by}} the set $A'_0$.

We will say that $A'$ is \textbf{\emph{stable under equivalences}} if the subcategory $(hA)'\subset hA$ above can be chosen to be stable under isomorphisms. Note that any inclusion $\into{A'}{A}$ of a subcategory that is stable under equivalences is a categorical fibration, i.e., a fibration for the Joyal model structure \cite[Cor. 2.4.6.5]{HTT}.
\end{rec}


\subsection*{Pairs of $\infty$-categories} Now we are prepared to introduce the notion of a \emph{pair} of $\infty$-categories.

\begin{dfn}\label{dfn:pair}
\begin{enumerate}[(\ref{dfn:pair}.1)]
\item By a \textbf{\emph{pair $(\mathscr{C},\mathscr{C}_{\dag})$ of $\infty$-categories}} (or simply a \textbf{\emph{pair}}), we shall mean an $\infty$-category $\mathscr{C}$ along with a subcategory \textup{(\ref{rec:subcats})} $\mathscr{C}_{\dag}\subset\mathscr{C}$ containing the maximal Kan complex $\iota\mathscr{C}\subset\mathscr{C}$. We shall call $\mathscr{C}$ the \textbf{\emph{underlying $\infty$-category}} of the pair $(\mathscr{C},\mathscr{C}_{\dag})$. A morphism of $\mathscr{C}_{\dag}$ will be said to be an \textbf{\emph{ingressive}} morphism. 
\item A \textbf{\emph{functor of pairs}} $\psi\colon\fromto{(\mathscr{C},\mathscr{C}_{\dag})}{(\mathscr{D},\mathscr{D}_{\dag})}$ is functor $\fromto{\mathscr{C}}{\mathscr{D}}$ that carries ingressive morphisms to ingressive morphisms; that is, it is a (strictly!) commutative diagram \addtocounter{equation}{2}
\begin{equation}\label{eqn:functorofpairs}
\begin{tikzpicture}[baseline]
\matrix(m)[matrix of math nodes, 
row sep=4ex, column sep=4ex, 
text height=1.5ex, text depth=0.25ex] 
{\mathscr{C}_{\dag}&\mathscr{D}_{\dag}\\ 
\mathscr{C}&\mathscr{D}\\}; 
\path[>=stealth,->,font=\scriptsize] 
(m-1-1) edge node[above]{$\psi_{\dag}$} (m-1-2) 
edge[right hook->] (m-2-1) 
(m-1-2) edge[left hook->] (m-2-2) 
(m-2-1) edge node[below]{$\psi$} (m-2-2); 
\end{tikzpicture}
\end{equation}
of $\infty$-categories.
\addtocounter{enumi}{1}
\item A functor of pairs $\fromto{\mathscr{C}}{\mathscr{D}}$ is said to be \textbf{\emph{strict}} if a morphism of $\mathscr{C}$ is ingressive just in case its image in $\mathscr{D}$ is so --- that is, if the diagram \eqref{eqn:functorofpairs} is a pullback diagram in $\Cat_{\infty}$.
\item\label{item:subpair} A \textbf{\emph{subpair}} of a pair $(\mathscr{C},\mathscr{C}_{\dag})$ is a subcategory \textup{(\ref{rec:subcats})} $\mathscr{D}\subset\mathscr{C}$ equipped with a pair structure $(\mathscr{D},\mathscr{D}_{\dag})$ such that the inclusion $\into{\mathscr{D}}{\mathscr{C}}$ is a strict functor of pairs. If the subcategory $\mathscr{D}\subset\mathscr{C}$ is full, then we'll say that $(\mathscr{D},\mathscr{D}_{\dag})$ is a \textbf{\emph{full subpair}} of $(\mathscr{C},\mathscr{C}_{\dag})$.
\end{enumerate}
\end{dfn}

Since a subcategory of an $\infty$-category is uniquely specified by a subcategory of its homotopy category, and since a morphism of an $\infty$-category is an equivalence if and only if the corresponding morphism of the homotopy category is an isomorphism \cite[Pr. 1.2.4.1]{HTT}, we deduce that a pair $(\mathscr{C},\mathscr{C}_{\dag})$ of $\infty$-categories may simply be described as an $\infty$-category $\mathscr{C}$ and a subcategory $(h\mathscr{C})_{\dag}\subset h\mathscr{C}$ of the homotopy category that contains all the isomorphisms. In particular, note that $\mathscr{C}_{\dag}$ contains all the objects of $\mathscr{C}$.

Note that pairs are a bit rigid: if $(\mathscr{C},\mathscr{C}_{\dag})$ and $(\mathscr{C},\mathscr{C}_{\dag\dag})$ are two pairs, then any equivalence of $\infty$-categories $\equivto{\mathscr{C}_{\dag}}{\mathscr{C}_{\dag\dag}}$ that is (strictly) compatible with the inclusions into $\mathscr{C}$ must be the identity. It follows that for any equivalence of $\infty$-categories $\equivto{C}{D}$, the set of pairs with underlying $\infty$-category $C$ is in bijection with the set of pairs with underlying $\infty$-category $D$.

Consequently, we will often identify a pair $(\mathscr{C},\mathscr{C}_{\dag})$ of $\infty$-categories by defining the underlying $\infty$-category $\mathscr{C}$ and then declaring which morphisms of $\mathscr{C}$ are ingressive. As long as the condition given holds for all equivalences and is stable under homotopies between morphisms and under composition, this will specify a well-defined pair of $\infty$-categories.

\begin{ntn} Suppose $(\mathscr{C},\mathscr{C}_{\dag})$ a pair. Then an ingressive morphism will frequently be denoted by an arrow with a tail: $\!\cofto{}{}\!$. We will often abuse notation by simply writing $\mathscr{C}$ for the pair $(\mathscr{C},\mathscr{C}_{\dag})$.
\end{ntn}

\begin{exm}\label{exm:minpairmaxpair} Any $\infty$-category $C$ can be given the structure of a pair in two ways: the \textbf{\emph{minimal pair}} $C^{\flat}\coloneq(C,\iota C)$ and the \textbf{\emph{maximal pair}} $C^{\sharp}\coloneq(C,C)$.
\end{exm}


\subsection*{The $\infty$-category of pairs} We describe an $\infty$-category $\Pair_{\infty}$ of pairs of $\infty$-categories in much the same manner as we described the $\infty$-category $\Cat_{\infty}$ of $\infty$-categories (Nt. \ref{ntn:catofcat}).

\begin{ntn}\label{ntn:pair} Suppose $\mathscr{C}=(\mathscr{C},\mathscr{C}_{\dag})$ and $\mathscr{D}=(\mathscr{D},\mathscr{D}_{\dag})$ two pairs of $\infty$-categories. Let us denote by $\Fun_{\Pair_{\infty}}(\mathscr{C},\mathscr{D})$ the the full subcategory of the $\infty$-category $\Fun(\mathscr{C},\mathscr{D})$ spanned by the functors $\fromto{\mathscr{C}}{\mathscr{D}}$ that carry ingressives to ingressives.

The large simplicial category $\Pair^{\Delta}_{\infty}$ is defined in the following manner. The objects of $\Pair^{\Delta}_{\infty}$ are small pairs of $\infty$-categories, and for any two pairs of $\infty$-categories $\mathscr{C}$ and $\mathscr{D}$, the morphism space $\Pair^{\Delta}_{\infty}(\mathscr{C},\mathscr{D})$ is interior $\infty$-groupoid (Nt. \ref{ntn:interior})
\begin{equation*}
\Pair^{\Delta}_{\infty}(\mathscr{C},\mathscr{D})\coloneq\iota\Fun_{\Pair_{\infty}}(\mathscr{C},\mathscr{D}).
\end{equation*}
Note that $\Pair^{\Delta}_{\infty}(\mathscr{C},\mathscr{D})$ is the union of connected components of $\Cat_{\infty}^{\Delta}(\mathscr{C},\mathscr{D})$ that correspond to functors of pairs.

Now the $\infty$-category $\Pair_{\infty}$ is the simplicial nerve of this simplicial category (Nt. \ref{ntn:superscriptsscats}).
\end{ntn}


\subsection*{Pair structures} It will be convenient to describe pairs of $\infty$-categories as certain functors between $\infty$-categories. This will permit us to exhibit $\Pair_{\infty}$ as a full subcategory of the arrow $\infty$-category $\Fun(\Delta^{1},\Cat_{\infty})$. This description will in fact imply (Pr. \ref{thm:pairasloc}) that the $\infty$-category $\Pair_{\infty}$ is presentable.

\begin{ntn} For any simplicial set $X$, write $\mathscr{O}(X)$ for the simplicial mapping space from $\Delta^1$ to $X$, whose $n$-simplices are given by
\begin{equation*}
\mathscr{O}(X)_n=\Mor(\Delta^1\times\Delta^n,X).
\end{equation*}
If $C$ is an $\infty$-category, then $\mathscr{O}(C)=\Fun(\Delta^1,C)$ is an $\infty$-category as well \cite[Pr. 1.2.7.3]{HTT}; this is the \textbf{\emph{arrow $\infty$-category}} of $C$. (In fact, $\mathscr{O}$ is a right Quillen functor for the Joyal model structure, since this model structure is cartesian.)
\end{ntn}

\begin{dfn} Suppose $C$ and $D$ $\infty$-categories. We say that a functor $\fromto{D}{C}$ \textbf{\emph{exhibits a pair structure on $C$}} if it factors as an equivalence $\equivto{D}{E}$ followed by an inclusion $\into{E}{C}$ of a subcategory \textup{(\ref{rec:subcats})} such that $(C,E)$ is a pair.
\end{dfn}

\begin{lem}\label{lem:pairstruct} Suppose $C$ and $D$ $\infty$-categories. Then a functor $\psi\colon\fromto{D}{C}$ exhibits a pair structure on $C$ if and only if the following conditions are satisfied.
\begin{enumerate}[(\ref{lem:pairstruct}.1)]
\item The functor $\psi$ induces an equivalence $\fromto{\iota D}{\iota C}$.
\item The functor $\psi$ is a (homotopy) monomorphism in the $\infty$-category $\Cat_{\infty}$; i.e., the diagonal morphism
\begin{equation*}
\fromto{D}{D\times^h_CD}
\end{equation*}
in $h\Cat_{\infty}$ is an isomorphism.
\end{enumerate}
\begin{proof} Clearly any equivalence of $\infty$-categories satisfies these criteria. If $\psi$ is an inclusion of a subcategory such that $(C,D)$ is a pair, then $\psi$, restricted to $\iota D$, is the identity map, and it is an inner fibration such that the diagonal map $\fromto{D}{D\times_CD}$ is an isomorphism. This shows that if $\psi$ exhibits a pair structure on $C$, then $\psi$ satisfies the conditions listed.

Conversely, suppose $\psi$ satisfies the conditions listed. Then it is hard not to show that for any objects $x,y\in D$, the functor $\psi$ induces a homotopy monomorphism
\begin{equation*}
\fromto{\Map_D(x,y)}{\Map_C(\psi(x),\psi(y))},
\end{equation*}
whence the natural map
\begin{equation*}
\fromto{\Map_D(x,y)}{\Map_{NhD}(x,y)\times^h_{\Map_{NhC}(\psi(x),\psi(y))}\Map_C(\psi(x),\psi(y))}
\end{equation*}
is a weak equivalence. This, combined with the fact that the map $\fromto{\iota D}{\iota C}$ is an equivalence, now implies that the natural map $\fromto{D}{NhD\times^h_{NhC}C}$ of $\Cat_{\infty}$ is an equivalence.

Since isomorphisms in $hC$ are precisely equivalences in $C$, the induced functor of homotopy categories $\fromto{hD}{hC}$ identifies $hD$ with a subcategory of $hC$ that contains all the isomorphisms. Denote by $hE\subset hC$ this subcategory. Now let $E$ be the corresponding subcategory of $C$; we thus have a diagram of $\infty$-categories
\begin{equation*}
\begin{tikzpicture} 
\matrix(m)[matrix of math nodes, 
row sep=4ex, column sep=4ex, 
text height=1.5ex, text depth=0.25ex] 
{D&E&C\\ 
NhD&NhE&NhC\\}; 
\path[>=stealth,->,font=\scriptsize,inner sep=0.75pt] 
(m-1-1) edge (m-1-2) 
edge (m-2-1)
(m-1-2) edge[right hook->] (m-1-3)
edge (m-2-2)
(m-1-3) edge (m-2-3)
(m-1-2) edge (m-2-2)
(m-2-1) edge node[below]{$\sim$} (m-2-2)
(m-2-2) edge[right hook->] (m-2-3); 
\end{tikzpicture}
\end{equation*}
in which the square on the right and the big rectangle are homotopy pullbacks (for the Joyal model structure). Thus the square on the left is a homotopy pullback as well, and so the functor $\fromto{D}{E}$ is an equivalence, giving our desired factorization.
\end{proof}
\end{lem}

\begin{cnstr} We now consider the following simplicial functor
\[U':\fromto{\Pair^{\Delta}_{\infty}}{\Fun(\mathbf{[1]},\Cat_{\infty}^\Delta)}.\]
On objects, $U'$ carries a pair $(\mathscr{C},\mathscr{C}_\dag)$ to the inclusion of $\infty$-categories $\into{\mathscr{C}_\dag}{\mathscr{C}}$. On mapping spaces, $U'$ is given by the obvious forgetful maps
\[\fromto{\iota\Fun_{\Pair_{\infty}}((\mathscr{C},\mathscr{C}_\dag),(\mathscr{D},\mathscr{D}_\dag))}{\iota\Fun(\mathscr{C},\mathscr{D})\times_{\iota\Fun(\mathscr{C}_\dag,\mathscr{D})}\iota\Fun(\mathscr{C}_\dag,\mathscr{D}_\dag).}\]

Now note that since $\into{\iota\Fun(\mathscr{C}_\dag,\mathscr{D}_\dag)}{\iota\Fun(\mathscr{C}_\dag,\mathscr{D})}$ is the inclusion of a union of connected components, it follows that the projection
\[\fromto{\iota\Fun(\mathscr{C},\mathscr{D})\times_{\iota\Fun(\mathscr{C}_\dag,\mathscr{D})}\iota\Fun(\mathscr{C}_\dag,\mathscr{D}_\dag)}{\iota\Fun(\mathscr{C},\mathscr{D})}\]
is an inclusion of a union of connected components as well; in particular, it is the inclusion of those connected components corresponding to those functors $\fromto{\mathscr{C}}{\mathscr{D}}$ that carry morphisms of $\mathscr{C}_\dag$ to morphisms of $\mathscr{D}_\dag$. That is, the inclusion
\[\into{\iota\Fun_{\Pair_{\infty}}((\mathscr{C},\mathscr{C}_\dag),(\mathscr{D},\mathscr{D}_\dag))}{\iota\Fun(\mathscr{C},\mathscr{D})}\]
factors through an equivalence
\[\equivto{\iota\Fun_{\Pair_{\infty}}((\mathscr{C},\mathscr{C}_\dag),(\mathscr{D},\mathscr{D}_\dag))}{\iota\Fun(\mathscr{C},\mathscr{D})\times_{\iota\Fun(\mathscr{C}_\dag,\mathscr{D})}\iota\Fun(\mathscr{C}_\dag,\mathscr{D}_\dag)}.\]
In other words, the functor $U'$ is fully faithful.
\end{cnstr}

We therefore conclude
\begin{prp}\label{prp:PairsubcatOcat} The functor
\[\fromto{\Pair_{\infty}}{N\Fun([1],\Cat_{\infty}^{\Delta})\simeq\mathscr{O}(\Cat_{\infty})}\]
induced by $U'$ exhibits an equivalence between $\Pair_{\infty}$ and the full subcategory of $\mathscr{O}(\Cat_{\infty})$ spanned by those functors $\fromto{D}{C}$ that exhibit a pair structure on $C$.
\end{prp}


\subsection*{The $\infty$-categories of pairs as a relative nerve} It will be convenient for us to have a description of $\Pair_{\infty}$ as a relative nerve (Df. \ref{dfn:relnerve}). First, we record the following trivial result.

\begin{prp}\label{lem:functorofpairsequiv} The following are equivalent for a functor of pairs $\psi\colon\fromto{\mathscr{C}}{\mathscr{D}}$.
\begin{enumerate}[(\ref{lem:functorofpairsequiv}.1)]
\item The functor of pairs $\psi$ is an equivalence in the $\infty$-category $\Pair_{\infty}$.
\item The underlying functor of $\infty$-categories is a categorical equivalence, and $\psi$ is strict.
\item The underlying functor of $\infty$-categories is a categorical equivalence that induces an equivalence $h\mathscr{C}_{\dag}\simeq h\mathscr{D}_{\dag}$.
\end{enumerate}
\begin{proof} The equivalence of the first two conditions follows from the equivalence between $\Pair_{\infty}$ and a full subcategory of $\mathscr{O}(\Cat_{\infty})$. The second condition clearly implies the third. To prove that the third condition implies the second, consider the commutative diagram
\begin{equation*}
\begin{tikzpicture}[cross line/.style={preaction={draw=white, -, 
line width=6pt}}]
\matrix(m)[matrix of math nodes, 
row sep=2ex, column sep=0.75ex, 
text height=1.5ex, text depth=0.25ex]
{&\mathscr{D}_{\dag}&&\mathscr{D}\\
\mathscr{C}_{\dag}&&\mathscr{C}&\\
&Nh\mathscr{D}_{\dag}&&Nh\mathscr{D}.\\
Nh\mathscr{C}_{\dag}&&Nh\mathscr{C}&\\
};
\path[>=stealth,->,font=\scriptsize]
(m-1-2) edge[<-] (m-2-1)
edge (m-3-2)
edge[right hook->] (m-1-4)
(m-3-2) edge[<-] (m-4-1)
edge[right hook->] (m-3-4)
(m-2-1) edge[cross line,right hook->] (m-2-3)
edge (m-4-1)
(m-1-4) edge[<-] (m-2-3)
edge (m-3-4)
(m-4-1) edge[right hook->] (m-4-3)
(m-3-4) edge[<-] (m-4-3)
(m-2-3) edge[cross line] (m-4-3);
\end{tikzpicture}
\end{equation*}
The front and back faces are pullback squares and therefore homotopy pullback squares. Since both $\equivto{Nh\mathscr{C}}{Nh\mathscr{D}}$ and $\equivto{Nh\mathscr{C}_{\dag}}{Nh\mathscr{D}_{\dag}}$ are equivalences, the bottom face is a homotopy pullback as well. Hence the top square is a homotopy pullback. But since $(\mathscr{C},\mathscr{C}_{\dag})$ is a pair, it must be an actual pullback; that is, $\psi$ is strict.
\end{proof}
\end{prp}

\noindent This proposition implies that the $\infty$-category of functors of pairs is compatible with equivalences of pairs.
\begin{cor} Suppose $\mathscr{A}$ a pair, and suppose $\equivto{\mathscr{C}}{\mathscr{D}}$ an equivalence of pairs of $\infty$-categories. Then the induced functor $\fromto{\Fun_{\Pair_{\infty}}(\mathscr{A},\mathscr{C})}{\Fun_{\Pair_{\infty}}(\mathscr{A},\mathscr{D})}$ is an equivalence of $\infty$-categories.
\begin{proof} The proposition implies that any homotopy inverse $\equivto{\mathscr{D}}{\mathscr{C}}$ of the equivalence $\equivto{\mathscr{C}}{\mathscr{D}}$ of underlying $\infty$-categories must carry ingressives to ingressives. This induces a homotopy inverse $\fromto{\Fun_{\Pair_{\infty}}(\mathscr{A},\mathscr{D})}{\Fun_{\Pair_{\infty}}(\mathscr{A},\mathscr{C})}$, completing the proof.
\end{proof}
\end{cor}

\noindent Furthermore, Pr. \ref{lem:functorofpairsequiv} may be combined with Pr. \ref{prp:PairsubcatOcat} and \ref{nul:Catrelnerve} to yield the following.
\begin{cor}\label{prp:Pairisrelnerve} Denote by $w(\Pair_{\infty}^{\Delta})_0\subset(\Pair_{\infty}^{\Delta})_0$ the subcategory of the ordinary category of pairs of $\infty$-categories \textup{(Nt. \ref{ntn:superscriptsscats})} consisting of those functors of pairs $\fromto{\mathscr{C}}{\mathscr{D}}$ whose underlying functor of $\infty$-categories is a categorical equivalence that induces an equivalence $h\mathscr{C}_{\dag}\simeq h\mathscr{D}_{\dag}$. Then the $\infty$-category $\Pair_{\infty}$ is a relative nerve \textup{(Df. \ref{dfn:relnerve})} of the relative category $((\Pair_{\infty}^{\Delta})_0,w(\Pair_{\infty}^{\Delta})_0)$.
\end{cor}


\subsection*{The dual picture} Let us conclude this section by briefly outlining the dual picture of $\infty$-categories with \emph{fibrations}.

\begin{dfn}\label{dfn:egressives} Suppose $(\mathscr{C}^{\op},(\mathscr{C}^{\op})_{\dag})$ a pair. Then write $\mathscr{C}^{\dag}$ for the subcategory
\begin{equation*}
((\mathscr{C}^{\op})_{\dag})^{\op}\subset\mathscr{C}.
\end{equation*}
We call the morphisms of $\mathscr{C}^{\dag}$ \textbf{\emph{egressive}} morphisms or \textbf{\emph{fibrations}}. The pair $(\mathscr{C},\mathscr{C}^{\dag})$ will be called the \textbf{\emph{opposite pair}} to $(\mathscr{C}^{\op},(\mathscr{C}^{\op})_{\dag})$. One may abuse terminology slightly by referring to $(\mathscr{C},\mathscr{C}^{\dag})$ as a \textbf{\emph{pair structure on $\mathscr{C}^{\op}$}}.
\end{dfn}

\begin{ntn} Suppose $(\mathscr{C}^{\op},(\mathscr{C}^{\op})_{\dag})$ a pair. Then a fibration of $\mathscr{C}$ will frequently be denoted by a double headed arrow: $\!\fibto{}{}\!$. We will often abuse notation by simply writing $\mathscr{C}$ for the opposite pair $(\mathscr{C},\mathscr{C}^{\dag})$.
\end{ntn}

\noindent We summarize this discussion with the following.
\begin{prp}\label{prp:opinvolonPair} The formation $\goesto{(\mathscr{C},\mathscr{C}_{\dag})}{(\mathscr{C}^{\op},(\mathscr{C}^{\op})^{\dag})}$ of the opposite pair defines an involution $(-)^{\op}$ of the $\infty$-category $\Pair_{\infty}$.
\end{prp}


\section{Waldhausen $\infty$-categories} In developing his abstract framework for $K$-theory, Waldhausen introduced first \cite[\S 1.1]{MR86m:18011} the notion of a \emph{category with cofibrations}, and then \cite[\S 1.2]{MR86m:18011} layered the added structure of a subcategory of weak equivalences satisfying some additional compatibilities to obtain what today is often called a \emph{Waldhausen category}. This added structure introduces homotopy theory, and Waldhausen required that the structure of a category with cofibrations interacts well with this homotopy theory.

The theory of \emph{Waldhausen $\infty$-categories}, which we introduce in this section, reverses these two priorities. The layer of homotopy theory is already embedded in the implementation of quasicategories. Then, because it is effectively impossible to formulate $\infty$-categorical notions that do \emph{not} interact well with the homotopy theory, we arrive at a suitable definition of Waldhausen $\infty$-categories by writing the quasicategorical analogues of the axioms for Waldhausen's categories with cofibrations. Consequently, a Waldhausen $\infty$-category will be a pair of $\infty$-categories that enjoys the following properties.
\begin{enumerate}[(A)]
\item The underlying $\infty$-category admits a \emph{zero object} $0$ such that the morphisms $\fromto{0}{X}$ are all ingressive.
\item \emph{Pushouts} of ingressives exist and are ingressives.
\end{enumerate}


\subsection*{Limits and colimits in $\infty$-categories} To work with these conditions effectively, it is convenient to fix some notations and terminology for the study of \emph{limits} and \emph{colimits} in $\infty$-categories, as defined in \cite[\S 1.2.13]{HTT}.

\begin{rec} Recall \cite[Df. 1.2.12.1]{HTT} that an object $X$ of an $\infty$-category $C$ is said to be \textbf{\emph{initial}} if for any object $Y$ of $C$, the mapping space $\Map(X,Y)$ is weakly contractible. Dually, $X$ is said to be \textbf{\emph{terminal}} if for any object $Y$ of $C$, the mapping space $\Map(Y,X)$ is weakly contractible.
\end{rec}

\begin{dfn} A \textbf{\emph{zero object}} of an $\infty$-category is an object that is both initial and terminal.
\end{dfn}

\begin{ntn} For any simplicial set $K$, one has \cite[Nt. 1.2.8.4]{HTT} the \textbf{\emph{right cone}} $K^{\rhd}\coloneq K\star\Delta^0$ and the \textbf{\emph{left cone}} $K^{\lhd}\coloneq\Delta^0\star K$; we write $+\infty$ for the cone point of $K^{\rhd}$, and we write $-\infty$ for the cone point of $K^{\lhd}$.
\end{ntn}

\begin{rec} Just as in ordinary category theory, a colimit and limit in an $\infty$-category can be described as an initial and terminal object of a suitable associated $\infty$-category. For any simplicial set $K$, a \textbf{\emph{limit diagram}} in an $\infty$-category $C$ is a diagram
\begin{equation*}
\overline{p}\colon\fromto{K^{\lhd}}{C}
\end{equation*}
that is a terminal object in the overcategory $C_{/p}$ \cite[\S 1.2.9]{HTT}, where $p=\overline{p}|K$. Dually, a \textbf{\emph{colimit diagram}} in an $\infty$-category $C$ is a diagram
\begin{equation*}
\overline{p}\colon\fromto{K^{\rhd}}{C}
\end{equation*}
that is a terminal object in the undercategory $C_{p/}$, where $p=\overline{p}|K$.

For any $\infty$-category $A$ and any $\infty$-category $C$, we denote by
\begin{equation*}
\Colim(A^{\rhd},C)\subset\Fun(A^{\rhd},C)
\end{equation*}
the full subcategory spanned by colimit diagrams $\fromto{A^{\rhd}}{C}$.
\end{rec}

\begin{dfn} A \textbf{\emph{pushout square}} in an $\infty$-category $C$ is a colimit diagram
\begin{equation*}
X\colon\fromto{(\Lambda^2_0)^{\rhd}\cong\Delta^1\times\Delta^1}{C}.
\end{equation*}
Such a diagram may be drawn
\begin{equation*}
\begin{tikzpicture} 
\matrix(m)[matrix of math nodes, 
row sep=4ex, column sep=4ex, 
text height=1.5ex, text depth=0.25ex] 
{X_{00}&X_{01}\\ 
X_{10}&X_{11};\\}; 
\path[>=stealth,->,font=\scriptsize] 
(m-1-1) edge (m-1-2) 
edge (m-2-1) 
(m-1-2) edge (m-2-2) 
(m-2-1) edge (m-2-2); 
\end{tikzpicture}
\end{equation*}
the edge $\fromto{X_{10}}{X_{11}}$ is called the \textbf{\emph{pushout}} of the edge $\fromto{X_{00}}{X_{01}}$.
\end{dfn}

\begin{rec}\label{rec:Ashapedcolim} A key result of Joyal \cite[Pr. 1.2.12.9]{HTT} states that for any functor $\psi\colon\fromto{A}{C}$, the fiber of the canonical restriction functor
\begin{equation*}
\rho\colon\fromto{\Colim(A^{\rhd},C)}{\Fun(A,C)}
\end{equation*}
over $\psi$ is either empty or a contractible Kan space. One says that $C$ \textbf{\emph{admits all $A$-shaped colimits}} if the fibers of the functor $\rho$ are all nonempty. In this case, $\rho$ is an equivalence of $\infty$-categories.

More generally, if $\mathscr{A}$ is a family of $\infty$-categories, then one says that $C$ \textbf{\emph{admits all $\mathscr{A}$-shaped colimits}} if the fibers of the functor $\fromto{\Colim(A^{\rhd},C)}{\Fun(A,C)}$ are all nonempty for every $A\in\mathscr{A}$.

Finally, if $\mathscr{A}$ is a family of $\infty$-categories, then a functor $f\colon\fromto{C'}{C}$ will be said to \textbf{\emph{preserve all $\mathscr{A}$-shaped colimits}} if for any $A\in\mathscr{A}$, the composite
\begin{equation*}
\Colim(A^{\rhd},C')\ \tikz[baseline]\draw[>=stealth,->](0,0.5ex)--(0.5,0.5ex);\ \Fun(A^{\rhd},C')\ \tikz[baseline]\draw[>=stealth,->](0,0.5ex)--(0.5,0.5ex);\ \Fun(A^{\rhd},C)
\end{equation*}
factors through $\Colim(A^{\rhd},C)\subset\Fun(A^{\rhd},C)$. We write $\Fun_{\mathscr{A}}(C',C)\subset\Fun(C',C)$ for the full subcategory spanned by those functors that preserve all $\mathscr{A}$-shaped colimits.
\end{rec}


\subsection*{Waldhausen $\infty$-categories} We are now introduce the notion of \emph{Waldhausen $\infty$-categories}, which are the primary objects of study in this work.

\begin{dfn}\label{dfn:preWald} A \textbf{\emph{Waldhausen $\infty$-category}} $(\mathscr{C},\mathscr{C}_{\dag})$ is a pair of essentially small $\infty$-categories such that the following axioms hold.
\begin{enumerate}[(\ref{dfn:preWald}.1)]
\item The $\infty$-category $\mathscr{C}$ contains a zero object.
\item\label{item:0toxingressive} For any zero object $0$, any morphism $\fromto{0}{X}$ is ingressive.
\item\label{item:pushcof} Pushouts of ingressive morphisms exist. That is, for any diagram $G\colon\fromto{\Lambda^2_0}{\mathscr{C}}$ represented as
\begin{equation*}
\begin{tikzpicture} 
\matrix(m)[matrix of math nodes, 
row sep=4ex, column sep=4ex, 
text height=1.5ex, text depth=0.25ex] 
{X&Y\\ 
X'&\\}; 
\path[>=stealth,->,font=\scriptsize] 
(m-1-1) edge[>->] (m-1-2) 
edge (m-2-1); 
\end{tikzpicture}
\end{equation*}
in which the morphism $\cofto{X}{Y}$ is ingressive, there exists a pushout square $\overline{G}\colon\fromto{(\Lambda^2_0)^{\rhd}\cong\Delta^1\times\Delta^1}{\mathscr{C}}$ extending $G$:
\begin{equation*}
\begin{tikzpicture} 
\matrix(m)[matrix of math nodes, 
row sep=4ex, column sep=4ex, 
text height=1.5ex, text depth=0.25ex] 
{X&Y\\ 
X'&Y'.\\}; 
\path[>=stealth,->,font=\scriptsize] 
(m-1-1) edge[>->] (m-1-2) 
edge (m-2-1) 
(m-1-2) edge (m-2-2) 
(m-2-1) edge (m-2-2); 
\end{tikzpicture}
\end{equation*}
\item\label{item:pushcofcof} Pushouts of ingressives are ingressives. That is, for any pushout square $\fromto{(\Lambda^2_0)^{\rhd}\cong\Delta^1\times\Delta^1}{\mathscr{C}}$ represented as
\begin{equation*}
\begin{tikzpicture} 
\matrix(m)[matrix of math nodes, 
row sep=4ex, column sep=4ex, 
text height=1.5ex, text depth=0.25ex] 
{X&Y\\ 
X'&Y',\\}; 
\path[>=stealth,->,font=\scriptsize] 
(m-1-1) edge[>->] (m-1-2) 
edge (m-2-1) 
(m-1-2) edge (m-2-2) 
(m-2-1) edge (m-2-2); 
\end{tikzpicture}
\end{equation*}
if the morphism $\cofto{X}{Y}$ is ingressive, then so is the morphism $\fromto{X'}{Y'}$.
\suspend{enumerate}

Call a functor of pairs $\psi\colon\fromto{\mathscr{C}}{\mathscr{D}}$ between two Waldhausen $\infty$-categories \textbf{\emph{exact}} if it satisfies the following conditions.
\resume{enumerate}[{[(\ref{dfn:preWald}.1)]}]
\item The underlying functor of $\psi$ carries zero objects of $\mathscr{C}$ to zero objects of $\mathscr{D}$.
\item\label{item:psipreservepushcof} For any pushout square $F\colon\fromto{(\Lambda^2_0)^{\rhd}\cong\Delta^1\times\Delta^1}{\mathscr{C}}$ represented as
\begin{equation*}
\begin{tikzpicture} 
\matrix(m)[matrix of math nodes, 
row sep=4ex, column sep=4ex, 
text height=1.5ex, text depth=0.25ex] 
{X&Y\\ 
X'&Y'\\}; 
\path[>=stealth,->,font=\scriptsize] 
(m-1-1) edge[>->] (m-1-2) 
edge (m-2-1) 
(m-1-2) edge (m-2-2) 
(m-2-1) edge[>->] (m-2-2); 
\end{tikzpicture}
\end{equation*}
in which $\cofto{X}{Y}$ and hence $\cofto{X'}{Y'}$ are ingressive, the induced square $\psi\circ F\colon\fromto{(\Lambda^2_0)^{\rhd}\cong\Delta^1\times\Delta^1}{\mathscr{D}}$ represented as
\begin{equation*}
\begin{tikzpicture} 
\matrix(m)[matrix of math nodes, 
row sep=4ex, column sep=4ex, 
text height=1.5ex, text depth=0.25ex] 
{\psi(X)&\psi(Y)\\ 
\psi(X')&\psi(Y')\\}; 
\path[>=stealth,->,font=\scriptsize] 
(m-1-1) edge[>->] (m-1-2) 
edge (m-2-1) 
(m-1-2) edge (m-2-2) 
(m-2-1) edge[>->] (m-2-2); 
\end{tikzpicture}
\end{equation*}
is a pushout as well.
\end{enumerate}

A \textbf{\emph{Waldhausen subcategory}} of a Waldhausen $\infty$-category $\mathscr{C}$ is a subpair $\mathscr{D}\subset\mathscr{C}$ such that $\mathscr{D}$ is a Waldhausen $\infty$-category, and the inclusion $\into{\mathscr{D}}{\mathscr{C}}$ is exact.
\end{dfn}

Let us repackage some of these conditions.

\begin{nul}\label{exm:LambdaDelta} Denote by $\Lambda_0\mathscr{Q}^2$ the pair $(\Lambda^2_0,\Delta^{\{0,1\}}\sqcup\Delta^{\{2\}})$, which may be represented as
\begin{equation*}
\begin{tikzpicture} 
\matrix(m)[matrix of math nodes, 
row sep=4ex, column sep=4ex, 
text height=1.5ex, text depth=0.25ex] 
{0&1\\ 
2.&\\}; 
\path[>=stealth,->,font=\scriptsize] 
(m-1-1) edge[>->] (m-1-2) 
edge (m-2-1); 
\end{tikzpicture}
\end{equation*} 
Denote by $\mathscr{Q}^2$ the pair
\begin{equation*}
((\Lambda^2_0)^{\rhd},\Delta^{\{0,1\}}\sqcup\Delta^{\{2,\infty\}})\cong(\Delta^1)^{\flat}\times(\Delta^1)^{\sharp}\cong(\Delta^1\times\Delta^1,(\Delta^{\{0\}}\sqcup\Delta^{\{1\}})\times\Delta^1)
\end{equation*}
(Ex. \ref{exm:minpairmaxpair}), which may be represented as
\begin{equation*}
\begin{tikzpicture} 
\matrix(m)[matrix of math nodes, 
row sep=4ex, column sep=4ex, 
text height=1.5ex, text depth=0.25ex] 
{0&1\\ 
2&\infty.\\}; 
\path[>=stealth,->,font=\scriptsize] 
(m-1-1) edge[>->] (m-1-2) 
edge (m-2-1) 
(m-1-2) edge (m-2-2) 
(m-2-1) edge[>->] (m-2-2); 
\end{tikzpicture}
\end{equation*}
There is an obvious strict inclusion of pairs $\into{\Lambda_0\mathscr{Q}^2}{\mathscr{Q}^2}$.

Conditions (\ref{dfn:preWald}.\ref{item:pushcof}) and (\ref{dfn:preWald}.\ref{item:pushcofcof}) can be rephrased as the single condition that the functor
\begin{equation*}
\fromto{\Fun_{\Pair_{\infty}}(\mathscr{Q}^2,\mathscr{C})}{\Fun_{\Pair_{\infty}}(\Lambda_0\mathscr{Q}^2,\mathscr{C})}
\end{equation*}
induces an equivalence of $\infty$-categories
\begin{equation*}
\equivto{\Colim_{\Pair_{\infty}}(\mathscr{Q}^2,\mathscr{C})}{\Fun_{\Pair_{\infty}}(\Lambda_0\mathscr{Q}^2,\mathscr{C})}
\end{equation*}
where $\Colim_{\Pair_{\infty}}(\mathscr{Q}^2,\mathscr{C})$ denotes the full subcategory of $\Fun_{\Pair_{\infty}}(\mathscr{Q}^2,\mathscr{C})$ spanned by those functors of pairs $\fromto{\mathscr{Q}^2}{\mathscr{C}}$ whose underlying functor $\fromto{(\Lambda^2_0)^{\rhd}}{\mathscr{C}}$ is a pushout square.

Condition (\ref{dfn:preWald}.\ref{item:psipreservepushcof}) on a functor of pairs $\psi\colon\fromto{\mathscr{C}}{\mathscr{D}}$ between Waldhausen $\infty$-categories is equivalent to the condition that the composite functor
\begin{equation*}
\Colim_{\Pair_{\infty}}(\mathscr{Q}^2,\mathscr{C})\subset\Fun_{\Pair_{\infty}}(\mathscr{Q}^2,\mathscr{C})\ {\tikz[baseline]\draw[>=stealth,->](0,0.5ex)--(0.5,0.5ex);}\ \Fun_{\Pair_{\infty}}(\mathscr{Q}^2,\mathscr{D})
\end{equation*}
factors through the full subcategory
\begin{equation*}
\Colim_{\Pair_{\infty}}(\mathscr{Q}^2,\mathscr{D})\subset\Fun_{\Pair_{\infty}}(\mathscr{Q}^2,\mathscr{D}).
\end{equation*}
\end{nul}


\subsection*{Some examples} To get a sense for how these axioms apply, let's give some examples of Waldhausen $\infty$-categories.

\begin{exm}\label{exm:whenisminandmaxWald} When equipped with the \emph{minimal} pair structure (Ex. \ref{exm:minpairmaxpair}), an $\infty$-category $C$ is a Waldhausen $\infty$-category $C^{\flat}$ if and only if $C$ is a contractible Kan complex.

Equipped with the \emph{maximal} pair structure (Ex. \ref{exm:minpairmaxpair}), any $\infty$-category $C$ that admits a zero object and all finite colimits can be regarded as a Waldhausen $\infty$-category $C^{\sharp}$.
\end{exm}

\begin{exm}\label{exm:inftytopoigiveWaldcats} As a special case of the above, suppose that $\mathscr{E}$ is an $\infty$-topos \cite[Df. 6.1.0.2]{HTT}. For example, one may consider the example $\mathscr{E}=\Fun(S,\Kan)$ for some simplicial set $S$. Then the $\infty$-category $\mathscr{E}^{\omega}_{\ast}$ of compact, pointed objects of $\mathscr{E}$, when equipped with its maximal pair structure, is a Waldhausen $\infty$-category. Its algebraic $K$-theory will be called the \emph{$A$-theory of $\mathscr{E}$}. For any Kan simplicial set $X$, the $A$-theory of the $\infty$-topos $\Fun(X,\Kan)$ agrees with Waldhausen's $A$-theory of $X$ (where one defines the latter via the category $\mathscr{R}_{\mathrm{df}}(X)$ of finitely dominated retractive spaces over $X$ \cite[p. 389]{MR86m:18011}). See Ex. \ref{exm:Atheoryofinfintopoi} for more.
\end{exm}

\begin{exm} Any stable $\infty$-category $\mathscr{A}$ \cite[Df. 1.1.1.9]{HA}, when equipped with its maximal pair structure, is a Waldhausen $\infty$-category. If $\mathscr{A}$ admits a t-structure \cite[Df. 1.2.1.4]{HA}, then one may define a pair structure on any of the $\infty$-categories $\mathscr{A}_{\leq n}$ by declaring that a morphism $\fromto{X}{Y}$ be ingressive just in case the induced morphism $\fromto{\pi_nX}{\pi_nY}$ is a monomorphism of the heart $\mathscr{A}^{\heartsuit}$. We study the relationship between the algebraic $K$-theory of these $\infty$-categories to the algebraic $K$-theory of $\mathscr{A}$ itself in a follow-up to this paper \cite{K21}.
\end{exm}

\begin{exm}\label{exm:catwithcofibsisWaldinftycat} If $(C,\cof C)$ is an ordinary \textbf{\emph{category with cofibrations}} in the sense of Waldhausen \cite[\S 1.1]{MR86m:18011}, then the pair $(NC,N(\cof C))$ is easily seen to be a Waldhausen $\infty$-category. If $(C,\cof C,wC)$ is a category with cofibrations and weak equivalences in the sense of Waldhausen \cite[\S 1.2]{MR86m:18011}, then one may endow a relative nerve (Df. \ref{dfn:relnerve}) $N(C,wC)$ of the relative category $(C,wC)$ with a pair structure by defining the subcategory $N(C,wC)_{\dag}\subset N(C,wC)$ as the smallest subcategory containing the equivalences and the images of the edges in $NC$ corresponding to cofibrations. In Pr. \ref{thm:classWaldareWald}, we will show that if $(C,wC)$ is a \emph{partial model category} in which the weak equivalences and trivial cofibrations are part of a three-arrow calculus of fractions, then any relative nerve of $(C,wC)$ is in fact a Waldhausen $\infty$-category with this pair structure.
\end{exm}


\subsection*{The $\infty$-category of Waldhausen $\infty$-categories} We now define the $\infty$-category of Waldhausen $\infty$-categories as a subcategory of the $\infty$-category of pairs.

\begin{ntn}\label{ntn:preWald}
\begin{enumerate}[(\ref{ntn:preWald}.1)]
\item Suppose $\mathscr{C}$ and $\mathscr{D}$ two Waldhausen $\infty$-categories. We denote by $\Fun_{\Wald}(\mathscr{C},\mathscr{D})$ the full subcategory of $\Fun_{\Pair_{\infty}}(\mathscr{C},\mathscr{D})$ spanned by the exact functors $\fromto{\mathscr{C}}{\mathscr{D}}$ of Waldhausen $\infty$-categories.
\item Define $\Wald^{\Delta}$ as the following simplicial subcategory of $\Pair^{\Delta}_{\infty}$. The objects of $\Wald^{\Delta}$ are small Waldhausen $\infty$-categories, and for any Waldhausen $\infty$-categories $\mathscr{C}$ and $\mathscr{D}$, the morphism space $\Wald^{\Delta}(\mathscr{C},\mathscr{D})$ is defined by the formula
\begin{equation*}
\Wald^{\Delta}(\mathscr{C},\mathscr{D})\coloneq\iota\Fun_{\Wald}(\mathscr{C},\mathscr{D}),
\end{equation*}
or, equivalently, $\Wald^{\Delta}(\mathscr{C},\mathscr{D})$ is the union of the connected components of $\Pair^{\Delta}_{\infty}(\mathscr{C},\mathscr{D})$ corresponding to the exact morphisms.
\item We now define the $\infty$-category $\Wald$ as the simplicial nerve of $\Wald^{\Delta}$ (Nt. \ref{ntn:superscriptsscats}), or, equivalently, as the subcategory of $\Pair_{\infty}$ whose objects are Waldhausen $\infty$-categories and whose morphisms are exact functors.
\end{enumerate}
\end{ntn}

\begin{lem} The subcategory $\Wald\subset\Pair_{\infty}$ is stable under equivalences.
\begin{proof} Suppose $\mathscr{C}$ a Waldhausen $\infty$-category, and suppose $\psi\colon\equivto{\mathscr{C}}{\mathscr{D}}$ an equivalence of pairs. The functor of pairs $\psi$ induces an equivalence of underlying $\infty$-categories, whence $\mathscr{D}$ admits a zero object as well. We also have, in the notation of Nt. \ref{exm:LambdaDelta}, a commutative square
\begin{equation*}
\begin{tikzpicture} 
\matrix(m)[matrix of math nodes, 
row sep=4ex, column sep=4ex, 
text height=1.5ex, text depth=0.25ex] 
{\Colim_{\Pair_{\infty}}(\mathscr{Q}^2,\mathscr{C})&\Fun_{\Pair_{\infty}}(\Lambda_0\mathscr{Q}^2,\mathscr{C})\\ 
\Colim_{\Pair_{\infty}}(\mathscr{Q}^2,\mathscr{D})&\Fun_{\Pair_{\infty}}(\Lambda_0\mathscr{Q}^2,\mathscr{D})\\}; 
\path[>=stealth,->,font=\scriptsize] 
(m-1-1) edge (m-1-2) 
edge (m-2-1) 
(m-1-2) edge (m-2-2) 
(m-2-1) edge (m-2-2); 
\end{tikzpicture}
\end{equation*}
in which the top functor is an equivalence since $\mathscr{C}$ is a Waldhausen $\infty$-category, and the vertical functors are equivalences since $\equivto{\mathscr{C}}{\mathscr{D}}$ is an equivalence of pairs. Hence the bottom functor is an equivalence of $\infty$-categories, whence $\mathscr{D}$ is a Waldhausen $\infty$-category.
\end{proof}
\end{lem}


\subsection*{Equivalences between maximal Waldhausen $\infty$-categories} Equivalences between Waldhausen $\infty$-categories with a \emph{maximal} pair structure (Ex \ref{exm:whenisminandmaxWald}) are often easy to detect, thanks to the following result.
\begin{prp}\label{prp:preapprox} Suppose $\mathscr{C}$ and $\mathscr{D}$ two $\infty$-categories that each contain zero objects and all finite colimits. Regard them as Waldhausen $\infty$-categories equipped with the maximal pair structure \textup{(Ex \ref{exm:whenisminandmaxWald})}. Assume that the suspension functor $\Sigma\colon\fromto{\mathscr{C}}{\mathscr{C}}$ is essentially surjective. Then an exact functor $\psi\colon\fromto{\mathscr{C}}{\mathscr{D}}$ is an equivalence if and only if it induces an equivalence of homotopy categories $\equivto{h\mathscr{C}}{h\mathscr{D}}$.
\begin{proof} We need only show that $\psi$ is fully faithful. Since $\psi$ preserves all finite colimits \cite[Cor. 4.4.2.5]{HTT}, it follows that $\psi$ preserves the tensor product with any finite Kan complex \cite[Cor. 4.4.4.9]{HTT}. Thus for any finite simplicial set $K$ and any objects $X$ and $Y$ of $\mathscr{C}$, the map
\begin{equation*}
\fromto{[K,\Map_{\mathscr{C}}(X,Y)]}{[K,\Map_{\mathscr{D}}(\psi(X),\psi(Y))]}
\end{equation*}
can be identified with the map
\begin{equation*}
\fromto{\pi_0\Map(X\otimes K,Y)}{\pi_0\Map(\psi(X\otimes K),\psi(Y))\cong\pi_0\Map(\psi(X)\otimes K,\psi(Y))}.
\end{equation*}
This map is a bijection for any finite simplicial set $K$. In particular, the map $\fromto{\Map(X,Y)}{\Map(\psi(X),\psi(Y))}$ is a weak homotopy equivalence on the connected components at $0$, whence $\equivto{\Map(\Sigma X,Y)}{\Map(\psi(\Sigma X),\psi(Y))}$ is an equivalence. Now since every object in $\mathscr{C}$ is a suspension, the functor $\psi$ is fully faithful.
\end{proof}
\end{prp}


\subsection*{The dual picture} Entirely dual to the theory of Waldhausen $\infty$-categories is the theory of \emph{coWaldhausen $\infty$-categories}. We record the definition here; clearly any result or construction in the theory of Waldhausen $\infty$-categories can be immediately dualized.

\begin{dfn}\label{dfn:coWald} \begin{enumerate}[(\ref{dfn:coWald}.1)]
\item A \textbf{\emph{coWaldhausen $\infty$-category}} $(\mathscr{C},\mathscr{C}^{\dag})$ is an opposite pair $(\mathscr{C},\mathscr{C}^{\dag})$ such that the opposite $(\mathscr{C}^{\op},(\mathscr{C}^{\op})_{\dag})$ is a Waldhausen $\infty$-category.
\item A functor of pairs $\psi\colon\fromto{\mathscr{C}}{\mathscr{D}}$ between two coWaldhausen $\infty$-categories is said to be \textbf{\emph{exact}} if its opposite $\psi^{\op}\colon\fromto{\mathscr{C}^{\op}}{\mathscr{D}^{\op}}$ is exact.
\end{enumerate}
\end{dfn}

\begin{ntn}\label{ntn:precoWald}
\begin{enumerate}[(\ref{ntn:precoWald}.1)]
\item Suppose $\mathscr{C}$ and $\mathscr{D}$ two coWaldhausen $\infty$-categories. Denote by $\Fun_{\coWald}(\mathscr{C},\mathscr{D})$ the full subcategory of $\Fun_{\Pair_{\infty}}(\mathscr{C},\mathscr{D})$ spanned by the exact morphisms of coWaldhausen $\infty$-categories.
\item Define $\coWald^{\Delta}_{\infty}$ as the following large simplicial subcategory of $\Pair^{\Delta}_{\infty}$. The objects of $\coWald^{\Delta}_{\infty}$ are small coWaldhausen $\infty$-categories, and for any coWaldhausen $\infty$-categories $\mathscr{C}$ and $\mathscr{D}$, the morphism space is defined by the formula
\begin{equation*}
\coWald^{\Delta}_{\infty}(\mathscr{C},\mathscr{D})\coloneq\iota\Fun_{\coWald}(\mathscr{C},\mathscr{D}),
\end{equation*}
or equivalently, $\coWald^{\Delta}_{\infty}(\mathscr{C},\mathscr{D})$ is the union of the connected components of $\Pair^{\Delta}_{\infty}(\mathscr{C},\mathscr{D})$ corresponding to the exact morphisms.
\item We then define an $\infty$-category $\coWald$ as the simplicial nerve (Df. \ref{dfn:relnerve}) of the simplicial category $\coWald^{\Delta}_{\infty}$. 
\end{enumerate}
\end{ntn}

\noindent We summarize these constructions with the following.
\begin{prp} The opposite involution on $\Pair_{\infty}$ \textup{(Pr. \ref{prp:opinvolonPair})} restricts to an equivalence between $\Wald$ and $\coWald$.
\end{prp}


\section{Waldhausen fibrations}\label{sect:Waldfib} A key component of Waldhausen's algebraic $K$-theory of spaces is his $S_{\bullet}$ construction \cite[\S 1.3]{MR86m:18011}. In effect, this is a diagram of categories
\begin{equation*}
S\colon\fromto{\Delta^{op}}{\Cat}
\end{equation*}
such that for any object $\mathbf{m}\in\Delta$, the category $S_m$ is the category of filtered spaces
\begin{equation*}
\ast=X_0\subset X_1\subset\cdots\subset X_m
\end{equation*}
of length $m$, and, for any simplicial operator $[\phi\colon\fromto{\mathbf{n}}{\mathbf{m}}]\in\Delta$, the induced functor $\phi_!\colon\fromto{S_m}{S_n}$ carries a filtered space $\ast=X_0\subset X_1\subset\cdots\subset X_m$ to a filtered space
\begin{equation*}
\ast=X_{\phi(0)}/X_{\phi(0)}\subset X_{\phi(1)}/X_{\phi(0)}\subset\cdots\subset X_{\phi(n)}/X_{\phi(0)}.
\end{equation*}
We will want to construct an $\infty$-categorical variant of $S_{\bullet}$, but there is a little wrinkle here: as written, this is not a functor on the nose. Rather, it is a \emph{pseudofunctor}, because quotients are defined only up to (canonical) isomorphism. To rectify this, Waldhausen constructs \cite[\S 1.3]{MR86m:18011} an honest functor by replacing each category $S_m$ with a fattening thereof, in which an object is a filtered space
\begin{equation*}
\ast=X_0\subset X_1\subset\cdots\subset X_m
\end{equation*}
along with compatible choices of all the quotient spaces $X_s/X_t$.

If one wishes to pass to a more homotopical variant of the $S_{\bullet}$ construction, matters become even more complicated. After all, any sequence of simplicial sets
\begin{equation*}
\ast\simeq X_0\ {\tikz[baseline]\draw[>=stealth,->](0,0.5ex)--(0.5,0.5ex);}\ X_1\ {\tikz[baseline]\draw[>=stealth,->](0,0.5ex)--(0.5,0.5ex);}\ \cdots\ {\tikz[baseline]\draw[>=stealth,->](0,0.5ex)--(0.5,0.5ex);}\ X_m
\end{equation*}
can, up to homotopy, be regarded as a filtered space. To extend the $S_{\bullet}$ construction to accept these objects, a simplicial operator should then induce functor that carries such a sequence to a corresponding sequence of \emph{homotopy} quotients, in which each map is replaced by a cofibration, and the suitable quotients are formed. This now presents not only a functoriality problem but also a homotopy coherence problem, which is precisely solved for Waldhausen categories satisfying a technical hypothesis (functorial factorizations of weak $w$-cofibrations) by means of Blumberg--Mandell's $S'_{\bullet}$-construction \cite[Df. 2.7]{BM}.

Unfortunately, these homotopy coherence problems grow less tractable as $K$-theoretic constructions become more involved. For example, if one seeks multiplicative structures on algebraic $K$-theory spectra, it becomes a challenge to perform all the necessary rectifications to turn a suitable pairing of Waldhausen categories into an $E_k$ multipciation on the $K$-theory. The work of Elmendorf and Mandell \cite{MR2254311} manages the case $k=\infty$ by using different (and quite rigid) inputs for the $K$-theory functor. More generally, Blumberg and Mandell \cite[Th. 2.6]{MR2805994} generalize this by providing, for any (colored) operad $O$ in categories, an $O$-algebra structure on the $K$-theory of any $O$-algebra in  Waldhausen categories.

However, the theory $\infty$-categories provides a powerful alternative to such explicit solutions to homotopy coherence problems. Namely, the theory of \emph{cartesian} and \emph{cocartesian} fibrations allows one, in effect, to leave the homotopy coherence problems \emph{unsolved} yet, at the same time, to work effectively with the resulting objects. For this reason, these concepts play a central role in our work here. (For fully general solutions to the problem of finding $O$ structures on $K$-theory spectra using machinery of the kind developed here, see either Blumberg--Gepner--Tabuada \cite{BGT2} or \cite{K3}.)


\subsection*{Cocartesian fibrations} The idea goes back at least to Grothendieck (and probably further). If $X\colon\fromto{C}{\Cat}$ is an (honest) diagram of ordinary categories, then one can define the \emph{Grothendieck construction} of $X$. This is a category $G(X)$ whose objects are pairs $(c,x)$ consisting of an object $c\in C$ and an object $x\in X(c)$, in which a morphism $(f,\phi)\colon\fromto{(d,y)}{(c,x)}$ is a morphism $f\colon\fromto{d}{c}$ of $C$ and a morphism
\begin{equation*}
\phi\colon\fromto{X(f)(y)}{x}
\end{equation*}
of $X(c)$. There is an obvious forgetful functor $p\colon\fromto{G(X)}{C}$.

One may now attempt to reverse-engineer the Grothendieck construction by trying to extract the salient features of the forgetful functor $p$ that ensures that it ``came from'' a diagram of categories. What we may notice is that for any morphism $f\colon\fromto{d}{c}$ of $C$ and any object $y\in X(d)$ there is a special morphism
\begin{equation*}
F=(f,\phi)\colon\fromto{(d,y)}{(c,X(f)(y))}
\end{equation*}
of $G(X)$ in which
\begin{equation*}
\phi\colon\fromto{X(f)(y)}{X(f)(y)}
\end{equation*}
is simply the identity morphism. This morphism is \emph{initial} among all the morphisms $F'$ of $G(X)$ such that $p(F')=f$; that is, for any morphism $F'$ of $G(X)$ such that $p(F')=f$, there exists a morphism $H$ of $G(X)$ such that $p(H)=\id_c$ such that $F'=H\circ F$.

We call morphisms of $G(X)$ that are initial in this sense \emph{$p$-cocartesian}. Since a $p$-cocartesian edge lying over a morphism $\fromto{d}{c}$ is defined by a universal property, it is uniquely specified up to a unique isomorphism lying over $\id_{c}$. The key condition that we are looking for is then that \emph{for any morphism of $C$ and any lift of its source, there is a $p$-cocartesian morphism with that source lying over it}. A functor $p$ satisfying this condition is called a \emph{Grothendieck opfibration}.

Now for \emph{any} Grothendieck opfibration $q\colon\fromto{D}{C}$, let us attempt to extract a functor $Y\colon\fromto{C}{\Cat}$ whose Grothendieck construction $G(Y)$ is equivalent (as a category over $C$) to $D$. We proceed in the following manner. To any object $c\in C$ assign the fiber $D_c$ of $q$ over $c$. To any morphism $f\colon\fromto{d}{c}$ assign a functor $Y(f)\colon\fromto{D_d}{D_c}$ that carries any object $y\in D_d$ to the target $Y(f)(y)\in D_c$ of ``the'' $q$-cocartesian edge lying over $f$. However, the problem is already apparent in the scare quotes around the word ``the.'' These functors will not be strictly compatible with composition; rather, one will obtain natural isomorphisms
\begin{equation*}
Y(g\circ f)\simeq Y(g)\circ Y(f)
\end{equation*}
that will satisfy a secondary layer of coherences that make $Y$ into a \emph{pseudofunctor}.

It is in fact possible to rectify any pseudofunctor to an equivalent honest functor, and this gives an honest functor whose Grothendieck construction is equivalent to our original $D$.

In light of all this, three options present themselves for contending with weak diagrams of ordinary categories:
\begin{enumerate}[(1)]
\item Rectify all pseudofunctors, and keep track of the rectifications as constructions become more involved.
\item Work systematically with pseudofunctors, verifying all the coherence laws as needed.
\item Work directly with Grothendieck opfibrations.
\end{enumerate}
Which of these one selects is largely a matter of taste. When we pass to diagrams of higher categories, however, the first two options veer sharply into the realm of impracticality. A pseudofunctor $\fromto{S}{\Cat_{\infty}}$ has not only a secondary level of coherences, but also an infinite progression of coherences between witnesses of lower-order coherences. Though rectifications of these pseudofunctors do exist (see \ref{rec:straighten} below), they are usually not terribly explicit, and it would be an onerous task to keep them all straight.

Fortunately, the last option generalizes quite comfortably to the context of quasicategories, yielding the theory of \emph{cocartesian fibrations}.
\begin{rec}\label{rec:cocart} Suppose $p\colon\fromto{X}{S}$ an inner fibration of simplicial sets. Recall \cite[Rk. 2.4.1.4]{HTT} that an edge $f\colon\fromto{\Delta^{1}}{X}$ is \textbf{\emph{$p$-cocartesian}} just in case, for each integer $n\geq 2$, any extension
\begin{equation*}
\begin{tikzpicture} 
\matrix(m)[matrix of math nodes, 
row sep=4ex, column sep=4ex, 
text height=1.5ex, text depth=0.25ex] 
{\Delta^{\{0,1\}}&X,\\ 
\Lambda^n_0&\\}; 
\path[>=stealth,->,font=\scriptsize] 
(m-1-1) edge node[above]{$f$} (m-1-2) 
edge[right hook->] (m-2-1) 
(m-2-1) edge node[below]{$F$} (m-1-2); 
\end{tikzpicture}
\end{equation*}
and any solid arrow commutative diagram
\begin{equation*}
\begin{tikzpicture} 
\matrix(m)[matrix of math nodes, 
row sep=4ex, column sep=4ex, 
text height=1.5ex, text depth=0.25ex] 
{\Lambda^n_0&X\\ 
\Delta^n&S,\\}; 
\path[>=stealth,->,font=\scriptsize] 
(m-1-1) edge node[above]{$F$} (m-1-2) 
edge[right hook->] (m-2-1) 
(m-1-2) edge node[right]{$p$} (m-2-2) 
(m-2-1) edge (m-2-2)
(m-2-1) edge[dotted] node[below]{$\overline{F}$} (m-1-2);
\end{tikzpicture}
\end{equation*}
the dotted arrow $\overline{F}$ exists, rendering the diagram commutative.

We say that $p$ is a \textbf{\emph{cocartesian fibration}} \cite[Df. 2.4.2.1]{HTT} if, for any edge $\eta\colon\fromto{s}{t}$ of $S$ and for every vertex $x\in X_0$ such that $p(x)=s$, there exists a $p$-cocartesian edge $f\colon\fromto{x}{y}$ such that $\eta=p(f)$.

\textbf{\emph{Cartesian edges}} and \textbf{\emph{cartesian fibrations}} are defined dually, so that an edge of $X$ is $p$-cartesian just in case the corresponding edge of $X^{\op}$ is cocartesian for the inner fibration $p^{\op}\colon\fromto{X^{\op}}{S^{\op}}$, and $p$ is a cartesian fibration just in case $p^{\op}$ is a cocartesian fibration.
\end{rec}

\begin{exm}\label{exm:nerveofGrothfib} A functor $p\colon\fromto{D}{C}$ between ordinary categories is a Grothendieck opfibration if and only if the induced functor $N(p)\colon\fromto{ND}{NC}$ on nerves is a cocartesian fibration \cite[Rk 2.4.2.2]{HTT}.
\end{exm}

\begin{exm}\label{exm:scocartfib} Recall that for any $\infty$-category $C$, we write $\mathscr{O}(C)\coloneq\Fun(\Delta^1,C)$. By \cite[Cor. 2.4.7.12]{HTT}, evaluation at $0$ defines a cartesian fibration $s\colon\fromto{\mathscr{O}(C)}{C}$, and evaluation at $1$ defines a cocartesian fibration $t\colon\fromto{\mathscr{O}(C)}{C}$.

One can ask whether the functor $s\colon\fromto{\mathscr{O}(C)}{C}$ is also a \emph{cocartesian} fibration. One may observe \cite[Lm. 6.1.1.1]{HTT} that an edge $\fromto{\Delta^1}{\mathscr{O}(C)}$ is $s$-cocartesian just in case the correponding diagram $\fromto{(\Lambda^2_0)^{\rhd}\cong\Delta^1\times\Delta^1}{C}$ is a pushout square.
\end{exm}

\begin{rec}\label{rec:straighten} Suppose $S$ a simplicial set. Then the collection of cocartesian fibrations to $S$ with small fibers is naturally organized into an $\infty$-category $\Cat_{\infty/S}^{\cocart}$. To construct it, let $\Cat_{\infty}^{\cocart}$ be the following subcategory of $\mathscr{O}(\Cat_{\infty})$: an object $\fromto{X}{U}$ of $\mathscr{O}(\Cat_{\infty})$ lies in $\Cat_{\infty}^{\cocart}$ if and only if it is a cocartesian fibration, and a morphism $\fromto{p}{q}$ in $\mathscr{O}(\Cat_{\infty})$ between cocartesian fibrations represented as a square
\begin{equation*}
\begin{tikzpicture} 
\matrix(m)[matrix of math nodes, 
row sep=4ex, column sep=4ex, 
text height=1.5ex, text depth=0.25ex] 
{X&Y\\ 
U&V\\}; 
\path[>=stealth,->,font=\scriptsize] 
(m-1-1) edge node[above]{$f$} (m-1-2) 
edge node[left]{$p$} (m-2-1) 
(m-1-2) edge node[right]{$q$} (m-2-2) 
(m-2-1) edge (m-2-2); 
\end{tikzpicture}
\end{equation*}
lies in in $\Cat_{\infty}^{\cocart}$ if and only if $f$ carries $p$-cocartesian edges to $q$-cocartesian edges. We now define $\Cat_{\infty/S}^{\cocart}$ as the fiber over $S$ of the target functor
\begin{equation*}
t\colon\Cat_{\infty}^{\cocart}\subset\mathscr{O}(\Cat_{\infty})\ {\tikz[baseline]\draw[>=stealth,->](0,0.5ex)--(0.5,0.5ex);}\ \Cat_{\infty}.
\end{equation*}
Equivalently \cite[Pr. 3.1.3.7]{HTT}, one may describe $\Cat_{\infty/S}^{\cocart}$ as the simplicial nerve (Nt. \ref{ntn:superscriptsscats}) of the (fibrant) simplicial category of marked simplicial sets \cite[Df. 3.1.0.1]{HTT} over $S$ that are fibrant for the \emph{cocartesian model structure} --- i.e., of the form $\fromto{X^{\natural}}{S}$ for $\fromto{X}{S}$ a cocartesian fibration \cite[Df. 3.1.1.8]{HTT}.

The straightening/unstraightening Quillen equivalence of \cite[Th. 3.2.0.1]{HTT} now yields an equivalence of $\infty$-categories
\begin{equation*}
\Cat_{\infty/S}^{\cocart}\simeq\Fun(S,\Cat_{\infty}).
\end{equation*}
So the dictionary between Grothendieck opfibrations and diagrams of categories generalizes gracefully to a dictionary between cocartesian fibrations $p\colon\fromto{X}{S}$ with small fibers and functors $\XX\colon\fromto{S}{\Cat_{\infty}}$. As for ordinary categories, for any vertex $s\in S_0$, the value $\XX(s)$ is equivalent to the fiber $X_s$, and for any edge $\eta\colon\fromto{s}{t}$, the functor $\fromto{h\XX(s)}{h\XX(t)}$ assigns to any object $x\in X_s$ an object $y\in X_t$ with the property that there is a cocartesian edge $\fromto{x}{y}$ that covers $\eta$. We say that $\XX$ \textbf{\emph{classifies}} $p$ \cite[Df. 3.3.2.2]{HTT}, and we will abuse terminology slightly by speaking of \textbf{\emph{the functor $\eta_{!}\colon\fromto{X_s}{X_t}$ induced by}} an edge $\eta\colon\fromto{s}{t}$ of $S$, even though $\eta_!$ is defined only up to canonical equivalence.

Dually, the collection of cartesian fibrations to $S$ with small fibers is naturally organized into an $\infty$-category $\Cat_{\infty/S}^{\cart}$, and the straightening/unstraightening Quillen equivalence yields an equivalence of $\infty$-categories
\begin{equation*}
\Cat_{\infty/S}^{\cart}\simeq\Fun(S^{\op},\Cat_{\infty}).
\end{equation*}
\end{rec}

\begin{exm} For any $\infty$-category $C$, the functor $\fromto{C^{\op}}{\Cat_{\infty}}$ that classifies the cartesian fibration $s\colon\fromto{\mathscr{O}(C)}{C}$ is the functor that carries any object $X$ of $C$ to the undercategory $C_{X/}$ and any morphism $f\colon\fromto{Y}{X}$ to the forgetful functor $f^{\star}\colon\fromto{C_{X/}}{C_{Y/}}$.

If $C$ admits all pushouts, then the cocartesian fibration $s\colon\fromto{\mathscr{O}(C)}{C}$ is classified by a functor $\fromto{C}{\Cat_{\infty}}$ that carries any object $X$ of $C$ to the undercategory $C_{X/}$ and any morphism $f\colon\fromto{Y}{X}$ to the functor $f_{!}\colon\fromto{C_{Y/}}{C_{X/}}$ that is given by pushout along $f$.
\end{exm}

\begin{rec}\label{rec:leftfib} A cocartesian fibration with the special property that each fiber is a Kan complex --- or equivalently, with the special property that the functor that classifies it factors through the full subcategory $\Kan\subset\Cat_{\infty}$ --- is called a \textbf{\emph{left fibration}}. These are more efficiently described as maps that satisfy the right lifting property with respect to horn inclusions $\into{\Lambda^n_k}{\Delta^n}$ such that $n\geq1$ and $0\leq k\leq n-1$ \cite[Pr. 2.4.2.4]{HTT}.

For any cocartesian fibration $p\colon\fromto{X}{S}$, one may consider the smallest simplicial subset $\iota_SX\subset X$ that contains the $p$-cocartesian edges. The restriction $\iota_S(p)\colon\fromto{\iota_SX}{S}$ of $p$ to $\iota_SX$ is a left fibration. The functor $\fromto{S}{\Kan}$ that classifies $\iota_Sp$ is then the functor given by the composition
\begin{equation*}
S\ \tikz[baseline]\draw[>=stealth,->,font=\scriptsize](0,0.5ex)--node[above]{$F$}(0.5,0.5ex);\ \Cat_{\infty}\ \tikz[baseline]\draw[>=stealth,->,font=\scriptsize](0,0.5ex)--node[above]{$\iota$}(0.5,0.5ex);\ \Kan,
\end{equation*}
where $F$ is the functor that classifies $p$.
\end{rec}

Let us recall a particularly powerful construction with cartesian and cocartesian fibrations, which will form the cornerstone for our study of filtered objects of Waldhausen $\infty$-categories.
\begin{rec}\label{rec:htt32213} Suppose $S$ a simplicial set, and suppose $\XX\colon\fromto{S^{\op}}{\Cat_{\infty}}$ and $\YY\colon\fromto{S}{\Cat_{\infty}}$ two diagrams of $\infty$-categories. Then one may define a functor
\begin{equation*}
\Fun(\XX,\YY)\colon\fromto{S}{\Cat_{\infty}}
\end{equation*}
that carries a vertex $s$ of $S$ to the $\infty$-category $\Fun(\XX(s),\YY(s))$ and an edge $\eta\colon\fromto{s}{t}$ of $S$ to the functor
\begin{equation*}
\fromto{\Fun(\XX(s),\YY(s))}{\Fun(\XX(t),\YY(t))}
\end{equation*}
given by the assignment $\goesto{F}{\YY(\eta)\circ F\circ\XX(\eta)}$.

If one wishes to work instead with the cartesian and cocartesian fibrations classified by $\XX$ and $\YY$, the following construction provides an elegant way of writing explicitly the cocartesian fibration classified by the functor $\Fun(\XX,\YY)$. If $p\colon\fromto{X}{S}$ is the cartesian fibration classified by $\XX$ and if $q\colon\fromto{Y}{S}$ is the cocartesian fibration classified by $\YY$, one may define a map $r\colon\fromto{T}{S}$ defined by the following universal property: for any map $\sigma\colon\fromto{K}{S}$, one has a bijection
\begin{equation*}
\Mor_{S}(K,T)\cong\Mor_{S}(X\times_SK,Y),
\end{equation*}
functorial in $\sigma$. It is then shown in \cite[Cor. 3.2.2.13]{HTT} that $p$ is a cocartesian fibration, and an edge $g\colon\fromto{\Delta^1}{T}$ is $r$-cocartesian just in case the induced map $\fromto{X\times_S\Delta^1}{Y}$ carries $p$-cartesian edges to $q$-cocartesian edges. The fiber of the map $\fromto{T}{S}$ over a vertex $s$ is the $\infty$-category $\Fun(X_s,Y_s)$, and for any edge $\eta\colon\fromto{s}{t}$ of $S$, the functor $\eta_!\colon\fromto{T_s}{T_t}$ induced by $\eta$ is equivalent to the functor $\goesto{F}{\YY(\eta)\circ F\circ\XX(\eta)}$ described above.
\end{rec}


\subsection*{Pair cartesian and cocartesian fibrations} Just as cartesian and cocartesian fibrations are well adapted to the study of weak diagrams of $\infty$-categories, so we will introduce the theory of \emph{Waldhausen cartesian} and \emph{cocartesian fibrations}, which make available a robust notion of weak diagrams of Waldhausen $\infty$-categories. In order to introduce this notion, we first discuss \emph{pair cartesian} and \emph{cocartesian fibations} in some detail. These will provide a notion of weak diagrams of pairs of $\infty$-categories.

\begin{dfn}\label{dfn:paircartfib} Suppose $S$ an $\infty$-category. Then a \textbf{\emph{pair cartesian fibration}} $\fromto{\mathscr{X}}{S}$ is a pair $\mathscr{X}$ and a morphism of pairs $p\colon\fromto{\mathscr{X}}{S^{\flat}}$ (where the target is the minimal pair $(S,\iota S)$ --- see Ex. \ref{exm:minpairmaxpair}) such that the following conditions are satisfied.
\begin{enumerate}[(\ref{dfn:paircartfib}.1)]
\item The underlying functor of $p$ is a cartesian fibration.
\item\label{item:etastarisfunctorofpairs} For any edge $\eta\colon\fromto{s}{t}$ of $S$, the induced functor $\eta^{\star}\colon\fromto{\mathscr{X}_t}{\mathscr{X}_s}$ carries ingressive morphisms to ingressive morphisms. 
\end{enumerate}

Dually, a \textbf{\emph{pair cocartesian fibration}} $\fromto{\mathscr{X}}{S}$ is a pair $\mathscr{X}$ and a morphism of pairs $p\colon\fromto{\mathscr{X}}{S^{\flat}}$ such that $p^{\op}\colon\fromto{\mathscr{X}^{\op}}{S^{\op}}$ is a pair cartesian fibration.
\end{dfn}

\begin{prp}\label{cor:pairfibsarefunsintopair} If $S$ is an $\infty$-category and $p\colon\fromto{\mathscr{X}}{S}$ is a pair cartesian fibration \textup{[}respectively, a pair cocartesian fibration\textup{]} with small fibers, then the functor $\fromto{S^{\op}}{\Cat_{\infty}}$ \textup{[}resp., the functor $\fromto{S}{\Cat_{\infty}}$\textup{]} that classifies $p$ lifts to a functor $\fromto{S^{\op}}{\Pair_{\infty}}$ \textup{[}resp., $\fromto{S}{\Pair_{\infty}}$\textup{]}.
\begin{proof} We employ the adjunction $(\mathfrak{C},N)$ of \cite[\S 1.1.5]{HTT}. Since $\Pair_{\infty}$ and $\Cat_{\infty}$ are both defined as simplicial nerves, the data of a lift $\fromto{S^{\op}}{\Pair_{\infty}}$ of $\fromto{S^{\op}}{\Cat_{\infty}}$ is tantamount to the data of a lift $\overline{\XX}\colon\fromto{\mathfrak{C}[S]^{\op}}{\Pair^{\Delta}_{\infty}}$ of the corresponding simplicial functor $\XX\colon\fromto{\mathfrak{C}[S]^{\op}}{\Cat^{\Delta}_{\infty}}$. Now for any object $s$ of $\mathfrak{C}[S]$, the categories $\XX(s)$ inherits a pair structure via the canonical equivalence $\XX(s)\simeq\mathscr{X}_s$. For any two objects $s$ and $t$ of $\mathfrak{C}[S]$, condition (\ref{dfn:paircartfib}.\ref{item:etastarisfunctorofpairs}) ensures that the map
\begin{equation*}
\fromto{\mathfrak{C}[S](t,s)}{\Cat_{\infty}^{\Delta}(\XX(s),\XX(t))}
\end{equation*}
factors through the simplicial subset (Nt. \ref{ntn:pair})
\begin{equation*}
\Pair_{\infty}^{\Delta}(\XX(s),\XX(t))\subset\Cat_{\infty}^{\Delta}(\XX(s),\XX(t)).
\end{equation*}
This now defines the desired simplicial functor $\overline{\XX}$.
\end{proof}
\end{prp}

\begin{dfn} In the situation of Pr. \ref{cor:pairfibsarefunsintopair}, we will say that the lifted functor $\fromto{S^{\op}}{\Pair_{\infty}}$ \textup{[}respectively, the lifted functor $\fromto{S}{\Pair_{\infty}}$\textup{]} \textbf{\emph{classifies}} the cartesian \textup{[}resp., cocartesian\textup{]} fibration $p$.
\end{dfn}

\begin{prp}\label{cor:paircartfibsbasechange} The classes of pair cartesian fibrations and pair cocartesian fibrations are each stable under base change. That is, for any pair cartesian [respectively, cocartesian] fibration $\fromto{\mathscr{X}}{S}$ and for any functor $f\colon\fromto{S'}{S}$, if the pullback $\mathscr{X}'\coloneq\mathscr{X}\times_SS'$ is endowed with the pair structure in which a morphism is ingressive just in case it is carried to an equivalence in $S'$ and to an ingressive morphism of $\mathscr{X}$, then $\fromto{\mathscr{X}'}{S'}$ is a pair cartesian [resp., cocartesian] fibration.
\begin{proof} We treat the case of pair cartesian fibrations. Cartesian fibrations are stable under pullbacks \cite[Pr. 2.4.2.3(2)]{HTT}, so it remains to note that for any morphism $\eta\colon\fromto{s}{t}$ of $S'$, the induced functor
\begin{equation*}
\eta^{\star}\simeq f(\eta)^{\star}\colon\fromto{\mathscr{X}'_t\cong\mathscr{X}_{f(t)}}{\mathscr{X}_{f(s)}\cong\mathscr{X}'_s}
\end{equation*}
carries ingressive morphisms to ingressive morphisms.
\end{proof}
\end{prp}


\subsection*{The $\infty$-categories of pair (co)cartesian fibrations} The collection of all pair cocartesian fibrations are organized into an $\infty$-category $\Pair_{\infty}^{\cocart}$, which is analogous to the $\infty$-category $\Cat_{\infty}^{\cocart}$ of \ref{rec:straighten}. Furthermore, pair cocartesian fibrations with a fixed base $\infty$-category $S$ organize themselves into an $\infty$-category $\Pair_{\infty/S}^{\cocart}$.

\begin{ntn} Denote by
\begin{equation*}
\Pair_{\infty}^{\cart}\textrm{\qquad[respectively, by\quad}\Pair_{\infty}^{\cocart}\textrm{\quad]}
\end{equation*}
the following subcategory of $\mathscr{O}(\Pair_{\infty})$. The objects of $\Pair_{\infty}^{\cart}$ [resp., $\Pair_{\infty}^{\cocart}$] are pair cartesian fibrations (resp., pair cocartesian fibrations) $\fromto{\mathscr{X}}{S}$. For any pair cartesian (resp., cocartesian) fibrations $p\colon\fromto{\mathscr{X}}{S}$ and $q\colon\fromto{\mathscr{Y}}{T}$,  a commutative square
\begin{equation*}
\begin{tikzpicture} 
\matrix(m)[matrix of math nodes, 
row sep=4ex, column sep=4ex, 
text height=1.5ex, text depth=0.25ex] 
{\mathscr{X}&\mathscr{Y}\\ 
S^{\flat}&T^{\flat}\\}; 
\path[>=stealth,->,font=\scriptsize] 
(m-1-1) edge node[above]{$\psi$} (m-1-2) 
edge node[left]{$p$} (m-2-1) 
(m-1-2) edge node[right]{$q$} (m-2-2) 
(m-2-1) edge (m-2-2); 
\end{tikzpicture}
\end{equation*}
of pairs of $\infty$-categories is a morphism $\fromto{p}{q}$ of $\Pair_{\infty}^{\cart}$ [resp., of $\Pair_{\infty}^{\cocart}$] if and only if $\psi$ carries $p$-cartesian (resp. $p$-cocartesian) edges to $q$-cartesian (resp. $q$-cocartesian) edges.

By an abuse of notation, we will denote by $(\mathscr{X}/S)$ an object $\fromto{\mathscr{X}}{S}$ of $\Pair_{\infty}^{\cart}$ [resp., of $\Pair_{\infty}^{\cocart}$].
\end{ntn}

The following is immediate from Pr. \ref{cor:paircartfibsbasechange} and \cite[Lm. 6.1.1.1]{HTT}.

\begin{lem} The target functors
\begin{equation*}
\fromto{\Pair_{\infty}^{\cart}}{\Cat_{\infty}}\textrm{\quad and\quad}\fromto{\Pair_{\infty}^{\cocart}}{\Cat_{\infty}}
\end{equation*}
induced by the inclusion $\{1\}\subset\Delta^1$ are both cartesian fibrations.
\end{lem}

\begin{ntn} The fibers of the cartesian fibrations
\begin{equation*}
\fromto{\Pair_{\infty}^{\cart}}{\Cat_{\infty}}\textrm{\quad and\quad}\fromto{\Pair_{\infty}^{\cocart}}{\Cat_{\infty}}
\end{equation*}
over an object $\{S\}\subset\Cat_{\infty}$ will be denoted $\Pair_{\infty,/S}^{\cart}$ and $\Pair_{\infty,/S}^{\cocart}$, respectively.

By an abuse of notation, denote by
\begin{equation*}
(\Pair_{\infty/S}^{\cart})_0\textrm{\qquad[respectively, by\quad}(\Pair_{\infty/S}^{\cocart})_0\textrm{\quad]}
\end{equation*}
the subcategory of the ordinary category $((\Pair_{\infty}^{\Delta})_{0}\downarrow S^{\flat})$ whose objects are pair cartesian fibrations [resp., pair cocartesian fibrations] $\fromto{\mathscr{X}}{S}$ and whose morphisms are functors of pairs $\fromto{\mathscr{X}}{\mathscr{Y}}$ over $S$ that carry cartesian morphisms to cartesian morphisms [resp., that carry cocartesian morphisms to cocartesian morphisms]. Denote by
\begin{equation*}
w(\Pair_{\infty/S}^{\cart})_0\subset(\Pair_{\infty/S}^{\cart})_0\textrm{\qquad [resp., by\quad}w(\Pair_{\infty/S}^{\cocart})_0\subset(\Pair_{\infty/S}^{\cocart})_0\textrm{\quad]}
\end{equation*}
the subcategory consisting of those morphisms $\fromto{\mathscr{X}}{\mathscr{Y}}$ over $S$ that are \emph{fiberwise equivalences of pairs} --- i.e., such that for any vertex $s\in S_0$, the induced functor $\fromto{\mathscr{X}_s}{\mathscr{Y}_s}$ is a weak equivalence of pairs. Equivalently, $w(\Pair_{\infty/S}^{\cart})_0$ is the collection of those equivalences of pairs $\equivto{\mathscr{X}}{\mathscr{Y}}$ over $S$ that are fiberwise equivalences of $\infty$-categories --- i.e., such that for any vertex $s\in S_0$, the induced functor $\fromto{\mathscr{X}_s}{\mathscr{Y}_s}$ is an equivalence of underyling $\infty$-categories.
\end{ntn}

\begin{lem}\label{lem:Paircocartisarelnerve} For any $\infty$-category $S$, the $\infty$-category $\Pair_{\infty/S}^{\cart}$ \textup{[}respectively, the $\infty$-category $\Pair_{\infty/S}^{\cocart}$\textup{]} is a relative nerve \textup{(Df. \ref{dfn:relnerve})} of
\begin{equation*}
((\Pair_{\infty/S}^{\cart})_0,\ w(\Pair_{\infty/S}^{\cart})_0)\textrm{\qquad\textup{[}resp., of\quad}((\Pair_{\infty/S}^{\cocart})_0,\ w(\Pair_{\infty/S}^{\cocart})_0)\textrm{\quad\textup{]}.}
\end{equation*}
\begin{proof} To show that $\Pair_{\infty/S}^{\cart}$ is a relative nerve of $((\Pair_{\infty/S}^{\cart})_0,w(\Pair_{\infty/S}^{\cart})_0)$, we first note that the analogous result for $\infty$-categories of cartesian fibrations $\fromto{X}{S}$ holds. More precisely, recall (\ref{rec:straighten}) that $\Cat_{\infty/S}^{\cart}$ may be identified with the nerve of the cartesian simplicial model category of marked simplicial sets over $S$, whence it is a relative nerve of the category $(\Cat_{\infty/S}^{\cart})_0$ of cartesian fibrations over $S$, equipped with the subcategory $w(\Cat_{\infty/S}^{\cart})_0$ consisting of fiberwise equivalences.

To extend this result to a characterization of $\Pair_{\infty/S}^{\cart}$ as a relative nerve, let us contemplate the square
\begin{equation*}
\begin{tikzpicture}[baseline]
\matrix(m)[matrix of math nodes,
row sep=4ex, column sep=4ex,
text height=1.5ex, text depth=0.25ex]
{N((\Pair_{\infty/S}^{\cart})_0,W) & N((\Cat_{\infty/S}^{\cart})_0\times_{(\Cat_{\infty})_0}(\Pair_{\infty})_0,W) \\
\Pair_{\infty/S}^{\cart} & \Cat_{\infty/S}^{\cart}\times_{\Cat_{\infty}}\Pair_{\infty}, \\ };
\path[>=stealth,->,font=\scriptsize]
(m-1-1) edge node[above]{} (m-1-2)
edge node[left]{} (m-2-1)
(m-1-2) edge node[right]{} (m-2-2)
(m-2-1) edge node[below]{} (m-2-2);
\end{tikzpicture}
\end{equation*}
where we have written $W$ for the obvious classes of weak equivalences. The horizontal maps are the forgetful functors, and the vertical maps are the ones determined by the universal property of the relative nerve. The vertical functor on the right is an equivalence, and the vertical functor on the left is essentially surjective. It therefore remains only to note that the horizontal functors are fully faithful.
\end{proof}
\end{lem}

We may now employ this lemma to lift the equivalence of $\infty$-categories
\begin{equation*}
\Cat_{\infty/S}^{\cart}\simeq\Fun(S^{\op},\Cat_{\infty})
\end{equation*}
of \cite[\S 3.2]{HTT} to an equivalence of $\infty$-categories
\begin{equation*}
\Pair_{\infty/S}^{\cart}\simeq\Fun(S^{\op},\Pair_{\infty}).
\end{equation*}
\begin{prp}\label{prp:WaldcocartisFunSWald} For any $\infty$-category $S$, the $\infty$-category $\Fun(S^{\op},\Pair_{\infty})$ \textup{[}respectively, the $\infty$-category $\Fun(S,\Pair_{\infty})$\textup{]} is a relative nerve \textup{(Df. \ref{dfn:relnerve})} of
\begin{equation*}
((\Pair_{\infty/S}^{\cart})_0,\ w(\Pair_{\infty/S}^{\cart})_0)\textrm{\qquad\textup{[}resp., of\quad}((\Pair_{\infty/S}^{\cocart})_0,\ w(\Pair_{\infty/S}^{\cocart})_0)\textrm{\quad\textup{]}.}
\end{equation*}
\begin{proof} The unstraightening functor of \cite[\S 3.2]{HTT} is a weak equivalence-preserving functor
\begin{equation*}
\mathrm{Un}^+\colon\fromto{(\Cat_{\infty}^{\Delta})^{\mathfrak{C}[S]^{\op}}}{(\Cat_{\infty/S}^{\cart})_0}
\end{equation*}
that induces an equivalence of relative nerves. (Here, $(\Cat_{\infty}^{\Delta})^{\mathfrak{C}[S]^{\op}}$ denotes the relative category of simplicial functors $\fromto{\mathfrak{C}[S]^{\op}}{\Cat_{\infty}^{\Delta}}$.) For any simplicial functor
\begin{equation*}
\XX\colon\fromto{\mathfrak{C}[S]^{\op}}{\Pair_{\infty}^{\Delta}},
\end{equation*}
endow the unstraightening $\mathrm{Un}^+(\XX)$ with a pair structure by letting $\mathrm{Un}^+(\XX)_{\dag}\subset\mathrm{Un}^+(\XX)$ be the smallest subcategory containing all the equivalences as well as any cofibration of any fiber $\mathrm{Un}^+(\XX)_s\cong\XX(s)$. With this definition, we obtain a weak equivalence-preserving functor
\begin{equation*}
\mathrm{Un}^+\colon\fromto{(\Pair_{\infty}^{\Delta})^{\mathfrak{C}[S]^{\op}}}{(\Pair_{\infty/S}^{\cart})_0}.
\end{equation*}

This functor induces a functor on relative nerves, which is essentially surjective by Pr. \ref{cor:pairfibsarefunsintopair}. Moreover, for any simplicial functors
\begin{equation*}
\XX,\YY\colon\fromto{\mathfrak{C}[S]^{\op}}{\Pair_{\infty}^{\Delta}},
\end{equation*}
the simplicial set
\begin{equation*}
\Map_{N((\Pair_{\infty}^{\Delta})^{\mathfrak{C}[S]^{\op}})}(\XX,\YY)
\end{equation*}
may be identified with the simplicial subset of
\begin{equation*}
\Map_{N((\Cat_{\infty}^{\Delta})^{\mathfrak{C}[S]^{\op}})}(\XX,\YY)
\end{equation*}
given by the union of the connected components corresponding to natural transformations $\fromto{\XX}{\YY}$ such that for any $s\in S_0$, the functor $\fromto{\XX(s)}{\YY(s)}$ is a functor of pairs. Similarly, the simplicial set
\begin{equation*}
\Map_{\Pair_{\infty/S}^{\cart}}(\mathrm{Un}^+(\XX),\mathrm{Un}^+(\YY))
\end{equation*}
may be identified with the subspace of
\begin{equation*}
\Map_{\Cat_{\infty/S}^{\cart}}(\mathrm{Un}^+(\XX),\mathrm{Un}^+(\YY))
\end{equation*}
given by the union of the connected components corresponding to functors
\begin{equation*}
\fromto{\mathrm{Un}^+(\XX)}{\mathrm{Un}^+(\YY)}
\end{equation*}
over $S$ that send cartesian edges to cartesian edges with the additional property that for any $s\in S_0$, the functor
\begin{equation*}
\fromto{\mathrm{Un}^+(\XX)_s\cong\XX(s)}{\YY(s)\cong\mathrm{Un}^+(\YY)_s}
\end{equation*}
is a functor of pairs. We thus conclude that $\mathrm{Un}^+$ is fully faithful. 
\end{proof}
\end{prp}
\noindent Armed with this, we may characterize colimits of pair cartesian fibrations fiberwise.
\begin{cor}\label{cor:colimitsofpairfibs} Suppose $S$ a small $\infty$-category, $K$ a small simplicial set. A functor $\mathscr{X}\colon\fromto{K^{\rhd}}{\Pair_{\infty/S}^{\cart}}$ \textup{[}respectively, a functor $\mathscr{X}\colon\fromto{K^{\rhd}}{\Pair_{\infty/S}^{\cocart}}$\textup{]} is a colimit diagram if and only if, for every vertex $s\in S_0$, the induced functor
\begin{equation*}
\mathscr{X}_s\colon\fromto{K^{\rhd}}{\Pair_{\infty}}
\end{equation*}
is a colimit diagram.
\end{cor}
\noindent Of course the same characterization of limits holds, but it will not be needed. We will take up the question of the \emph{existence} of colimits in the $\infty$-category $\Pair_{\infty}$ in Cor. \ref{cor:colimsofpairs} below.


\subsection*{A pair version of \protect{\ref{rec:htt32213}}} The theory of pair cartesian and cocartesian fibrations is a relatively mild generalization of the theory of cartesian and cocartesian fibrations, and many of the results extend to this setting. In particular, we now set about proving a pair version of \ref{rec:htt32213} (i.e., of \cite[Cor. 3.2.2.13]{HTT}).

In effect, the objective is to give a fibration-theoretic version of the following observation. For any $\infty$-category $S$, any diagram $\XX\colon\fromto{S^{\op}}{\Pair_{\infty}}$, and any diagram $\YY\colon\fromto{S}{\Pair_{\infty}}$, there is a functor
\begin{equation*}
\Fun_{\Pair_{\infty}}(\XX,\YY)\colon\fromto{S}{\Cat_{\infty}}
\end{equation*}
that carries any object $s$ of $S$ to the $\infty$-category $\Fun_{\Pair_{\infty}}(\XX(s),\YY(s))$. 

\begin{ntn}\label{ntn:ordcatofpairs} Consider the ordinary category $s\Set(2)$ of pairs $(V,U)$ consisting of a small simplicial set $U$ and a simplicial subset $U\subset V$.
\end{ntn}

\begin{prp}\label{prp:htt32213} Suppose $p\colon\fromto{\mathscr{X}}{S}$ a pair cartesian fibration, and suppose $q\colon\fromto{\mathscr{Y}}{S}$ a pair cocartesian fibration. Let $r\colon\fromto{T_p\mathscr{Y}}{S}$ be the map defined by the following universal property. We require, for any simplicial set $K$ and any map $\sigma\colon\fromto{K}{S}$, a bijection
\begin{equation*}
\Mor_{S}(K,T_p\mathscr{Y})\cong\Mor_{s\Set(2)/(S,\iota S)}((K\times_S\mathscr{X},K\times_S\mathscr{X}_{\dag}),(\mathscr{Y},\mathscr{Y}_{\dag}))
\end{equation*}
\textup{(Nt. \ref{ntn:ordcatofpairs})}, functorial in $\sigma$. Then $r$ is a cocartesian fibration.
\begin{proof} We may use \cite[Cor. 3.2.2.13]{HTT} to define a cocartesian fibration $r'\colon\fromto{T'_p\mathscr{Y}}{S}$ with the universal property
\begin{equation*}
\Mor_{S}(K,T'_p\mathscr{Y})\cong\Mor_{S}(K\times_S\mathscr{X},\mathscr{Y}).
\end{equation*}
Thus $T'_p\mathscr{Y}$ is an $\infty$-category whose objects are pairs $(s,\phi)$ consisting of an object $s\in S_0$ and a functors $\phi\colon\fromto{\mathscr{X}_s}{\mathscr{Y}_s}$, and $T_p\mathscr{Y}\subset T'_p\mathscr{Y}$ is the full subcategory spanned by those pairs $(s,\phi)$ such that $\phi$ is a functor of pairs. An edge $\fromto{(s,\phi)}{(t,\psi)}$ in $T'_p\mathscr{Y}$ over an edge $\eta\colon\fromto{s}{t}$ of $S$ is $r'$-cocartesian if and only if the corresponding natural transformation $\fromto{\eta_{\mathscr{Y},!}\circ\phi\circ\eta_{\mathscr{X}}^{\star}}{\psi}$ is an equivalence. Since composites of functors of pairs are again functors of pairs, it follows that if $(s,\phi)$ is an object of $T_p\mathscr{Y}$, then so is $(t,\psi)$, whence it follows that $r$ is a cocartesian fibration.
\end{proof}
\end{prp}

Suppose that $\XX$ classifies $p$ and that $\YY$ classifies $q$. Since $\Fun_{\Pair_{\infty}}(\XX(s),\YY(s))$ is a full subcategory of $\Fun(\XX(s),\YY(s))$, it follows from \ref{rec:htt32213} that $T_p\mathscr{Y}$ is in fact classified by $\Fun_{\Pair_{\infty}}(\XX,\YY)$.

Suppose $S$ an $\infty$-category, and suppose $p\colon\fromto{\mathscr{X}}{S}$ a pair cartesian fibration. The construction $T_p$ is visibly a functor
\begin{equation*}
\fromto{(\Pair_{\infty/S}^{\cocart})_0}{(\Cat_{\infty/S}^{\cocart})_0}.
\end{equation*}
To show that $T_p$ defines a functor of $\infty$-categories $\fromto{\Pair_{\infty/S}^{\cocart}}{\Cat_{\infty/S}^{\cocart}}$, it suffices by Lm. \ref{lem:Paircocartisarelnerve} just to observe that the functor $T_p$ carries weak equivalences of $\Pair_{\infty/S}^{\cocart,0}$ to cocartesian equivalences. Hence we have the following.
\begin{prp}\label{nul:Tpafunctor} Suppose $p\colon\fromto{\mathscr{X}}{S}$ a cartesian fibration; then the assignment $\goesto{\mathscr{Y}}{T_p\mathscr{Y}}$ defines a functor
\begin{equation*}
\fromto{\Pair_{\infty/S}^{\cocart}}{\Cat_{\infty/S}^{\cocart}}.
\end{equation*}
\end{prp}


\subsection*{Waldhausen cartesian and cocartesian fibrations} Now we have laid the groundwork for our theory of \emph{Waldhausen cartesian} and \emph{cocartesian} fibrations.

\begin{dfn}\label{dfn:Waldcartfib} Suppose $S$ an $\infty$-category. A \textbf{\emph{Waldhausen cartesian fibration}} $p\colon\fromto{\mathscr{X}}{S}$ is a pair cartesian fibration satisfying the following conditions.
\begin{enumerate}[(\ref{dfn:Waldcartfib}.1)]
\item For any object $s$ of $S$, the pair
\begin{equation*}
\mathscr{X}_s\coloneq(\mathscr{X}\times_S\{s\},\mathscr{X}_{\dag}\times_S\{s\})
\end{equation*}
is a Waldhausen $\infty$-category.
\item For any morphism $\eta\colon\fromto{s}{t}$, the corresponding functor of pairs
\begin{equation*}
\eta^{\star}\colon\fromto{\mathscr{X}_t}{\mathscr{X}_s}
\end{equation*}
is an exact functor of Waldhausen $\infty$-categories.
\end{enumerate}

Dually, a \textbf{\emph{Waldhausen cocartesian fibration}} $p\colon\fromto{\mathscr{X}}{S}$ is a pair cocartesian fibration satisfying the following conditions.
\begin{enumerate}[(\ref{dfn:Waldcartfib}.1)]\addtocounter{enumi}{2}
\item For any object $s$ of $S$, the pair
\begin{equation*}
\mathscr{X}_s\coloneq(\mathscr{X}\times_S\{s\},\mathscr{X}_{\dag}\times_S\{s\})
\end{equation*}
is a Waldhausen $\infty$-category.
\item For any morphism $\eta\colon\fromto{s}{t}$, the corresponding functor of pairs
\begin{equation*}
\eta_!\colon\fromto{\mathscr{X}_s}{\mathscr{X}_t}
\end{equation*}
is an exact functor of Waldhausen $\infty$-categories.
\end{enumerate}
\end{dfn}

As with pair cartesian fibrations, Waldhausen cartesian fibrations classify functors to $\Wald$. The following is an immediate consequence of the definition.

\begin{prp}\label{prp:Waldcartsclassify} Suppose $S$ an $\infty$-category. Then a pair cartesian \textup{[}respectively, cocartesian\textup{]} fibration $p\colon\fromto{\mathscr{X}}{S}$ is a Waldhausen cartesian fibration \textup{[}resp., a Waldhausen cocartesian fibration\textup{]} if and only if the functor $\fromto{S^{\op}}{\Pair_{\infty}}$ \textup{[}resp., the functor $\fromto{S}{\Pair_{\infty}}$\textup{]} that classifies $p$ factors through $\Wald\subset\Pair_{\infty}$.
\end{prp}

\begin{prp}\label{prp:Waldcartfibsbasechange} The classes of Waldhausen cartesian fibrations and Waldhausen cocartesian fibrations are each stable under base change. That is, for any Waldhausen cartesian \textup{[}respectively, cocartesian\textup{]} fibration $\fromto{\mathscr{X}}{S}$ and for any functor $f\colon\fromto{S'}{S}$, if the pullback $\mathscr{X}'\coloneq\mathscr{X}\times_SS'$ is endowed with the pair structure in which a morphism is ingressive just in case it is carried to an equivalence in $S'$ and to an ingressive morphism of $\mathscr{X}$, then $\fromto{\mathscr{X}'}{S'}$ is a Waldhausen cartesian \textup{[}resp., cocartesian\textup{]} fibration.
\begin{proof} We treat the case of Waldhausen cartesian fibrations. By Pr. \ref{cor:paircartfibsbasechange}, $\fromto{\mathscr{X}'}{S'}$ is a pair cartesian fibration, so it remains to note that for any morphism $\eta\colon\fromto{s}{t}$ of $S'$, the induced functor of pairs
\begin{equation*}
\eta^{\star}\simeq f(\eta)^{\star}\colon\fromto{\mathscr{X}'_t\cong\mathscr{X}_{f(t)}}{\mathscr{X}_{f(s)}\cong\mathscr{X}'_s}
\end{equation*}
is an exact functor.
\end{proof}
\end{prp}

\begin{ntn} Denote by
\begin{equation*}
\Wald^{\cart}\textrm{\qquad [respectively, by\quad}\Wald^{\cocart}\textrm{\quad]}
\end{equation*}
the following subcategory of
\begin{equation*}
\Pair_{\infty}^{\cart}\textrm{\qquad [resp., of\quad}\Pair_{\infty}^{\cocart}\textrm{\quad].}
\end{equation*}
The objects of $\Wald^{\cart}$ [resp., of $\Wald^{\cocart}$] are Waldhausen cartesian fibrations [resp., Waldhausen cocartesian fibrations] $\fromto{\mathscr{X}}{S}$. A morphism
\begin{equation*}
\begin{tikzpicture} 
\matrix(m)[matrix of math nodes, 
row sep=4ex, column sep=4ex, 
text height=1.5ex, text depth=0.25ex] 
{\mathscr{X}&\mathscr{Y}\\ 
S^{\flat}&T^{\flat}\\}; 
\path[>=stealth,->,font=\scriptsize] 
(m-1-1) edge node[above]{$\psi$} (m-1-2) 
edge node[left]{$p$} (m-2-1) 
(m-1-2) edge node[right]{$q$} (m-2-2) 
(m-2-1) edge node[below]{$\phi$} (m-2-2); 
\end{tikzpicture}
\end{equation*}
of $\Pair_{\infty}^{\cart}$ (resp., $\Pair_{\infty}^{\cocart}$) is a morphism $\fromto{p}{q}$ of the subcategory $\Wald^{\cart}$ [resp., of $\Wald^{\cocart}$] if and only if $\psi$ induces exact functors $\fromto{\mathscr{X}_s}{\mathscr{Y}_{\phi(s)}}$ for every vertex $s\in S_0$.
\end{ntn}

The following is again a consequence of Pr. \ref{prp:Waldcartfibsbasechange} and \cite[Lm. 6.1.1.1]{HTT}.
\begin{lem}\label{lem:waldcartdiawald} The target functors
\begin{equation*}
\fromto{\Wald^{\cart}}{\Cat_{\infty}}\textrm{\quad and\quad}\fromto{\Wald^{\cocart}}{\Cat_{\infty}}
\end{equation*}
induced by the inclusion $\{1\}\subset\Delta^1$ are both cartesian fibrations.
\end{lem}

\begin{ntn} The fibers of the cartesian fibrations
\begin{equation*}
\fromto{\Wald^{\cart}}{\Cat_{\infty}}\textrm{\quad and\quad}\fromto{\Wald^{\cocart}}{\Cat_{\infty}}
\end{equation*}
over an object $\{S\}\subset\Cat_{\infty}$ will be denoted $\mathbf{Wald}_{\infty/S}^{\cart}$ and $\mathbf{Wald}_{\infty/S}^{\cocart}$, respectively.
\end{ntn}

\begin{prp}\label{prp:Waldstraightening} The equivalence of $\infty$-categories $\Pair_{\infty/S}^{\cart}\simeq\Fun(S^{\op},\Pair_{\infty})$ \textup{[}respectively, the equivalence of $\infty$-categories $\Pair_{\infty/S}^{\cocart}\simeq\Fun(S,\Pair_{\infty})$\textup{]} of \textup{Pr. \ref{prp:WaldcocartisFunSWald}} restricts to an equivalence of $\infty$-categories
\begin{equation*}
\mathbf{Wald}_{\infty/S}^{\cart}\simeq\Fun(S^{\op},\Wald)\textrm{\qquad\textup{[}resp.,\quad}\mathbf{Wald}_{\infty/S}^{\cocart}\simeq\Fun(S,\Wald)\textrm{\quad\textup{]}.}
\end{equation*}
\begin{proof} We treat the cartesian case. Note that $\mathbf{Wald}_{\infty/S}^{\cart}$ is the subcategory of the $\infty$-category $\Pair_{\infty/S}^{\cart}$ consisting of those objects and morphisms whose image under the equivalence $\Pair_{\infty/S}^{\cart}\simeq\Fun(S^{\op},\Pair_{\infty})$, lies in the subcategory $\Fun(S^{\op},\Wald)\subset\Fun(S^{\op},\Pair_{\infty})$. So one may identify $\mathbf{Wald}_{\infty/S}^{\cart}$ as the pullback
\begin{equation*}
\begin{tikzpicture} 
\matrix(m)[matrix of math nodes, 
row sep=4ex, column sep=4ex, 
text height=1.5ex, text depth=0.25ex] 
{\mathbf{Wald}_{\infty/S}^{\cart}&\Fun(S^{\op},\Wald)\\ 
\Pair_{\infty/S}^{\cart}&\Fun(S^{\op},\Pair_{\infty}).\\}; 
\path[>=stealth,->,font=\scriptsize] 
(m-1-1) edge (m-1-2) 
edge[right hook->] (m-2-1) 
(m-1-2) edge[right hook->] (m-2-2) 
(m-2-1) edge node[below,inner sep=0.6pt]{$\sim$} (m-2-2); 
\end{tikzpicture}
\end{equation*}
The result now follows from the fact that because the right-hand vertical map is a categorical fibration (\ref{rec:subcats}), this square is a homotopy pullback for the Joyal model structure.
\end{proof}
\end{prp}
\noindent As with pair fibrations (Cor. \ref{cor:colimitsofpairfibs}), we employ this result to observe that colimits of Waldhausen cartesian fibrations may be characterized fiberwise.
\begin{cor}\label{cor:colimsinWaldScocart} Suppose $S$ a small $\infty$-category, $K$ a small simplicial set. A functor $\mathscr{X}\colon\fromto{K^{\rhd}}{\mathbf{Wald}_{\infty/S}^{\cart}}$ \textup{[}respectively, a functor $\mathscr{X}\colon\fromto{K^{\rhd}}{\mathbf{Wald}_{\infty/S}^{\cocart}}$\textup{]} is a colimit diagram if and only if, for every vertex $s\in S_0$, the induced functor
\begin{equation*}
\mathscr{X}_s\colon\fromto{K^{\rhd}}{\Wald}
\end{equation*}
is a colimit diagram.
\end{cor}


\section{The derived $\infty$-category of Waldhausen $\infty$-categories} So far, we have built up a language for talking about the $\infty$-categories of interest to $K$-theorists. Now we want to study the $\infty$-category $\Wald$ of all these objects in some detail. More importantly, in later sections we'll need an enlargement of $\Wald$ on which we can define suitable \emph{derived functors}.

We take our inspiration from the following construction. Let $V(k)$ denote the ordinary category of vector spaces over a field $k$, and let $\mathrm{D}_{\geq 0}(k)$ be the \emph{connective derived $\infty$-category} of $V(k)$. That is, $\mathrm{D}_{\geq 0}(k)$ is a relative nerve of the relative category of (homologically graded) chain complexes whose homology vanishes in negative degrees, where a weak equivalence is declared to be a quasi-isomorphism.

The connective derived $\infty$-category is the vehicle with which one may define \emph{left derived functors} of right exact functors: one very general way of formulating this is to characterize $\mathrm{D}_{\geq 0}(k)$ as the $\infty$-category obtained from $V(k)$ by adding formal \emph{geometric realizations} --- that is, homotopy colimits of simplicial diagrams. More precisely, for any $\infty$-category $C$ that admits all geometric realizations, the functor
\begin{equation*}
\fromto{\Fun(\mathrm{D}_{\geq 0}(k),C)}{\Fun(NV(k),C)}
\end{equation*}
induced by the inclusion $\into{NV(k)}{\mathrm{D}_{\geq0}(k)}$ restricts to an equivalence from the full subcategory of $\Fun(\mathrm{D}_{\geq 0}(k),C)$ spanned by those functors $\fromto{\mathrm{D}_{\geq0}(k)}{C}$ that preserve geometric realizations to $\Fun(NV(k),C)$. (This characterization follows from the Dold--Kan correspondence; see \cite[Pr. 1.3.3.8]{HA} for a proof.) The objects of $\mathrm{D}_{\geq0}(k)$ can be represented as presheaves (of spaces) on the nerve of the category of \emph{finite-dimensional} vector spaces that carry direct sums to products.

In this section, we wish to mimic this construction, treating the $\infty$-category $\Wald$ of Waldhausen $\infty$-categories as formally analogous to the category $V(k)$. We thus define $\VWald$ as the $\infty$-category presheaves (of spaces) on the nerve of the category of suitably finite Waldhausen $\infty$-categories that carry direct sums to products. We call these presheaves \emph{virtual Waldhausen $\infty$-categories}. As with $\mathrm{D}_{\geq0}(k)$, virtual Waldhausen $\infty$-categories can be viewed as formal geometric realizations of simplicial Waldhausen $\infty$-categories, and the $\infty$-category $\VWald$ enjoys the following universal property: for any $\infty$-category $C$ that admits all geometric realizations, the functor
\begin{equation*}
\fromto{\Fun(\VWald,C)}{\Fun(\Wald,C)}
\end{equation*}
induced by the Yoneda embedding $\into{\Wald}{\VWald}$ restricts to an equivalence from the full subcategory of $\Fun(\VWald,C)$ spanned by those functors $\fromto{\VWald}{C}$ that preserve geometric realizations to $\Fun(\Wald,C)$.

To get this idea off the ground, it is clear that we must analyze limits and colimits in $\Wald$. Along the way, we'll find that, indeed, $\Wald$ is rather a lot like $V(k)$.


\subsection*{Limits and colimits of pairs of $\infty$-categories} We first analyze limits and colimits in the $\infty$-category $\Pair_{\infty}$.

\begin{rec} Suppose $C$ a locally small $\infty$-category \cite[Df. 5.4.1.3]{HTT}. For a regular cardinal $\kappa<\kappa_0$, recall \cite[Df. 5.5.7.1]{HTT} that $C$ is said to be \textbf{\emph{$\kappa$-compactly generated}} (or simply \textbf{\emph{compactly generated}} if $\kappa=\omega$) if it is $\kappa$-accessible and admits all small colimits. From this it will follow that $C$ admits all small limits as well. It follows from Simpson's theorem \cite[Th. 5.5.1.1]{HTT} that $C$ is $\kappa$-compactly generated if and only if it is a $\kappa$-accessible localization of the $\infty$-category of presheaves $\mathscr{P}(C_0)=\Fun(C_0^{\op},\Kan)$ of small spaces on some small $\infty$-category $C_0$.
\end{rec}

\begin{prp}\label{thm:pairasloc} The $\infty$-category $\Pair_{\infty}$ is an $\omega$-accessible localization of the arrow $\infty$-category $\mathscr{O}(\Cat_{\infty})$.
\begin{proof} We use \ref{prp:PairsubcatOcat} to identify $\Pair_{\infty}$ with a full subcategory of $\mathscr{O}(\Cat_{\infty})$. Now the condition that an object $\fromto{C'}{C}$ of $\mathscr{O}(\Cat_{\infty})$ be a monomorphism is equivalent to the demand that the functors
\begin{equation*}
\fromto{\iota C'}{\iota C'\times^h_{\iota C}\iota C'}\textrm{\quad and\quad}\fromto{\iota\mathscr{O}(C')}{\iota\mathscr{O}(C')\times^h_{\iota\mathscr{O}(C)}\iota\mathscr{O}(C')}
\end{equation*}
be isomorphisms of $h\Cat_{\infty}$. This, in turn, is the requirement that the object $\fromto{C'}{C}$ be $S$-local, where $S$ is the set
\begin{equation*}
S\coloneq\left\{
\begin{tikzpicture}[baseline]
\matrix(m)[matrix of math nodes, 
row sep=6ex, column sep=6ex, 
text height=1.5ex, text depth=0.25ex] 
{\Delta^p\sqcup\Delta^p&\Delta^p\\ 
\Delta^p&\Delta^p\\}; 
\path[>=stealth,->,font=\scriptsize] 
(m-1-1) edge node[above]{$\nabla$} (m-1-2) 
edge node[left]{$\nabla$} (m-2-1) 
(m-1-2) edge[-,double distance=1.5pt] (m-2-2) 
(m-2-1) edge[-,double distance=1.5pt] (m-2-2); 
\end{tikzpicture}
\ \Big|\quad\mathbf{p}\in\Delta\right\}
\end{equation*}
of morphisms of $\mathscr{O}(\Cat_{\infty})$. The condition that an object $\fromto{C'}{C}$ of $\mathscr{O}(\Cat_{\infty})$ induce an equivalence $\fromto{\iota C'}{\iota C}$ is equivalent to the requirement that it be local with respect to the singleton
\begin{equation*}
\{\phi\colon\into{[\fromto{\varnothing}{\Delta^0}]}{[\fromto{\Delta^0}{\Delta^0}]}\}.
\end{equation*}
Hence $\Pair_{\infty}$ is equivalent to the full subcategory of the $S\cup\{\phi\}$-local objects of $\mathscr{O}(\Cat_{\infty})$. Now it is easy to see that the $S\cup\{\phi\}$-local objects of $\mathscr{O}(\Cat_{\infty})$ are closed under filtered colimits; hence by \cite[Pr. 5.5.3.6 and Cor. 5.5.7.3]{HTT}, the $\infty$-category $\Pair_{\infty}$ is an $\omega$-accessible localization.
\end{proof}
\end{prp}

\begin{cor}\label{cor:paircpctlygen} The $\infty$-category $\Pair_{\infty}$ is compactly generated.
\end{cor}

\begin{cor}\label{cor:limsofpairs} The $\infty$-category $\Pair_{\infty}$ admits all small limits, and the inclusion
\begin{equation*}
\into{\Pair_{\infty}}{\mathscr{O}(\Cat_{\infty})}
\end{equation*}
preserves them.
\end{cor}

\begin{cor}\label{cor:colimsofpairs} The $\infty$-category $\Pair_{\infty}$ admits all small colimits, and the inclusion
\begin{equation*}
\into{\Pair_{\infty}}{\mathscr{O}(\Cat_{\infty})}
\end{equation*}
preserves small filtered colimits.
\end{cor}

\begin{cor} Any pair $\mathscr{C}$ is the colimit of its compact subpairs.
\end{cor}

\begin{exm}\label{exm:finiteimpliescompactforpairs} Suppose $\mathscr{C}$ a pair such that $\mathscr{C}$ and $\mathscr{C}_{\dag}$ are each compact in $\Cat_{\infty}$. Then $\mathscr{C}$ is compact in $\Pair_{\infty}$. Indeed, suppose $\mathscr{D}\colon\fromto{\Lambda^{\rhd}}{\Pair_{\infty}}$ is a colimit of a filtered diagram of pairs. The compactness of $\mathscr{C}$ and $\mathscr{C}_{\dag}$ yields an equivalence
\[\Pair_{\infty}^{\Delta}(\mathscr{C},\mathscr{D}_{+\infty})\simeq\colim_{\alpha}\Cat_{\infty}^{\Delta}(\mathscr{C},\mathscr{D}_{\alpha})\times_{\colim_{\beta}\Cat_{\infty}^{\Delta}(\mathscr{C}_{\dag},\mathscr{D}_{\beta})}\colim_{\gamma}\Cat_{\infty}^{\Delta}(\mathscr{C}_{\dag},\mathscr{D}_{\gamma,\dag}).\]
Now since filtered colimits in spaces commute with finite limits, one has
\[\Pair_{\infty}^{\Delta}(\mathscr{C},\mathscr{D}_{+\infty})\simeq\colim_{\alpha}\Cat_{\infty}^{\Delta}(\mathscr{C},\mathscr{D}_{\alpha})\times_{\Cat_{\infty}^{\Delta}(\mathscr{C}_{\dag},\mathscr{D}_{\alpha})}\Cat_{\infty}^{\Delta}(\mathscr{C}_{\dag},\mathscr{D}_{\alpha,\dag}),\]
which implies that $\mathscr{C}$ is compact in $\Pair_{\infty}$.

In particular, any pair $\mathscr{C}$ in which both $\mathscr{C}$ and $\mathscr{C}_{\dag}$ are finite simplicial sets is compact.
\end{exm}


\subsection*{Limits and filtered colimits of Waldhausen $\infty$-categories} Now we construct limits and colimits in $\Wald$.

\begin{prp}\label{thm:Waldlims} The $\infty$-category $\Wald$ admits all small limits, and the inclusion functor $\fromto{\Wald}{\Pair_{\infty}}$ preserves them.
\begin{proof} We employ \cite[Pr. 4.4.2.6]{HTT} to reduce the problem to proving the existence of products and pullbacks in $\Wald$. To complete the proof, we make the following observations.
\begin{enumerate}[(\ref{thm:Waldlims}.1)]
\item Suppose $I$ a set, suppose $(\mathscr{C}_i)_{i\in I}$ an $I$-tuple of pairs of $\infty$-categories, and suppose $\mathscr{C}$ the product of these pairs. If for each $i\in I$, the pair $\mathscr{C}_i$ is a Waldhausen $\infty$-category, then so is $\mathscr{C}$. Moreover, if $\mathscr{D}$ is a Waldhausen $\infty$-category, then a functor of pairs $\fromto{\mathscr{D}}{\mathscr{C}}$ is exact if and only if the composite
\begin{equation*}
\mathscr{D}\ \tikz[baseline]\draw[>=stealth,->](0,0.5ex)--(0.5,0.5ex);\ \mathscr{C}\ \tikz[baseline]\draw[>=stealth,->](0,0.5ex)--(0.5,0.5ex);\ \mathscr{C}_i
\end{equation*}
is exact for any $i\in I$. This follows directly from the fact that limits and colimits of a product are computed objectwise \cite[Cor. 5.1.2.3]{HTT}.
\item Suppose
\begin{equation*}
\begin{tikzpicture} 
\matrix(m)[matrix of math nodes, 
row sep=4ex, column sep=4ex, 
text height=1.5ex, text depth=0.25ex] 
{\mathscr{E}'&\mathscr{F}'\\ 
\mathscr{E}&\mathscr{F}\\}; 
\path[>=stealth,->,font=\scriptsize] 
(m-1-1) edge node[above]{$q'$} (m-1-2) 
edge node[left]{$p'$} (m-2-1) 
(m-1-2) edge node[right]{$p$} (m-2-2) 
(m-2-1) edge node[below]{$q$} (m-2-2); 
\end{tikzpicture}
\end{equation*}
a pullback diagram of pairs of $\infty$-categories. Suppose moreover that $\mathscr{E}$, $\mathscr{F}$, and $\mathscr{F}'$ are all Waldhausen $\infty$-categories, and $p$ and $q$ are exact functors. Then by \cite[Lm. 5.4.5.2]{HTT} and its dual, $\mathscr{E}'$ admits both an initial object and a terminal object, each of which is preserved by $p'$ and $q'$, and they are equivalent since they are so in $\mathscr{E}$, $\mathscr{F}$, and $\mathscr{F}'$. It now follows from \cite[Lm. 5.4.5.5]{HTT} that $\mathscr{E}'$ is a Waldhausen $\infty$-category, and for any Waldhausen $\infty$-category $\mathscr{D}$, a functor of pairs $\psi\colon\fromto{\mathscr{D}}{\mathscr{E}}$ is exact if and only if the composites $p'\circ\psi$ and $q'\circ\psi$ are exact.\qedhere
\end{enumerate}
\end{proof}
\end{prp}

We obtain a similar characterization of filtered colimits in $\Wald$.

\begin{prp}\label{thm:Waldfiltcolims} The $\infty$-category $\Wald$ admits all small filtered colimits, and the inclusion functor $\fromto{\Wald}{\Pair_{\infty}}$ preserves them.
\begin{proof} Suppose $A$ a filtered $\infty$-category, and suppose $\fromto{A}{\Wald}$ a functor given by the assignment $\goesto{a}{\mathscr{C}_a}$, and suppose $\mathscr{C}$ the colimit of the composite functor
\begin{equation*}
A\ \tikz[baseline]\draw[>=stealth,->](0,0.5ex)--(0.5,0.5ex);\ \Wald\ \tikz[baseline]\draw[>=stealth,->](0,0.5ex)--(0.5,0.5ex);\ \Pair_{\infty}.
\end{equation*}
Pushouts of ingressive morphisms in $\mathscr{C}$ exist and are ingressive morphisms. Furthermore, the image of any zero object in any $\mathscr{C}_a$ is initial in both $\mathscr{C}$ and in $\mathscr{C}_{\dag}$. Both of these facts follow by precisely the same argument as \cite[Pr. 5.5.7.11]{HTT}. The dual argument ensures that this image is also terminal in $\mathscr{C}$, whence it is a zero object.
\end{proof}
\end{prp}


\subsection*{Direct sums of Waldhausen $\infty$-categories} The $\infty$-category $\Wald$ also admits finite \emph{direct sums}, i.e., that finite products in $\Wald$ are also finite coproducts.

\begin{dfn} \label{item:directsums} Suppose $C$ is an $\infty$-category. Then $C$ is said to \emph{admit finite direct sums} if the following conditions hold.
\begin{enumerate}[(\ref{item:directsums}.1)]
\item The $\infty$-category $C$ is pointed.
\item The $\infty$-category $C$ has all finite products and coproducts.
\item For any finite set $I$ and any $I$-tuple $(X_i)_{i\in I}$ of objects of $C$, the map
\begin{equation*}
\fromto{\coprod X_I}{\prod X_I}
\end{equation*}
in $hC$ --- given by the maps $\phi_{ij}\colon\fromto{X_i}{X_j}$, where $\phi_{ij}$ is zero unless $i=j$, in which case it is the identity --- is an isomorphism.
\end{enumerate}
If $C$ admits finite direct sums, then for any finite set $I$ and any $I$-tuple $(X_i)_{i\in I}$ of objects of $C$, we denote by $\bigoplus X_I$ the product (or, equivalently, the coproduct) of the $X_i$.

We will say that $C$ is \emph{additive} if it admits direct sums, and the resulting commutative monoids $\Mor_{h\mathscr{A}}(X,Y)$ are all abelian groups.
\end{dfn}

\begin{prp}\label{prp:Waldhassums} The $\infty$-category $\Wald$ admits finite direct sums.
\begin{proof} The Waldhausen $\infty$-category $\Delta^0$ is a zero object. To complete the proof, it suffices to show that for any finite set $I$ and any $I$-tuple of Waldhausen $\infty$-categories $(\mathscr{C}_i)_{i\in I}$ with product $\mathscr{C}$, the functors $\phi_i\colon\fromto{\mathscr{C}_i}{\mathscr{C}}$ --- given by the functors $\phi_{ij}\colon\fromto{\mathscr{C}_i}{\mathscr{C}_j}$, where $\phi_{ij}$ is zero unless $j=i$, in which case it is the identity --- are exact and exhibit $\mathscr{C}$ as the \emph{coproduct} of $(\mathscr{C}_i)_{i\in I}$. To prove this, we employ \cite[Th. 4.2.4.1]{HTT} to reduce the problem to showing that for any Waldhausen $\infty$-category $\mathscr{D}$, the map
\begin{equation*}
\fromto{\Wald^{\Delta}(\mathscr{C},\mathscr{D})}{\prod_{i\in I}\Wald^{\Delta}(\mathscr{C}_i,\mathscr{D})}
\end{equation*}
induced by the functor $\phi_i$ is a weak homotopy equivalence. We prove the stronger claim that the functor
\begin{equation*}
w\colon\fromto{\Fun_{\Wald}(\mathscr{C},\mathscr{D})}{\prod_{i\in I}\Fun_{\Wald}(\mathscr{C}_i,\mathscr{D})}
\end{equation*}
is an equivalence of $\infty$-categories. 

For this, consider the following composite
\begin{equation*}
\begin{tikzpicture} 
\matrix(m)[matrix of math nodes, 
row sep=4ex, column sep=-6ex, 
text height=1.5ex, text depth=0.25ex] 
{&\Fun(\mathscr{C},\Fun(NI,\mathscr{D}))&[14ex]\Fun(\mathscr{C},\Colim((NI)^{\rhd},\mathscr{D}))&\\ 
\prod_{i\in I}\Fun(\mathscr{C}_i,\mathscr{D})&&&\Fun(\mathscr{C},\mathscr{D})\\}; 
\path[>=stealth,->,font=\scriptsize] 
(m-1-2) edge node[above]{$r$} (m-1-3) 
edge[<-,inner sep=4pt] node[left]{$u$} (m-2-1) 
(m-1-3) edge[inner sep=4pt] node[right]{$e$} (m-2-4); 
\end{tikzpicture}
\end{equation*}
where $u$ is the functor corresponding to the functor
\begin{equation*}
\fromto{\mathscr{C}\times\prod_{i\in I}\Fun(\mathscr{C}_i,\mathscr{D})\cong\prod_{i\in I}(\mathscr{C}_i\times\Fun(\mathscr{C}_i,\mathscr{D}))}{\prod_{i\in I}\mathscr{D}},
\end{equation*}
where $r$ is a section of the trivial fibration
\begin{equation*}
\fromto{\Fun(\mathscr{C},\Colim((NI)^{\rhd},\mathscr{D}))}{\Fun(\mathscr{C},\Fun(NI,\mathscr{D}))},
\end{equation*}
and $e$ is the functor induced by the functor $\fromto{\Colim((NI)^{\rhd},\mathscr{D})}{\mathscr{D}}$ given by evaluation at the cone point $\infty$. This composite restricts to a functor
\begin{equation*}
v\colon\fromto{\prod_{i\in I}\Fun_{\Wald}(\mathscr{C}_i,\mathscr{D})}{\Fun_{\Wald}(\mathscr{C},\mathscr{D})};
\end{equation*}
indeed, one checks directly that if $(\psi_i\colon\fromto{\mathscr{C}_i}{\mathscr{D}})_{i\in I}$ is an $I$-tuple of exact functors, then a functor $\psi\colon\fromto{\mathscr{C}}{\mathscr{D}}$ that sends a simplex $\sigma=(\sigma_i)_{i\in I}$ to a coproduct $\coprod_{i\in I}\psi_i(\sigma_i)$ in $\mathscr{D}$ is exact, and the situation is similar for natural transformations of exact functors.

We claim that the functor $v$ is a homotopy inverse to $w$. A homotopy $w\circ v\simeq\id$ can be constructed directly from the canonical equivalences
\begin{equation*}
Y\simeq Y\sqcup\coprod_{i\in I-\{j\}}0_i
\end{equation*}
for any zero objects $0_i$ in $\mathscr{D}$. In the other direction, the existence of a homotopy $v\circ w\simeq\id$ follows from the observation that the natural transformations $\fromto{\phi_i\circ\pr_i}{\id}$ exhibit the identity functor on $\mathscr{C}$ as the coproduct $\coprod_{i\in I}\phi_i\circ\pr_i$.
\end{proof}
\end{prp}

Since any small coproduct can be written as a filtered colimit of finite coproducts, we deduce the following.
\begin{cor} The $\infty$-category $\Wald$ admits all small coproducts.
\end{cor}
\noindent Coproducts in $\Wald$ enjoy a description reminiscent of the description of coproducts in the category of vector spaces over a field: for any set $I$ and an $I$-tuple $(\mathscr{C}_i)_{i\in I}$ of Waldhausen $\infty$-categories, $\coprod_{i\in I}\mathscr{C}_i$ is equivalent to the full subcategory of $\prod_{i\in I}\mathscr{C}_i$ spanned by those objects $(X_i)_{i\in I}$ such that all but a finite number of the objects $X_i$ are zero objects.


\subsection*{Accessibility of \protect{$\Wald$}} Finally, we set about showing that $\Wald$ is an accessible $\infty$-category. In fact, we prove the following stronger result.
\begin{prp}\label{prp:Waldcpctlygen} The $\infty$-category $\Wald$ is compactly generated.
\begin{proof} The $\infty$-category $\Kan$ is compactly generated, as is the $\infty$-category $\Kan_{\ast}$ of pointed Kan complexes. We have already seen that $\Pair_{\infty}$ is compactly generated. Additionally, we may contemplate the full subcategory $\mathbf{Mono}\subset\Fun(\Delta^1,\Kan)$ spanned by those functors $\fromto{C}{D}$ that are monomorphisms. We claim that $\mathbf{Mono}$ is also compactly generated. Indeed, $\mathbf{Mono}$ is nothing more than the full subcategory of $\{\phi\}$-local objects, where $\phi$ denotes the map
\begin{equation*}
\begin{tikzpicture}[baseline]
\matrix(m)[matrix of math nodes, 
row sep=6ex, column sep=6ex, 
text height=1.5ex, text depth=0.25ex] 
{\partial\Delta^1&\Delta^0\\ 
\Delta^0&\Delta^0,\\}; 
\path[>=stealth,->,font=\scriptsize] 
(m-1-1) edge (m-1-2) 
edge (m-2-1) 
(m-1-2) edge[-,double distance=1.5pt] (m-2-2) 
(m-2-1) edge[-,double distance=1.5pt] (m-2-2); 
\end{tikzpicture}
\end{equation*}
and $\mathbf{Mono}\subset\Fun(\Delta^1,\Kan)$ is clearly stable under filtered colimits, whence it is an $\omega$-accessible localization by \cite[Cor. 5.5.7.3]{HTT}.

Now we define some functors among these $\infty$-categories. Denote by $\iota$ the interior functor $\fromto{\Pair_{\infty}}{\Kan}$ \ref{ntn:interior}. Write $F\colon\fromto{\Pair_{\infty}}{\Kan}$ for the functor $\goesto{\mathscr{C}}{\Map_{\Pair_{\infty}}(\mathscr{Q}^2,\mathscr{C})}$ corepresented by $\mathscr{Q}^2$. We also have the target functor $\fromto{\mathbf{Mono}}{\Kan}$ and the forgetful functor $\fromto{\Kan_{\ast}}{\Kan}$. It is easy to see that all of these functors preserve limits and filtered colimits. Therefore we may form the fiber product
\[\CC\coloneq\mathbf{Mono}\times_{\Kan,F}\Pair_{\infty}\times_{U,\Kan}\Kan_{\ast},\]
which by \cite[Pr. 5.5.7.6]{HTT} is thus compactly generated.

The objects of $\CC$ can thus be thought of as $4$-tuples $(\mathscr{C},\mathscr{C}_{\dag},I,M)$, where $(\mathscr{C},\mathscr{C}_{\dag})$ is a pair, $I$ is an object of $\mathscr{C}$, and $M\subset\Map_{\Pair_{\infty}}(\mathscr{Q}^2,\mathscr{C})$ is a collection of functors of pairs $\fromto{\mathscr{Q}^2}{\mathscr{C}}$. A morphism $\fromto{(\mathscr{C},\mathscr{C}_{\dag},I,M)}{(\mathscr{D},\mathscr{D}_{\dag},J,N)}$ is a functor of pairs $\fromto{(\mathscr{C},\mathscr{C}_{\dag})}{(\mathscr{D},\mathscr{D}_{\dag})}$ that carries $I$ to $J$ and carries any square in $M$ to a square in $N$. In particular, $\Wald$ can be identified with the full subcategory of $\CC$ spanned by those objects $(\mathscr{C},\mathscr{C}_{\dag},I,M)$ such that $(\mathscr{C},\mathscr{C}_{\dag})$ is a Waldhausen $\infty$-category, $I$ is a zero object of $\mathscr{C}$, and $M$ is the collection of pushout squares $\fromto{\mathscr{Q}^2}{\mathscr{C}}$.

Now we have already shown that the inclusion $\into{\Wald}{\CC}$ preserves limits and filtered colimits. We now intend to construct a left adjoint to this inclusion, whence $\Wald$ is compactly generated by \cite[Cor. 5.5.7.3]{HTT}.

In light of \cite[Pr. 5.2.7.8]{HTT}, it suffices, for any object $(\mathscr{C},\mathscr{C}_{\dag},I,M)$ of $\CC$, to give a localization $F\colon\fromto{(\mathscr{C},\mathscr{C}_{\dag},I,M)}{(\mathscr{D},\mathscr{D}_{\dag},J,N)}$ relative to $\Wald\subset\CC$. To do this, we present a kind of pair version of \cite[\S 5.3.6]{HTT}.

First, we form the $\infty$-category of presheaves of pointed spaces
\[\mathscr{P}_{\ast}(\mathscr{C})\coloneq\Fun(\mathscr{C}^{\op},\Kan_{\ast}),\]
and we write $j$ for the composite of the Yoneda embedding $\into{\mathscr{C}}{\mathscr{P}(\mathscr{C})}$ with the pointing functor $\fromto{\mathscr{P}(\mathscr{C})}{\mathscr{P}_{\ast}(\mathscr{C})}$.

Now for any square $p\colon\fromto{\mathscr{Q}^2}{\mathscr{C}}$ in $M$, select a colimit $x_p$ of $j\circ p|_{\Lambda_0\mathscr{Q}^2}$, and consider the natural map $f_p\colon\fromto{x_p}{j(p(1,1))}$ (which is unique up to a contractible choice). Now let $\phi$ be the canonical map $\fromto{j(I)}{0}$ from $j(I)$ to the zero object of $\mathscr{P}(\mathscr{C})$. Write $S$ for the set
\[\{f_p\ |\ p\in M\}\cup\{\phi\},\]
and form the $\infty$-category $L_S\mathscr{P}_{\ast}(\mathscr{C})$ of $S$-local objects is $\mathscr{P}_{\ast}(\mathscr{C})$. Write $L$ for the left adjoint to the inclusion $\into{L_S\mathscr{P}_{\ast}(\mathscr{C})}{\mathscr{P}_{\ast}(\mathscr{C})}$.

We define $L_S\mathscr{P}_{\ast}(\mathscr{C})_{\dag}$ as the smallest subcategory of $L_S\mathscr{P}_{\ast}(\mathscr{C})_{\dag}$ that contains all the equivalences, the image of any map of $M$ under $L\circ j$, and any map $\fromto{0}{x}$, and that is stable under pushouts.

Finally, we select the smallest full subcategory $\mathscr{D}\subset L_S\mathscr{P}_{\ast}(\mathscr{C})$ that contains the essential image of $L\circ j$ that is closed under pushouts along any morphism of $L_S\mathscr{P}_{\ast}(\mathscr{C})_{\dag}$, and we set
\[\mathscr{D}_{\dag}\coloneq\mathscr{D}\cap L_S\mathscr{P}_{\ast}(\mathscr{C})_{\dag}.\]
We set $F\coloneq L\circ j$, and we set $J\coloneq F(I)$, and we let $N$ be the collection of all pushout squares in $\mathscr{D}$ along a map of $\mathscr{D}_{\dag}$.

The claim is now threefold:
\begin{enumerate}[(\ref{prp:Waldcpctlygen}.1)]
\item The pair $(\mathscr{D},\mathscr{D}_{\dag})$ is a Waldhausen $\infty$-category, $J$ is a zero object, and $N$ consists of pushout squares $\fromto{\mathscr{Q}^2}{\mathscr{D}}$.
\item The functor $F$ carries $\mathscr{C}_{\dag}$ to $\mathscr{D}_{\dag}$, $I$ to $J$, and $M$ to $N$.
\item For any Waldhausen $\infty$-category $\mathscr{E}$, the functor $F$ induces an equivalence
\[\fromto{\Map_{\Wald}(\mathscr{D},\mathscr{E})}{\Map_{\CC}(\mathscr{C},\mathscr{E})}.\]
\end{enumerate}
The first two claims are now obvious from the construction. The last claim as in the proof of \cite[Pr. 5.3.6.2(2)]{HTT}.
\end{proof}
\end{prp}

This result shows that in fact the $\infty$-category $\Wald$ admits \emph{all} small colimits, not only the filtered ones. However, these other colimits are not preserved by the sorts of invariants in which we are interested, and so we will regard them as pathological. Nevertheless, we will have use for the following.

\begin{cor} The $\infty$-category $\Wald$ is $\omega$-accessible. 
\end{cor}

\begin{cor}\label{cor:WaldisInd} The $\infty$-category $\Wald$ may be identified with the $\mathrm{Ind}$-objects of the full subcategory $\Wald^{\omega}\subset\Wald$ spanned by the compact Waldhausen $\infty$-categories:
\begin{equation*}
\Wald\simeq\Ind(\Wald^{\omega}).
\end{equation*}
\end{cor}
\noindent We obtain a further corollary by combining Prs. \ref{prp:Waldcpctlygen}, \ref{thm:Waldlims}, and \ref{thm:Waldfiltcolims} together with the adjoint functor theorem \cite[Cor. 5.5.2.9]{HTT}.
\begin{cor}\label{cor:Wconstruction} The forgetful functor $\fromto{\Wald}{\Pair_{\infty}}$ admits a left adjoint $W\colon\fromto{\Pair_{\infty}}{\Wald}$.
\end{cor}

\begin{nul} Since the opposite functor $\fromto{\Wald}{\coWald}$ is an equivalence of $\infty$-categories, it follows that this whole crop of structural results also hold for $\coWald$. That is, $\coWald$ admits all small limits and all small filtered colimits, and the inclusion functor $\fromto{\coWald}{\Pair_{\infty}}$ preserves each of them. Similarly, $\coWald$ admits finite direct sums and all small coproducts, and it is compactly generated.
\end{nul}


\subsection*{Virtual Waldhausen $\infty$-categories} Now we are prepared to introduce a convenient enlargement of the $\infty$-category $\Wald$. In effect, we aim to ``correct'' the colimits of $\Wald$ that we regard as pathological. As with the formation of $\mathrm{D}_{\geq0}(k)$ from $NV(k)$ (see the introduction of this section) --- or indeed with the formation of the $\infty$-category of spaces from the nerve of the category of sets ---, we will add to $\Wald$ formal geometric realizations and nothing more. The result is the \emph{derived $\infty$-category of Waldhausen $\infty$-categories}, whose homotopy theory forms the basis of our work here.

The definition is exactly as for $\mathrm{D}_{\geq0}(k)$:
\begin{dfn}\label{dfn:virtWald} A \textbf{\emph{virtual Waldhausen $\infty$-category}} is a presheaf
\begin{equation*}
\mathscr{X}\colon\fromto{(\Wald^{\omega})^{\op}}{\Kan}
\end{equation*}
that preserves products.
\end{dfn}

\begin{ntn} Denote by
\begin{equation*}
\VWald\subset\Fun(\Wald^{\omega,\op},\Kan)
\end{equation*}
the full subcategory spanned by the virtual Waldhausen $\infty$-categories. In other words, $\VWald$ is the \emph{nonabelian derived $\infty$-category} of $\Wald^{\omega}$ \cite[\S 5.5.8]{HTT}. We simply call $\VWald$ the \textbf{\emph{derived $\infty$-category of Waldhausen $\infty$-categories}}.
\end{ntn}

\begin{ntn} For any $\infty$-category $C$, we shall write $\mathscr{P}(C)$ for the $\infty$-category $\Fun(C^{\op},\Kan)$ of presheaves of small spaces on $C$. If $C$ is locally small, then there exists a Yoneda embedding \cite[Pr. 5.1.3.1]{HTT}
\begin{equation*}
j\colon\into{C}{\mathscr{P}(C)}.
\end{equation*}
\end{ntn}

\begin{rec}\label{rec:PABofC} Suppose $\mathscr{A}\subset\mathscr{B}$ two classes of small simplicial sets, and suppose $C$ an $\infty$-category that admits all $\mathscr{A}$-shaped colimits (\ref{rec:Ashapedcolim}). Recall \cite[\S 5.3.6]{HTT} that there exist an $\infty$-category $\mathscr{P}_{\mathscr{A}}^{\mathscr{B}}(C)$ and a fully faithful functor $j\colon\into{C}{\mathscr{P}_{\mathscr{A}}^{\mathscr{B}}(C)}$ such that for any $\infty$-category $D$ with all $\mathscr{B}$-shaped colimits, $j$ induces an equivalance of $\infty$-categories (\ref{rec:Ashapedcolim})
\begin{equation*}
\equivto{\Fun_{\mathscr{B}}(\mathscr{P}_{\mathscr{A}}^{\mathscr{B}}(C),D)}{\Fun_{\mathscr{A}}(C,D)}.
\end{equation*}

Recall also \cite[Nt. 6.1.2.12]{HTT} that for any $\infty$-category $C$, the colimit of a simplicial diagram $X\colon\fromto{N\Delta^{\op}}{C}$ will be called the \textbf{\emph{geometric realization}} of $X$.
\end{rec}

\begin{nul}\label{nul:VWaldasPDKWald} In the notation of \ref{rec:PABofC}, the $\infty$-category $\VWald$ can be identified with any of the following $\infty$-categories:
\begin{enumerate}[(\ref{nul:VWaldasPDKWald}.1)]
\item the $\infty$-category $\mathscr{P}_{\varnothing}^{\{N\Delta^{\op}\}}\Wald$,
\item the $\infty$-category $\mathscr{P}_{\mathscr{R}}^{\mathscr{S}}\Wald$, where $\mathscr{R}$ is the collection of small, filtered simplicial sets and $\mathscr{S}$ is the collection of small, sifted simplicial sets,
\item the $\infty$-category $\mathscr{P}_{\varnothing}^{\mathscr{S}}\Wald^{\omega}$, and
\item the $\infty$-category $\mathscr{P}_{\mathscr{D}}^{\mathscr{K}}\Wald^{\omega}$, where $\mathscr{D}$ is the collection of finite discrete simplicial sets, and $\mathscr{K}$ is the collection of small simplicial sets.
\end{enumerate}
The equivalence of these characterizations follow directly from Cor. \ref{cor:WaldisInd}, the description of the nonabelian derived $\infty$-category of \cite[Pr. 5.5.8.16]{HTT}, the fact that sifted colimits can be decomposed as geometric realizations of filtered colimits \cite[Pr. 5.5.8.15]{HTT}, and the transitivity assertion of \cite[Pr. 5.3.6.11]{HTT}.

We may summarize these characterizations by saying that the Yoneda embedding is a fully faithful functor
\begin{equation*}
j\colon\into{\Wald}{\VWald}
\end{equation*}
that induces, for any $\infty$-category $E$ that admits geometric realizations, any $\infty$-category $E'$ that admits all sifted colimits, and any $\infty$-category that admits all small colimits, equivalences (\ref{rec:Ashapedcolim})
\begin{eqnarray}
&\equivto{\Fun_{\{N\Delta^{\op}\}}(\VWald,E)}{\Fun(\Wald,E)};&\nonumber\\
&\equivto{\Fun_{\mathscr{J}}(\VWald,E')}{\Fun_{\mathscr{I}}(\Wald,E')};&\nonumber\\
&\equivto{\Fun_{\mathscr{J}}(\VWald,E')}{\Fun_{\mathscr{I}}(\Wald^{\omega},E')};&\nonumber\\
&\equivto{\Fun_{\mathscr{K}}(\VWald,E'')}{\Fun_{\mathscr{D}}(\Wald^{\omega},E'')}.&\nonumber
\end{eqnarray}
\end{nul}

\begin{dfn}\label{dfn:lderivWald} Suppose $E$ an $\infty$-category that admits all sifted colimits. Then a functor
\begin{equation*}
\Phi\colon\fromto{\VWald}{E}
\end{equation*}
that preserves all sifted colimits will be said to be the \textbf{\emph{left derived functor}} of the corresponding $\omega$-continuous functor $\phi=\Phi\circ j\colon\fromto{\Wald}{E}$ (which preserves filtered colimits) or of the further restriction $\fromto{\Wald^{\omega}}{E}$ of $\phi$ to $\Wald^{\omega}$.
\end{dfn}

\begin{prp} The $\infty$-category $\VWald$ is compactly generated. Moreover, it admits all direct sums, and the inclusion $j$ preserves them.
\begin{proof} The first statement is \cite[Pr. 5.5.8.10(6)]{HTT}. To see that $\VWald$ admits direct sums, we use the fact that we may exhibit any object of $\VWald$ as a sifted colimit of compact Waldhausen $\infty$-categories in $\mathscr{P}(\Wald^{\omega})$ \cite[Lm. 5.5.8.14]{HTT}; now since sifted colimits commute with both finite products \cite[Lm. 5.5.8.11]{HTT} and coproducts, and since $j$ preserves products and finite coproducts \cite[Lm. 5.5.8.10(2)]{HTT}, the proof is complete.
\end{proof}
\end{prp}


\subsection*{Realizations of Waldhausen cocartesian fibrations} We now give an explicit construction of colimits in $\VWald$ of sifted diagrams of Waldhausen $\infty$-categories when they are exhibited as Waldhausen cocartesian fibrations.

The idea behind our construction comes from the following observation.

\begin{rec}\label{rec:colimitastotalspace} For any left fibration $p\colon\fromto{X}{S}$ (\ref{rec:leftfib}), the total space $X$ is a model for the colimit of the functor $\fromto{S}{\Kan}$ that classifies $p$ \cite[Cor. 3.3.4.6]{HTT}.
\end{rec}

If $S$ is an $\infty$-category and $\XX\colon\fromto{S}{\Wald}$ is a diagram of Waldhausen $\infty$-categories, then the colimit of the composite 
\begin{equation*}
S\ \tikz[baseline]\draw[>=stealth,->,font=\scriptsize](0,0.5ex)--node[above]{$\XX$}(0.5,0.5ex);\ \Wald\ \tikz[baseline]\draw[>=stealth,right hook->,font=\scriptsize](0,0.5ex)--node[above]{$j$}(0.5,0.5ex);\ \mathscr{P}(\Wald)
\end{equation*}
is computed objectwise \cite[Cor. 5.1.2.3]{HTT}. If $S$ is sifted, then since $\VWald\subset\mathscr{P}(\Wald)$ is stable under sifted colimits, the colimit of the composite 
\begin{equation*}
S\ \tikz[baseline]\draw[>=stealth,->,font=\scriptsize](0,0.5ex)--node[above]{$\XX$}(0.5,0.5ex);\ \Wald\ \tikz[baseline]\draw[>=stealth,right hook->,font=\scriptsize](0,0.5ex)--node[above]{$j$}(0.5,0.5ex);\ \VWald
\end{equation*}
is also computed objectwise. That is, for any compact Waldhausen $\infty$-category $\mathscr{C}$, one has
\begin{equation*}
(\colim_{s\in S}\XX(s))(\mathscr{C})\simeq\colim_{s\in S}\iota\Fun_{\Wald}(\mathscr{C},\XX(s)).
\end{equation*}
Suppose that $\XX$ classifies a Waldhausen cartesian fibration $\fromto{\mathscr{X}}{S}$; then we aim to produce a left fibration (\ref{rec:leftfib})
\begin{equation*}
\fromto{\mathrm{H}(\mathscr{C},(\mathscr{X}/S))\coloneq\iota_S\mathscr{H}(\mathscr{C},(\mathscr{X}/S))}{S}
\end{equation*}
that classifies the colimit of the composite 
\begin{equation*}
S\ \tikz[baseline]\draw[>=stealth,->,font=\scriptsize](0,0.5ex)--node[above]{$\XX$}(0.5,0.5ex);\ \Wald\ \tikz[baseline]\draw[>=stealth,right hook->,font=\scriptsize](0,0.5ex)--node[above]{$j$}(0.5,0.5ex);\ \mathscr{P}(\Wald)\ \tikz[baseline]\draw[>=stealth,->,font=\scriptsize](0,0.5ex)--node[above]{$\ev_{\mathscr{C}}$}(0.75,0.5ex);\ \Kan.
\end{equation*}
We can avoid choosing a straightening of the Waldhausen cocartesian fibration by means of the following.
\begin{cnstr}\label{ntn:geomrealinWald} Suppose $S$ a sifted $\infty$-category, and suppose $\fromto{\mathscr{X}}{S}$ a Waldhausen cocartesian fibration. Then for any compact Waldhausen $\infty$-category $\mathscr{C}$, define a simplicial set $\mathscr{H}'(\mathscr{C},(\mathscr{X}/S))$ over $S$ via the universal property
\begin{equation*}
\Mor_{S}(K,\mathscr{H}'(\mathscr{C},(\mathscr{X}/S)))\cong\Mor_{S}(\mathscr{C}\times K,\mathscr{X}),
\end{equation*}
functorial in simplicial sets $K$ over $S$. The resulting map
\begin{equation*}
\fromto{\mathscr{H}'(\mathscr{C},(\mathscr{X}/S))}{S}
\end{equation*}
is a cocartesian fibration by \ref{rec:htt32213} and \cite[Cor. 3.2.2.13]{HTT}. Denote by $\mathscr{H}(\mathscr{C},(\mathscr{X}/S))$ the full subcategory of $\mathscr{H}'(\mathscr{C},(\mathscr{X}/S))$ spanned by those functors $\fromto{\mathscr{C}}{\mathscr{X}_s}$ that are exact functors of Waldhausen $\infty$-categories; here too the canonical functor
\begin{equation*}
p\colon\fromto{\mathscr{H}(\mathscr{C},(\mathscr{X}/S))}{S}
\end{equation*}
is a cocartesian fibration. Now denote by $\mathrm{H}(\mathscr{C},(\mathscr{X}/S))$ the subcategory
\begin{equation*}
\iota_S\mathscr{H}(\mathscr{C},(\mathscr{X}/S))\subset\mathscr{H}(\mathscr{C},(\mathscr{X}/S))
\end{equation*}
consisting of the $p$-cocartesian morphisms (\ref{rec:leftfib}). The functor
\begin{equation*}
\iota_S(p)\colon\fromto{\mathrm{H}(\mathscr{C},(\mathscr{X}/S))}{S}
\end{equation*}
is now a left fibration.
\end{cnstr}

Of course, we may simply realize the assignment $\goesto{(\mathscr{C},(\mathscr{X}/S))}{\mathrm{H}(\mathscr{C},(\mathscr{X}/S))}$ as a functor
\begin{equation*}
\mathrm{H}\colon\fromto{\Wald^{\omega,\op}\times\mathbf{Wald}_{\infty/S}^{\cocart}}{\Kan}
\end{equation*}
by choosing both an equivalence $\equivto{\mathbf{Wald}_{\infty/S}^{\cocart}}{\Fun(S,\Wald)}$ and a colimit functor $\fromto{\Fun(S,\Kan)}{\Kan}$. We have given this explicit construction of the values of this functor in terms of Waldhausen cocartesian fibrations for later use.

In the meantime, since virtual Waldhausen $\infty$-categories are closed under sifted colimits in $\mathscr{P}(\Wald^{\omega})$, we have the following.
\begin{prp} If $S$ is a small sifted $\infty$-category and if $\fromto{\mathscr{X}}{S}$ is a Waldhausen cocartesian fibration in which $\mathscr{X}$ is small, then the corresponding functor $\mathrm{H}(-,(\mathscr{X}/S))\colon\fromto{\Wald^{\op}}{\Kan}$ is a virtual Waldhausen $\infty$-category.
\end{prp}

\begin{cor} If $S$ is a small sifted $\infty$-category, the functor
\begin{equation*}
\mathrm{H}\colon\fromto{\mathbf{Wald}_{\infty,/S}^{\cocart}}{\mathscr{P}(\Wald^{\omega})}
\end{equation*}
factors through the $\infty$-category of virtual Waldhausen $\infty$-categories:
\begin{equation*}
|\cdot|_S\colon\fromto{\mathbf{Wald}_{\infty/S}^{\cocart}}{\VWald}.
\end{equation*}
\end{cor}
\noindent A presheaf on $\Wald^{\omega}$ lies in the nonabelian derived $\infty$-category just in case it can be written as the geometric realization of a diagram of Ind-objects of $\Wald^{\omega}$ \cite[Lm. 5.5.8.14]{HTT}. In other words, we have the following.
\begin{cor}\label{cor:everythingsarealization} Suppose $\mathscr{X}$ a virtual Waldhausen $\infty$-category. Then there exists a Waldhausen cocartesian fibration $\fromto{\mathscr{Y}}{N\Delta^{\op}}$ and an equivalence $\mathscr{X}\simeq|\mathscr{Y}|_{N\Delta^{\op}}$.
\end{cor}

\begin{dfn}\label{dfn:realization} For any small sifted simplicial set and any Waldhausen cocartesian fibration $\mathscr{X}/S$, the virtual Waldhausen $\infty$-category $|\mathscr{X}|_S$ will be called the \textbf{\emph{realization}} of $\mathscr{X}/S$.
\end{dfn}


\part{Filtered objects and additive theories} In this part, we study reduced and finitary functors from $\Wald$ to the $\infty$-category of pointed objects of an $\infty$-topos, which we simply call \emph{theories}. We begin by studying the virtual Waldhausen $\infty$-categories of filtered and totally filtered objects of a Waldhausen $\infty$-category. Using these, we study the class of \emph{fissile} virtual Waldhausen $\infty$-categories; these form a localization of $\VWald$, and we show that suspension in this $\infty$-category is given by the formation of the virtual Waldhausen $\infty$-category of totally filtered objects, which is in turn an $\infty$-categorical analogue of Waldhausen's $S_{\bullet}$ construction. We then show that suitable excisive functors on the $\infty$-category of fissile virtual Waldhausen $\infty$-categories correspond to additive theories that satisfy the consequences of an $\infty$-categorical analogue of Waldhausen's additivity theorem, and we construct an \emph{additivization} as a Goodwillie derivative, employing our newly minted suspension functor.


\section{Filtered objects of Waldhausen $\infty$-categories}\label{sect:filt} The phenomenon behind additivity is the interaction between a filtered object and its various quotients. For example, for a category $\mathscr{C}$ with cofibrations in the sense of Waldhausen, the universal property of $K_0(\mathscr{C})$ ensures that it regards an object with a filtration of finite length
\begin{equation*}
X_0\ \tikz[baseline]\draw[>=stealth,>->](0,0.5ex)--(0.5,0.5ex);\ X_1\ \tikz[baseline]\draw[>=stealth,>->](0,0.5ex)--(0.5,0.5ex);\ \cdots\ \tikz[baseline]\draw[>=stealth,>->](0,0.5ex)--(0.5,0.5ex);\ X_n
\end{equation*}
as the sum of the first term $X_0$ and the filtered object obtained by quotienting by $X_0$:
\begin{equation*}
0\ \tikz[baseline]\draw[>=stealth,>->](0,0.5ex)--(0.5,0.5ex);\ X_1/X_0\ \tikz[baseline]\draw[>=stealth,>->](0,0.5ex)--(0.5,0.5ex);\ \cdots\ \tikz[baseline]\draw[>=stealth,>->](0,0.5ex)--(0.5,0.5ex);\ X_n/X_0,
\end{equation*}
or, by induction, as the sum of $X_0$, $X_1/X_0$, \dots, $X_n/X_{n-1}$. In order to formulate this condition properly for the entire $K$-theory \emph{space}, it is necessary to study $\infty$-categories of filtered objects in a Waldhausen $\infty$-category and the various quotient functors all as suitable inputs for algebraic $K$-theory. This is the subject of this section.

In particular, for any integer $m\geq 0$ and any Waldhausen $\infty$-category $\mathscr{C}$, we construct a Waldhausen $\infty$-category $\mathscr{F}_m(\mathscr{C})$ of filtered objects of length $m$, and we define not only the exact functors between these Waldhausen $\infty$-categories corresponding to changing the length of the filtration (given by morphisms of $\Delta$), but also sundry quotient functors. Since quotient functors are only defined up to coherent equivalences, we employ the language of Waldhausen (co)cartesian fibrations (\S \ref{sect:Waldfib}).

After we pass to suitable colimits in $\VWald$, we end up with two functors $\fromto{\VWald}{\VWald}$. The first of these, which we denote $\mathscr{F}$, is a model for the \emph{cone} in $\VWald$ (Pr. \ref{prp:FCisacone}). The second, which we will denote $\mathscr{S}$, will be a \emph{suspension}, not quite in $\VWald$, but in a suitable localization of $\VWald$ (Cor. \ref{cor:Sdotisreallysuspension}). The study of these functors is thus central to our interpretation of additive functors as excisive functors (Th. \ref{thm:additiveequiv}). 

\subsection*{The cocartesian fibration of filtered objects} Filtered objects are defined in the familiar manner.
\begin{dfn} A \textbf{\emph{filtered object of length $m$}} of a pair of $\infty$-categories $\mathscr{C}$ is a sequence of ingressive morphisms
\begin{equation*}
X_0\ \tikz[baseline]\draw[>=stealth,>->](0,0.5ex)--(0.5,0.5ex);\ X_1\ \tikz[baseline]\draw[>=stealth,>->](0,0.5ex)--(0.5,0.5ex);\ \cdots\ \tikz[baseline]\draw[>=stealth,>->](0,0.5ex)--(0.5,0.5ex);\ X_m;
\end{equation*}
that is, it is a functor of pairs $X\colon\fromto{(\Delta^m)^{\sharp}}{\mathscr{C}}$ (Ex. \ref{exm:minpairmaxpair}).
\end{dfn}

For any morphism $\eta\colon\fromto{\mathbf{m}}{\mathbf{n}}$ of $\Delta$ and any filtered object $X$ of length $n$, one may precompose $X$ with the induced functor of pairs $\fromto{(\Delta^m)^{\sharp}}{(\Delta^n)^{\sharp}}$ to obtain a filtered object $\psi^{\star}X$ of length $m$:
\begin{equation*}
X_{\eta(0)}\ \tikz[baseline]\draw[>=stealth,>->](0,0.5ex)--(0.5,0.5ex);\ X_{\eta(1)}\ \tikz[baseline]\draw[>=stealth,>->](0,0.5ex)--(0.5,0.5ex);\ \cdots\ \tikz[baseline]\draw[>=stealth,>->](0,0.5ex)--(0.5,0.5ex);\ X_{\eta(m)}.
\end{equation*}
One thus obtains a functor $\fromto{N\Delta^{\op}}{\Cat_{\infty}}$ that assigns to any object $\mathbf{m}\in\Delta$ the $\infty$-category $\Fun_{\Pair_{\infty}}((\Delta^m)^{\sharp},\mathscr{C})$. This is all simple enough.

But we will soon be forced to make things more complicated: if $\mathscr{C}$ is a Waldhausen $\infty$-category, we will below have to contemplate not only filtered objects but also \emph{totally filtered objects} in $\mathscr{C}$; these are filtered objects $X$ such that the object $X_0$ is a zero object. The $\infty$-category of totally filtered objects of length $m$ is also functorial in $\mathbf{m}$: for any morphism $\eta\colon\fromto{\mathbf{m}}{\mathbf{n}}$ of $\Delta$ and any totally filtered object $X$ of length $n$, one may still precompose $X$ with the induced functor of pairs $\fromto{(\Delta^m)^{\sharp}}{(\Delta^n)^{\sharp}}$ to obtain the filtered object $\eta^{\star}X$ of length $m$, and then one may get a \emph{totally} filtered object by forming a quotient by the object $X_{\eta(0)}$:
\begin{equation*}
0\ \tikz[baseline]\draw[>=stealth,>->](0,0.5ex)--(0.5,0.5ex);\ X_{\eta(1)}/X_{\eta(0)}\ \tikz[baseline]\draw[>=stealth,>->](0,0.5ex)--(0.5,0.5ex);\ \cdots\ \tikz[baseline]\draw[>=stealth,>->](0,0.5ex)--(0.5,0.5ex);\ X_{\eta(m)}/X_{\eta(0)}.
\end{equation*}
As we noted just before \ref{rec:cocart}, this does not specify a functor $\fromto{N\Delta^{\op}}{\Cat_{\infty}}$ on the nose, because the formation of quotients is only unique up to canonical equivalences.

This can be repaired in a variety of ways; for example, one may follow in Waldhausen's footsteps \cite[\S 1.3]{MR86m:18011} and rectify this construction by choosing all the compatible homotopy quotients at once. (For example, Lurie makes use of Waldhausen's idea in \cite[\S 1.2.2]{HA}.) But this is overkill: the theory of $\infty$-categories is precisely designed to finesse these homotopy coherence problems, and there is a genuine technical advantage in doing so. (For example, the total space of a left fibration is a ready-to-wear model for the homotopy colimit of the functor that classifies it; see \ref{rec:colimitastotalspace} or \cite[Cor. 3.3.4.6]{HTT}.) More specifically, the theory of cocartesian fibrations allows us to work effectively with this construction without solving homotopy coherence problems like this.

To that end, let's first use Pr. \ref{prp:htt32213} to access the cocartesian fibration
\begin{equation*}
\fromto{\mathscr{F}(\mathscr{C})}{N\Delta^{\op}}
\end{equation*}
classified by the functor
\begin{equation*}
\goesto{\mathbf{m}}{\Fun_{\Pair_{\infty}}((\Delta^m)^{\sharp},\mathscr{C})}.
\end{equation*}
At no extra cost, for any Waldhausen cocartesian fibration $\fromto{\mathscr{X}}{S}$ classified by a functor $\XX\colon\fromto{S}{\Wald}$, we can actually write down the cocartesian fibration
\begin{equation*}
\fromto{\mathscr{F}(\mathscr{X}/S)}{N\Delta^{\op}\times S}
\end{equation*}
classified by the functor
\begin{equation*}
\goesto{(\mathbf{m},s)}{\Fun_{\Pair_{\infty}}((\Delta^m)^{\sharp},\XX(s))}.
\end{equation*}
Once this has been done, we'll be in a better position to define a Waldhausen cocartesian fibration of totally filtered objects.

The first step to using Pr. \ref{prp:htt32213} is to identify the pair cartesian fibration (Df. \ref{dfn:paircartfib}) that is classified by the functor $\goesto{\mathbf{m}}{(\Delta^m)^{\sharp}}$.

\begin{ntn}\label{ntn:M} Denote by $\mathrm{M}$ the ordinary category whose objects are pairs $(\mathbf{m},i)$ consisting of an object $\mathbf{m}\in\Delta$ and an element $i\in\mathbf{m}$ and whose morphisms $\fromto{(\mathbf{n},j)}{(\mathbf{m},i)}$ are maps $\phi\colon\fromto{\mathbf{m}}{\mathbf{n}}$ of $\Delta$ such that $j\leq\phi(i)$. This category comes equipped with a natural projection $\fromto{\mathrm{M}}{\Delta^{\op}}$.

It is easy to see that the projection $\fromto{\mathrm{M}}{\Delta^{\op}}$ is a Grothendieck fibration, and so the projection $\pi\colon\fromto{N\mathrm{M}}{N\Delta^{\op}}$ is a cartesian fibration. In fact, the category $\mathrm{M}$ is nothing more than the Grothendieck construction applied to the natural inclusion $\into{(\Delta^{\op})^{\op}\cong\Delta}{\Cat}$. So the functor $\fromto{\Delta}{\Cat_{\infty}}$ that classifies $\pi$ is given by the assignment $\goesto{\mathbf{m}}{\Delta^m}$.

The nerve $N\mathrm{M}$ can be endowed with a pair structure by setting
\begin{equation*}
(N\mathrm{M})_{\dag}\coloneq N\mathrm{M}\times_{N\Delta^{\op}}\iota N\Delta^{\op}.
\end{equation*}
Put differently, an edge of $\mathrm{M}$ is ingressive just in case it covers an equivalence of $\Delta^{\op}$.  Consequently, $\pi$ is automatically a pair cartesian fibration (Df. \ref{dfn:paircartfib}); the functor $\fromto{N\Delta}{\Pair_{\infty}}$ classified by $\pi$ is given by the assignment $\goesto{\mathbf{m}}{(\Delta^m)^{\sharp}}$.
\end{ntn}
\noindent Now it is no problem to use the technology from Pr. \ref{prp:htt32213} to define the cocartesian fibration $\fromto{\mathscr{F}(\mathscr{X}/S)}{N\Delta^{\op}\times S}$ that we seek.
\begin{cnstr} For any pair cocartesian fibration $\fromto{\mathscr{X}}{S}$, define a map $\fromto{\mathscr{F}(\mathscr{X}/S)}{N\Delta^{\op}\times S}$, using the notation of Pr. \ref{prp:htt32213} and Ex. \ref{exm:minpairmaxpair}, as
\begin{equation*}
\mathscr{F}(\mathscr{X}/S)\coloneq T_{\pi\times\id_S}((N\Delta^{\op})^{\flat}\times\mathscr{X}).
\end{equation*}
Equivalently, we require, for any simplicial set $K$ and any map $\sigma\colon\fromto{K}{N\Delta^{\op}\times S}$, a bijection between the set $\Mor_{N\Delta^{\op}\times S}(K,\mathscr{F}(\mathscr{X}/S))$ and the set
\begin{equation*}
\Mor_{s\Set(2)/(S,\iota S)}((K\times_{N\Delta^{\op}}N\mathrm{M},K\times_{N\Delta^{\op}}(N\mathrm{M})_{\dag}),(\mathscr{X},\mathscr{X}_{\dag}))
\end{equation*}
(Nt. \ref{ntn:ordcatofpairs}), functorial in $\sigma$.
\end{cnstr}

With this definition, Pr. \ref{prp:htt32213} now implies the following.
\begin{prp} Suppose $p\colon\fromto{\mathscr{X}}{S}$ a pair cocartesian fibration. Then the functor
\begin{equation*}
\fromto{\mathscr{F}(\mathscr{X}/S)}{N\Delta^{\op}\times S}
\end{equation*}
is a cocartesian fibration.
\end{prp}
Furthermore, the functor $\fromto{N\Delta^{\op}\times S}{\Cat_{\infty}}$ that classifies the cocartesian fibration $\fromto{\mathscr{F}(\mathscr{X}/S)}{N\Delta^{\op}\times S}$ is indeed the functor
\begin{equation*}
\goesto{(\mathbf{m},s)}{\Fun_{\Pair_{\infty}}((\Delta^m)^{\sharp},\XX(s))},
\end{equation*}
where $\XX\colon\fromto{S}{\Pair_{\infty}}$ is the functor that classifies $p$.

\begin{ntn} When $S=\Delta^0$, write $\mathscr{F}(\mathscr{C})$ for $\mathscr{F}(\mathscr{C}/S)$, and for any integer $m\geq 0$, write $\mathscr{F}_m(\mathscr{C})$ for the fiber $\Fun_{\Pair_{\infty}}((\Delta^m)^{\sharp},\mathscr{C})$ of the cocartesian fibration $\fromto{\mathscr{F}(\mathscr{C})}{N\Delta^{\op}}$ over $\mathbf{m}$.

Hence for any Waldhausen cocartesian fibration $\fromto{\mathscr{X}}{S}$, the fiber of the cocartesian fibration $\fromto{\mathscr{F}(\mathscr{X}/S)}{N\Delta^{\op}\times S}$ over a vertex $(\mathbf{m},s)$ is the Waldhausen $\infty$-category $\mathscr{F}_m(\mathscr{X}_s)$.
\end{ntn}


\subsection*{A Waldhausen structure on filtered objects of a Waldhausen $\infty$-category} We may endow the $\infty$-categories $\mathscr{F}(\mathscr{X}/S)$ of filtered objects with a pair structure in a variety of ways, but we wish to focus on one pair structure that will retain good formal properties when we pass to quotients.

More specifically, suppose $\mathscr{C}$ a Waldhausen $\infty$-category. A morphism $f\colon\fromto{X}{Y}$ of $\mathscr{F}_m(\mathscr{C})$ can be represented as a diagram
\begin{equation*}
\begin{tikzpicture} 
\matrix(m)[matrix of math nodes, 
row sep=4ex, column sep=4ex, 
text height=1.5ex, text depth=0.25ex] 
{X_0&X_1&\cdots&X_m\\ 
Y_0&Y_1&\cdots&Y_m.\\}; 
\path[>=stealth,->,font=\scriptsize] 
(m-1-1) edge[>->] (m-1-2) 
edge (m-2-1) 
(m-1-2) edge[>->] (m-1-3)
edge (m-2-2)
(m-1-3) edge[>->] (m-1-4)
(m-1-4) edge (m-2-4)
(m-2-1) edge[>->] (m-2-2)
(m-2-2) edge[>->] (m-2-3)
(m-2-3) edge[>->] (m-2-4); 
\end{tikzpicture}
\end{equation*}
What should it mean to say that $f$ is ingressive? It is natural to demand, first and foremost, that each morphism $f_i\colon\fromto{X_i}{Y_i}$ is ingressive, but this will not be enough to ensure that the morphisms $\fromto{X_j/X_i}{Y_j/Y_i}$ are all ingressive. Guaranteeing this turns out to be equivalent to the claim that in each of the squares
\begin{equation*}
\begin{tikzpicture} 
\matrix(m)[matrix of math nodes, 
row sep=4ex, column sep=4ex, 
text height=1.5ex, text depth=0.25ex] 
{X_i&X_j\\ 
Y_i&Y_j,\\}; 
\path[>=stealth,->,font=\scriptsize] 
(m-1-1) edge[>->] (m-1-2) 
edge[>->] (m-2-1) 
(m-1-2) edge[>->] (m-2-2) 
(m-2-1) edge[>->] (m-2-2); 
\end{tikzpicture}
\end{equation*}
the morphism from the pushout $X_j\cup^{X_i}Y_i$ to $Y_j$ is a cofibration as well. This was noted by Waldhausen \cite[Lm. 1.1.2]{MR86m:18011}.

Our approach is thus to define a pair structure in such a concrete manner on $\mathscr{F}_1(\mathscr{C})$, and then to declare that a morphism $f$ of $\mathscr{F}_m(\mathscr{C})$ is ingressive just in case $\eta^{\star}(f)$ is so for any $\eta\colon\fromto{\Delta^1}{\Delta^m}$.

\begin{dfn}\label{dfn:injpairstruct} Suppose $\mathscr{C}$ a Waldhausen $\infty$-category. We now endow the $\infty$-category $\mathscr{F}_1(\mathscr{C})$ with a pair structure by letting $\mathscr{F}_1(\mathscr{C})_{\dag}\subset\mathscr{F}_1(\mathscr{C})$ be the smallest subcategory containing the following classes of edges of $\mathscr{C}$:
\begin{enumerate}[(\ref{dfn:injpairstruct}.1)]
\item\label{item:degensourceingresstarg} any edge $\fromto{X}{Y}$ represented as a square
\begin{equation*}
\begin{tikzpicture} 
\matrix(m)[matrix of math nodes, 
row sep=4ex, column sep=4ex, 
text height=1.5ex, text depth=0.25ex] 
{X_0&X_1\\ 
Y_0&Y_1\\}; 
\path[>=stealth,->,font=\scriptsize] 
(m-1-1) edge[>->] (m-1-2) 
edge[inner sep=0.5pt] node[left]{$\sim$} (m-2-1) 
(m-1-2) edge[>->] (m-2-2) 
(m-2-1) edge[>->] (m-2-2); 
\end{tikzpicture}
\end{equation*}
in which $\equivto{X_0}{Y_0}$ is an equivalence and $\cofto{X_1}{Y_1}$ is ingressive, and
\item\label{item:cocartedgeovercof} any edge $\fromto{X}{Y}$ represented as a pushout square
\begin{equation*}
\begin{tikzpicture} 
\matrix(m)[matrix of math nodes, 
row sep=4ex, column sep=4ex, 
text height=1.5ex, text depth=0.25ex] 
{X_0&X_1\\ 
Y_0&Y_1\\}; 
\path[>=stealth,->,font=\scriptsize] 
(m-1-1) edge[>->] (m-1-2) 
edge[>->] (m-2-1) 
(m-1-2) edge[>->] (m-2-2) 
(m-2-1) edge[>->] (m-2-2);
\end{tikzpicture}
\end{equation*}
in which $\cofto{X_0}{Y_0}$ and thus also $\cofto{X_1}{Y_1}$ are ingressive.
\end{enumerate}
\end{dfn}

Let's compare this definition to our more concrete one outlined above it. To this end, we need a bit of notation.
\begin{ntn} Let us denote by $\mathscr{R}$ the pair of $\infty$-categories whose underlying $\infty$-category is
\begin{equation*}
(\Delta^1\times(\Lambda^2_0)^{\rhd})/(\Delta^1\times\Lambda^2_0),
\end{equation*}
which may be drawn
\begin{equation*}
\begin{tikzpicture} 
\matrix(m)[matrix of math nodes, 
row sep=5ex, column sep=5ex, 
text height=1.5ex, text depth=0.25ex] 
{0&1&[-3ex]\\ 
2&\infty'&\\
[-3ex]&&\infty,\\}; 
\path[>=stealth,->,font=\scriptsize] 
(m-1-1) edge (m-1-2) 
edge (m-2-1) 
(m-1-2) edge (m-2-2)
edge (m-3-3)
(m-2-1) edge (m-2-2)
edge (m-3-3)
(m-2-2) edge (m-3-3); 
\end{tikzpicture}
\end{equation*}
in which only the edges $\cofto{0}{1}$, $\cofto{2}{\infty}$ and $\cofto{2}{\infty'}$ are ingressive. In the notation of \ref{exm:LambdaDelta}, there is an obvious strict inclusion of pairs $\into{\Lambda_0\mathscr{Q}^2}{\mathscr{R}}$, and there are two strict inclusions of pairs
\begin{equation*}
\into{\mathscr{Q}^2\cong\mathscr{Q}^2\times\Delta^{\{0\}}}{\mathscr{R}}\textrm{\quad and\quad}\into{\mathscr{Q}^2\cong\mathscr{Q}^2\times\Delta^{\{1\}}}{\mathscr{R}}.
\end{equation*} 
\end{ntn}

\begin{lem}\label{lem:concreteingressiveF1C} Suppose $\mathscr{C}$ a Waldhausen $\infty$-category. Then a morphism $f\colon\fromto{X}{Y}$ of $\mathscr{F}_1(\mathscr{C})$ is ingressive just in case the morphism $\fromto{X_0}{Y_0}$ is ingressive and the corresponding square
\begin{equation*}
F\colon\fromto{\mathscr{Q}^2\cong(\Delta^1)^{\flat}\times(\Delta^1)^{\sharp}}{\mathscr{C}}
\end{equation*}
has the property that for any diagram $\overline{F}\colon\fromto{\mathscr{R}}{\mathscr{C}}$ such that $\overline{F}|_{\mathscr{Q}^2\times\Delta^{\{0\}}}$ is a pushout square, and $F=\overline{F}|_{\mathscr{Q}^2\times\Delta^{\{1\}}}$, the edge $\fromto{\overline{F}(\infty')}{\overline{F}(\infty)}$ is ingressive.
\begin{proof} An easy argument shows that morphisms with this property form a subcategory of $\mathscr{F}_1(\mathscr{C})$, and it is clear that morphisms either of type (\ref{dfn:injpairstruct}.\ref{item:degensourceingresstarg}) or of type (\ref{dfn:injpairstruct}.\ref{item:cocartedgeovercof}) enjoy this property. Consequently, every ingressive morphism enjoys this property. On the other hand, a morphism $\fromto{X}{Y}$ that enjoys this property can clearly be factored as $X\ \tikz[baseline]\draw[>=stealth,->](0,0.5ex)--(0.5,0.5ex);\ Y'\ \tikz[baseline]\draw[>=stealth,->](0,0.5ex)--(0.5,0.5ex);\ Y$, where $\fromto{X}{Y'}$ is of type (\ref{dfn:injpairstruct}.\ref{item:cocartedgeovercof}), and $\fromto{Y'}{Y}$ is of type (\ref{dfn:injpairstruct}.\ref{item:degensourceingresstarg}), viz.:
\begin{equation*}
\begin{tikzpicture} 
\matrix(m)[matrix of math nodes, 
row sep=4ex, column sep=4ex, 
text height=1.5ex, text depth=0.25ex] 
{X_0&X_1\\
Y_0&Y_{01}\\
Y_0&Y_1,\\}; 
\path[>=stealth,->,font=\scriptsize] 
(m-1-1) edge[>->] (m-1-2) 
edge[>->] (m-2-1) 
(m-1-2) edge[>->] (m-2-2) 
(m-2-1) edge[>->] (m-2-2)
edge[-,double distance=1.5pt] (m-3-1)
(m-2-2) edge[>->] (m-3-2)
(m-3-1) edge[>->] (m-3-2); 
\end{tikzpicture}
\end{equation*}
where the top square is a pushout square and $\cofto{Y_{01}}{Y_1}$ is ingressive.
\end{proof}
\end{lem}

\begin{dfn}\label{dfn:ingressinFmC} Now suppose $\fromto{\mathscr{X}}{S}$ a Waldhausen cocartesian fibration. We endow the $\infty$-category $\mathscr{F}(\mathscr{X}/S)$ with the following pair structure. Let $\mathscr{F}(\mathscr{X}/S)_{\dag}\subset\mathscr{F}(\mathscr{X}/S)$ be the smallest pair structure containing any edge $f\colon\fromto{\Delta^1}{\mathscr{F}(\mathscr{X}/S)}$ covering a degenerate edge $\id_{(\mathbf{m},s)}$ of $N\Delta^{\op}\times S$ such that for any edge $\eta\colon\fromto{\Delta^1}{\Delta^m}$, the edge
\begin{equation*}
\Delta^1\ \tikz[baseline]\draw[>=stealth,->,font=\scriptsize](0,0.5ex)--node[above]{$f$}(0.5,0.5ex);\mathscr{F}_m(\mathscr{X}_s)\ \tikz[baseline]\draw[>=stealth,->,font=\scriptsize](0,0.5ex)--node[above]{$\eta^{\star}$}(0.5,0.5ex);\ \mathscr{F}_1(\mathscr{X}_s)
\end{equation*}
is ingressive in the sense of Df. \ref{dfn:injpairstruct}.
\end{dfn}

\begin{lem}\label{lem:checkingressonneighbors} Suppose $\mathscr{C}$ a Waldhausen $\infty$-category. Then a morphism $f\colon\fromto{X}{Y}$ of $\mathscr{F}_m(\mathscr{C})$ is ingressive just in case, for any integer $1\leq i\leq m$, the restricted morphism $\fromto{X|_{(\Delta^{\{i-1,i\}})^{\sharp}}}{Y|_{(\Delta^{\{i-1,i\}})^{\sharp}}}$ is ingressive in $\mathscr{F}_{\{i-1,i\}}(\mathscr{X}_s)$.
\begin{proof} Suppose $f$ satisfies this condition. It is immediate that every morphism $\cofto{X_i}{Y_i}$ is ingressive, so we can regard $f$ as an $m$-simplex $\sigma\colon\fromto{\Delta^m}{\mathscr{F}_1(\mathscr{C})}$. By Lm. \ref{lem:concreteingressiveF1C}, this condition is equivalent to the condition that each edge $\sigma|_{\Delta^{i-1,i}}$ is ingressive, and since ingressive edges are closed under composition, it follows that every edge $\sigma|_{\Delta^{\{i,j\}}}$ is ingressive.
\end{proof}
\end{lem}

\begin{prp}\label{prp:FXSisWaldcocart} Suppose $p\colon\fromto{\mathscr{X}}{S}$ a Waldhausen cocartesian fibration. Then with the pair structure of \textup{Df. \ref{dfn:injpairstruct}}, the functor
\begin{equation*}
\fromto{\mathscr{F}(\mathscr{X}/S)}{N\Delta^{\op}\times S}
\end{equation*}
is a Waldhausen cocartesian fibration.
\begin{proof} It is easy to see that $\fromto{\mathscr{F}(\mathscr{X}/S)}{N\Delta^{\op}\times S}$ a pair cocartesian fibration.

We claim that for any vertex $(\mathbf{m},s)\in N\Delta^{\op}\times S$, the pair $\mathscr{F}_m(\mathscr{X}_s)$ is a Waldhausen $\infty$-category. Note that since $\mathscr{X}_s$ admits a zero object, so does $\mathscr{F}_m(\mathscr{X}_s)$. For the remaining two axioms, one reduces immediately to the case where $m=1$. Then (\ref{dfn:preWald}.\ref{item:0toxingressive}) follows from the presence of (\ref{dfn:injpairstruct}.\ref{item:degensourceingresstarg}) among ingressive morphisms. To prove (\ref{dfn:preWald}.\ref{item:pushcof}), one may note that cofibations of $\mathscr{F}_1(\mathscr{X}_s)$ are in particular ingressive morphisms of $\mathscr{O}(\mathscr{C})$, for which the existence of pushouts is clear. Finally, to prove (\ref{dfn:preWald}.\ref{item:pushcofcof}), it suffices to see that a pushout of any edge of either of the classes (\ref{dfn:injpairstruct}.\ref{item:degensourceingresstarg}) or (\ref{dfn:injpairstruct}.\ref{item:cocartedgeovercof}) is of the same class. For the class (\ref{dfn:injpairstruct}.\ref{item:degensourceingresstarg}), this follows from the fact that pushouts in $\mathscr{F}_1(\mathscr{X}_s)$ are computed pointwise. A pushout of a morphism of the class (\ref{dfn:injpairstruct}.\ref{item:cocartedgeovercof}) is a cube
\begin{equation*}
X\colon\fromto{(\Delta^1)^{\flat}\times(\Delta^1)^{\sharp}\times(\Delta^1)^{\sharp}}{\mathscr{X}_s}
\end{equation*}
in which the faces
\begin{equation*}
X|_{\Delta^{\{0\}}\times(\Delta^1)^{\sharp}\times(\Delta^1)^{\sharp}},\quad X|_{(\Delta^1)^{\flat}\times\Delta^{\{0\}}\times(\Delta^1)^{\sharp}},\textrm{\quad and\quad}X|_{(\Delta^1)^{\flat}\times\Delta^{\{1\}}\times(\Delta^1)^{\sharp}}
\end{equation*}
are all pushouts. If $X$ is represented by the commutative diagram
\begin{equation*}
\begin{tikzpicture}[cross line/.style={preaction={draw=white, -, 
line width=6pt}}]
\matrix(m)[matrix of math nodes, 
row sep=2ex, column sep=0.75ex, 
text height=1.5ex, text depth=0.25ex]
{&X_{100}&&X_{101}\\
X_{000}&&X_{001}&\\
&X_{110}&&X_{111}\\
X_{010}&&X_{011},&\\
};
\path[>=stealth,->,font=\scriptsize]
(m-1-2) edge[<-] (m-2-1)
edge[>->] (m-3-2)
edge[>->] (m-1-4)
(m-3-2) edge[<-] (m-4-1)
edge[>->] (m-3-4)
(m-2-1) edge[cross line,>->] (m-2-3)
edge[>->] (m-4-1)
(m-1-4) edge[<-] (m-2-3)
edge[>->](m-3-4)
(m-4-1) edge[>->] (m-4-3)
(m-3-4) edge[<-] (m-4-3)
(m-2-3) edge[cross line,>->] (m-4-3);
\end{tikzpicture}
\end{equation*}
then the front face, the top face, and the bottom face are all pushouts. By Quetzalcoatl (e.g., by \cite[Lm. 4.4.2.1]{HTT}), the back face $X|_{\Delta^{\{1\}}\times(\Delta^1)^{\sharp}\times(\Delta^1)^{\sharp}}$ must be a pushout as well; this is precisely the claim that the pushout is of the class (\ref{dfn:injpairstruct}.\ref{item:cocartedgeovercof}).

For any $\mathbf{m}\in\Delta$ and any edge $f\colon\fromto{s}{t}$ of $S$, since the functor $f_{\mathscr{X},!}\colon\fromto{\mathscr{X}_s}{\mathscr{X}_t}$ is exact, it follows directly that the functor
\begin{equation*}
f_{\mathscr{F},!}\colon\fromto{\mathscr{F}_m(\mathscr{X}_s)}{\mathscr{F}_m(\mathscr{X}_t)}
\end{equation*}
is exact as well. Now for any fixed vertex $s\in S_0$ and any simplicial operator $\phi\colon\fromto{\mathbf{n}}{\mathbf{m}}$ of $\Delta$, the functor
\begin{equation*}
\phi_{\mathscr{F},!}\colon\fromto{\mathscr{F}_m(\mathscr{X}_s)}{\mathscr{F}_n(\mathscr{X}_s)}
\end{equation*}
visibly carries ingressive morphisms to ingressive morphisms, and it preserves zero objects as well as any pushouts that exist, since limits and colimits are formed pointwise. 
\end{proof}
\end{prp}

Thanks to \ref{nul:Tpafunctor}, we have:
\begin{cor} The assignment $\goesto{(\mathscr{X}/S)}{\mathscr{F}(\mathscr{X}/S)}$ defines a functor
\begin{equation*}
\mathscr{F}\colon\fromto{\Wald^{\cocart}}{\Wald^{\cocart}}
\end{equation*}
covering the endofunctor $\goesto{S}{N\Delta^{\op}\times S}$ of $\Cat_{\infty}$.
\end{cor}

\begin{prp}\label{prp:Fisacategory} Suppose $\fromto{\mathscr{X}}{S}$ a Waldhausen cocartesian fibration, and suppose
\begin{equation*}
\FF_{\ast}(\mathscr{X}/S)\colon\fromto{N\Delta^{\op}}{\Fun(S,\Wald)}
\end{equation*}
a functor that classifies the Waldhausen cocartesian fibration
\begin{equation*}
\fromto{\mathscr{F}(\mathscr{X}/S)}{N\Delta^{\op}\times S}.
\end{equation*}
Then $\FF_{\ast}(\mathscr{X}/S)$ is a \textbf{\emph{category object}} \cite[Df. 1.1.1]{G}; that is, the morphisms of $\Delta$ of the form $\into{\{i-1,i\}}{\mathbf{m}}$ induce morphisms that exhibit $\FF_m(\mathscr{X}/S)$ as the limit in $\Fun(S,\Wald)$ of the diagram
\begin{equation*}
\begin{tikzpicture} 
\matrix(m)[matrix of math nodes, 
row sep=4ex, column sep=-6ex, 
text height=1.5ex, text depth=0.25ex] 
{\FF_{\{0,1\}}(\mathscr{X}/S)&&\FF_{\{1,2\}}(\mathscr{X}/S)&[6ex]&[6ex]\FF_{\{m-2,m-1\}}(\mathscr{X}/S)&&\FF_{\{m-1,m\}}(\mathscr{X}/S).\\
&\FF_{\{1\}}(\mathscr{X}/S)&&\cdots&&\FF_{\{m-1\}}(\mathscr{X}/S)&\\}; 
\path[>=stealth,->,font=\scriptsize] 
(m-1-1) edge (m-2-2) 
(m-1-3) edge (m-2-2)
edge (m-2-4) 
(m-1-5) edge (m-2-4)
edge (m-2-6)
(m-1-7) edge (m-2-6); 
\end{tikzpicture}
\end{equation*}
\begin{proof} Since limits in $\Fun(S,\Wald)$ are computed objectwise, it suffices to assume that $S=\Delta^0$. It is easy to see that $(\Delta^m)^{\sharp}$ decomposes in $\Pair_{\infty}$ as the pushout of the diagram
\begin{equation*}
\begin{tikzpicture} 
\matrix(m)[matrix of math nodes, 
row sep=4ex, column sep=0ex, 
text height=1.5ex, text depth=0.25ex] 
{&(\Delta^{\{1\}})^{\sharp}&&[2ex]\cdots&&(\Delta^{\{m-1\}})^{\sharp}&&\\ 
(\Delta^{\{0,1\}})^{\sharp}&&(\Delta^{\{1,2\}})^{\sharp}&&(\Delta^{\{m-2,m-1\}})^{\sharp}&&(\Delta^{\{m-1,m\}})^{\sharp},\\}; 
\path[>=stealth,->,font=\scriptsize] 
(m-1-2) edge (m-2-1)
edge (m-2-3) 
(m-1-4) edge (m-2-3)
edge (m-2-5) 
(m-1-6) edge (m-2-5)
edge (m-2-7); 
\end{tikzpicture}
\end{equation*}
since the analogous statement is true in $\Cat_{\infty}$. Thus $\FF_m(\mathscr{X})$ is the desired limit in $\Cat_{\infty}$, and it follows immediately from Lm. \ref{lem:checkingressonneighbors} that $\FF_m(\mathscr{X})$ is the desired limit in the $\infty$-category $\Pair_{\infty}$ and thus also in the $\infty$-category $\Wald$.
\end{proof}
\end{prp}


\subsection*{Totally filtered objects} Now we are in a good position to study the functoriality of filtered objects $X$ that are \emph{separated} in the sense that $X_0$ is a zero object. We call these \emph{totally filtered} objects.
\begin{dfn}\label{dfn:totallyfilteredobject} Suppose $\mathscr{C}$ a Waldhausen $\infty$-category. Then a filtered object $X\colon\fromto{(\Delta^m)^{\sharp}}{\mathscr{C}}$ will be said to be \textbf{\emph{totally filtered}} if $X_0$ is a zero object.
\end{dfn}

\begin{ntn} Suppose $p\colon\fromto{\mathscr{X}}{S}$ a Waldhausen cocartesian fibration. Denote by $\mathscr{S}(\mathscr{X}/S)$ the full subpair (\ref{dfn:pair}.\ref{item:subpair}) of $\mathscr{F}(\mathscr{X}/S)$ spanned by those filtered objects $X$ such that $X$ is totally filtered in $\mathscr{X}_{p(X)}$. When $S=\Delta^0$, write $\mathscr{S}(\mathscr{X})$ for $\mathscr{S}(\mathscr{X}/S)$, and for any integer $m\geq 0$, write $\mathscr{S}_m(\mathscr{X})$ for the fiber of $\fromto{\mathscr{S}(\mathscr{X})}{N\Delta^{\op}}$ over the object $\mathbf{m}\in N\Delta^{\op}$.
\end{ntn}

\begin{prp}\label{prp:FmisS1plusm} Suppose $\mathscr{C}$ a Waldhausen $\infty$-category. For any integer $m\geq0$, the $0$-th face map defines an equivalence of $\infty$-categories
\begin{equation*}
\equivto{\mathscr{S}_{1+m}(\mathscr{C})}{\mathscr{F}_m(\mathscr{C})},
\end{equation*}
and the map $\fromto{\mathscr{S}_0(\mathscr{C})}{\Delta^0}$ is an equivalence.
\begin{proof} It follows from Joyal's theorem \cite[Pr. 1.2.12.9]{HTT} that the natural functor $\fromto{\mathscr{S}_{1+m}(\mathscr{C})}{\mathscr{F}_m(\mathscr{C})}$ is a left fibration whose fibers are contractible Kan complexes --- hence a trivial fibration.
\end{proof}
\end{prp}

As $\mathbf{m}$ varies, the functoriality of $\mathscr{S}_m(\mathscr{X})$ is, as we have observed, traditionally a matter of some consternation, as the functors involve various (homotopy) quotients, which are not uniquely defined on the nose. We all share the intuition that the uniqueness of these quotients is good enough for all practical purposes and that the coherence issues that appear to arise are mere technical issues. The theory of cocartesian fibrations allows us to make this intuition honest.

Below (Th. \ref{thm:SXiscocartfib}), we'll show that for any Waldhausen cocartesian fibration $\fromto{\mathscr{X}}{S}$, the functor $\fromto{\mathscr{S}(\mathscr{X}/S)}{N\Delta^{\op}\times S}$ is a Waldhausen cocartesian fibration. Let's reflect on what this means when $S=\Delta^0$; in this case, $\mathscr{X}$ is just a Waldhausen $\infty$-category. An edge $\fromto{X}{Y}$ of $\mathscr{S}(\mathscr{X})$ that covers an edge given by a morphism $\eta\colon\fromto{\mathbf{m}}{\mathbf{n}}$ of $\Delta$ is by definition a commutative diagram
\begin{equation*}
\begin{tikzpicture} 
\matrix(m)[matrix of math nodes, 
row sep=4ex, column sep=4ex, 
text height=1.5ex, text depth=0.25ex] 
{X_{\eta(0)}&X_{\eta(1)}&\cdots&X_{\eta(m)}\\ 
0&Y_{1}&\cdots&Y_{m}.\\}; 
\path[>=stealth,->,font=\scriptsize] 
(m-1-1) edge[>->] (m-1-2) 
edge (m-2-1) 
(m-1-2) edge[>->] (m-1-3)
edge (m-2-2)
(m-1-3) edge[>->] (m-1-4)
(m-1-4) edge (m-2-4)
(m-2-1) edge (m-2-2)
(m-2-2) edge (m-2-3)
(m-2-3) edge (m-2-4); 
\end{tikzpicture}
\end{equation*}
To say that $\fromto{X}{Y}$ is a cocartesian edge over $\eta$ is to say that $Y$ is initial among totally filtered objects under $\eta^{\star}X$. This is equivalent to the demand that each of the squares above must be pushout squares, i.e., that $Y_k\simeq X_{\eta(k)}/X_{\eta(0)}$. So if $\fromto{\mathscr{S}(\mathscr{C})}{N\Delta^{\op}}$ is a Waldhausen cocartesian fibration, then the functor $\SSS_{\ast}\colon\fromto{N\Delta^{\op}}{\Wald}$ that classifies it works exactly as Waldhausen's $S_{\bullet}$ construction: it carries an object $\mathbf{m}\in N\Delta^{\op}$ to the Waldhausen $\infty$-category $\SSS_m(\mathscr{C})$ of totally filtered objects of length, and it carries a morphism $\eta\colon\fromto{\mathbf{m}}{\mathbf{n}}$ of $\Delta$ to the exact functor $\fromto{\SSS_n(\mathscr{C})}{\SSS_m(\mathscr{C})}$ given by
\begin{equation*}
\goesto{[X_0\ \tikz[baseline]\draw[>=stealth,>->](0,0.5ex)--(0.5,0.5ex);\ X_1\ \tikz[baseline]\draw[>=stealth,>->](0,0.5ex)--(0.5,0.5ex);\ \cdots\ \tikz[baseline]\draw[>=stealth,>->](0,0.5ex)--(0.5,0.5ex);\ X_n]}{[0\ \tikz[baseline]\draw[>=stealth,>->](0,0.5ex)--(0.5,0.5ex);\ X_{\eta(1)}/X_{\eta(0)}\ \tikz[baseline]\draw[>=stealth,>->](0,0.5ex)--(0.5,0.5ex);\ \cdots\ \tikz[baseline]\draw[>=stealth,>->](0,0.5ex)--(0.5,0.5ex);\ X_{\eta(m)}/X_{\eta(0)}].}
\end{equation*}
In other words, the \emph{data} of the $\infty$-categorical $S_{\bullet}$ construction is already before us; we just need to confirm that it works as desired.  

To prove Th. \ref{thm:SXiscocartfib}, it turns out to be convenient to study the ``mapping cylinder'' $\mathscr{M}(\mathscr{X}/S)$ of the inclusion functor $\into{\mathscr{S}(\mathscr{X}/S)}{\mathscr{F}(\mathscr{X}/S)}$. We will discover that this inclusion admits a left adjoint, and then we will use this left adjoint to complete the proof of Th. \ref{thm:SXiscocartfib}.

\begin{ntn} For any Waldhausen cocartesian fibration $\fromto{\mathscr{X}}{S}$, let us write $\mathscr{M}(\mathscr{X}/S)$ for the full subcategory of $\Delta^1\times\mathscr{F}(\mathscr{X}/S)$ spanned by those pairs $(i,X)$ such that $X$ is totally filtered if $i=1$. This $\infty$-category comes equipped with an inner fibration
\begin{equation*}
\fromto{\mathscr{M}(\mathscr{X}/S)}{\Delta^1\times N\Delta^{\op}\times S}.
\end{equation*}
Define a pair structure on $\mathscr{M}(\mathscr{X}/S)$ so that it is a subpair of $(\Delta^1)^{\flat}\times\mathscr{F}(\mathscr{X}/S)$; that is, let $\mathscr{M}(\mathscr{X}/S)_{\dag}\subset\mathscr{M}(\mathscr{X}/S)$ be the subcategory whose edges are maps $\fromto{(i,X)}{(j,Y)}$ such that $i=j$ and $\fromto{X}{Y}$ is an ingressive morphism of $\mathscr{F}(\mathscr{X}/S)$.
\end{ntn}

Our first lemma is obvious by construction.
\begin{lem} For any Waldhausen cocartesian fibration $\fromto{\mathscr{X}}{S}$, the natural projection $\fromto{\mathscr{M}(\mathscr{X}/S)}{\Delta^1}$ is a pair cartesian fibration.
\end{lem}
\noindent Our next lemma, however, is subtler.
\begin{lem}\label{lem:MXScocartfib} For any Waldhausen cocartesian fibration $\fromto{\mathscr{X}}{S}$, the natural projection $\fromto{\mathscr{M}(\mathscr{X}/S)}{\Delta^1}$ is a pair cocartesian fibration.
\begin{proof} By \cite[2.4.1.3(3)]{HTT}, it suffices to show that for any vertex $(\mathbf{m},s)\in(N\Delta^{\op}\times S)_0$, the inner fibration
\begin{equation*}
q\colon\fromto{\mathscr{M}_m(\mathscr{X}_s)}{\Delta^1}
\end{equation*}
is a pair cocartesian fibration. Note that an edge $\fromto{X}{Y}$ of $\mathscr{M}_m(\mathscr{X}_s)$ covering the nondegenerate edge $\sigma$ of $\Delta^1$ is $q$-cocartesian if and only if it is an initial object of the fiber $\mathscr{M}_m(\mathscr{X}_s)_{X/}\times_{\Delta^1_{0/}}\{\sigma\}$. If $m=0$, then the map
\begin{equation*}
\fromto{\mathscr{M}_0(\mathscr{X}_s)_{X/}}{\Delta^1_{0/}}
\end{equation*}
is a trivial fibration \cite[Pr. 1.2.12.9]{HTT}, so the fiber over $\sigma$ is a contractible Kan complex. Let us now induct on $m$; assume that $m>0$ and that the functor $p\colon\fromto{\mathscr{M}_{m-1}(\mathscr{X}_s)}{\Delta^1}$ is a cocartesian fibration. It is easy to see that the inclusion $\into{\{0,1,\dots,m-1\}}{\mathbf{m}}$ induces an inner fibration $\phi\colon\fromto{\mathscr{M}_m(\mathscr{X}_s)}{\mathscr{M}_{m-1}(\mathscr{X}_s)}$ such that $q=p\circ\phi$. Again by \cite[2.4.1.3(3)]{HTT}, it suffices to show that for any object $X$ of $\mathscr{M}_m(\mathscr{X}_s)$ and any $p$-cocartesian edge $\eta\colon\fromto{\phi(X)}{Y'}$ covering $\sigma$, there exists a $\phi$-cocartesian edge $\fromto{X}{Y}$ of $\mathscr{M}_m(\mathscr{X}_s)$ covering $\eta$. But this follows directly from (\ref{dfn:preWald}.\ref{item:pushcof}).

We now show that $q$ is a \emph{pair} cocartesian fibration. Suppose
\begin{equation*}
\begin{tikzpicture} 
\matrix(m)[matrix of math nodes, 
row sep=4ex, column sep=4ex, 
text height=1.5ex, text depth=0.25ex] 
{X'&X\\ 
Y'&Y\\}; 
\path[>=stealth,->,font=\scriptsize] 
(m-1-1) edge (m-1-2) 
edge (m-2-1) 
(m-1-2) edge (m-2-2) 
(m-2-1) edge (m-2-2); 
\end{tikzpicture}
\end{equation*}
is a square of $\mathscr{M}_m(\mathscr{X}_s)$ in which $\fromto{X'}{X}$ and $\fromto{Y'}{Y}$ are $q$-cocartesian morphisms and $\fromto{X}{Y}$ is ingressive. We aim to show that for any edge $\eta\colon\fromto{\Delta^{\{p,q\}}}{\Delta^m}$, the morphism $\fromto{X'|\Delta^{\{p,q\}}}{Y'|\Delta^{\{p,q\}}}$ is ingressive. For this, we may factor $\fromto{X}{Y}$ as
\begin{equation*}
X\ \tikz[baseline]\draw[>=stealth,>->](0,0.5ex)--(0.5,0.5ex);\ Z\ \tikz[baseline]\draw[>=stealth,>->](0,0.5ex)--(0.5,0.5ex);\ Y,
\end{equation*}
where $Z|\Delta^{\{0,\dots,p\}}=Y|\Delta^{\{0,\dots,p\}}$, and for any $r>p$, the edge $\fromto{X|\Delta^{\{p,r\}}}{Z|\Delta^{\{p,r\}}}$ is cocartesian. Now choose a cocartesian morphism $\fromto{Z'}{Z}$ as well. The proof is now completed by the following observations.
\begin{enumerate}[(\ref{lem:MXScocartfib}.1)]
\item Since the morphism $\fromto{X|\Delta^{\{p,q\}}}{Z|\Delta^{\{p,q\}}}$ is of type (\ref{dfn:injpairstruct}.\ref{item:cocartedgeovercof}), it follows by Quetzalcoatl that the morphism $\fromto{X'|\Delta^{\{p,q\}}}{Z'|\Delta^{\{p,q\}}}$ is of type (\ref{dfn:injpairstruct}.\ref{item:cocartedgeovercof}) as well.
\item The morphism $\fromto{Z|\Delta^{\{p,q\}}}{Y|\Delta^{\{p,q\}}}$ is of type (\ref{dfn:injpairstruct}.\ref{item:degensourceingresstarg}) and the morphism $\fromto{Z'_p}{X'_p}$ is an equivalence; so again by Quetzalcoatl, the morphism $\fromto{Z'|\Delta^{\{p,q\}}}{Y'|\Delta^{\{p,q\}}}$ is of type (\ref{dfn:injpairstruct}.\ref{item:degensourceingresstarg}).\qedhere
\end{enumerate}
\end{proof}
\end{lem}

\begin{ntn}\label{ntn:JandF} Together, these lemmas state that for any Waldhausen cocartesian fibration $\fromto{\mathscr{X}}{S}$, the functor $\fromto{\mathscr{M}(\mathscr{X}/S)}{\Delta^1}$ exhibits an adjunction of $\infty$-categories \cite[Df. 5.2.2.1]{HTT}
\begin{equation*}
\adjunct{F}{\mathscr{F}(\mathscr{X}/S)}{\mathscr{S}(\mathscr{X}/S)}{J}
\end{equation*}
over $N\Delta^{\op}\times S$ in which both $F$ and $J$ are functors of pairs. In particular, for any integer $m\geq 0$ and any vertex $s\in S_0$, the fiber $\fromto{\mathscr{M}_m(\mathscr{X}_s)}{\Delta^1}$ over $(\mathbf{m},s)$ also exhibits an adjunction
\begin{equation*}
\adjunct{F_m}{\mathscr{F}_m(\mathscr{X}_{s})}{\mathscr{S}(\mathscr{X}_{s})}{J_m}.
\end{equation*}

Let's unravel this a bit. Assume $S=\Delta^0$. The functor $J$ is the functor of pairs specified by the edge $\fromto{\Delta^1}{\Pair_{\infty}}$ that classifies the \emph{cartesian} fibration $\fromto{\mathscr{M}(\mathscr{X})}{\Delta^1}$. By construction, this is a forgetful functor: it carries a totally filtered object of $\mathscr{X}$ to its underlying filtered object. The functor $F$ is the functor of pairs specified by the edge $\fromto{\Delta^1}{\Pair_{\infty}}$ that classifies the \emph{cocartesian} fibration $\fromto{\mathscr{M}(\mathscr{X})}{\Delta^1}$, and it is much more interesting: it carries a filtered object $X$ represented as
\begin{equation*}
X_0\ \tikz[baseline]\draw[>=stealth,>->](0,0.5ex)--(0.5,0.5ex);\ X_1\ \tikz[baseline]\draw[>=stealth,>->](0,0.5ex)--(0.5,0.5ex);\ \cdots\ \tikz[baseline]\draw[>=stealth,>->](0,0.5ex)--(0.5,0.5ex);\ X_m
\end{equation*}
to the totally filtered object $FX$ that is initial among all totally filtered objects under $X$; in other words, $FX$ is the quotient of $X$ by $X_0$:
\begin{equation*}
0\simeq X_0/X_0\ \tikz[baseline]\draw[>=stealth,>->](0,0.5ex)--(0.5,0.5ex);\ X_1/X_0\ \tikz[baseline]\draw[>=stealth,>->](0,0.5ex)--(0.5,0.5ex);\ \cdots\ \tikz[baseline]\draw[>=stealth,>->](0,0.5ex)--(0.5,0.5ex);\ X_m/X_0.
\end{equation*}
\end{ntn}

This functor $F$ is a cornerstone for the following result:
\begin{thm}\label{thm:SXiscocartfib} Suppose $\fromto{\mathscr{X}}{S}$ a Waldhausen cocartesian fibration. Then the functor
\begin{equation*}
\fromto{\mathscr{S}(\mathscr{X}/S)}{N\Delta^{\op}\times S}
\end{equation*}
is a Waldhausen cocartesian fibration.
\begin{proof} We first show that the functor $\fromto{\mathscr{S}(\mathscr{X}/S)}{N\Delta^{\op}\times S}$ is a cocartesian fibration by proving the stronger assertion that the inner fibration
\begin{equation*}
p\colon\fromto{\mathscr{M}(\mathscr{X}/S)}{\Delta^1\times N\Delta^{\op}\times S}
\end{equation*}
is a cocartesian fibration. By \ref{prp:FXSisWaldcocart}, the map
\begin{equation}\label{eqn:MXSover0}
\fromto{\Delta^{\{0\}}\times_{\Delta^1}\mathscr{M}(\mathscr{X}/S)}{\Delta^{\{0\}}\times N\Delta^{\op}\times S}
\end{equation}
is a cocartesian fibration. By \ref{lem:MXScocartfib}, for any vertex $(\mathbf{m},s)\in(N\Delta^{\op}\times S)_0$, the map
\begin{equation}\label{eqn:MXSoverms}
\fromto{\mathscr{M}(\mathscr{X}/S)\times_{N\Delta^{\op}\times S}\{(\mathbf{m},s)\}}{\Delta^1\times\{(\mathbf{m},s)\}}
\end{equation}
is a cocartesian fibration. Finally, for any $\mathbf{m}\in\Delta$ and any edge $f\colon\fromto{s}{t}$ of $S$, the functor $f_!\colon\fromto{\mathscr{X}_s}{\mathscr{X}_t}$ carries zero objects to zero objects; consequently, any cocartesian edge of $\mathscr{F}(\mathscr{X}/S)$ that covers $(\id_{\mathbf{m}},f)$ lies in $\mathscr{S}(\mathscr{X}/S)$ if and only if its source does. Thus the map
\begin{equation*}
\fromto{(\Delta^{\{1\}}\times\{\mathbf{m}\})\times_{\Delta^1\times N\Delta^{\op}}\mathscr{M}(\mathscr{X}/S)}{\Delta^{\{1\}}\times \{\mathbf{m}\}\times S}
\end{equation*}
is a cocartesian fibration.

Now to complete the proof that $p$ is a cocartesian fibration, thanks to \cite[2.4.1.3(3)]{HTT} it remains to show that for any vertex $s\in S_0$, any simplicial operator $\phi\colon\fromto{\mathbf{n}}{\mathbf{m}}$, and any totally $m$-filtered object $X$ of $\mathscr{X}_s$, there exists a $p$-cartesian morphism $\fromto{(1,X)}{(1,Y)}$ of $\mathscr{F}(\mathscr{X}/S)$ covering $(\id_1,\phi,\id_s)$. Write $\sigma$ for the nondegenerate edge of $\Delta^1$. The $p$-cartesian edge $e\colon\fromto{(0,X)}{(1,X)}$ covering $(\sigma,\id_{\mathbf{m}},\id_s)$ is also $p$-cocartesian. Since \eqref{eqn:MXSover0} is a cocartesian fibration, there exists a $p$-cocartesian edge $\eta'\colon\fromto{(0,X)}{(0,Y')}$ covering $(\id_0,\phi,\id_s)$. Since \eqref{eqn:MXSoverms} is a cocartesian fibration, there exists a $p$-cocartesian edge $e'\colon\fromto{(0,Y')}{(1,Y)}$ covering $(\sigma,\id_{\mathbf{n}},\id_s)$. Since $e$ is $p$-cocartesian, we have a diagram
\begin{equation*}
\fromto{\Delta^1\times\Delta^1}{\mathscr{M}(\mathscr{X}/S)\times_S\{s\}}
\end{equation*}
of the form
\begin{equation*}
\begin{tikzpicture} 
\matrix(m)[matrix of math nodes, 
row sep=4ex, column sep=4ex, 
text height=1.5ex, text depth=0.25ex] 
{(0,X)&(0,Y')\\ 
(1,X)&(1,Y).\\}; 
\path[>=stealth,->,font=\scriptsize] 
(m-1-1) edge node[above]{$\eta'$} (m-1-2) 
edge node[left]{$e$} (m-2-1) 
(m-1-2) edge node[right]{$e'$} (m-2-2) 
(m-2-1) edge node[below]{$\eta$} (m-2-2); 
\end{tikzpicture}
\end{equation*}
It follows from \cite[2.4.1.7]{HTT} that $\eta$ is $p$-cocartesian.

From \ref{prp:FmisS1plusm} and \ref{prp:FXSisWaldcocart} it follows that the fibers of $\fromto{\mathscr{S}(\mathscr{X}/S)}{N\Delta^{\op}\times S}$ are all Waldhausen $\infty$-categories. For any $\mathbf{m}\in\Delta$ and any edge $f\colon\fromto{s}{t}$ of $S$, the functor $f_{\mathscr{X},!}\colon\fromto{\mathscr{X}_s}{\mathscr{X}_t}$ is exact, whence it follows by \ref{prp:FmisS1plusm} that the functor
\begin{equation*}
f_{\mathscr{S},!}\colon\fromto{\mathscr{S}_m(\mathscr{X}_s)\simeq\Fun_{\Pair_{\infty}}((\Delta^{m-1})^{\sharp},\mathscr{X}_s)}{\Fun_{\Pair_{\infty}}((\Delta^{m-1})^{\sharp},\mathscr{X}_t)\simeq\mathscr{S}_m(\mathscr{X}_t)}
\end{equation*}
is exact, just as in the proof of \ref{prp:FXSisWaldcocart}. Now for any fixed vertex $s\in S_0$ and any simplicial operator $\phi\colon\fromto{\mathbf{n}}{\mathbf{m}}$ of $\Delta$, the functor $\phi_{\mathscr{S},!}\colon\fromto{\mathscr{S}_m(\mathscr{X}_s)}{\mathscr{S}_n(\mathscr{X}_s)}$ is by construction the composite
\begin{equation*}
\mathscr{S}_m(\mathscr{X}_s)\ \tikz[baseline]\draw[>=stealth,right hook->,font=\scriptsize](0,0.5ex)--node[above]{$J_{m,s}$}(0.75,0.5ex);\ \mathscr{F}_m(\mathscr{X}_s)\ \tikz[baseline]\draw[>=stealth,->,font=\scriptsize](0,0.5ex)--node[above]{$\phi_{\mathscr{F},!}$}(0.5,0.5ex);\ \mathscr{F}_n(\mathscr{X}_s)\ \tikz[baseline]\draw[>=stealth,->,font=\scriptsize](0,0.5ex)--node[above]{$F_{n,s}$}(0.75,0.5ex);\ \mathscr{S}_n(\mathscr{X}_s),
\end{equation*}
and as $\phi_{\mathscr{F},!}$ is an exact functor (\ref{prp:FXSisWaldcocart}), we are reduced to checking that the functors of pairs $J_{m,s}$ and $F_{n,s}$ are each exact functors.

For this, it is clear that $J_{m,s}$ and $F_{n,s}$ each carry zero objects to zero objects, and as $F_{n,s}$ is a left adjoint, it preserves any pushout squares that exist in $\mathscr{F}_n(\mathscr{X}_s)$. Moreover, a pushout square in $\mathscr{S}_m(\mathscr{X}_s)$ is nothing more than a pushout square in $\mathscr{F}_m(\mathscr{X}_s)$ of totally $m$-filtered objects; hence $J_{m,s}$ preserves pushouts along ingressive morphisms.
\end{proof}
\end{thm}

For any Waldhausen cocartesian fibration $\fromto{\mathscr{X}}{S}$, write
\begin{equation*}
\mathbf{S}_{\ast}(\mathscr{X}/S)\colon\fromto{N\Delta^{\op}\times S}{\Wald}
\end{equation*}
for the diagram of Waldhausen $\infty$-categories that classifies the Waldhausen cocartesian fibration $\fromto{\mathscr{S}(\mathscr{X}/S)}{N\Delta^{\op}\times S}$, and, similarly, write
\begin{equation*}
\FF_{\ast}(\mathscr{X}/S)\colon\fromto{N\Delta^{\op}\times S}{\Wald}
\end{equation*}
for the diagram of Waldhausen $\infty$-categories that classifies the Waldhausen cocartesian fibration $\fromto{\mathscr{F}(\mathscr{X}/S)}{N\Delta^{\op}\times S}$. An instant consequence of the construction of the functoriality of $\mathscr{S}$ in the proof above is the following.
\begin{cor}\label{cor:Jisanattrans} The functors $F_m\colon\fromto{\mathscr{F}_m(\mathscr{X}/S)}{\mathscr{S}_m(\mathscr{X}/S)}$ assemble to a morphism $F\colon\fromto{\mathscr{F}(\mathscr{X}/S)}{\mathscr{S}(\mathscr{X}/S)}$ of $\mathbf{Wald}_{\infty/N\Delta^{\op}\times S}^{\cocart}$, or, equivalently, a natural transformation
\begin{equation*}
F\colon\fromto{\FF_{\ast}(\mathscr{X}/S)}{\SSS_{\ast}(\mathscr{X}/S)}.
\end{equation*}
\end{cor}
\noindent Note, however, that it is \emph{not} the case that the functors $J_m$ assemble to a natural transformation of this kind.


\subsection*{Virtual Waldhausen $\infty$-categories of filtered objects} Thanks to \ref{nul:Tpafunctor}, the assignments
\begin{equation*}
\goesto{(\mathscr{X}/S)}{(\mathscr{F}(\mathscr{X}/S)/(N\Delta^{\op}\times S))}\textrm{\quad and\quad}\goesto{(\mathscr{X}/S)}{(\mathscr{S}(\mathscr{X}/S)/(N\Delta^{\op}\times S))}
\end{equation*}
define endofunctors of $\Wald^{\cocart}$ over the endofunctor $\goesto{S}{N\Delta^{\op}\times S}$ of $\Cat_{\infty}$. We now aim to descend these functors to endofunctors of the $\infty$-category of virtual Waldhausen $\infty$-categories.

\begin{lem}\label{lem:FandScontinuous} The functors $\fromto{\Wald}{\mathbf{Wald}_{\infty/N\Delta^{\op}}^{\cocart}}$ given by
\begin{equation*}
\goesto{\mathscr{C}}{(\mathscr{F}(\mathscr{C})/N\Delta^{\op})}\textrm{\qquad and\qquad}\goesto{\mathscr{C}}{(\mathscr{S}(\mathscr{C})/N\Delta^{\op})}
\end{equation*}
each preserve filtered colimits.
\begin{proof} By Cor. \ref{cor:colimsinWaldScocart}, it is enough to check the claim fiberwise. The assignment $\goesto{\mathscr{C}}{\mathscr{S}_0(\mathscr{C})}$ is an essentially constant functor whose values are all terminal objects; hence since filtered simplicial sets are weakly contractible, this functor preserves filtered colimits. We are now reduced to the claim that for any natural number $m$, the assignment $\goesto{\mathscr{C}}{\mathscr{F}_m(\mathscr{C})}$ defines a functor $\fromto{\Wald}{\Wald}$ that preserves filtered colimits.

Suppose now that $\Lambda$ is a filtered simplicial set; by \cite[Pr. 5.3.1.16]{HTT}, we may assume that $\Lambda$ is the nerve of a filtered poset. Suppose $\mathscr{C}\colon\fromto{\Lambda^{\rhd}}{\Wald}$ a colimit digram of Waldhausen $\infty$-categories. Suppose $\widetilde{\mathscr{F}}_m(\mathscr{C})\colon\fromto{\Lambda^{\rhd}}{\Pair_{\infty}}$ be a colimit diagram such that $\widetilde{\mathscr{F}}_m(\mathscr{C})|\Lambda=\mathscr{F}_m(\mathscr{C}|\Lambda)$. By \ref{thm:Waldfiltcolims}, we are reduced to showing that the natural functor of pairs
\begin{equation*}
\nu\colon\fromto{\widetilde{\mathscr{F}}_m(\mathscr{C})_{\infty}}{\mathscr{F}_m(\mathscr{C}_{\infty})}
\end{equation*}
is an equivalence. Indeed, $\nu$ induces an equivalence of the underlying $\infty$-categories, since $(\Delta^m)^{\sharp}\times(\Delta^n)^{\flat}$ is a compact object of $\Pair_{\infty}$ (Ex. \ref{exm:finiteimpliescompactforpairs}); hence it remains to show that $\nu$ is a strict functor of pairs. For this it suffices to show that for any ingressive morphism $\psi\colon\cofto{X}{Y}$ of $\mathscr{F}_m(\mathscr{C}_{\infty})$, there exists a vertex $\alpha\in\Lambda$ and an edge $\overline{\psi}\colon\fromto{\overline{X}}{\overline{Y}}$ of $\mathscr{F}_m(\mathscr{C}_{\alpha})$ lifting $\psi$. It is enough to assume that $m=1$ and to show that $\psi$ is either of type (\ref{dfn:injpairstruct}.\ref{item:degensourceingresstarg}) or of type (\ref{dfn:injpairstruct}.\ref{item:cocartedgeovercof}). That is, we may assume that $\psi$ is represented by a square
\begin{equation}\label{eqn:morphinF1C}
\begin{tikzpicture}[baseline]
\matrix(m)[matrix of math nodes, 
row sep=4ex, column sep=4ex, 
text height=1.5ex, text depth=0.25ex] 
{X&Y\\ 
X'&Y'\\}; 
\path[>=stealth,>->,font=\scriptsize] 
(m-1-1) edge (m-1-2) 
edge (m-2-1) 
(m-1-2) edge (m-2-2) 
(m-2-1) edge (m-2-2); 
\end{tikzpicture}
\end{equation}
of ingressive morphisms such that either $\cofto{X}{X'}$ is an equivalence or else the square \eqref{eqn:morphinF1C} is a pushout. Since $(\Delta^1)^{\sharp}\times(\Delta^1)^{\sharp}$ is compact in $\Pair_{\infty}$ (Ex. \ref{exm:finiteimpliescompactforpairs}), a square of ingressive morphisms of the form \eqref{eqn:morphinF1C} must lift to a square of ingressive morphisms
\begin{equation}\label{eqn:morphinF1Calpha}
\begin{tikzpicture}[baseline]
\matrix(m)[matrix of math nodes, 
row sep=4ex, column sep=4ex, 
text height=1.5ex, text depth=0.25ex] 
{\overline{X}&\overline{Y}\\ 
\overline{X}'&\overline{Y}'\\}; 
\path[>=stealth,>->,font=\scriptsize] 
(m-1-1) edge (m-1-2) 
edge (m-2-1) 
(m-1-2) edge (m-2-2) 
(m-2-1) edge (m-2-2); 
\end{tikzpicture}
\end{equation}
of $\mathscr{C}_{\alpha}$ for some vertex $\alpha\in\Lambda$. Now the argument is completed by the following brace of observations.
\begin{enumerate}[(\ref{lem:FandScontinuous}.1)]\addtocounter{enumi}{2}
\item If $\cofto{X}{X'}$ is an equivalence, then, increasing $\alpha$ if necessary, we may assume that its lift $\cofto{\overline{X}}{\overline{X}'}$ in $\mathscr{C}_{\alpha}$ is an equivalence as well, since for example the pushout
\begin{equation*}
\Delta^3\cup^{(\Delta^{\{0,2\}}\sqcup\Delta^{\{1,3\}})}(\Delta^0\sqcup\Delta^0)
\end{equation*}
is compact in the Joyal model structure; hence it represents an ingressive morphism of type (\ref{dfn:injpairstruct}.\ref{item:degensourceingresstarg}) of $\mathscr{F}_1(\mathscr{C}_{\alpha})$.
\item If \eqref{eqn:morphinF1C} is a pushout, then one may form the pushout of $\overline{X}'\ \tikz[baseline]\draw[>=stealth,<-<](0,0.5ex)--(0.5,0.5ex);\ \overline{X}\ \tikz[baseline]\draw[>=stealth,>->](0,0.5ex)--(0.5,0.5ex);\ \overline{Y}$ in $\mathscr{C}_{\alpha}$. Since $\fromto{\mathscr{C}_{\alpha}}{\mathscr{C}_{\infty}}$ preserves such pushouts, we may assume that \eqref{eqn:morphinF1Calpha} is a pushout square in $\mathscr{C}_{\alpha}$; hence it represents an ingressive morphism of type (\ref{dfn:injpairstruct}.\ref{item:cocartedgeovercof}) of $\mathscr{F}_1(\mathscr{C}_{\alpha})$.\qedhere
\end{enumerate}
\end{proof}
\end{lem}

\begin{cnstr} One may compose the functors
\begin{equation*}
\mathscr{F}\colon\fromto{\Wald}{\mathbf{Wald}_{\infty,/N\Delta^{\op}}}\textrm{\qquad and\qquad}\mathscr{S}\colon\fromto{\Wald}{\mathbf{Wald}_{\infty,/N\Delta^{\op}}}
\end{equation*}
with the realization functor $|\cdot|_{N\Delta^{\op}}$ of Df. \ref{dfn:realization}; the results are models for the functors $\fromto{\Wald}{\VWald}$ that assign to any Waldhausen $\infty$-category $\mathscr{C}$ the formal geometric realizations of the simplicial Waldhausen $\infty$-categories $\FF_{\ast}(\mathscr{C})$ and $\SSS_{\ast}(\mathscr{C})$ that classify $\mathscr{F}(\mathscr{C})$ and $\mathscr{S}(\mathscr{C})$. In particular, these composites
\begin{equation*}
|\FF_{\ast}|_{N\Delta^{\op}},|\SSS_{\ast}|_{N\Delta^{\op}}\colon\fromto{\Wald}{\VWald}
\end{equation*}
each preserve filtered colimits, whence one may form their left derived functors (Df. \ref{dfn:lderivWald}), which we will abusively also denote $\mathscr{F}$ and $\mathscr{S}$. These are the essentially unique endofunctors of $\VWald$ that preserve sifted colimits such that the squares
\begin{equation*}
\begin{tikzpicture}[baseline]
\matrix(m)[matrix of math nodes, 
row sep=6ex, column sep=4ex, 
text height=1.5ex, text depth=0.25ex] 
{\Wald&\mathbf{Wald}_{\infty,/N\Delta^{\op}}^{\cocart}\\ 
\VWald&\VWald\\}; 
\path[>=stealth,->,font=\scriptsize] 
(m-1-1) edge node[above]{$\mathscr{F}$} (m-1-2) 
edge node[left]{$j$} (m-2-1) 
(m-1-2) edge node[right]{$|\cdot|_{N\Delta^{\op}}$} (m-2-2) 
(m-2-1) edge node[below]{$\mathscr{F}$} (m-2-2); 
\end{tikzpicture}
\textrm{\qquad and\qquad}
\begin{tikzpicture}[baseline]
\matrix(m)[matrix of math nodes, 
row sep=6ex, column sep=4ex, 
text height=1.5ex, text depth=0.25ex] 
{\Wald&\mathbf{Wald}_{\infty,/N\Delta^{\op}}^{\cocart}\\ 
\VWald&\VWald\\}; 
\path[>=stealth,->,font=\scriptsize] 
(m-1-1) edge node[above]{$\mathscr{S}$} (m-1-2) 
edge node[left]{$j$} (m-2-1) 
(m-1-2) edge node[right]{$|\cdot|_{N\Delta^{\op}}$} (m-2-2) 
(m-2-1) edge node[below]{$\mathscr{S}$} (m-2-2); 
\end{tikzpicture}
\end{equation*}
commute via a specified homotopy.

Also note that the natural transformation $F$ from Cor. \ref{cor:Jisanattrans} descends further to a natural transformation $F\colon\fromto{\mathscr{F}}{\mathscr{S}}$ of endofunctors of $\VWald$.
\end{cnstr}

As it happens, the functor $\mathscr{F}\colon\fromto{\VWald}{\VWald}$ is not particularly exciting:
\begin{prp}\label{prp:FCisacone} For any virtual Waldhausen $\infty$-category $\mathscr{X}$, the virtual Waldhausen $\infty$-category $\mathscr{F}(\mathscr{X})$ is the zero object.
\begin{proof} For any Waldhausen $\infty$-category $\mathscr{C}$, the virtual Waldhausen $\infty$-category $|\mathscr{F}(\mathscr{C})|_{N\Delta^{\op}}$ is by definition a functor $\fromto{\Wald^{\omega}}{\Kan}$ that assigns to any compact Waldhausen $\infty$-category $\mathscr{Y}$ the geometric realization of the simplicial space
\begin{equation*}
\goesto{\mathbf{m}}{\Wald^{\Delta}(\mathscr{Y},\mathscr{F}_m(\mathscr{C}))}.
\end{equation*}
By Pr. \ref{prp:FmisS1plusm}, this simplicial space is the path space of the simplicial space
\begin{equation*}
\goesto{\mathbf{m}}{\Wald^{\Delta}(\mathscr{Y},\mathscr{S}_m(\mathscr{C}))}.\qedhere
\end{equation*}
\end{proof}
\end{prp}

For any Waldhausen $\infty$-category $\mathscr{C}$, we have a natural morphism $\fromto{\mathscr{C}}{\mathscr{F}(\mathscr{C})}$ in $\VWald$, which is induced by the inclusion of the fiber over $0$. The previous result now entitles us to regard the virtual Waldhausen $\infty$-category $\mathscr{F}(\mathscr{C})$ as a \emph{cone} on $\mathscr{C}$. With this perspective, in the next section we will end up thinking of the induced morphism $F\colon\fromto{\mathscr{F}(\mathscr{C})}{\mathscr{S}(\mathscr{C})}$ induced by the functor $F$ as the quotient of $\mathscr{F}(\mathscr{C})$ by $\mathscr{C}$, thereby identifying $\mathscr{S}(\mathscr{C})$ as a \emph{suspension} of $\mathscr{C}$ in a suitable localization of $\VWald$.

The fact that the extensions $\mathscr{F}$ and $\mathscr{S}$ to $\VWald$ preserve sifted colimits now easily implies the following.
\begin{prp} If $S$ is a small sifted $\infty$-category, then the squares
\begin{equation*}
\begin{tikzpicture}[baseline]
\matrix(m)[matrix of math nodes, 
row sep=6ex, column sep=3ex, 
text height=1.5ex, text depth=0.25ex] 
{\mathbf{Wald}_{\infty,/S}^{\cocart}&\mathbf{Wald}_{\infty,/N\Delta^{\op}\times S}^{\cocart}\\ 
\VWald&\VWald\\}; 
\path[>=stealth,->,font=\scriptsize] 
(m-1-1) edge node[above]{$\mathscr{F}$} (m-1-2) 
edge node[left]{$|\cdot|_{S}$} (m-2-1) 
(m-1-2) edge node[right]{$|\cdot|_{N\Delta^{\op}\times S}$} (m-2-2) 
(m-2-1) edge node[below]{$\mathscr{F}$} (m-2-2); 
\end{tikzpicture}
\textrm{\quad and\quad}
\begin{tikzpicture}[baseline]
\matrix(m)[matrix of math nodes, 
row sep=6ex, column sep=3ex, 
text height=1.5ex, text depth=0.25ex] 
{\mathbf{Wald}_{\infty,/S}^{\cocart}&\mathbf{Wald}_{\infty,/N\Delta^{\op}\times S}^{\cocart}\\ 
\VWald&\VWald\\}; 
\path[>=stealth,->,font=\scriptsize] 
(m-1-1) edge node[above]{$\mathscr{S}$} (m-1-2) 
edge node[left]{$|\cdot|_{S}$} (m-2-1) 
(m-1-2) edge node[right]{$|\cdot|_{N\Delta^{\op}\times S}$} (m-2-2) 
(m-2-1) edge node[below]{$\mathscr{S}$} (m-2-2); 
\end{tikzpicture}
\end{equation*}
commute via a specified homotopy.
\end{prp}
\noindent Of course this is no surprise for $\mathscr{F}\colon\fromto{\VWald}{\VWald}$, as we have already seen that $\mathscr{F}$ is constant at zero.


\section{The fissile derived $\infty$-category of Waldhausen $\infty$-categories} A functor $\phi\colon\fromto{\Wald}{\Kan}$ may be described and studied through its left derived functor (Df. \ref{dfn:lderivWald})
\begin{equation*}
\Phi\colon\fromto{\VWald}{\Kan}.
\end{equation*}
In this section, we construct a somewhat peculiar localization $\VaddWald$ of the $\infty$-category $\VWald$ on which the functor $\mathscr{S}\colon\fromto{\VWald}{\VWald}$ constructed in the previous section can be identified as the suspension (Cor. \ref{cor:Sdotisreallysuspension}). In the next section we will use this to show that $\phi$ is \emph{additive} in the sense of Waldhausen just in case $\Phi$ factors through an excisive functor on $\VaddWald$ (Th. \ref{thm:additiveequiv}).


\subsection*{Fissile virtual Waldhausen $\infty$-categories} In Df. \ref{dfn:virtWald}, we defined a virtual Waldhausen $\infty$-category as a presheaf $\mathscr{X}\colon\fromto{\Wald^{\omega,\op}}{\Kan}$ such that the natural maps
\begin{equation*}
\equivto{\mathscr{X}(\mathscr{C}\oplus\mathscr{D})}{\mathscr{X}(\mathscr{C})\times\mathscr{X}(\mathscr{D})}
\end{equation*}
are equivalences. This condition implies in particular that the value of $\mathscr{X}$ on the Waldhausen $\infty$-category of \emph{split} cofiber sequences in a Waldhausen $\infty$-category $\mathscr{C}$ agrees with the product $\mathscr{X}(\mathscr{C})\times\mathscr{X}(\mathscr{C})$. We can ask for more: we can demand that $\mathscr{X}$ be able split even those cofiber sequences that are not already split. That is, we can ask that $\mathscr{X}$ regard the Waldhausen $\infty$-categories of split exact sequences and that of all exact sequences in $\mathscr{C}$ as indistinguishable. This is obviously very closely related to Waldhausen's additivity, and it is what we will mean by a \emph{fissile} virtual Waldhausen $\infty$-category, and the $\infty$-category of these will be called the \emph{fissile derived $\infty$-category} of Waldhausen $\infty$-categories. (The word ``fissile'' in geology and nuclear physics means, in essence, ``easily split.'' The intuition is that when we pass to the fissile derived $\infty$-category, filtered objects can be identified with the sum of their layers.)

But this is asking a lot of our presheaf $\mathscr{X}$. For example, while Waldhausen $\infty$-categories always represent virtual Waldhausen $\infty$-categories, they are almost never fissile. Nevertheless, any virtual Waldhausen $\infty$-category has a best fissile approximation. In other words, the inclusion of fissile virtual Waldhausen $\infty$-categories into virtual Waldhausen $\infty$-categories actually admits a left adjoint, which exhibits the $\infty$-category of fissile virtual Waldhausen $\infty$-categories as a localization of the $\infty$-category of all virtual Waldhausen $\infty$-categories.

\begin{cnstr}\label{cnstr:Em} Suppose $\mathscr{C}$ a Waldhausen $\infty$-category. Then for any integer $m\geq 0$, we may define a fully faithful functor
\begin{equation*}
E_m\colon\into{\mathscr{C}\simeq\mathscr{F}_0(\mathscr{C})}{\mathscr{F}_m(\mathscr{C})}
\end{equation*}
that carries an object $X$ of $\mathscr{C}$ to the constant filtration of length $m$:
\begin{equation*}
X\ \tikz[baseline]\draw[>=stealth,-,double distance=1.5pt](0,0.5ex)--(0.5,0.5ex);\ X\ \tikz[baseline]\draw[>=stealth,-,double distance=1.5pt](0,0.5ex)--(0.5,0.5ex);\ \cdots\ \tikz[baseline]\draw[>=stealth,-,double distance=1.5pt](0,0.5ex)--(0.5,0.5ex);\ X.
\end{equation*}
This is the functor induced by the simplicial operator $\fromto{\mathbf{0}}{\mathbf{m}}$. One has a similiar functor
\begin{equation*}
E_m'\colon\into{\Delta^0\simeq\mathscr{S}_0(\mathscr{C})}{\mathscr{S}_m(\mathscr{C})},
\end{equation*}
which is of course just the inclusion of a contractible Kan complex of zero objects into $\mathscr{S}_m(\mathscr{C})$.

We will also need to have a complete picture of how these functors transform as $\mathbf{m}$ and $\mathscr{C}$ each vary, so we give the following abstract description of them. Since there is an equivalence of $\infty$-categories
\begin{equation*}
\mathbf{Wald}_{\infty/N\Delta^{\op}}^{\cocart}\simeq\Fun(N\Delta^{\op},\Wald)
\end{equation*}
(Pr. \ref{prp:Waldstraightening}), and since $\mathbf{0}$ is an initial object of $N\Delta^{\op}$, it is easy to see that there is an adjunction
\begin{equation*}
\adjunct{C}{\Wald}{\mathbf{Wald}_{\infty/N\Delta^{\op}}^{\cocart}}{R},
\end{equation*}
where $C$ is the functor $\goesto{\mathscr{C}}{\mathscr{C}\times N\Delta^{\op}}$, which represents the contstant functor $\fromto{\Wald}{\Fun(N\Delta^{\op},\Wald)}$, and $R$ is the functor $\goesto{(\mathscr{X}/N\Delta^{\op})}{\mathscr{X}_0}$, which represents evaluation at zero $\fromto{\Fun(N\Delta^{\op},\Wald)}{\Wald}$. The counit $\fromto{CR}{\id}$ of this adjunction can now be composed with the the natural transformation $F\colon\fromto{\mathscr{F}}{\mathscr{S}}$ (which we regard as a morphism of $\Fun(\Wald,\mathbf{Wald}_{\infty/N\Delta^{\op}}^{\cocart})$) to give a commutative square
\begin{equation*}
\begin{tikzpicture} 
\matrix(m)[matrix of math nodes, 
row sep=4ex, column sep=8ex, 
text height=1.5ex, text depth=0.25ex] 
{CR\circ\mathscr{F}\times N\Delta^{\op}&CR\circ\mathscr{S}\\ 
\mathscr{F}&\mathscr{S}\\}; 
\path[>=stealth,->,font=\scriptsize] 
(m-1-1) edge node[above]{$CR\circ F$} (m-1-2) 
edge node[left]{$E$} (m-2-1) 
(m-1-2) edge node[right]{$E'$} (m-2-2) 
(m-2-1) edge node[below]{$F$} (m-2-2); 
\end{tikzpicture}
\end{equation*}
in the $\infty$-category $\Fun(\Wald,\mathbf{Wald}_{\infty/N\Delta^{\op}}^{\cocart})$.

Forming the fiber over an object $\mathbf{m}\in N\Delta^{\op}$, we obtain a commutative square
\begin{equation*}
\begin{tikzpicture} 
\matrix(m)[matrix of math nodes, 
row sep=4ex, column sep=4ex, 
text height=1.5ex, text depth=0.25ex] 
{\mathscr{F}_0&\mathscr{S}_0\\ 
\mathscr{F}_m&\mathscr{S}_m\\}; 
\path[>=stealth,->,font=\scriptsize] 
(m-1-1) edge node[above]{$F_0$} (m-1-2) 
edge node[left]{$E_m$} (m-2-1) 
(m-1-2) edge node[right]{$E'_m$} (m-2-2) 
(m-2-1) edge node[below]{$F_m$} (m-2-2); 
\end{tikzpicture}
\end{equation*}
in the $\infty$-category $\Fun(\Wald,\Wald)$. We see that $E_m$ and $E_m'$ are the functors we identified above.

On the other hand, applying the realization functor $|\cdot|_{N\Delta^{\op}}$ (Df. \ref{dfn:realization}), and noting that
\begin{equation*}
|CR\circ\mathscr{F}|_{N\Delta^{\op}}\simeq|\mathscr{F}_0\times N\Delta^{\op}|_{N\Delta^{\op}}\simeq\id
\end{equation*}
and
\begin{equation*}
|CR\circ\mathscr{S}|_{N\Delta^{\op}}\simeq|\mathscr{S}_0\times N\Delta^{\op}|_{N\Delta^{\op}}\simeq\Delta^0,
\end{equation*}
we obtain a commutative square
\begin{equation}\label{eqn:Sdotissuspension}
\begin{tikzpicture}[baseline]
\matrix(m)[matrix of math nodes, 
row sep=4ex, column sep=4ex, 
text height=1.5ex, text depth=0.25ex] 
{\id&0\\ 
\mathscr{F}&\mathscr{S},\\}; 
\path[>=stealth,->,font=\scriptsize] 
(m-1-1) edge (m-1-2) 
edge node[left]{$E$} (m-2-1) 
(m-1-2) edge (m-2-2) 
(m-2-1) edge node[below]{$F$} (m-2-2); 
\end{tikzpicture}
\end{equation}
in the $\infty$-category $\Fun(\VWald,\VWald)$. When we pass to the fissile derived $\infty$-category, we will actually force this square to be pushout. Since $\mathscr{F}$ is the zero functor (Pr. \ref{prp:FCisacone}), this will exhibit $\mathscr{S}$ as a suspension.
\end{cnstr}

Before we give our definition of \emph{fissibility}, we need a spot of abusive notation.
\begin{ntn} Recall (Cor. \ref{cor:WaldisInd}) that we have an equivalence of $\infty$-categories $\Wald\simeq\Ind(\Wald^{\omega})$. Consequently, we may use the transitivity result of \cite[Pr. 5.3.6.11]{HTT} to conclude that, in the notation of \ref{rec:PABofC}, we also have an equivalence $\mathscr{P}(\Wald^{\omega})\simeq\mathscr{P}_{\mathscr{I}}^{\mathscr{K}}(\Wald)$, where $\mathscr{I}$ is the class of all small filtered simplicial sets, and $\mathscr{K}$ is the class of all small simplicial sets.

In particular, every presheaf $\mathscr{X}\colon\fromto{\Wald^{\omega,\op}}{\Kan}$ extends to an essentially unique presheaf $\fromto{\Wald^{\op}}{\Kan}$ with the property that it carries filtered colimits in $\Wald$ to the corresponding limits in $\Kan$. We will abuse notation by denoting this extended functor by $\mathscr{X}$ as well. This entitles us to speak of the value of a presheaf $\mathscr{X}\colon\fromto{\Wald^{\omega,\op}}{\Kan}$ even on Waldhausen $\infty$-categories that may not be compact.
\end{ntn}

\begin{dfn}\label{dfn:fissilepresheaf} A presheaf $\mathscr{X}\colon\fromto{\Wald^{\omega,\op}}{\Kan}$ will be said to be \textbf{\emph{fissile}} if for every Waldhausen $\infty$-category $\mathscr{C}$ and every integer $m\geq0$, the exact functors $E_m$ and $J_m$ (Cnstr. \ref{cnstr:Em} and Nt. \ref{ntn:JandF}) induce functors
\begin{equation*}
\fromto{\mathscr{X}(\mathscr{F}_m(\mathscr{C}))}{\mathscr{X}(\mathscr{C})}\textrm{\qquad and\qquad}\fromto{\mathscr{X}(\mathscr{F}_m(\mathscr{C}))}{\mathscr{X}(\mathscr{S}_m(\mathscr{C}))}
\end{equation*}
that together exhibit $\mathscr{X}(\mathscr{F}_m(\mathscr{C}))$ as the product of $\mathscr{X}(\mathscr{C})$ and $\mathscr{X}(\mathscr{S}_m(\mathscr{C}))$:
\begin{equation*}
(E_m^{\star},J_m^{\star})\colon\equivto{\mathscr{X}(\mathscr{F}_m(\mathscr{C}))}{\mathscr{X}(\mathscr{C})\times\mathscr{X}(\mathscr{S}_m(\mathscr{C}))}.
\end{equation*}
\end{dfn}

An induction using Pr. \ref{prp:FmisS1plusm} demonstrates that the value of a fissile presheaf $\mathscr{X}\colon\fromto{\Wald^{\omega,\op}}{\Kan}$ on the Waldhausen $\infty$-category of filtered objects $\mathscr{F}_m(\mathscr{C})$ of length $m$ is split into $1+m$ copies of $\mathscr{X}(\mathscr{C})$. That is, the $1+m$ different functors $\fromto{\mathscr{C}}{\mathscr{F}_m(\mathscr{C})}$ of the form 
\begin{equation*}
\goesto{X}{[0\ \tikz[baseline]\draw[>=stealth,-,double distance=1.5pt](0,0.5ex)--(0.5,0.5ex);\ \cdots\ \tikz[baseline]\draw[>=stealth,-,double distance=1.5pt](0,0.5ex)--(0.5,0.5ex);\ 0\ \tikz[baseline]\draw[>=stealth,>->](0,0.5ex)--(0.5,0.5ex);\ X\ \tikz[baseline]\draw[>=stealth,-,double distance=1.5pt](0,0.5ex)--(0.5,0.5ex);\ \cdots\ \tikz[baseline]\draw[>=stealth,-,double distance=1.5pt](0,0.5ex)--(0.5,0.5ex);\ X]}
\end{equation*}
induce an equivalence
\begin{equation*}
\equivto{\mathscr{X}(\mathscr{F}_m(\mathscr{C}))}{\mathscr{X}(\mathscr{C})^{1+m}}.
\end{equation*}

We began our discussion of fissile presheaves by thinking of them as special examples of virtual Waldhausen $\infty$-categories. That wasn't wrong:
\begin{lem} A presheaf $\mathscr{X}\in\mathscr{P}(\Wald^{\omega})$ is fissile only if $\mathscr{X}$ carries direct sums in $\Wald^{\omega}$ to products --- that is, only if $\mathscr{X}$ is a virtual Waldhausen $\infty$-category.
\begin{proof} Suppose $\mathscr{C}$ and $\mathscr{D}$ two compact Waldhausen $\infty$-categories. Consider the retract diagrams
\begin{equation*}
\begin{tikzpicture}[baseline]
\matrix(m)[matrix of math nodes, 
row sep=4ex, column sep=6ex, 
text height=1.5ex, text depth=0.25ex] 
{\mathscr{C}&\mathscr{C}\oplus\mathscr{D}&\mathscr{C}\\
\mathscr{C}\oplus\mathscr{D}&\mathscr{F}_1(\mathscr{C}\oplus\mathscr{D})&\mathscr{C}\oplus\mathscr{D}\\}; 
\path[>=stealth,->,font=\scriptsize] 
(m-1-1) edge (m-1-2) 
edge (m-2-1)
(m-1-2) edge (m-1-3)
edge node[right]{$E_1$} (m-2-2)
(m-1-3) edge (m-2-3)
(m-2-1) edge node[below]{$E_1\oplus J_1$} (m-2-2)
(m-2-2) edge node[below]{$I_{1,0}\oplus F_1$} (m-2-3); 
\end{tikzpicture}
\end{equation*}
and
\begin{equation*}
\begin{tikzpicture}[baseline]
\matrix(m)[matrix of math nodes, 
row sep=4ex, column sep=6ex, 
text height=1.5ex, text depth=0.25ex] 
{\mathscr{D}&\mathscr{C}\oplus\mathscr{D}&\mathscr{D}\\
\mathscr{C}\oplus\mathscr{D}&\mathscr{F}_1(\mathscr{C}\oplus\mathscr{D})&\mathscr{C}\oplus\mathscr{D}.\\}; 
\path[>=stealth,->,font=\scriptsize] 
(m-1-1) edge (m-1-2) 
edge (m-2-1)
(m-1-2) edge (m-1-3)
edge node[right]{$J_1$} (m-2-2)
(m-1-3) edge (m-2-3)
(m-2-1) edge node[below]{$E_1\oplus J_1$} (m-2-2)
(m-2-2) edge node[below]{$I_{1,0}\oplus F_1$} (m-2-3); 
\end{tikzpicture}
\end{equation*}
Here $I_{1,0}$ is the functor induced by the morphism $\goesto{0}{0}$. For any fissile virtual Waldhausen $\infty$-category $\mathscr{X}$, we have an induced retract diagram
\begin{equation}\label{eqn:retfordecompsarevirts}
\begin{tikzpicture}[baseline]
\matrix(m)[matrix of math nodes, 
row sep=4ex, column sep=4ex, 
text height=1.5ex, text depth=0.25ex] 
{\mathscr{X}(\mathscr{C}\oplus\mathscr{D})&\mathscr{X}(\mathscr{F}_1(\mathscr{C}\oplus\mathscr{D}))&\mathscr{X}(\mathscr{C}\oplus\mathscr{D})\\ 
\mathscr{X}(\mathscr{C})\times\mathscr{X}(\mathscr{D})&\mathscr{X}(\mathscr{C}\oplus\mathscr{D})\times\mathscr{X}(\mathscr{C}\oplus\mathscr{D})&\mathscr{X}(\mathscr{C})\times\mathscr{X}(\mathscr{D}).\\}; 
\path[>=stealth,->,font=\scriptsize] 
(m-1-1) edge (m-1-2) 
edge (m-2-1)
(m-1-2) edge (m-1-3)
edge (m-2-2)
(m-1-3) edge (m-2-3)
(m-2-1) edge (m-2-2)
(m-2-2) edge (m-2-3); 
\end{tikzpicture}
\end{equation}
Since the center vertical map is an equivalence, and since equivalences are closed under retracts, so are the outer vertical maps.
\end{proof}
\end{lem}

\begin{ntn} Denote by
\begin{equation*}
\VaddWald\subset\VWald
\end{equation*}
the full subcategory spanned by the fissile functors. We'll call this the \textbf{\emph{fissile derived $\infty$-category}} of Waldhausen $\infty$-categories.
\end{ntn}

Since sifted colimits in $\VWald$ commute with products \cite[Lm. 5.5.8.11]{HTT}, we deduce the following.
\begin{lem}\label{cor:VaddWaldstableundersiftedcolims} The subcategory $\VaddWald\subset\VWald$ is stable under sifted colimits.
\end{lem}


\subsection*{Fissile approximations to virtual Waldhausen $\infty$-categories} Note that representable presheaves are typically \emph{not} fissile. Consequently, the obvious fully faithful inclusion $\into{\Wald^{\omega}}{\VWald}$ does not factor through $\VaddWald\subset\VWald$. Instead, in order to make a representable presheaf fissile, we'll have to form a fissile approximation to it. Fortunately, there's a universal way to do that.

\begin{prp}\label{prp:Laddexists} The inclusion functor admits a left adjoint
\begin{equation*}
L^{\fiss}\colon\fromto{\VWald}{\VaddWald},
\end{equation*}
which exhibits $\VaddWald$ as an accessible localization of $\VWald$.
\begin{proof} For any compact Waldhausen $\infty$-category $\mathscr{C}$ and every integer $m\geq 0$, consider the exact functor
\begin{equation*}
E_m\oplus J_m\colon\fromto{\mathscr{C}\oplus\mathscr{S}_m(\mathscr{C})}{\mathscr{F}_m(\mathscr{C})};
\end{equation*}
let $S$ be the set of morphisms of $\VWald$ of this form; let $\overline{S}$ be the strongly saturated class it generates. Since $\Wald^{\omega}$ is essentially small, the class $\overline{S}$ is of small generation. Hence we may form the accessible localization $S^{-1}\VWald$. Since virtual Waldhausen $\infty$-categories are functors $\mathscr{X}\colon\fromto{\Wald^{\omega,\op}}{\Kan}$ that preserve products, one sees that $S^{-1}\VWald$ coincides with the full subcategory $\VaddWald\subset\VWald$.
\end{proof}
\end{prp}
\noindent The fully faithful inclusion $\into{\VaddWald}{\VWald}$ preserve finite products, and its left adjoint $L^{\fiss}$ preserve finite coproducts, whence we deduce the following.
\begin{cor} The $\infty$-category $\VaddWald$ is compactly generated and admits finite direct sums, which are preserved by the inclusion
\begin{equation*}
\into{\VaddWald}{\VWald}.
\end{equation*}
\end{cor}
\noindent Combining this with Lm. \ref{cor:VaddWaldstableundersiftedcolims} and \cite[Lm. 1.3.2.9]{HA}, we deduce the following somewhat surprising fact.
\begin{cor} The subcategory $\VaddWald\subset\VWald$ is stable under all small colimits.
\end{cor}


\subsection*{Suspension of fissile virtual Waldhausen $\infty$-categories} We now show that the suspension in the fissile derved $\infty$-category is essentially given by the functor $\mathscr{S}$. This is the key to showing that Waldhausen's additivity is essentially equivalent to excision on the fissile derived $\infty$-category (Th. \ref{thm:additiveequiv}). As a first step, we have the following observation.
\begin{prp}\label{prp:Sdotissuspension} The diagram
\begin{equation*}
\begin{tikzpicture} 
\matrix(m)[matrix of math nodes, 
row sep=5ex, column sep=5ex, 
text height=1.5ex, text depth=0.25ex]
{\VWald&\VWald\\ 
\VaddWald&\VaddWald\\}; 
\path[>=stealth,->,font=\scriptsize] 
(m-1-1) edge node[above]{$\mathscr{S}$} (m-1-2) 
edge node[left]{$L^{\fiss}$} (m-2-1)
(m-1-2) edge node[right]{$L^{\fiss}$} (m-2-2)
(m-2-1) edge node[below]{$\Sigma$} (m-2-2); 
\end{tikzpicture}
\end{equation*}
commutes (up to homotopy), where $\Sigma$ is the suspension endofunctor on the fissile derived $\infty$-category $\VaddWald$.
\begin{proof} Apply $L^{\fiss}$ to the square \eqref{eqn:Sdotissuspension} to obtain a square
\begin{equation}\label{eqn:Laddpushoutsquare}
\begin{tikzpicture}[baseline]
\matrix(m)[matrix of math nodes, 
row sep=4ex, column sep=4ex, 
text height=1.5ex, text depth=0.25ex] 
{L^{\fiss}&0\\ 
L^{\fiss}\circ\mathscr{F}&L^{\fiss}\circ\mathscr{S}\\}; 
\path[>=stealth,->,font=\scriptsize] 
(m-1-1) edge (m-1-2) 
edge (m-2-1) 
(m-1-2) edge (m-2-2) 
(m-2-1) edge node[below]{$F$} (m-2-2); 
\end{tikzpicture}
\end{equation}
of natural transformations between functors $\fromto{\VWald}{\VaddWald}$. Since $\mathscr{F}$ is essentially constant with value the zero object, this gives rise to a natural transformation $\fromto{\Sigma\circ L^{\fiss}}{L^{\fiss}\circ\mathscr{S}}$. To see that this natural transformation is an equivalence, it suffices to consider its value on a compact Waldhausen $\infty$-category $\mathscr{C}$. Now for any $\mathbf{m}\in N\Delta^{op}$, we have a diagram
\begin{equation*}
\begin{tikzpicture} 
\matrix(m)[matrix of math nodes, 
row sep=4ex, column sep=4ex, 
text height=1.5ex, text depth=0.25ex] 
{L^{\fiss}\mathscr{S}_0(\mathscr{C})&L^{\fiss}\mathscr{F}_0(\mathscr{C})&L^{\fiss}\mathscr{S}_0(\mathscr{C})\\ 
L^{\fiss}\mathscr{S}_m(\mathscr{C})&L^{\fiss}\mathscr{F}_m(\mathscr{C})&L^{\fiss}\mathscr{S}_m(\mathscr{C})\\}; 
\path[>=stealth,->,font=\scriptsize]
(m-1-1) edge node[above]{$J_0$} (m-1-2)
edge node[left]{$E'_m$} (m-2-1)
(m-2-1) edge node[below]{$J_m$} (m-2-2)
(m-1-2) edge node[above]{$F_0$} (m-1-3) 
edge node[left]{$E_m$} (m-2-2) 
(m-1-3) edge node[right]{$E'_m$} (m-2-3) 
(m-2-2) edge node[below]{$F_m$} (m-2-3); 
\end{tikzpicture}
\end{equation*}
of Waldhausen $\infty$-categories in which the horizontal composites are equivalences. Since $\mathscr{S}_0(\mathscr{C})$ is a zero object, the left-hand square is a pushout by definition; hence the right-hand square is as well. The geometric realization of the right-hand square is precisely the value of the square \eqref{eqn:Laddpushoutsquare} on $\mathscr{C}$.
\end{proof}
\end{prp}

The observation that $\Sigma\circ L^{\fiss}\simeq L^{\fiss}\circ\mathscr{S}$, nice though it is, doesn't quite cut it: we want an even closer relationship between $\mathscr{S}$ and the suspension in the fissile derived $\infty$-category. More precisely, we'd like to know that it isn't necessary to apply $L^{\fiss}$ to $\mathscr{S}(\mathscr{C})$ in order to get $\Sigma L^{\fiss}\mathscr{C}$. So we conclude this section with a proof that the functor $\mathscr{S}\colon\fromto{\VWald}{\VWald}$ already takes values in the fissile derived $\infty$-category $\VaddWald$.
\begin{prp}\label{prp:SXisdistrib} For any virtual Waldhausen $\infty$-category $\mathscr{X}$, the virtual Waldhausen $\infty$-category $\mathscr{S}\mathscr{X}$ is fissile.
\begin{proof} We may write $\mathscr{X}$ as a geometric realization of a simplicial diagram $\mathscr{Y}_{\ast}$ of Waldhausen $\infty$-categories. So our claim is that for any compact Waldhausen $\infty$-category $\mathscr{C}$ and any integer $m\geq 0$, the map
\begin{equation*}
\fromto{(\colim\mathscr{S}(\mathscr{Y}_{\ast}))(\mathscr{F}_m(\mathscr{C}))}{(\colim\mathscr{S}(\mathscr{Y}_{\ast}))(\mathscr{C})\times(\colim\mathscr{S}(\mathscr{Y}_{\ast}))(\mathscr{S}_m(\mathscr{C}))}
\end{equation*}
induced by $(E_m,J_m)$ is an equivalence. Since geometric realization commutes with products, we reduce to the case in which $\mathscr{Y}_{\ast}$ is constantat a Waldhausen $\infty$-category $\mathscr{Y}$. Now our claim is that for any compact Waldhausen $\infty$-category $\mathscr{C}$ and any integer $m\geq 0$, the map
\begin{equation*}
\begin{tikzpicture} 
\matrix(m)[matrix of math nodes, 
row sep=4ex, column sep=4ex, 
text height=1.5ex, text depth=0.25ex] 
{\mathrm{H}(\mathscr{F}_m(\mathscr{C}),(\mathscr{S}(\mathscr{Y})/N\Delta^{\op}))\\ 
\mathrm{H}(\mathscr{C},(\mathscr{S}(\mathscr{Y})/N\Delta^{\op}))\times\mathrm{H}(\mathscr{S}_m(\mathscr{C}),(\mathscr{S}(\mathscr{Y})/N\Delta^{\op}))\\}; 
\path[>=stealth,->,font=\scriptsize] 
(m-1-1) edge (m-2-1); 
\end{tikzpicture}
\end{equation*}
(Cnstr. \ref{ntn:geomrealinWald}) induced by $(E_m,J_m)$ is a weak homotopy equivalence.  To simplify notation, we write $\mathrm{H}(-,\mathscr{S}(\mathscr{Y}))$ for $\mathrm{H}(-,(\mathscr{S}(\mathscr{Y})/N\Delta^{\op}))$ in what follows.

Let's use Joyal's $\infty$-categorical variant of Quillen's Theorem A \cite[Th. 4.1.3.1]{HTT}. Fix an object
\begin{equation*}
((\mathbf{p},\alpha),(\mathbf{q},\beta))\in\mathrm{H}(\mathscr{C},\mathscr{S}(\mathscr{Y}))\times\mathrm{H}(\mathscr{S}_m(\mathscr{C}),\mathscr{S}(\mathscr{Y})).
\end{equation*}
So $\mathbf{p}$ and $\mathbf{q}$ are objects of $N\Delta^{op}$, $\alpha$ is an exact functor $\fromto{\mathscr{C}}{\mathscr{S}_p(\mathscr{Y})}$, and $\beta$ is an exact functor $\fromto{\mathscr{S}_m(\mathscr{C})}{\mathscr{S}_p(\mathscr{Y})}$. Write $J((\mathbf{p},\alpha),(\mathbf{q},\beta))$ for the pullback
\begin{equation*}
\begin{tikzpicture} 
\matrix(m)[matrix of math nodes, 
row sep=4ex, column sep=4ex, 
text height=1.5ex, text depth=0.25ex] 
{J((\mathbf{p},\alpha),(\mathbf{q},\beta))&\mathrm{H}(\mathscr{C},\mathscr{S}(\mathscr{Y}))\times\mathrm{H}(\mathscr{S}_m(\mathscr{C}),\mathscr{S}(\mathscr{Y}))\\ 
\mathrm{H}(\mathscr{F}_m(\mathscr{C}),\mathscr{S}(\mathscr{Y}))&(\mathrm{H}(\mathscr{C},\mathscr{S}(\mathscr{Y}))\times\mathrm{H}(\mathscr{S}_m(\mathscr{C}),\mathscr{S}(\mathscr{Y})))_{((\mathbf{p},\alpha),(\mathbf{q},\beta))/}\\}; 
\path[>=stealth,->,font=\scriptsize] 
(m-1-1) edge (m-1-2) 
edge (m-2-1) 
(m-1-2) edge (m-2-2) 
(m-2-1) edge (m-2-2); 
\end{tikzpicture}
\end{equation*}
We may identify $J((\mathbf{p},\alpha),(\mathbf{q},\beta))$ with a quasicategory whose objects are tuples $(\mathbf{r},\gamma,\mu,\nu,\sigma,\tau)$ consisting of:
\begin{itemize}
\item[---] $\mathbf{r}$ is an object of $\Delta$,
\item[---] $\gamma\colon\fromto{\mathscr{S}_m(\mathscr{C})}{\mathscr{S}_r\mathscr{Y}}$ is an exact functor,
\item[---] $\mu\colon\fromto{\mathbf{r}}{\mathbf{p}}$ and $\nu\colon\fromto{\mathbf{r}}{\mathbf{q}}$ are morphisms of $\Delta$, and
\item[---] $\sigma\colon\equivto{\mu^{\ast}\alpha}{\gamma|_{\mathscr{C}}}$ and $\tau\colon\equivto{\nu^{\ast}\beta}{\gamma|_{\mathscr{S}_m(\mathscr{C})}}$ are equivalences of exact functors.
\end{itemize}

Denote by $\kappa$ the constant functor $\fromto{J((\mathbf{p},\alpha),(\mathbf{q},\beta))}{J((\mathbf{p},\alpha),(\mathbf{q},\beta))}$ at the object
\begin{equation*}
\left(0,0,\into{\{0\}}{\mathbf{p}},\into{\{0\}}{\mathbf{q}},0,0\right).
\end{equation*}
To prove that $J((\mathbf{p},\alpha),(\mathbf{q},\beta))$ is contractible, we construct an endofunctor $\lambda$ and natural transformations
\begin{equation*}
\id\ \tikz[baseline]\draw[>=stealth,<-](0,0.5ex)--(0.5,0.5ex);\ \lambda\ \tikz[baseline]\draw[>=stealth,->](0,0.5ex)--(0.5,0.5ex);\ \kappa.
\end{equation*}
We define the functor $\lambda$ by
\begin{equation*}
\lambda(\mathbf{r},\gamma,\mu,\nu,\sigma,\tau)\coloneq(\mathbf{r}^{\lhd},s_0\circ\gamma,\mu',\nu',\sigma',\tau'),
\end{equation*}
where $\mu'|_{\mathbf{r}}=\mu$ and $\mu'(-\infty)=0$, $\nu'|_{\mathbf{r}}=\nu$ and $\nu'(-\infty)=0$, and $\sigma'$ and $\tau'$ are the obvious extensions of $\sigma$ and $\tau$. The inclusion $\into{\mathbf{r}}{\mathbf{r}^{\lhd}}$ induces a natural transformation $\fromto{\lambda}{\id}$, and the inclusion $\into{\{-\infty\}}{\mathbf{r}^{\lhd}}$ induces a natural transformation $\fromto{\lambda}{\kappa}$.
\end{proof}
\end{prp}
\noindent We thus have the following enhancement of Pr. \ref{prp:Sdotissuspension}.
\begin{cor}\label{cor:Sdotisreallysuspension} The diagram
\begin{equation*}
\begin{tikzpicture} 
\matrix(m)[matrix of math nodes, 
row sep=4ex, column sep=4ex, 
text height=1.5ex, text depth=0.25ex]
{\VWald&\\
&\VaddWald\\
\VaddWald&\\}; 
\path[>=stealth,->,font=\scriptsize] 
(m-1-1) edge node[above right]{$\mathscr{S}$} (m-2-2) 
edge node[left]{$L^{\fiss}$} (m-3-1)
(m-3-1) edge node[below right]{$\Sigma$} (m-2-2); 
\end{tikzpicture}
\end{equation*}
commutes (up to homotopy), where $\Sigma$ is the suspension endofunctor on the fissile derived $\infty$-category $\VaddWald$.
\end{cor}


\section{Additive theories} In this section we introduce the $\infty$-categorical analogue of Waldhausen's notion of additivity, and we prove our Structure Theorem (Th. \ref{thm:additiveequiv}), which identifies the homotopy theory of additive functors $\fromto{\Wald}{\Kan}$ with the homotopy theory of certain excisive functors $\fromto{\VaddWald}{\Kan}$ on the fissile derived $\infty$-category of the previous section. Using this, we can find the best additive approximation to any functor $\phi\colon\fromto{\Wald}{\Kan}$ as a Goodwillie differential. Since suspension in this $\infty$-category is given by the functor $\mathscr{S}$, this best excisive approximation $D\phi$ can be exhibited by a formula
\begin{equation*}
\goesto{\mathscr{C}}{\colim_n\Omega^n|\phi(S_{\ast}^n(\mathscr{C}))|}.
\end{equation*}
If $\phi$ preserves finite products, the colimit turns out to be unnecessary, and $D\phi$ can be given by an even simpler formula:
\begin{equation*}
\goesto{\mathscr{C}}{\Omega|\phi(S_{\ast}(\mathscr{C}))|}.
\end{equation*}

In the next section, we'll use this perspective on additivity to prove some fundamental things, such as the Eilenberg Swindle and Waldhausen's Fibration Theorem, for general additive functors. In \S \ref{sect:univpropKthy}, we'll apply our additive approximation to the ``moduli space of objects'' functor $\iota$ to give a universal description of algebraic $K$-theory of Waldhausen $\infty$-categories, and the formula above shows that our algebraic $K$-theory extends Waldhausen's.


\subsection*{Theories and additive theories} The kinds of functors we're going to be thinking about are called \emph{theories}. What we'll show is that among theories, one can isolate the class of \emph{additive theories}, which split all exact sequences.

\begin{dfn}\label{dfn:reducedtheories} Suppose $C$ and $D$ $\infty$-categories, and suppose that $C$ is pointed. Recall (\cite[p. 1]{MR1076523} or \cite[Df. 1.4.2.1(ii)]{HA}) that a functor $\fromto{C}{D}$ is \textbf{\emph{reduced}} if it carries the zero object of $C$ to the terminal object of $D$. We write $\Fun^{\ast}(C,D)\subset\Fun(C,D)$ for the full subcategory spanned by the reduced functors, and if $\mathscr{A}$ is a collection of simplicial sets, then we write $\Fun_{\mathscr{A}}^{\ast}(C,D)\subset\Fun(C,D)$ for the full subcategory spanned by the reduced functors that preserve $\mathscr{A}$-shaped colimits (\ref{rec:Ashapedcolim}).

Similarly, recall thet a functor $\fromto{C}{D}$ is \textbf{\emph{excisive}} if it carries pushout squares in $C$ to pullback squares in $D$.

Suppose $\mathscr{E}$ an $\infty$-topos. By an \textbf{\emph{$\mathscr{E}$-valued theory}}, we shall here mean a reduced functor $\fromto{\Wald}{\mathscr{E}}$ that preserves filtered colimits. We write $\Thy(\mathscr{E})$ for the full subcategory of $\Fun(\Wald,\mathscr{E})$ spanned by $\mathscr{E}$-valued theories.
\end{dfn}

Those who grimace the prospect of contemplating general $\infty$-topoi can enjoy a complete picture of what's going on by thinking only of examples of the form $\mathscr{E}=\Fun(S,\Kan)$. The extra generality comes at no added expense, but we won't get around to using it here.

Note that a theory $\phi\colon\fromto{\Wald}{\mathscr{E}}$ may be uniquely identified in different ways. On one hand, $\phi$ is (Cor. \ref{cor:WaldisInd}) the left Kan extension of its restriction
\begin{equation*}
\phi|_{\Wald^{\omega}}\colon\fromto{\Wald^{\omega}}{\mathscr{E}};
\end{equation*}
on the other, we can extend $\phi$ to its left derived functor (Df. \ref{dfn:lderivWald})
\begin{equation*}
\Phi\colon\fromto{\VWald}{\mathscr{E}},
\end{equation*}
which is the unique extension of $\phi$ that preserves all sifted colimits.

Many examples of theories that arise in practice have the property that the natural morphism $\fromto{\phi(\mathscr{C}\oplus\mathscr{D})}{\phi(\mathscr{C})\times\phi(\mathscr{D})}$ is an equivalence. We'll look at these theories more closely below (Df. \ref{dfn:preadditive}). In any case, when this happens, the sum functor $\fromto{\mathscr{C}\oplus\mathscr{C}}{\mathscr{C}}$ defines a monoid structure on $\pi_0\phi(\mathscr{C})$. For invariants like $K$-theory, we'll want to demand that this monoid actually be a group. We thus make the following definition, which is sensible for \emph{any} theory.
\begin{dfn} A theory $\phi\in\Thy(\mathscr{E})$ will be said to be \textbf{\emph{grouplike}} if, for any Waldhausen $\infty$-category $\mathscr{C}$, the \emph{shear functor} $\fromto{\mathscr{C}\oplus\mathscr{C}}{\mathscr{C}\oplus\mathscr{C}}$ defined by the assignment $\goesto{(X,Y)}{(X,X\vee Y)}$ induces an equivalence $\equivto{\pi_0\phi(\mathscr{C}\oplus\mathscr{C})}{\pi_0\phi(\mathscr{C}\oplus\mathscr{C})}$.
\end{dfn}

To formulate our Structure Theorem, we need to stare at a few functors between various Waldhausen $\infty$-categories of filtered objects.
\begin{cnstr} Suppose $m\geq 0$ an integer, and suppose $0\leq k\leq m$. We consider the morphism $i_k\colon\into{\mathbf{0}\cong\{k\}}{\mathbf{m}}$ of $\Delta$. For any Waldhausen $\infty$-category $\mathscr{C}$, write $I_{m,k}$ for the induced functor $\fromto{\mathscr{F}_m(\mathscr{C})}{\mathscr{F}_0(\mathscr{C})}$, and write $I'_{m,k}$ for the induced functor $\fromto{\mathscr{S}_m(\mathscr{C})}{\mathscr{S}_0(\mathscr{C})}$. Of course $\mathscr{F}_0(\mathscr{C})\simeq\mathscr{C}$ and $\mathscr{S}_0(\mathscr{C})\simeq0$. So the functor $I_{m,k}$ extracts from a filtered object
\begin{equation*}
X_0\ \tikz[baseline]\draw[>=stealth,>->](0,0.5ex)--(0.5,0.5ex);\ X_1\ \tikz[baseline]\draw[>=stealth,>->](0,0.5ex)--(0.5,0.5ex);\ \cdots\ \tikz[baseline]\draw[>=stealth,>->](0,0.5ex)--(0.5,0.5ex);\ X_m
\end{equation*}
its $k$-th filtered piece $X_k$, and the functor $I'_{m,k}$ is, by necessity, the trivial functor.

We may now contemplate a square of retract diagrams
\begin{equation*}
\fromto{(\Delta^2/\Delta^{\{0,2\}})\times(\Delta^2/\Delta^{\{0,2\}})}{\Wald}
\end{equation*}
given by
\begin{equation}\label{eqn:pushpullsquares}
\begin{tikzpicture}[baseline]
\matrix(m)[matrix of math nodes, 
row sep=4ex, column sep=4ex, 
text height=1.5ex, text depth=0.25ex] 
{\mathscr{S}_0(\mathscr{C})&\mathscr{F}_0(\mathscr{C})&\mathscr{S}_0(\mathscr{C})\\ 
\mathscr{S}_m(\mathscr{C})&\mathscr{F}_m(\mathscr{C})&\mathscr{S}_m(\mathscr{C})\\
\mathscr{S}_0(\mathscr{C})&\mathscr{F}_0(\mathscr{C})&\mathscr{S}_0(\mathscr{C}).\\}; 
\path[>=stealth,->,font=\scriptsize] 
(m-1-1) edge node[above]{$J_0$} (m-1-2) 
edge node[left]{$E'_m$} (m-2-1) 
(m-1-2) edge node[above]{$F_0$} (m-1-3)
edge node[left]{$E_m$} (m-2-2) 
(m-1-3) edge node[right]{$E'_m$} (m-2-3)
(m-2-1) edge node[below]{$J_m$} (m-2-2) 
edge node[left]{$I'_{m,k}$} (m-3-1) 
(m-2-2) edge node[above]{$F_m$} (m-2-3)
edge node[right]{$I_{m,k}$} (m-3-2) 
(m-2-3) edge node[right]{$I_{m,k}$} (m-3-3)
(m-3-1) edge node[below]{$J_0$} (m-3-2)
(m-3-2) edge node[below]{$F_0$} (m-3-3); 
\end{tikzpicture}
\end{equation}
Only the upper right square of \eqref{eqn:pushpullsquares} is (by Cnstr. \ref{cnstr:Em}) functorial in $\mathbf{m}$.

We may now apply the localization functor $L^{\fiss}$ to \eqref{eqn:pushpullsquares}. In the resulting diagram
\begin{equation}\label{eqn:Laddpushpullsquares}
\begin{tikzpicture}[baseline]
\matrix(m)[matrix of math nodes, 
row sep=4ex, column sep=4ex, 
text height=1.5ex, text depth=0.25ex] 
{L^{\fiss}\mathscr{S}_0(\mathscr{C})&L^{\fiss}\mathscr{F}_0(\mathscr{C})&L^{\fiss}\mathscr{S}_0(\mathscr{C})\\ 
L^{\fiss}\mathscr{S}_m(\mathscr{C})&L^{\fiss}\mathscr{F}_m(\mathscr{C})&L^{\fiss}\mathscr{S}_m(\mathscr{C})\\
L^{\fiss}\mathscr{S}_0(\mathscr{C})&L^{\fiss}\mathscr{F}_0(\mathscr{C})&L^{\fiss}\mathscr{S}_0(\mathscr{C}),\\}; 
\path[>=stealth,->,font=\scriptsize] 
(m-1-1) edge node[above]{$J_0$} (m-1-2) 
edge node[left]{$E'_m$} (m-2-1) 
(m-1-2) edge node[above]{$F_0$} (m-1-3)
edge node[left]{$E_m$} (m-2-2) 
(m-1-3) edge node[right]{$E'_m$} (m-2-3)
(m-2-1) edge node[below]{$J_m$} (m-2-2) 
edge node[left]{$I'_{m,k}$} (m-3-1) 
(m-2-2) edge node[above]{$F_m$} (m-2-3)
edge node[right]{$I_{m,k}$} (m-3-2) 
(m-2-3) edge node[right]{$I_{m,k}$} (m-3-3)
(m-3-1) edge node[below]{$J_0$} (m-3-2)
(m-3-2) edge node[below]{$F_0$} (m-3-3); 
\end{tikzpicture}
\end{equation}
the square in the upper left corner is a pushout, whence every square is a pushout.
\end{cnstr}

Now we are ready to state the Structure Theorem.
\begin{thm}[Structure Theorem for Additive Theories]\label{thm:additiveequiv} Suppose $\mathscr{E}$ an $\infty$-topos. Suppose $\phi$ an $\mathscr{E}$-valued theory. Then the following are equivalent.
\begin{enumerate}[(\ref{thm:additiveequiv}.1)]
\item\label{item:add} For any Waldhausen $\infty$-category $\mathscr{C}$, any integer $m\geq1$, and any integer $0\leq k\leq m$, the functors 
\begin{equation*}
\phi(F_m)\colon\fromto{\phi(\mathscr{F}_m(\mathscr{C}))}{\phi(\mathscr{S}_m(\mathscr{C}))}\textrm{\quad and\quad}\phi(I_{m,k})\colon\fromto{\phi(\mathscr{F}_m(\mathscr{C}))}{\phi(\mathscr{F}_0(\mathscr{C}))}
\end{equation*}
exhibit $\phi(\mathscr{F}_m(\mathscr{C}))$ as a product of $\phi(\mathscr{S}_m(\mathscr{C}))$ and $\phi(\mathscr{F}_0(\mathscr{C}))$.
\item\label{item:group} For any Waldhausen $\infty$-category $\mathscr{C}$ and for any functor
\begin{equation*}
\SSS_{\ast}(\mathscr{C})\colon\fromto{N\Delta^{\op}}{\Wald}
\end{equation*}
that classifies the Waldhausen cocartesian fibration $\fromto{\mathscr{S}(\mathscr{C})}{N\Delta^{\op}}$, the induced functor $\phi\circ\SSS_{\ast}(\mathscr{C})\colon\fromto{N\Delta^{\op}}{\mathscr{E}_{\ast}}$ is a group object \cite[Df. 7.2.2.1]{HTT}.
\item\label{item:stableadd} The theory $\phi$ is grouplike, and for any Waldhausen $\infty$-category $\mathscr{C}$ and any integer $m\geq1$, the functors 
\begin{equation*}
\phi(F_m)\colon\fromto{\phi(\mathscr{F}_m(\mathscr{C}))}{\phi(\mathscr{S}_m(\mathscr{C}))}\textrm{\quad and\quad}\phi(I_{m,0})\colon\fromto{\phi(\mathscr{F}_m(\mathscr{C}))}{\phi(\mathscr{F}_0(\mathscr{C}))}
\end{equation*}
exhibit $\phi(\mathscr{F}_m(\mathscr{C}))$ as a product of $\phi(\mathscr{S}_m(\mathscr{C}))$ and $\phi(\mathscr{F}_0(\mathscr{C}))$.
\item\label{item:stableaddmisone} The theory $\phi$ is grouplike, and for any Waldhausen $\infty$-category $\mathscr{C}$, the functors 
\begin{equation*}
\phi(F_1)\colon\fromto{\phi(\mathscr{F}_1(\mathscr{C}))}{\phi(\mathscr{S}_1(\mathscr{C}))}\textrm{\quad and\quad}\phi(I_{1,0})\colon\fromto{\phi(\mathscr{F}_1(\mathscr{C}))}{\phi(\mathscr{F}_0(\mathscr{C}))}
\end{equation*}
exhibit $\phi(\mathscr{F}_1(\mathscr{C}))$ as a product of $\phi(\mathscr{S}_1(\mathscr{C}))$ and $\phi(\mathscr{F}_0(\mathscr{C}))$.
\item\label{item:condonhocat} The theory $\phi$ is grouplike, it carries direct sums to products, and, for any Waldhausen $\infty$-category $\mathscr{C}$, the images of $\phi(I_{1,1})$ and $\phi(I_{1,0}\oplus F_1)$ in the set $\Mor_{h\mathscr{E}_{\ast}}(\mathscr{F}_1(\mathscr{C}),\mathscr{C})$ are equal.
\item\label{item:stablecat} The theory $\phi$ is grouplike, and for any Waldhausen $\infty$-category $\mathscr{C}$ and any functor $\SSS_{\ast}(\mathscr{C})\colon\fromto{N\Delta^{\op}}{\Wald}$ that classifies the Waldhausen cocartesian fibration $\fromto{\mathscr{S}(\mathscr{C})}{N\Delta^{\op}}$, the induced functor $\phi\circ\SSS_{\ast}(\mathscr{C})\colon\fromto{N\Delta^{\op}}{\mathscr{E}_{\ast}}$ is a category object \textup{(}see \textup{Pr. \ref{prp:Fisacategory}} or \cite[Df. 1.1.1]{G}\textup{)}.
\item\label{item:excisive} The left derived functor $\Phi\colon\fromto{\VWald}{\mathscr{E}}$ of $\phi$ factors through an excisive functor
\begin{equation*}
\Phi_{\add}\colon\fromto{\VaddWald}{\mathscr{E}}.
\end{equation*}
\end{enumerate}
\begin{proof} The equivalence of conditions (\ref{thm:additiveequiv}.\ref{item:add}) and (\ref{thm:additiveequiv}.\ref{item:group}) follows from Pr. \ref{prp:FmisS1plusm} and the proof of \cite[Pr. 6.1.2.6]{HTT}. (Also see \cite[Rk. 6.1.2.8]{HTT}.) Conditions (\ref{thm:additiveequiv}.\ref{item:stableadd}) and (\ref{thm:additiveequiv}.\ref{item:stablecat}) are clearly special cases of (\ref{thm:additiveequiv}.\ref{item:add}) and (\ref{thm:additiveequiv}.\ref{item:group}), respectively, and condition (\ref{thm:additiveequiv}.\ref{item:stableaddmisone}) is a special case of (\ref{thm:additiveequiv}.\ref{item:stableadd}). The equivalence of (\ref{thm:additiveequiv}.\ref{item:stableadd}) and (\ref{thm:additiveequiv}.\ref{item:stablecat}) also follows directly from Pr. \ref{prp:FmisS1plusm}.

Let us show that (\ref{thm:additiveequiv}.\ref{item:stableaddmisone}) implies (\ref{thm:additiveequiv}.\ref{item:condonhocat}). We begin by noting that we have an analogue of the commutative diagram \eqref{eqn:retfordecompsarevirts}:
\begin{equation*}
\begin{tikzpicture}[baseline]
\matrix(m)[matrix of math nodes, 
row sep=4ex, column sep=4ex, 
text height=1.5ex, text depth=0.25ex] 
{\phi(\mathscr{C}\oplus\mathscr{D})&\phi(\mathscr{F}_1(\mathscr{C}\oplus\mathscr{D}))&\phi(\mathscr{C}\oplus\mathscr{D})\\ 
\phi(\mathscr{C})\times\phi(\mathscr{D})&\phi(\mathscr{C}\oplus\mathscr{D})\times\phi(\mathscr{C}\oplus\mathscr{D})&\phi(\mathscr{C})\times\phi(\mathscr{D}),\\}; 
\path[>=stealth,->,font=\scriptsize] 
(m-1-1) edge (m-1-2) 
edge (m-2-1)
(m-1-2) edge (m-1-3)
edge (m-2-2)
(m-1-3) edge (m-2-3)
(m-2-1) edge (m-2-2)
(m-2-2) edge (m-2-3); 
\end{tikzpicture}
\end{equation*}
and once again it is a retract diagram in $\mathscr{E}$. Since $\mathscr{E}$ admits filtered colimits, equivalences therein are closed under retracts, so since the center vertical morphism is an equivalence, the outer vertical morphisms are as well. Hence $\phi $ carries direct sums to products. Now the exact functor $I_{1,0}\oplus F_1$ admits a (homotopy) section $\sigma\colon\fromto{\mathscr{C}\oplus\mathscr{C}}{\mathscr{F}_1(\mathscr{C})}$ such that $I_{1,1}\circ\sigma\simeq\nabla$. Hence if $\phi$ satisfies (\ref{thm:additiveequiv}.\ref{item:stableaddmisone}), then $\phi(I_{1,0}\oplus F_1)$ is an equivalence with homotopy inverse $\phi(\sigma)$, whence $\phi(I_{1,1})$ and $\phi(I_{1,0}\oplus F_1)$ are equal in $\Mor_{h\mathscr{E}}(\mathscr{F}_1(\mathscr{C}),\mathscr{C})$.

It is now easy to see that (\ref{thm:additiveequiv}.\ref{item:stablecat}) implies (\ref{thm:additiveequiv}.\ref{item:group}).

We now show that (\ref{thm:additiveequiv}.\ref{item:condonhocat}) implies (\ref{thm:additiveequiv}.\ref{item:stableadd}). For any natural number $m$, suppose the images of $\phi(I_{1,1})$ and $\phi(I_{1,0}\oplus F_1)$ in $\Mor_{h\mathscr{E}}(\mathscr{F}_1(\mathscr{F}_m(\mathscr{C})),\mathscr{F}_m(\mathscr{C}))$ are equal; we must show that $\phi(I_{m,0}\oplus F_m)$ is an equivalence. Compose $I_{1,1}$ and $I_{1,0}\oplus F_1$ with the exact functor $\fromto{\mathscr{F}_m(\mathscr{C})}{\mathscr{F}_1(\mathscr{F}_m(\mathscr{C}))}$ that sends a filtered object
\begin{equation*}
X_0\ \tikz[baseline]\draw[>=stealth,>->](0,0.5ex)--(0.75,0.5ex);\ X_1\ \tikz[baseline]\draw[>=stealth,>->](0,0.5ex)--(0.75,0.5ex);\ X_2\ \tikz[baseline]\draw[>=stealth,>->](0,0.5ex)--(0.75,0.5ex);\ \cdots\ \tikz[baseline]\draw[>=stealth,>->](0,0.5ex)--(0.75,0.5ex);\ X_m
\end{equation*}
to the ingressive morphism of filtered objects given by the diagram
\begin{equation*}
\begin{tikzpicture} 
\matrix(m)[matrix of math nodes, 
row sep=4ex, column sep=4ex, 
text height=1.5ex, text depth=0.25ex] 
{X_0&X_0&X_0&\cdots&X_0\\ 
X_0&X_1&X_2&\cdots&X_m;\\}; 
\path[>=stealth,font=\scriptsize] 
(m-1-1) edge[-,double distance=1.5pt] (m-1-2) 
edge[-,double distance=1.5pt] (m-2-1) 
(m-1-2) edge[-,double distance=1.5pt] (m-1-3)
edge[>->] (m-2-2)
(m-1-3) edge[-,double distance=1.5pt] (m-1-4)
edge[>->] (m-2-3)
(m-1-4) edge[-,double distance=1.5pt] (m-1-5)
(m-1-5) edge[>->] (m-2-5)
(m-2-1) edge[>->] (m-2-2)
(m-2-2) edge[>->] (m-2-3)
(m-2-3) edge[>->] (m-2-4)
(m-2-4) edge[>->] (m-2-5); 
\end{tikzpicture}
\end{equation*}
the exact functor $I_{m,0}\oplus F_m$ also admits a section $\sigma\colon\fromto{\mathscr{C}\oplus\mathscr{S}_m(\mathscr{C})}{\mathscr{F}_m(\mathscr{C})}$ (up to homotopy) such that $I_{m,1}\circ\sigma\simeq\nabla$, and applying our condition on $\phi$, we find that $\phi(\sigma\circ(I_{m,0}\oplus F_m))\simeq\phi(\id)$.

We now set about showing that (\ref{thm:additiveequiv}.\ref{item:stableadd}) implies (\ref{thm:additiveequiv}.\ref{item:excisive}). First, we show that $\Phi$ factors through a functor
\begin{equation*}
\Phi_{\add}\colon\fromto{\VaddWald}{\mathscr{E}}.
\end{equation*}
As above, we find that $\Phi$ carries direct sums to products, and from this we deduce that $\Phi$ carries morphisms of the class $S$ described in Pr. \ref{prp:Laddexists} to equivalences. We further claim that the family $T$ of those morphisms of $\VWald$ that are carried to equivalences by $\Phi$ is a strongly saturated class. Since $\Phi$ sends direct sums to products, it carries any finite coproduct of elements of $T$ to equivalences. Moreover, since $\Phi$ preserves sifted colimits, it preserves any morphism that can be exhibited as a small sifted colimit of elements of $T$. Hence the full subcategory of $\mathscr{O}(\VWald)$ spanned by the elements of $T$ is closed under all small colimits. Finally, to prove that any pushout $\psi'\colon\fromto{\mathscr{X}'}{\mathscr{Y}'}$ of an element $\psi\colon\fromto{\mathscr{X}}{\mathscr{Y}}$ of $T$ (along any morphism $\fromto{\mathscr{X}}{\mathscr{X}'}$) lies again in $T$, we note that we may exhibit $\psi'$ as the natural morphism of geometric realizations\footnote{We are grateful to Jacob Lurie for this observation.}
\begin{equation*}
\fromto{|B_{\ast}(\mathscr{X}',\mathscr{X},\mathscr{X})|}{|B_{\ast}(\mathscr{X}',\mathscr{X},\mathscr{Y})|},
\end{equation*}
where the simplicial objects $B_{\ast}(\mathscr{X}',\mathscr{X},\mathscr{X})$ and $B_{\ast}(\mathscr{X}',\mathscr{X},\mathscr{Y})$ are two-sided bar constructions defined by
\begin{equation*}
B_n(\mathscr{X}',\mathscr{X},\mathscr{X})\coloneq\mathscr{X}'\oplus\mathscr{X}^{\oplus n}\oplus\mathscr{X}\textrm{\quad and\quad}B_n(\mathscr{X}',\mathscr{X},\mathscr{Y})\coloneq\mathscr{X}'\oplus\mathscr{X}^{\oplus n}\oplus\mathscr{Y}.
\end{equation*}
Since $T$ is closed under formation of products, each map
\begin{equation*}
\fromto{B_n(\mathscr{X}',\mathscr{X},\mathscr{X})}{B_n(\mathscr{X}',\mathscr{X},\mathscr{Y})}
\end{equation*}
is an element of $T$, and since $T$ is closed under geometric realizations, the morphism $\fromto{\mathscr{X}'}{\mathscr{Y}'}$ is an element of $T$. Hence $T$ is strongly saturated and therefore contains $\overline{S}$; thus $\Phi$ factors through a functor $\Phi_{\add}\colon\fromto{\VaddWald}{\mathscr{E}}$.

We now show that $\Phi_{\add}$ is excisive. For any nonnegative integer $m$, apply $\phi$ to the diagram \eqref{eqn:pushpullsquares} with $k=0$. The lower right corner of the resulting diagram is a pullback. Hence the upper right corner of the diagram resulting from applying $\phi$ to the diagram \eqref{eqn:pushpullsquares} is also a pullback. Now we may form the geometric realization of this simplicial diagram of squares to obtain a square
\begin{equation*}
\begin{tikzpicture} 
\matrix(m)[matrix of math nodes, 
row sep=4ex, column sep=4ex, 
text height=1.5ex, text depth=0.25ex] 
{\Phi(\mathscr{F}_0(\mathscr{C}))&\Phi(\mathscr{S}_0(\mathscr{C}))\\ 
\Phi(\mathscr{F}(\mathscr{C}))&\Phi(\mathscr{S}(\mathscr{C})).\\}; 
\path[>=stealth,->,font=\scriptsize] 
(m-1-1) edge (m-1-2) 
edge (m-2-1) 
(m-1-2) edge (m-2-2) 
(m-2-1) edge (m-2-2); 
\end{tikzpicture}
\end{equation*}
It follows from the Segal delooping machine (\cite{MR50:5782} and \cite[Lm. 7.2.2.11]{HTT}) that this square is a pullback as well, since for any functor $S_{\ast}(\mathscr{C})\colon\fromto{N\Delta^{\op}}{\Wald}$ classified by the Waldhausen cocartesian fibration $\fromto{\mathscr{S}(\mathscr{C})}{N\Delta^{\op}}$, the simplicial object $\Phi\circ S_{\ast}(\mathscr{C})$ is a group object, and $\mathscr{F}(\mathscr{C})$ and $\mathscr{S}_0(\mathscr{C})$ are zero objects. Since $\mathscr{S}$ is a suspension functor in $\VaddWald$, we find that the natural transformation $\fromto{\Phi_{\add}}{\Omega_{\mathscr{E}}\circ\Phi_{\add}\circ\Sigma}$ is an equivalence, whence $F_{\add}$ is excisive \cite[Pr. 1.4.2.13]{HA}.

To complete the proof, it remains to show that (\ref{thm:additiveequiv}.\ref{item:excisive}) implies (\ref{thm:additiveequiv}.\ref{item:add}). It follows from (\ref{thm:additiveequiv}.\ref{item:excisive}) that for any nonnegative integer $m$ and any integer $0\leq k\leq m$, applying $\Phi$ to \eqref{eqn:pushpullsquares} yields the same result as applying $\Phi_{\add}$ to \eqref{eqn:Laddpushpullsquares}. Since the lower right square of the latter diagram is a pushout in $\VaddWald$, the excisive functor $F_{\add}$ carries it to a pullback square in $\mathscr{E}$, whence we obtain the first condition.
\end{proof}
\end{thm}

\begin{dfn}\label{dfn:additivefunctor} Suppose $\mathscr{E}$ an $\infty$-topos. An $\mathscr{E}$-valued theory $\phi$ will be said to be \textbf{\emph{additive}} just in case it satisfies any of the equivalent conditions of Th. \ref{thm:additiveequiv}. We denote by $\Add(\mathscr{E})$ the full subcategory of $\Thy(\mathscr{E})$ spanned by the additive theories.
\end{dfn}

Our Structure Theorem (Th. \ref{thm:additiveequiv}) yields an identification of additive theories and excisive functors on fissile virtual Waldhausen $\infty$-categories.
\begin{thm}\label{thm:addfuncsareexc} Suppose $\mathscr{E}$ an $\infty$-topos. The functor $L^{\fiss}\circ j$ induces an equivalence of $\infty$-categories
\begin{equation*}
\equivto{\Exc_{\mathscr{G}}(\VaddWald,\mathscr{E})}{\Add(\mathscr{E})},
\end{equation*}
where $\Exc_{\mathscr{G}}(\VaddWald,\mathscr{E})\subset\Fun^{\ast}(\VaddWald,\mathscr{E})$ is the full subcategory spanned by the reduced excisive functors that preserve small sifted colimits.
\begin{proof} It follows from Th. \ref{thm:additiveequiv} that composition with $L^{\fiss}\circ j$ defines an essentially surjective functor
\begin{equation*}
\fromto{\Exc_{\mathscr{G}}(\VaddWald,\mathscr{E})}{\Add(\mathscr{E})}.
\end{equation*}
To see that this functor is fully faithful, it suffices to note that we have a commutative diagram
\begin{equation*}
\begin{tikzpicture} 
\matrix(m)[matrix of math nodes, 
row sep=4ex, column sep=4ex, 
text height=1.5ex, text depth=0.25ex] 
{\Exc_{\mathscr{G}}(\VaddWald,\mathscr{E})&\Add(\mathscr{E})\\ 
\Fun(\VaddWald,\mathscr{E})&\Fun(\VWald,\mathscr{E})\\}; 
\path[>=stealth,->,font=\scriptsize] 
(m-1-1) edge (m-1-2) 
edge[right hook->] (m-2-1) 
(m-1-2) edge[right hook->] (m-2-2) 
(m-2-1) edge[right hook->] (m-2-2); 
\end{tikzpicture}
\end{equation*}
in which the vertical functors are fully faithful by definition, and the bottom functor is fully faithful because the $\infty$-category $\VaddWald$ is a localization of $\VWald$ (Pr. \ref{prp:Laddexists}).
\end{proof}
\end{thm}
\noindent By virtue of \cite[Pr. 1.4.2.22]{HA}, this result now yields a \emph{canonical} delooping of any additive functor.
\begin{cor}\label{cor:additivevaluedinstab} Suppose $\mathscr{E}$ an $\infty$-topos. Then composition with the canonical functor $\Omega^{\infty}\colon\fromto{\Sp(\mathscr{E})}{\mathscr{E}_{\ast}}$ induces an equivalence of $\infty$-categories
\begin{equation*}
\fromto{\Fun^{\mathrm{rex}}_{\mathscr{G}}(\VaddWald,\Sp(\mathscr{E}))}{\Add(\mathscr{E})},
\end{equation*}
where $\Fun^{\mathrm{rex}}_{\mathscr{G}}(\VaddWald,\Sp(\mathscr{E}))\subset\Fun(\VaddWald,\Sp(\mathscr{E}))$ denotes the full subcategory spanned by the right exact functors $\mathbf{\Phi}\colon\fromto{\VaddWald}{\Sp(\mathscr{E})}$ such that $\Omega^{\infty}\circ\mathbf{\Phi}\colon\fromto{\VaddWald}{\mathscr{E}}$ preserves sifted colimits.
\end{cor}


\subsection*{Additivization} We now find that any theory admits an additive approximation given by a Goodwillie differential. The nature of colimits computed in $\VaddWald$ will then permit us to describe this additive approximation as an $\infty$-categorical $S_{\bullet}$ construction. As a result, we find that any such theory deloops to a \emph{connective} spectrum.

We first need the following well-known lemma, which follows from \cite[Lm. 5.3.6.17]{HA} or, alternately, from a suitable generalization of \cite[Cor. 5.1.3.7]{HA}.
\begin{lem}\label{lem:loopspreservesiftedcolims} For any $\infty$-topos $\mathscr{E}$, the loop functor $\Omega_{\mathscr{E}}\colon\fromto{\mathscr{E}_{\ast}}{\mathscr{E}_{\ast}}$ preserves sifted colimits of connected objects.
\end{lem}

\begin{thm}\label{thm:additivizationexists} Suppose $\mathscr{E}$ an $\infty$-topos. The inclusion functor
\begin{equation*}
\into{\Add(\mathscr{E})}{\Thy(\mathscr{E})}
\end{equation*}
admits a left adjoint $D$ given by a Goodwillie differential \cite{MR1076523,MR1162445,MR2026544}
\begin{equation*}
D\phi\simeq\underset{n\to\infty}{\colim}\ \Omega_{\mathscr{E}}^n\circ\Phi\circ\mathscr{S}^n\circ j,
\end{equation*}
where $\Phi\colon\fromto{\VWald}{\mathscr{E}}$ is the left derived functor of $\phi$.
\begin{proof} Let us write $\mathscr{F}$ for the class of small filtered colimits. By \cite[Th. 1.8]{MR2026544} or \cite[Cor. 7.1.1.10]{HA}, the inclusion
\begin{equation*}
\into{\Exc_{\mathscr{F}}(\VaddWald,\mathscr{E})}{\Fun^{\ast}_{\mathscr{F}}(\VaddWald,\mathscr{E})}
\end{equation*}
(Df. \ref{dfn:reducedtheories}) admits a left adjoint given by the assignment
\begin{equation*}
\goesto{\Phi}{\underset{n\to\infty}{\colim}\ \Omega_{\mathscr{E}}^n\circ\Phi\circ\Sigma_{\VaddWald}^n}.
\end{equation*}
Now the inclusion $i\colon\into{\VaddWald}{\VWald}$ induces a left adjoint
\begin{equation*}
\fromto{\Fun^{\ast}_{\mathscr{F}}(\VWald,\mathscr{E})}{\Fun^{\ast}_{\mathscr{F}}(\VaddWald,\mathscr{E})}
\end{equation*}
to the forgetful functor induced by $L^{\fiss}$. By composing these adjoints, we thus obtain a left adjoint $\DD$ to the forgetful functor
\begin{equation*}
\into{\Exc_{\mathscr{F}}(\VaddWald,\mathscr{E})}{\Fun^{\ast}_{\mathscr{F}}(\VWald,\mathscr{E})}.
\end{equation*}
The left adjoint $\DD$ is given by the assignment
\begin{equation*}
\goesto{\Phi}{\underset{n\to\infty}{\colim}\ \Omega_{\mathscr{E}}^n\circ\Phi\circ i\circ\Sigma_{\VaddWald}^n}.
\end{equation*}
By Cor. \ref{cor:Sdotisreallysuspension}, if $n\geq 1$, then one may rewrite the functor $\Omega_{\mathscr{E}}^n\circ F\circ i\circ\Sigma_{\VaddWald}^n$ as
\begin{equation*}
\Omega_{\mathscr{E}}^n\circ\Phi\circ i\circ\Sigma_{\VaddWald}^n\circ L^{\fiss}\circ i\simeq\Omega_{\mathscr{E}}^n\circ\Phi\circ i\circ\mathscr{S}^n.
\end{equation*}
Now if $\Phi\colon\fromto{\VWald}{\mathscr{E}}$ is the left derived functor of a theory, then for any virtual Waldhausen $\infty$-category $\mathscr{Y}$, since $\Phi$ is reduced, and since $\mathscr{S}(\mathscr{Y})$ is the colimit of a simplicial virtual Waldhausen $\infty$-category $\SSS_{\ast}(\mathscr{Y})$ with $\SSS_0(\mathscr{Y})\simeq 0$, the object $\Phi(\mathscr{S}(\mathscr{Y}))$ is connected as well. By Lm \ref{lem:loopspreservesiftedcolims}, $\Omega_{\mathscr{E}}$ commutes with sifted colimits of connected objects of $\mathscr{E}$, whence it follows that the restriction of $\DD\colon\fromto{\Fun_{\mathscr{F}}^{\star}(\VWald,\mathscr{E})}{\Exc_{\mathscr{F}}(\VaddWald,\mathscr{E})}$ to
\begin{equation*}
\Thy(\mathscr{E})\simeq\Fun_{\mathscr{G}}^{\star}(\VWald,\mathscr{E})\subset\Fun_{\mathscr{F}}^{\star}(\VWald,\mathscr{E})
\end{equation*}
in fact factors through the full subcategory
\begin{equation*}
\Exc_{\mathscr{G}}(\VaddWald,\mathscr{E})\subset\Exc_{\mathscr{F}}(\VaddWald,\mathscr{E}).
\end{equation*}
Thanks to Th. \ref{thm:addfuncsareexc}, the functor $\DD$ consequently descends to a functor
\begin{equation*}
D\colon\fromto{\Thy(\mathscr{E})}{\Add(\mathscr{E})}
\end{equation*}
given by the assignment
\begin{equation*}
\goesto{\Phi}{\underset{n\to\infty}{\colim}\ \Omega_{\mathscr{E}}^n\circ\Phi\circ\Sigma_{\VaddWald}^n\circ L^{\fiss}\circ j}.
\end{equation*}
Now another application of Cor. \ref{cor:Sdotisreallysuspension} completes the proof.
\end{proof}
\end{thm}

\begin{dfn} The left adjoint
\begin{equation*}
D\colon\fromto{\Thy(\mathscr{E})}{\Add(\mathscr{E})}
\end{equation*}
of the previous corollary will be called the \emph{additivization}.
\end{dfn}

Suppose $\phi\colon\fromto{\Wald}{\mathscr{E}}$ a theory; denote by $\Phi$ its left derived functor. For any virtual Waldhausen $\infty$-category $\mathscr{Y}$ and any natural number $n$, since the virtual Waldhausen $\infty$-category $\mathscr{S}^n(\mathscr{Y})$ is the colimit of a reduced $n$-simplicial diagram $\SSS_{\ast}(\SSS_{\ast}(\cdots\SSS_{\ast}(\mathscr{Y})\cdots))$, it follows that the object $\Phi(\mathscr{S}^n(\mathscr{Y}))$ is $n$-connected. This proves the following.
\begin{prp}\label{prp:deloopingofdfisconnective} The canonical delooping \textup{(Cor. \ref{cor:additivevaluedinstab})} of the additivization $D\phi$ of a theory $\phi\colon\fromto{\Wald}{\mathscr{E}_{\ast}}$ is valued in connective spectra:
\begin{equation*}
\fromto{\Wald}{\Sp(\mathscr{E})_{\geq0}}.
\end{equation*}
\end{prp}


\subsection*{Pre-additive theories} We have already mentioned that many of the theories that arise in practice have the property that they carry direct sums of Waldhausen $\infty$-categories to products. What's really useful about theories $\phi$ that enjoy this property is that the colimit
\begin{equation*}
\colim[\phi\ \tikz[baseline]\draw[>=stealth,->](0,0.5ex)--(0.5,0.5ex);\ \Omega\circ\Phi\circ\mathscr{S}\circ j\ \tikz[baseline]\draw[>=stealth,->](0,0.5ex)--(0.5,0.5ex);\ \ \Omega^{2}\circ\Phi\circ\mathscr{S}^{2}\circ j\ \tikz[baseline]\draw[>=stealth,->](0,0.5ex)--(0.5,0.5ex);\ \cdots]
\end{equation*}
that appears in the formula for the additivization (Th. \ref{thm:additivizationexists}) stabilizes after the first term; that is, only one loop is necessary to get an additive theory.
\begin{dfn}\label{dfn:preadditive} Suppose $\mathscr{E}$ an $\infty$-topos. Then a theory $\phi\in\Thy(\mathscr{E})$ is said to be \textbf{\emph{pre-additive}} if it carries direct sums of Waldhausen $\infty$-categories to products in $\mathscr{E}$.
\end{dfn}

\begin{prp} Suppose $\mathscr{E}$ an $\infty$-topos, and suppose $\phi\in\Thy(\mathscr{E})$ a pre-additive theory with left derived functor $\Phi$. Then the morphisms
\begin{equation*}
\fromto{\Phi(\mathscr{S}(\mathscr{F}_m(\mathscr{C})))}{\Phi(\mathscr{S}(\mathscr{C}))}\textrm{\quad and\quad}\fromto{\Phi(\mathscr{S}(\mathscr{F}_m(\mathscr{C})))}{\Phi(\mathscr{S}(\mathscr{S}_m(\mathscr{C})))}
\end{equation*}
induced by $I_{m,0}$ and $F_m$ together exhibit $\Phi(\mathscr{S}(\mathscr{F}_m(\mathscr{C})))$ as a product of $\Phi(\mathscr{S}(\mathscr{C}))$ and $\Phi(\mathscr{S}(\mathscr{S}_m(\mathscr{C})))$.
\begin{proof} Since $\phi$ is pre-additive, the morphism from $\Phi(\mathscr{S}(\mathscr{F}_m(\mathscr{C})))$ to the desired product may be identified with the morphism
\begin{equation*}
\fromto{\Phi(\mathscr{S}(\mathscr{F}_m(\mathscr{C})))}{\Phi(\mathscr{S}(\mathscr{C})\oplus\mathscr{S}(\mathscr{S}_m(\mathscr{C})))},
\end{equation*}
which can in turn be identified with the natural morphism
\begin{equation*}
\fromto{\Phi(i\circ\Sigma_{\VaddWald}\circ L^{\fiss}(\mathscr{F}_m(\mathscr{C})))}{\Phi(i\circ\Sigma_{\VaddWald}\circ L^{\fiss}(\mathscr{C}\oplus\mathscr{S}_m(\mathscr{C})))}
\end{equation*}
by Cor. \ref{cor:Sdotisreallysuspension}. The upper right corner of \eqref{eqn:Laddpushpullsquares} is a pushout, and since $E_m\oplus J_m$ is a section of $I_{m,0}\oplus F_m$, the natural morphism $\fromto{L^{\fiss}(\mathscr{F}_m(\mathscr{C}))}{L^{\fiss}(\mathscr{C}\oplus\mathscr{S}_m(\mathscr{C}))}$ is an equivalence.
\end{proof}
\end{prp}
\noindent By Th. \ref{thm:additiveequiv}, we obtain the following repackaging of Waldhausen's Additivity Theorem.
\begin{cor}\label{for:preaddsareeasy} Suppose $\mathscr{E}$ an $\infty$-topos, and suppose $\phi\in\Thy(\mathscr{E})$ a pre-additive theory with left derived functor $\Phi$. Then the additivization is given by
\begin{equation*}
D\phi\simeq\Omega\circ\Phi\circ\mathscr{S}\circ j.
\end{equation*}
\end{cor} 

Suppose $\mathscr{E}$ an $\infty$-topos, and suppose $\phi\in\Thy(\mathscr{E})$ a pre-additive theory. Then the counit $\fromto{\phi}{D\phi}$ is the initial object of the $\infty$-category $\Add(\mathscr{E})\times_{\Thy(\mathscr{E})}\Thy(\mathscr{E})_{\phi/}$. By Th. \ref{thm:additiveequiv}, this means that $D\phi$ is the initial object of the full subcategory of $\Thy(\mathscr{E})_{\phi/}$ spanned by those natural transformations $\fromto{\phi}{\phi'}$ such that for any Waldhausen $\infty$-category $\mathscr{C}$ and for any functor $\SSS_{\ast}(\mathscr{C})\colon\fromto{N\Delta^{\op}}{\Wald}$ that classifies the Waldhausen cocartesian fibration $\fromto{\mathscr{S}(\mathscr{C})}{N\Delta^{\op}}$, the induced functor $\phi'\circ\SSS_{\ast}(\mathscr{C})\colon\fromto{N\Delta^{\op}}{\mathscr{E}_{\ast}}$ is a group object.

Motivated by this, we may now note that the inclusion of the full subcategory $\Grp(\mathscr{E})$ of $\Fun(N\Delta^{\op},\mathscr{E})$ spanned by the group objects admits a left adjoint $L$. (It is an straightforward matter to note that $\Grp(\mathscr{E})\subset\Fun(N\Delta^{\op},\mathscr{E})$ is stable under arbitrary limits and filtered colimits; alternatively, one may find a small set $S$ of morphisms of $\Fun(N\Delta^{\op},\mathscr{E})$ such that a simplicial object $X$ of $\mathscr{E}$ is a group object if and only if $X$ is $S$-local.) Hence one may consider the following composite functor $L^{\phi}_{\ast}$:
\begin{equation*}
\Wald\ \tikz[baseline]\draw[>=stealth,->,font=\scriptsize](0,0.5ex)--node[above]{$\SSS_{\ast}$}(0.75,0.5ex);\ \Fun(N\Delta^{\op},\Wald)\ \tikz[baseline]\draw[>=stealth,->,font=\scriptsize](0,0.5ex)--node[above]{$\phi$}(0.75,0.5ex);\ \Fun(N\Delta^{\op},\mathscr{E})\ \tikz[baseline]\draw[>=stealth,->,font=\scriptsize](0,0.5ex)--node[above]{$L$}(0.75,0.5ex);\ \Grp(\mathscr{E}).
\end{equation*}
If $\ev_1\colon\fromto{\Grp(\mathscr{E})}{\mathscr{E}}$ is the functor given by evaluation at $1$, then the functor $\ev_1\circ L$ may be identified with the functor $\Omega_{\mathscr{E}}\circ\colim_{N\Delta^{\op}}$. (This is Segal's delooping machine.) It therefore follows from the previous corollary that the functor $L^{\phi}_1=\ev_1\circ L^{\phi}_{\ast}$ can be identified with the additivization of $\phi$. This provides us with a \emph{local} recognition principle for $D\phi$.
\begin{prp}\label{prp:localrecogofDf} Suppose $\mathscr{E}$ an $\infty$-topos, suppose $\phi\in\Thy(\mathscr{E})$ a pre-additive theory, and suppose $\mathscr{C}$ a Waldhausen $\infty$-category. Write 
\begin{equation*}
\SSS_{\ast}(\mathscr{C})\colon\fromto{N\Delta^{\op}}{\Wald}
\end{equation*}
for the functor that classifies the Waldhausen cocartesian fibration $\fromto{\mathscr{S}(\mathscr{C})}{N\Delta^{\op}}$. Then the object $D\phi(\mathscr{C})$ is canonically equivalent to underlying object of the group object that is initial in the $\infty$-category
\begin{equation*}
\Grp(\mathscr{E})\times_{\Fun(N\Delta^{\op},\mathscr{E})}\Fun(N\Delta^{\op},\mathscr{E})_{\phi\circ\SSS_{\ast}(\mathscr{C})/}.
\end{equation*}
\end{prp}

\begin{nul} One may hope to study the rest of the Taylor tower of a theory. In particular, for any positive integer $n$ and any theory $\phi\in\Thy(\mathscr{E})$, one may define a symmetric ``multi-additive'' theory $D^{(n)}\phi$ via a formula
\begin{equation*}
D^{(n)}\phi(\mathscr{C}_1,\dots,\mathscr{C}_n)=\underset{(j_1,\dots,j_n)}{\colim}\ \Omega_{\mathscr{E}}^{j_1+\cdots+j_n}\creff_n\Phi(\mathscr{S}^{j_1}\mathscr{C}_1,\dots,\mathscr{S}^{j_n}\mathscr{C}_n),
\end{equation*}
where $\Phi$ is the left derived functor of $\phi$, and $\creff_n\Phi$ is the $n$-th cross-effect functor of the restriction of $\Phi$ to $\VaddWald$. However, if $\phi$ is pre-additive, then for $n\geq 2$, the cross-effect functor $\creff_n\Phi$ vanishes, whence $D^{(n)}\phi$ vanishes as well. As a result, the Taylor tower for $\Phi$ is constant above the first level. More informally, the best polynomial approximation to $\Phi$ is linear. Consequently, if $\phi\colon\fromto{\Wald}{\mathscr{E}_{\ast}}$ is pre-additive, then $\Phi$ factors through an $n$-excisive functor $\fromto{\VaddWald}{\mathscr{E}_{\ast}}$ for some $n\geq 1$ if and only if $\phi$ is an additive theory, in which case $n$ may be allowed to be $1$. This seems to suggest a rather peculiar dichotomy: a pre-additive theory is either additive or staunchly non-analytic.
\end{nul}


\section{Easy consequences of additivity} Additive theories, which we introduced in the last section, are quite special. In this section, we'll prove some simple results that will illustrate just how special they really are. We'll show that additive theories vanish on any Waldhausen $\infty$-category that is ``too large'' (Pr. \ref{prp:Eilenbergswindle}), and we'll show that additive functors do not distinguish between Waldhausen $\infty$-categories whose pair structure is maximal and suitable stable $\infty$-categories extracted from them. As a side note, we'll remark that, rather curiously, the fissile derived $\infty$-category is only one loop away from being stable. Finally, and most importantly, we'll prove our $\infty$-categorical variant of Waldhausen's Fibration Theorem. In the next section, we'll introduce a richer structure into this story, to prove a more useful variant of this result.


\subsection*{The Eilenberg Swindle} We now show that Waldhausen $\infty$-categories with ``too many'' coproducts are invisible to additive theories.

\begin{prp}[Eilenberg Swindle]\label{prp:Eilenbergswindle} Suppose $\mathscr{E}$ an $\infty$-topos, and suppose $\phi\in\Add(\mathscr{E})$. Then for any Waldhausen $\infty$-category $\mathscr{C}$ that admits countable coproducts, $\phi(\mathscr{C})$ is terminal in $\mathscr{E}$.
\begin{proof} Denote by $I$ the set of natural numbers, regarded as a discrete $\infty$-category, and denote by $\psi\colon\fromto{\mathscr{C}}{\mathscr{C}}$ the composite of the constant functor $\fromto{\mathscr{C}}{\Fun(I,\mathscr{C})}$ followed by its left adjoint $\fromto{\Fun(I,\mathscr{C})}{\mathscr{C}}$. The inclusion $\into{\{0\}}{I}$ and the successor bijection $\sigma\colon\equivto{I}{I-\{0\}}$ together specify a natural ingressive $\cofto{\id}{\psi}$. This defines an exact functor $\fromto{\mathscr{C}}{\mathscr{F}_1(\mathscr{C})}$. Applying $I_{1,1}$ and $I_{1,0}\oplus F_1$ to this functor, we find that $\phi(\psi)=\phi(\id)+\phi(\psi)$, whence $\phi(\id)=0$.
\end{proof}
\end{prp}

 
\subsection*{Stabilization and approximation} We prove that the value of an additive theory on a Waldhausen $\infty$-category whose pair structure is maximal agrees with its value on a certain stable $\infty$-category. Using this, we show that for these Waldhausen $\infty$-categories, equivalences on the homotopy category suffice to give equivalences under any additive theory.

\begin{prp}[Suspension Theorem]\label{prp:suspension} Suppose $\mathscr{A}$ a Waldhausen $\infty$-category whose pair structure is maximal. Then for any additive theory $\phi\in\Add(\mathscr{E})$, the suspension functor $\Sigma\colon\fromto{\mathscr{A}}{\mathscr{A}}$ induces multiplication by $-1$ on the group object $\phi(\mathscr{A})$.
\begin{proof} This follows directly from the existence of the pushout square of endofunctors of $\mathscr{A}$
\begin{equation*}
\begin{tikzpicture} 
\matrix(m)[matrix of math nodes, 
row sep=4ex, column sep=4ex, 
text height=1.5ex, text depth=0.25ex] 
{\id&0\\ 
0&\Sigma.\\}; 
\path[>=stealth,->,font=\scriptsize] 
(m-1-1) edge (m-1-2) 
edge (m-2-1) 
(m-1-2) edge (m-2-2) 
(m-2-1) edge (m-2-2); 
\end{tikzpicture}\qedhere
\end{equation*}
\end{proof}
\end{prp}

\begin{cor}\label{cor:stabilization} Suppose $\mathscr{A}$ a Waldhausen $\infty$-category whose pair structure is maximal. Write $\widetilde{\Sp}(\mathscr{A})$ for the colimit
\begin{equation*}
\mathscr{A}\ \tikz[baseline]\draw[>=stealth,->,font=\scriptsize](0,0.5ex)--node[above]{$\Sigma$}(0.5,0.5ex);\ \mathscr{A}\ \tikz[baseline]\draw[>=stealth,->,font=\scriptsize](0,0.5ex)--node[above]{$\Sigma$}(0.5,0.5ex);\ \cdots\ \tikz[baseline]\draw[>=stealth,->,font=\scriptsize](0,0.5ex)--node[above]{$\Sigma$}(0.5,0.5ex);\ \mathscr{A}\ \tikz[baseline]\draw[>=stealth,->,font=\scriptsize](0,0.5ex)--node[above]{$\Sigma$}(0.5,0.5ex);\ \cdots
\end{equation*}
in $\Wald$. Then for any additive theory $\phi\in\Add(\mathscr{E})$, the canonical functor
\begin{equation*}
\Sigma^{\infty}\colon\fromto{\mathscr{A}}{\widetilde{\Sp}(\mathscr{A})}
\end{equation*}
induces an equivalence $\fromto{\phi(\mathscr{A})}{\phi(\widetilde{\Sp}(\mathscr{A}))}$.
\end{cor}
\noindent We now obtain the following corolary, which we can regard as a version of Waldhausen's Approximation Theorem. Very similar results appear in work of Cisinski \cite[Th. 2.15]{MR2746284} and Blumberg--Mandell \cite[Th. 1.3]{MR2764905}, and an interesting generalization has recently appeared in a preprint of Fiore \cite{fiore}.
\begin{cor}[Approximation]\label{prp:approx} Suppose $\mathscr{C}$ and $\mathscr{D}$ two $\infty$-categories that each contain zero objects and all finite colimits, and regard them as Waldhausen $\infty$-categories equipped with the maximal pair structure \textup{(Ex \ref{exm:whenisminandmaxWald})}. Then any exact functor $\psi\colon\fromto{\mathscr{C}}{\mathscr{D}}$ that induces an equivalence of homotopy categories $\equivto{h\mathscr{C}}{h\mathscr{D}}$ also induces an equivalence $\phi(\psi)\colon\equivto{\phi(\mathscr{C})}{\phi(\mathscr{D})}$ for any additive theory $\phi\in\Add(\mathscr{E})$.
\begin{proof} We note that since the homotopy category functor $\goesto{\mathscr{C}}{h\mathscr{C}}$ preserves colimits, the induced functor $\fromto{h\widetilde{\Sp}(\mathscr{C})}{h\widetilde{\Sp}(\mathscr{D})}$ is an equivalence. Now we combine Prs. \ref{cor:stabilization} and \ref{prp:preapprox}.
\end{proof}
\end{cor}
\noindent The $\infty$-category $\widetilde{\Sp}(\mathscr{A})$ is not always the stabilization of $\mathscr{A}$, but when $\mathscr{A}$ is idempotent complete, it is.
\begin{prp}\label{prp:costabisstab} Suppose $\mathscr{A}$ an idempotent complete $\infty$-category that contains a zero object and all finite colimits. Regard $\mathscr{A}$ as a Waldhausen $\infty$-category with its maximal pair structure. Then $\widetilde{\Sp}(\mathscr{A})$ is equivalent to the stabilization $\Sp(\mathscr{A})$ of $\mathscr{A}$.
\begin{proof} The colimit of the sequence
\begin{equation*}
\mathscr{A}\ \tikz[baseline]\draw[>=stealth,->,font=\scriptsize](0,0.5ex)--node[above]{$\Sigma$}(0.5,0.5ex);\ \mathscr{A}\ \tikz[baseline]\draw[>=stealth,->,font=\scriptsize](0,0.5ex)--node[above]{$\Sigma$}(0.5,0.5ex);\ \cdots\ \tikz[baseline]\draw[>=stealth,->,font=\scriptsize](0,0.5ex)--node[above]{$\Sigma$}(0.5,0.5ex);\ \mathscr{A}\ \tikz[baseline]\draw[>=stealth,->,font=\scriptsize](0,0.5ex)--node[above]{$\Sigma$}(0.5,0.5ex);\ \cdots
\end{equation*}
in $\Wald$ agrees with the same colimit taken in $\Cat_{\infty}(\kappa_1)^{\mathrm{Rex}}$ by \cite[Pr. 5.5.7.11]{HTT} and Pr. \ref{thm:Waldfiltcolims}. Since $\Ind$ is a left adjoint \cite[Pr. 5.5.7.10]{HTT}, the colimit of the sequence
\begin{equation*}
\Ind\mathscr{A}\ \tikz[baseline]\draw[>=stealth,->,font=\scriptsize](0,0.5ex)--node[above]{$\Sigma$}(0.5,0.5ex);\ \Ind\mathscr{A}\ \tikz[baseline]\draw[>=stealth,->,font=\scriptsize](0,0.5ex)--node[above]{$\Sigma$}(0.5,0.5ex);\ \cdots\ \tikz[baseline]\draw[>=stealth,->,font=\scriptsize](0,0.5ex)--node[above]{$\Sigma$}(0.5,0.5ex);\ \Ind\mathscr{A}\ \tikz[baseline]\draw[>=stealth,->,font=\scriptsize](0,0.5ex)--node[above]{$\Sigma$}(0.5,0.5ex);\ \cdots
\end{equation*}
in $\Pr^{\mathrm{L}}_{\omega}$ is $\Ind(\widetilde{\Sp}(\mathscr{A}))$. By \cite[Nt. 5.5.7.7]{HTT}, there is an equivalence between $\Pr^{\mathrm{L}}_{\omega}$ and $(\Pr^{\mathrm{R}}_{\omega})^{\op}$, whence $\Ind(\widetilde{\Sp}(\mathscr{A}))$ can be identified with the limit of the sequence
\begin{equation*}
\cdots\ \tikz[baseline]\draw[>=stealth,->,font=\scriptsize](0,0.5ex)--node[above]{$\Omega$}(0.5,0.5ex);\ \Ind\mathscr{A}\ \tikz[baseline]\draw[>=stealth,->,font=\scriptsize](0,0.5ex)--node[above]{$\Omega$}(0.5,0.5ex);\ \cdots\ \tikz[baseline]\draw[>=stealth,->,font=\scriptsize](0,0.5ex)--node[above]{$\Omega$}(0.5,0.5ex);\ \Ind\mathscr{A}\ \tikz[baseline]\draw[>=stealth,->,font=\scriptsize](0,0.5ex)--node[above]{$\Omega$}(0.5,0.5ex);\ \Ind\mathscr{A}
\end{equation*}
in $\Pr^{\mathrm{R}}_{\omega}$. Since the inclusion $\into{\Pr^{\mathrm{R}}_{\omega}}{\Cat_{\infty}(\kappa_1)}$ preserves limits \cite[Pr. 5.5.7.6]{HTT}, it follows that $\Ind(\widetilde{\Sp}(\mathscr{A}))\simeq\Sp(\Ind(\mathscr{A}))$. Now the functor $\goesto{C}{C^{\omega}}$ is an equivalence of $\infty$-categories between $\Pr^{\mathrm{R}}_{\omega}$ and the full subcategory of $\Cat_{\infty}(\kappa_1)^{\mathrm{Lex}}$ spanned by the essentially small, idempotent complete $\infty$-categories, whence it follows that
\begin{equation*}
\widetilde{\Sp}(\mathscr{A})\simeq\Ind(\widetilde{\Sp}(\mathscr{A}))^{\omega}\simeq\Sp(\Ind(\mathscr{A}))^{\omega}\simeq\Sp(\Ind(\mathscr{A})^{\omega})\simeq\Sp(\mathscr{A}).\qedhere
\end{equation*}
\end{proof}
\end{prp}

\begin{exm}\label{exm:stabinftytopoi} Suppose $\mathscr{E}$ an $\infty$-topos. (One may, again, think of $\Fun(X,\Kan)$ for a simplicial set $X$.) For any additive theory $\phi$, the results above show that one has an equivalence
\begin{equation*}
\phi(\mathscr{E}^{\omega}_{\ast})\simeq\phi(\Sp(\mathscr{E}^{\omega})).
\end{equation*}
\end{exm}


\subsection*{Digression: the near-stability of the fissile derived $\infty$-category} By analyzing the additivization of the Yoneda embedding, we now find that a fissile virtual Waldhausen $\infty$-category is one step away from being an infinite loop object. This implies that the $\infty$-category $\VaddWald$ can be said to admit a much stronger form of the Blakers--Massey excision theorem than the $\infty$-category of spaces. Armed with this, we give an easy necessary and sufficient criterion for a morphism of virtual Waldhausen $\infty$-categories to induce an equivalence on every additive theory.

\begin{dfn} We shall call a theory $\phi\in\Thy(\mathscr{E})$ \textbf{\emph{left exact}} just in case its left derived functor $\Phi$ preserves finite limits.
\end{dfn}
\noindent Clearly every left exact theory is pre-additive. Moreover, the best excisive approximation $P_1(G\circ F)$ to the composite $G\circ F$ of a suitable functor $F\colon\fromto{C}{D}$ with a functor $G\colon\fromto{D}{D'}$ that preserves finite limits is simply the composite $G\circ P_1(F)$. Accordingly, we have the following.

\begin{lem}\label{prp:OmegaSigmaisadditive} Suppose $\phi\in\Thy(\mathscr{E})$ a left exact theory. Then
\begin{equation*}
D\phi\simeq\Phi\circ\Omega_{\VaddWald}\circ\mathscr{S}.
\end{equation*}
\end{lem}

\begin{exm}\label{exm:Yoneda} The Yoneda embedding $y\colon\fromto{\Wald}{\mathscr{P}(\Wald^{\omega})}$ is a left exact theory; its left derived functor $Y\colon\into{\VWald}{\mathscr{P}(\Wald^{\omega})}$ is simply the canonical inclusion. Consequently, thanks to Cor. \ref{for:preaddsareeasy}, the additivization of $y$ is now given by the formula
\begin{equation*}
Dy\simeq\Omega\circ\mathscr{S}\circ j.
\end{equation*}

Let's give some equivalent descriptions of the functor $Dy$. Since $\mathscr{F}(\mathscr{C})$ is contractible, one may write
\begin{equation*}
Dy(\mathscr{C})\simeq\mathscr{F}(\mathscr{C})\times_{\mathscr{S}(\mathscr{C})}\mathscr{F}(\mathscr{C}).
\end{equation*}
Alternately, since suspension in $\VaddWald$ is given by $\mathscr{S}$, the functor
\begin{equation*}
Dy(\mathscr{C})\colon\fromto{\Wald^{\omega,\op}}{\Kan}
\end{equation*}
can be described by the formula
\begin{equation*}
Dy(\mathscr{C})(\mathscr{D})\simeq\Map_{\VWald}(\mathscr{S}(\mathscr{D}),\mathscr{S}(\mathscr{C})).
\end{equation*}

In other words, $\Omega\Sigma\simeq\Omega\mathscr{S}$ is the Goodwillie differential of the identity on $\VaddWald$.
\end{exm}


\subsection*{Waldhausen's Generic Fibration Theorem} Let's now examine the circumstances under which a sequence of virtual Waldhausen $\infty$-categories gives rise to a fiber sequence under any additive functor. In this direction we have Pr. \ref{prp:fibthmi} below, which is an analogue of Waldhausen's \cite[Pr. 1.5.5 and Cor. 1.5.7]{MR86m:18011}. We will deduce from this a necessary and sufficient condition for an exact functor to induce an equivalence under \emph{every} additive theory (Pr. \ref{cor:invertedbyeveryaddthy}).

\begin{ntn} Suppose $\psi\colon\fromto{\mathscr{B}}{\mathscr{A}}$ an exact functor of Waldhausen $\infty$-categories. Write
\begin{equation*}
\mathscr{K}(\psi)\coloneq|\mathscr{F}(\mathscr{A})\times_{\mathscr{S}(\mathscr{A})}\mathscr{S}(\mathscr{B})|_{N\Delta^{\op}}
\end{equation*}
for the realization (Df. \ref{dfn:realization}) of the Waldhausen cocartesian fibration
\begin{equation*}
\fromto{\mathscr{F}(\mathscr{A})\times_{\mathscr{S}(\mathscr{A})}\mathscr{S}(\mathscr{B})}{N\Delta^{\op}}.
\end{equation*}
\end{ntn}

In other words, the virtual Waldhausen $\infty$-category $\mathscr{K}(\psi)$ is the geometric realization of the simplicial Waldhausen $\infty$-category whose $m$-simplices consist of a totally filtered object
\begin{equation*}
0\ \tikz[baseline]\draw[>=stealth,>->](0,0.5ex)--(0.75,0.5ex);\ U_1\ \tikz[baseline]\draw[>=stealth,>->](0,0.5ex)--(0.75,0.5ex);\ U_2\ \tikz[baseline]\draw[>=stealth,>->](0,0.5ex)--(0.75,0.5ex);\ \dots\ \tikz[baseline]\draw[>=stealth,>->](0,0.5ex)--(0.75,0.5ex);\ U_m
\end{equation*}
of $\mathscr{B}$, a filtered object
\begin{equation*}
X_0\ \tikz[baseline]\draw[>=stealth,>->](0,0.5ex)--(0.75,0.5ex);\ X_1\ \tikz[baseline]\draw[>=stealth,>->](0,0.5ex)--(0.75,0.5ex);\ X_2\ \tikz[baseline]\draw[>=stealth,>->](0,0.5ex)--(0.75,0.5ex);\ \dots\ \tikz[baseline]\draw[>=stealth,>->](0,0.5ex)--(0.75,0.5ex);\ X_m
\end{equation*}
of $\mathscr{A}$, and a diagram
\begin{equation*}
\begin{tikzpicture} 
\matrix(m)[matrix of math nodes, 
row sep=4ex, column sep=4ex, 
text height=1.5ex, text depth=0.25ex] 
{X_0&X_1&X_2&\dots&X_m\\ 
0&\psi(U_1)&\psi(U_2)&\dots&\psi(U_m)\\}; 
\path[>=stealth,->,font=\scriptsize] 
(m-1-1) edge[>->] (m-1-2) 
edge (m-2-1)
(m-1-2) edge[>->] (m-1-3)
edge (m-2-2)
(m-1-3) edge[>->] (m-1-4)
edge (m-2-3)
(m-1-4) edge[>->] (m-1-5)
(m-1-5) edge (m-2-5)
(m-2-1) edge[>->] (m-2-2) 
(m-2-2) edge[>->] (m-2-3)
(m-2-3) edge[>->] (m-2-4)
(m-2-4) edge[>->] (m-2-5); 
\end{tikzpicture}
\end{equation*}
of $\mathscr{A}$ in which every square is a pushout.

The object $\mathscr{K}(\psi)$ is not itself the corresponding fiber product of virtual Waldhausen $\infty$-categories; however, for any additive functor $\phi\colon\fromto{\Wald}{\mathscr{E}}$ with left derived functor $\Phi$, we shall now show that $\Phi(\mathscr{K}(\psi))$ is in fact the fiber of the induced morphism $\fromto{\Phi(\mathscr{S}(\mathscr{B}))}{\Phi(\mathscr{S}(\mathscr{A}))}$.

\begin{thm}[Generic Fibration Theorem I]\label{prp:fibthmi} Suppose $\psi\colon\fromto{\mathscr{B}}{\mathscr{A}}$ an exact functor of Waldhausen $\infty$-categories. Then for any additive theory $\phi\colon\fromto{\Wald}{\mathscr{E}}$ with left derived functor $\Phi$, there is a diagram
\begin{equation*}
\begin{tikzpicture} 
\matrix(m)[matrix of math nodes, 
row sep=4ex, column sep=4ex, 
text height=1.5ex, text depth=0.25ex] 
{\phi(\mathscr{B})&\phi(\mathscr{A})&\ast\\ 
\ast&\Phi(\mathscr{K}(\psi))&\Phi(\mathscr{S}(\mathscr{B}))\\
&\ast&\Phi(\mathscr{S}(\mathscr{A}))\\}; 
\path[>=stealth,->,font=\scriptsize] 
(m-1-1) edge (m-1-2) 
edge (m-2-1)
(m-1-2) edge (m-2-2)
edge (m-1-3)
(m-1-3) edge (m-2-3)
(m-2-1) edge (m-2-2)
(m-2-2) edge (m-2-3)
edge (m-3-2)
(m-2-3) edge (m-3-3)
(m-3-2) edge (m-3-3); 
\end{tikzpicture}
\end{equation*}
of $\mathscr{E}$ in which each square is a pullback.
\begin{proof} For any vertex $\mathbf{m}\in N\Delta^{\op}$, there exist functors 
\begin{equation*}
s\coloneq(E_m\oplus\mathscr{S}_m(\psi),\pr_2)\colon\fromto{\mathscr{F}_0(\mathscr{A})\oplus\mathscr{S}_m(\mathscr{B})}{\mathscr{F}_m(\mathscr{A})\times_{\mathscr{S}_m(\mathscr{A})}\mathscr{S}_m(\mathscr{B})}
\end{equation*}
and
\begin{equation*}
p\coloneq(I_{m,0}\circ\pr_1)\oplus\pr_2\colon\fromto{\mathscr{F}_m(\mathscr{A})\times_{\mathscr{S}_m(\mathscr{A})}\mathscr{S}_m(\mathscr{B})}{\mathscr{F}_0(\mathscr{A})\oplus\mathscr{S}_m(\mathscr{B})}.
\end{equation*}
Clearly $p\circ s\simeq\id$; we claim that $\phi(s\circ p)\simeq\phi(\id)$ in $\mathscr{E}_{\ast}$. This follows from additivity applied to the functor
\begin{equation*}
\fromto{\mathscr{F}_m(\mathscr{A})\times_{\mathscr{S}_m(\mathscr{A})}\mathscr{S}_m(\mathscr{B})}{\mathscr{F}_1(\mathscr{F}_m(\mathscr{A})\times_{\mathscr{S}_m(\mathscr{A})}\mathscr{S}_m(\mathscr{B}))}
\end{equation*}
given by the ingressive morphism of functors $\cofto{(E_m\circ I_{m,0}\circ\pr_1,0)}{\id}$. Thus the value $\phi(\mathscr{F}_m(\mathscr{A})\times_{\mathscr{S}_m(\mathscr{A})}\mathscr{S}_m(\mathscr{B}))$ is exhibited as the product $\phi(\mathscr{F}_0(\mathscr{A}))\times\phi(\mathscr{S}_m(\mathscr{B}))$.

We may therefore consider the following commutative diagram of $\mathscr{E}_{\ast}$:
\begin{equation*}
\begin{tikzpicture}[baseline]
\matrix(m)[matrix of math nodes, 
row sep=4ex, column sep=5ex, 
text height=1.5ex, text depth=0.25ex] 
{\phi(\mathscr{F}_0(\mathscr{B}))&\phi(\mathscr{F}_0(\mathscr{A}))&\phi(\mathscr{S}_0(\mathscr{B}))\\ 
\phi(\mathscr{F}_m(\mathscr{B}))&\phi(\mathscr{F}_m(\mathscr{A})\times_{\mathscr{S}_m(\mathscr{A})}\mathscr{S}_m(\mathscr{B}))&\phi(\mathscr{S}_m(\mathscr{B}))\\
&\phi(\mathscr{F}_m(\mathscr{A}))&\phi(\mathscr{S}_m(\mathscr{A}))\\
&\phi(\mathscr{F}_0(\mathscr{A}))&\phi(\mathscr{S}_0(\mathscr{A})).\\}; 
\path[>=stealth,->,font=\scriptsize] 
(m-1-1) edge node[above]{$\mathscr{F}_0(\psi)$} (m-1-2) 
edge node[left]{$E_m$} (m-2-1) 
(m-1-2) edge node[above]{$F_0$} (m-1-3)
edge (m-2-2) 
(m-1-3) edge node[right]{$E'_m$} (m-2-3)
(m-2-1) edge node[below]{$(0,F_m)$} (m-2-2) 
(m-2-2) edge node[above]{$\pr_1$} (m-2-3)
edge node[left]{$\pr_2$} (m-3-2) 
(m-2-3) edge node[right]{$\mathscr{S}_m(\psi)$} (m-3-3)
(m-3-2) edge node[below]{$F_m$} (m-3-3)
edge node[left]{$I_{m,0}$} (m-4-2)
(m-3-3) edge node[right]{$I'_{m,0}$} (m-4-3)
(m-4-2) edge node[below]{$F_0$} (m-4-3); 
\end{tikzpicture}
\end{equation*}
The lower right-hand square is a pullback square by additivity; hence, in light of the identification above, all the squares on the right hand side are pullbacks as well. Again by additivity the wide rectangle of the top row is carried to a pullback square under $\phi$, whence all the squares of this diagram are carried to pullback squares.

Since $\phi$ is additive, so is $\Phi\circ\mathscr{S}$. Hence we obtain a commutative diagram in $\mathscr{E}_{\ast}$:
\begin{equation*}
\begin{tikzpicture}[baseline]
\matrix(m)[matrix of math nodes, 
row sep=4ex, column sep=4ex, 
text height=1.5ex, text depth=0.25ex] 
{\Phi(\mathscr{S}(\mathscr{F}_0(\mathscr{B})))&\Phi(\mathscr{S}(\mathscr{F}_0(\mathscr{A})))&\Phi(\mathscr{S}(\mathscr{S}_0(\mathscr{B})))\\ 
\Phi(\mathscr{S}(\mathscr{F}_m(\mathscr{B})))&\Phi(\mathscr{S}(\mathscr{F}_m(\mathscr{A})\times_{\mathscr{S}_m(\mathscr{A})}\mathscr{S}_m(\mathscr{B})))&\Phi(\mathscr{S}(\mathscr{S}_m(\mathscr{B})))\\
&\Phi(\mathscr{S}(\mathscr{F}_m(\mathscr{A})))&\Phi(\mathscr{S}(\mathscr{S}_m(\mathscr{A}))),\\}; 
\path[>=stealth,->,font=\scriptsize] 
(m-1-1) edge (m-1-2) 
edge (m-2-1) 
(m-1-2) edge (m-1-3)
edge (m-2-2) 
(m-1-3) edge (m-2-3)
(m-2-1) edge (m-2-2) 
(m-2-2) edge (m-2-3)
edge (m-3-2) 
(m-2-3) edge (m-3-3)
(m-3-2) edge (m-3-3); 
\end{tikzpicture}
\end{equation*}
in which every square is a pullback. All the squares in this diagram are functorial in $\mathbf{m}$, and since the objects that appear are all connected, it follows from \cite[Lm. 5.3.6.17]{HA} that the squares of the colimit diagram
\begin{equation*}
\begin{tikzpicture}[baseline]
\matrix(m)[matrix of math nodes, 
row sep=4ex, column sep=4ex, 
text height=1.5ex, text depth=0.25ex] 
{\Phi(\mathscr{S}(\mathscr{F}_0(\mathscr{B})))&\Phi(\mathscr{S}(\mathscr{F}_0(\mathscr{A})))&\Phi(\mathscr{S}(\mathscr{S}_0(\mathscr{B})))\\ 
\Phi(\mathscr{S}(\mathscr{F}(\mathscr{B})))&\Phi(\mathscr{S}(\mathscr{K}(\psi)))&\Phi(\mathscr{S}(\mathscr{S}(\mathscr{B})))\\
&\Phi(\mathscr{S}(\mathscr{F}(\mathscr{A})))&\Phi(\mathscr{S}(\mathscr{S}(\mathscr{A}))),\\}; 
\path[>=stealth,->,font=\scriptsize] 
(m-1-1) edge (m-1-2) 
edge (m-2-1) 
(m-1-2) edge (m-1-3)
edge (m-2-2) 
(m-1-3) edge (m-2-3)
(m-2-1) edge (m-2-2) 
(m-2-2) edge (m-2-3)
edge (m-3-2) 
(m-2-3) edge (m-3-3)
(m-3-2) edge (m-3-3); 
\end{tikzpicture}
\end{equation*}
are all pullbacks. Applying the loopspace functor $\Omega_{\mathscr{E}}$ to this diagram now produces a diagram equivalent to the diagram
\begin{equation*}
\begin{tikzpicture}[baseline]
\matrix(m)[matrix of math nodes, 
row sep=4ex, column sep=4ex, 
text height=1.5ex, text depth=0.25ex] 
{\phi(\mathscr{F}_0(\mathscr{B}))&\phi(\mathscr{F}_0(\mathscr{A}))&\phi(\mathscr{S}_0(\mathscr{B}))\\ 
\Phi(\mathscr{F}(\mathscr{B}))&\Phi(\mathscr{K}(\psi))&\Phi(\mathscr{S}(\mathscr{B}))\\
&\Phi(\mathscr{F}(\mathscr{A}))&\Phi(\mathscr{S}(\mathscr{A})),\\}; 
\path[>=stealth,->,font=\scriptsize] 
(m-1-1) edge (m-1-2) 
edge (m-2-1) 
(m-1-2) edge (m-1-3)
edge (m-2-2) 
(m-1-3) edge (m-2-3)
(m-2-1) edge (m-2-2) 
(m-2-2) edge (m-2-3)
edge (m-3-2) 
(m-2-3) edge (m-3-3)
(m-3-2) edge (m-3-3); 
\end{tikzpicture}
\end{equation*}
in which every square again is a pullback.
\end{proof}
\end{thm}

\begin{prp}\label{cor:invertedbyeveryaddthy} The following are equivalent for an exact functor $\psi\colon\fromto{\mathscr{B}}{\mathscr{A}}$ of Waldhausen $\infty$-categories.
\begin{enumerate}[(\ref{cor:invertedbyeveryaddthy}.1)]
\item\label{item:addequiv} For any $\infty$-topos $\mathscr{E}$ and any $\phi\in\Add(\mathscr{E})$ with left derived functor
\begin{equation*}
\Phi\colon\fromto{\VWald}{\mathscr{E}},
\end{equation*}
the induced morphism $\Phi(\psi)\colon\fromto{\Phi(\mathscr{B})}{\Phi(\mathscr{A})}$ is an equivalence of $\mathscr{E}$.
\item\label{item:gpfiber} For any $\infty$-topos $\mathscr{E}$ and any $\phi\in\Add(\mathscr{E})$ with left derived functor
\begin{equation*}
\Phi\colon\fromto{\VWald}{\mathscr{E}},
\end{equation*}
the object $\Phi(\mathscr{K}(\psi))$ is contractible.
\item\label{item:SKpsiiszero} The virtual Waldhausen $\infty$-category $\mathscr{S}(\mathscr{K}(\psi))$ is contractible.
\end{enumerate}
\begin{proof} In light of Lm. \ref{prp:OmegaSigmaisadditive} and Ex. \ref{exm:Yoneda}, if (\ref{cor:invertedbyeveryaddthy}.\ref{item:addequiv}) holds, then the induced morphism
\begin{equation*}
\Omega\mathscr{S}(\psi)\colon\fromto{\Omega\mathscr{S}(\mathscr{B})}{\Omega\mathscr{S}(\mathscr{A})}
\end{equation*}
is an equivalence of virtual Waldhausen $\infty$-categories. Since $\mathscr{S}(\mathscr{B})$ and $\mathscr{S}(\mathscr{A})$ are connected objects of $\mathscr{P}(\Wald^{\omega})$, this in turn implies (using, say, \cite[Cor. 5.1.3.7]{HA}) that the induced morphism of virtual Waldhausen $\infty$-categories
\begin{equation*}
\mathscr{S}(\psi)\colon\fromto{\mathscr{S}(\mathscr{B})}{\mathscr{S}(\mathscr{A})}
\end{equation*}
is an equivalence and therefore by Pr. \ref{prp:fibthmi} that (\ref{cor:invertedbyeveryaddthy}.\ref{item:gpfiber}) holds.

Now if (\ref{cor:invertedbyeveryaddthy}.\ref{item:gpfiber}) holds, then in particular, $\Omega\mathscr{S}(\mathscr{K}(\psi))$ is contractible. Since the virtual Waldhausen $\infty$-category $\mathscr{S}(\mathscr{K}(\psi))$ is connected, it is contractible, yielding (\ref{cor:invertedbyeveryaddthy}.\ref{item:SKpsiiszero}).

That the last condition implies the first now follows immediately from Pr. \ref{prp:fibthmi}.
\end{proof}
\end{prp}


\section{Labeled Waldhausen $\infty$-categories and Waldhausen's Fibration Theorem}\label{subsect:labeledWald} We have remarked (Ex. \ref{exm:catwithcofibsisWaldinftycat}) that nerves of Waldhausen's categories with cofibrations are natural examples of Waldhausen $\infty$-categories. But Waldhausen's categories with cofibrations and weak equivalences don't fit so easily into this story. One may attempt to form the relative nerve (Df. \ref{dfn:relnerve}) of the underlying relative category and to endow the resulting $\infty$-category with a suitable pair structure, but part of the point of Waldhausen's set-up was precisely that one didn't need to assume things such as the two-out-of-three axiom. For example, Waldhausen considers (\cite[Part 3]{MR86m:18011} or \cite{WJR:stabparhcob}) categories of spaces in which the weak equivalences are chosen to be the simple maps. In these situations the $K$-theory of the relative nerve will not correctly encode the Waldhausen $K$-theory.

The time has come to address this issue. Fortunately, the machinery we have developed provides a useful alternative. Namely, we introduce the notion of a \emph{labeled Waldhausen $\infty$-category} (Df. \ref{dfn:labeledWaldcats}), which is a Waldhausen $\infty$-category equipped with a subcategory of \emph{labeled edges} that satisfy the analogue of Waldhausen's axioms for a category with cofibrations and weak equivalences. There is a relative form of this, too, as an example, we show how to label Waldhausen cocartesian fibrations of filtered objects.

It is possible to extract from these categories with cofibrations and weak equivalences useful \emph{virtual Waldhausen $\infty$-categories} (Nt. \ref{ntn:BasavirtWaldcat}). These virtual Waldhausen $\infty$-categories are constructed as realizations of certain Waldhausen cocartesian fibrations over $N\Delta^{\op}$; they are not Waldhausen $\infty$-categories, but they are ``close'' (Pr. \ref{prp:BCwCiswFunblankC}). We also discuss the relationship between the virtual Waldhausen $\infty$-categories attached to a labeled Waldhausen $\infty$-category and the result from formally inverting (in the $\infty$-categorical sense, of course) the labeled edges (Nt. \ref{ntn:gammaC}).

The main result of this section is a familiar case of the Generic Fibration Theorem I (Th. \ref{prp:fibthmi}). This result (Th. \ref{thm:fibration}) gives, for any labeled Waldhausen $\infty$-category $(\mathscr{A},w\mathscr{A})$ satisfying a certain compatibility between the ingressives and the labeled edges (Df. \ref{dfn:enoughcofibrations}) a fiber sequence
\begin{equation*}
\phi(\mathscr{A}^w)\ \tikz[baseline]\draw[>=stealth,->](0,0.5ex)--(0.5,0.5ex);\ \phi(\mathscr{A})\ \tikz[baseline]\draw[>=stealth,->](0,0.5ex)--(0.5,0.5ex);\ \Phi(\mathscr{B}(\mathscr{A},w\mathscr{A}))
\end{equation*}
for any additive theory $\phi$ with left derived functor $\Phi$. This result is the foundation of virtually all fiber sequences that arise in $K$-theory.


\subsection*{Labeled Waldhausen $\infty$-categories} In analogy with Waldhausen's theory of categories with cofibrations and weak equivalences, we study here Waldhausen $\infty$-categories with certain compatible classes of \emph{labeled morphisms}.

\begin{dfn}\label{dfn:labeledWaldcats} Suppose $\mathscr{C}$ a Waldhausen $\infty$-category. Then a \textbf{\emph{gluing diagram}} in $\mathscr{C}$ is a functor of pairs
\begin{equation*}
X\colon\fromto{\mathscr{Q}^2\times(\Delta^1)^{\flat}}{\mathscr{C}}
\end{equation*}
(Ex. \ref{exm:minpairmaxpair} and \ref{exm:LambdaDelta}) such that the squares $X|_{(\mathscr{Q}^2\times\Delta^{\{0\}})}$ and $X|_{(\mathscr{Q}^2\times\Delta^{\{1\}})}$ are pushouts. We may depict such gluing diagrams as cubes
\begin{equation*}
\begin{tikzpicture}[cross line/.style={preaction={draw=white, -, 
line width=6pt}}]
\matrix(m)[matrix of math nodes, 
row sep=2ex, column sep=1.5ex, 
text height=1.5ex, text depth=0.25ex]
{&X_{00}&&X_{10}\\
X_{20}&&X_{\infty0}&\\
&X_{01}&&X_{11}\\
X_{21}&&X_{\infty1}&\\
};
\path[>=stealth,->,font=\scriptsize]
(m-1-2) edge (m-2-1)
edge (m-3-2)
edge[>->] (m-1-4)
(m-3-2) edge (m-4-1)
edge[>->] (m-3-4)
(m-2-1) edge[cross line,>->] (m-2-3)
edge (m-4-1)
(m-1-4) edge (m-2-3)
edge (m-3-4)
(m-4-1) edge[>->] (m-4-3)
(m-3-4) edge (m-4-3)
(m-2-3) edge[cross line] (m-4-3);
\end{tikzpicture}
\end{equation*}
in which the top and bottom faces are pushout squares.
\end{dfn}

\begin{dfn}\label{dfn:labeling} A \textbf{\emph{labeling}} of a Waldhausen $\infty$-category is a subcategory $w\mathscr{C}$ of $\mathscr{C}$ that contains $\iota\mathscr{C}$ (i.e., a a pair structure on $\mathscr{C}$) such that for any gluing diagram $X$ of $\mathscr{C}$ in which the morphisms
\begin{equation*}
\fromto{X_{00}}{X_{01}}\textrm{,\qquad}\fromto{X_{10}}{X_{11}}\textrm{,\quad and\quad}\fromto{X_{20}}{X_{21}}
\end{equation*}
lie in $w\mathscr{C}$, the morphism $\fromto{X_{\infty0}}{X_{\infty1}}$ lies in $w\mathscr{C}$ as well. In this case, the edges of $w\mathscr{C}$ will be called \textbf{\emph{labeled edges}}, and the pair $(\mathscr{C},w\mathscr{C})$ is called a \textbf{\emph{labeled Waldhausen $\infty$-category}}.

A \textbf{\emph{labeled exact functor}} between two labeled Waldhausen $\infty$-categories $\mathscr{C}$ and $\mathscr{D}$ is an exact functor $\fromto{\mathscr{C}}{\mathscr{D}}$ that carries labeled edges to labeled edges.
\end{dfn}

Note that a labeled Waldhausen $\infty$-category has two pair structures: the ingressives and the labeled edges.

\begin{exm}\label{exm:catwithcofibsandwesislabWaldinftycat} We have remarked (Ex. \ref{exm:catwithcofibsisWaldinftycat}) that the nerve of an ordinary \emph{category with cofibrations} in the sense of Waldhausen is a Waldhausen $\infty$-category. Similarly, if $(C,\cof C,wC)$ is a \emph{category with cofibrations and weak equivalences} in the sense of Waldhausen \cite[\S 1.2]{MR86m:18011}, then $(NC,N\cof C,NwC)$ is a labeled Waldhausen $\infty$-category.
\end{exm}

Suppose $(\mathscr{C},w\mathscr{C})$ a labeled Waldhausen $\infty$-category. For gluing diagrams $X$ of $\mathscr{C}$ in which the edges
\begin{eqnarray}
\fromto{X_{00}}{X_{20}},&&\fromto{X_{00}}{X_{01}},\nonumber\\
\fromto{X_{10}}{X_{\infty0}},&&\fromto{X_{10}}{X_{11}}\nonumber
\end{eqnarray}
are all degenerate, the condition above reduces to a guarantee that pushouts of labeled morphism along ingressive morphisms are labeled. For gluing diagrams $X$ of $\mathscr{C}$ in which the edges
\begin{eqnarray}
\fromto{X_{00}}{X_{10}},&&\fromto{X_{00}}{X_{01}},\nonumber\\
\fromto{X_{20}}{X_{\infty0}},&&\fromto{X_{20}}{X_{21}}\nonumber
\end{eqnarray}
are all degenerate, the condition above reduces to a guarantee that the pushout of any labeled ingressive morphism along any morphism exists and is again a labeled ingressive morphism.

\begin{ntn} Denote by
\begin{equation*}
\lWald\subset\Wald\times_{\Cat_{\infty}}\Pair_{\infty}
\end{equation*}
the full subcategory spanned by the labeled Waldhausen $\infty$-categories.
\end{ntn}

\begin{prp}\label{prp:lWaldpres} The $\infty$-category $\lWald$ is presentable.
\begin{proof} The inclusion
\begin{equation*}
\into{\lWald}{\Wald\times_{\Cat_{\infty}}\Pair_{\infty}}
\end{equation*}
admits a left adjoint, which assigns to any object $(\mathscr{C},\mathscr{C}_{\dag},w\mathscr{C})$ the labeled Waldhausen $\infty$-category $(\mathscr{C},\mathscr{C}_{\dag},\overline{w}\mathscr{C})$, where $\overline{w}\mathscr{C}$ is the smallest labeling containing $w\mathscr{C}$. It is easy to see that $\lWald$ is stable under filtered colimits in $\Wald\times_{\Cat_{\infty}}\Pair_{\infty}$; hence $\lWald$ is an accessible localization of $\Wald\times_{\Cat_{\infty}}\Pair_{\infty}$. Since the latter $\infty$-category is locally presentable by \cite[Pr. 5.5.7.6]{HTT}, the proof is complete.
\end{proof}
\end{prp}


\subsection*{The Waldhausen cocartesian fibration attached to a labeled Waldhausen $\infty$-category} In \S \ref{sect:filt}, we defined the virtual Waldhausen $\infty$-category of filtered objects of a Waldhausen $\infty$-category $\mathscr{C}$. We did this by first using Pr. \ref{prp:htt32213} to write down a cocartesian fibration that is classified by the simplicial $\infty$-category
\begin{equation*}
\FF_{\ast}(\mathscr{C})\colon\fromto{N\Delta^{\op}}{\Cat}
\end{equation*}
such that for any integer $m\geq 0$, the $\infty$-category $\FF_m(\mathscr{C})$ has as objects sequences of ingressive morphisms
\begin{equation*}
X_0\ \tikz[baseline]\draw[>=stealth,>->](0,0.5ex)--(0.5,0.5ex);\ X_1\ \tikz[baseline]\draw[>=stealth,>->](0,0.5ex)--(0.5,0.5ex);\ \cdots\ \tikz[baseline]\draw[>=stealth,>->](0,0.5ex)--(0.5,0.5ex);\ X_m.
\end{equation*}
Then we defined the virtual Waldhausen $\infty$-category we were after by forming the formal geometric realization of the diagram $\FF_{\ast}(\mathscr{C})$.

Here, we introduce an analogous construction when $\mathscr{C}$ admits a labeling, in which the role of the cofibrations is played instead by the labeled edges. That is, we will define a cocartesian fibration $\fromto{\mathscr{B}(\mathscr{C},w\mathscr{C})}{N\Delta^{\op}}$ that is classified by the simplicial $\infty$-category
\begin{equation*}
\BB_{\ast}(\mathscr{C},w\mathscr{C})\colon\fromto{N\Delta^{\op}}{\Cat_{\infty}}
\end{equation*}
such that for any integer $m\geq 0$, the $\infty$-category $\BB_m(\mathscr{C},w\mathscr{C})$ has as objects sequences of labeled edges
\begin{equation*}
X_0\ \tikz[baseline]\draw[>=stealth,->,font=\scriptsize,inner sep=0.5pt](0,0.5ex)--(0.5,0.5ex);\ X_1\ \tikz[baseline]\draw[>=stealth,->,font=\scriptsize,inner sep=0.5pt](0,0.5ex)--(0.5,0.5ex);\ \cdots\ \tikz[baseline]\draw[>=stealth,->,font=\scriptsize,inner sep=0.5pt](0,0.5ex)--(0.5,0.5ex);\ X_m.
\end{equation*}
The pair structure will be simpler than in \S \ref{sect:filt}, but once again we will define the virtual Waldhausen $\infty$-category we're after by forming the formal geometric realization of the diagram $\BB_{\ast}(\mathscr{C},w\mathscr{C})$.
\begin{cnstr}\label{cnstr:BXSwXs} Suppose $\mathscr{C}$ a Waldhausen $\infty$-category, and suppose $w\mathscr{C}\subset\mathscr{C}$ a labeling thereof. Define a map $\fromto{\mathscr{B}(\mathscr{C},w\mathscr{C})}{N\Delta^{\op}}$, using the notation of Pr. \ref{prp:htt32213}, Ex. \ref{exm:minpairmaxpair}, and Nt. \ref{ntn:M}, as
\begin{equation*}
\mathscr{B}(\mathscr{C},w\mathscr{C})\coloneq T_{\pi}((N\Delta^{\op})^{\flat}\times(\mathscr{C},w\mathscr{C})).
\end{equation*}
Equivalently, we require, for any simplicial set $K$ and any map $\sigma\colon\fromto{K}{N\Delta^{\op}}$, a bijection between the set $\Mor_{N\Delta^{\op}}(K,\mathscr{B}(\mathscr{C},w\mathscr{C}))$ and the set
\begin{equation*}
\Mor_{s\Set(2)}((K\times_{N\Delta^{\op}}N\mathrm{M},K\times_{N\Delta^{\op}}(N\mathrm{M})_{\dag}),(\mathscr{C},w\mathscr{C}))
\end{equation*}
(Nt. \ref{ntn:ordcatofpairs}), functorial in $\sigma$.

In other words, $\mathscr{B}(\mathscr{C},w\mathscr{C})$ is the simplicial set $\mathscr{F}(\mathscr{C},w\mathscr{C})$, where $\mathscr{C}$ is regarded as a pair with its subcategory of \emph{labeled edges}, rather than its subcategory of cofibrations.
\end{cnstr}

It follows from \ref{prp:htt32213} that $\fromto{\mathscr{B}(\mathscr{C},w\mathscr{C})}{N\Delta^{\op}}$ is a cocartesian fibration.

\begin{nul} For any Waldhausen $\infty$-category $\mathscr{C}$ and any labeling $w\mathscr{C}\subset\mathscr{C}$ thereof, we endow the $\infty$-category $\mathscr{B}(\mathscr{C},w\mathscr{C})$ with a pair structure in the following manner. We let $\mathscr{B}_{\dag}(\mathscr{C},w\mathscr{C})$ be the smallest pair structure containing morphisms of the form $(\id,\psi)\colon\fromto{(\mathbf{m},Y)}{(\mathbf{m},X)}$, where for any integer $0\leq k\leq m$, the induced morphism $\fromto{Y_k}{X_k}$ is ingressive.
\end{nul}

\begin{lem} For any Waldhausen $\infty$-category $\mathscr{C}$ and any labeling $w\mathscr{C}\subset\mathscr{C}$ thereof, the cocartesian fibration $p\colon\fromto{\mathscr{B}(\mathscr{C},w\mathscr{C})}{N\Delta^{\op}}$ is a Waldhausen cocartesian fibration.
\begin{proof} It is plain to see that $p$ is a pair cocartesian fibration.

Now suppose $m\geq0$ an integer. Since limits and colimits of the $\infty$-category $\Fun(\Delta^m,\mathscr{C})$ are computed pointwise, a zero object in $\Fun(\Delta^m,\mathscr{C})$ is an essentially constant functor whose value at any point of $\Delta^m$ is a zero object. Since any equivalence of $\mathscr{C}$ is contained in $w\mathscr{C}$, this zero object is contained in $\mathscr{B}(\mathscr{C},w\mathscr{C})_m$ as well. Again since pushouts in $\Fun(\Delta^m,\mathscr{C})$ are formed objectwise, a pushout square in $\Fun(\Delta^m,\mathscr{C})$ is a functor
\begin{equation*}
X\colon\fromto{\Delta^1\times\Delta^1\times\Delta^{\{0,k\}}}{\mathscr{C}}
\end{equation*}
such that for any integer $0\leq k\leq m$, the restriction $X|_{(\Delta^1\times\Delta^1\times\Delta^{\{0,k\}})}$ is a pushout square; now if $X$ is in addition a functor of pairs $\fromto{\mathscr{Q}^2\times(\Delta^m)^{\flat}}{\mathscr{C}}$, then it follows from the gluing axiom that if $X|_{(\{0\}\times\Delta^{\{0,k\}})}$, $X|_{(\{1\}\times\Delta^{\{0,k\}})}$, and $X|_{(\{2\}\times\Delta^{\{0,k\}})}$ all factor through $w\mathscr{C}\subset\mathscr{C}$, then so does $X|_{(\{\infty\}\times\Delta^{\{0,k\}})}$. Hence the fibers $\mathscr{B}_m(\mathscr{C},w\mathscr{C})$ of $p$ are Waldhausen $\infty$-categories, and, again using the fact that colimits and limits are computed objectwise, we conclude that $p$ is a Waldhausen cocartesian fibration.
\end{proof} 
\end{lem}


\subsection*{The virtual Waldhausen $\infty$-category attached to a labeled Waldhausen $\infty$-category} It follows from \ref{nul:Tpafunctor} that the assignment
\begin{equation*}
\goesto{(\mathscr{C},w\mathscr{C})}{\mathscr{B}(\mathscr{C},w\mathscr{C})}
\end{equation*}
defines a functor
\begin{equation*}
\mathscr{B}\colon\fromto{\lWald}{\mathbf{Wald}_{\infty/N\Delta^{\op}}^{\cocart}}.
\end{equation*}
By composing with the realization functor (Df. \ref{dfn:realization}), we find a functorial construction of virtual Waldhausen $\infty$-categories from labeled Waldhausen $\infty$-categories:
\begin{ntn}\label{ntn:BasavirtWaldcat} By a small abuse of notation, we denote also as $\mathscr{B}$ the composite functor
\begin{equation*}
\lWald\ \tikz[baseline]\draw[>=stealth,->,font=\scriptsize](0,0.5ex)--node[above]{$\mathscr{B}$}(0.5,0.5ex);\ \mathbf{Wald}_{\infty/N\Delta^{\op}}^{\cocart}\ \tikz[baseline]\draw[>=stealth,->,font=\scriptsize](0,0.5ex)--node[above]{$|\cdot|_{N\Delta^{\op}}$}(1.5,0.5ex);\ \VWald.
\end{equation*}
\end{ntn}

\begin{exm} One deduces from Ex. \ref{exm:catwithcofibsandwesislabWaldinftycat} that a category $(C,\cof C,wC)$ with cofibrations and weak equivalences gives rise to a virtual Waldhausen $\infty$-category $\mathscr{B}(NC,N\cof C,NwC)$.
\end{exm}

\begin{ntn}\label{ntn:sigmastar} Note that the pair cartesian fibration $\pi\colon\fromto{N\mathrm{M}}{N\Delta^{\op}}$ of Nt. \ref{ntn:M} admits a section $\sigma$ that assigns to any object $\mathbf{m}\in\Delta$ the pair $(\mathbf{m},0)\in\mathrm{M}$. For any labeled Waldhausen $\infty$-category $(\mathscr{C},w\mathscr{C})$, this section induces a functor of pairs over $N\Delta^{\op}$
\begin{equation*}
\sigma^{\star}_{(\mathscr{C},w\mathscr{C})}\colon\fromto{\mathscr{B}(\mathscr{C},w\mathscr{C})}{(N\Delta^{\op})^{\flat}\times\mathscr{C}},
\end{equation*}
which carries any object $(\mathbf{m},X)$ of $\mathscr{B}(\mathscr{C},w\mathscr{C})$ to the pair $(\mathbf{m},X_0)$ and any morphism $(\phi,\psi)\colon\fromto{(\mathbf{n},Y)}{(\mathbf{m},X)}$ to the composite
\begin{equation*}
Y_0\ \tikz[baseline]\draw[>=stealth,->,font=\scriptsize](0,0.5ex)--(0.5,0.5ex);\ Y_{\phi(0)}\ \tikz[baseline]\draw[>=stealth,->,font=\scriptsize](0,0.5ex)--node[above]{$\psi_0$}(0.5,0.5ex);\ X_0.
\end{equation*}

The section $\sigma$ induces a map of simplicial sets
\begin{equation*}
\fromto{\mathrm{H}(\mathscr{D},\mathscr{B}(\mathscr{C},w\mathscr{C}))}{w\Fun_{\Wald}(\mathscr{D},\mathscr{C})},
\end{equation*}
natural in $\mathscr{D}$, where $w\Fun_{\Wald}(\mathscr{D},\mathscr{C})\subset\Fun_{\Wald}(\mathscr{D},\mathscr{C})$ denotes the subcategory containing all exact functors $\fromto{\mathscr{D}}{\mathscr{C}}$ and those natural transformations that are pointwise labeled.
\end{ntn}

\begin{lem}\label{lem:iotaBCwCiswC} For any labeled Waldhausen $\infty$-category $(\mathscr{C},w\mathscr{C})$ and any compact Waldhausen $\infty$-category $\mathscr{D}$, the map $\fromto{\mathrm{H}(\mathscr{D},\mathscr{B}(\mathscr{C},w\mathscr{C}))}{w\Fun_{\Wald}(\mathscr{D},\mathscr{C})}$ induced by $\sigma$ is a weak homotopy equivalence.
\begin{proof}[Proof A] Using (the dual of) Joyal's $\infty$-categorical version of Quillen's Theorem A \cite[Th. 4.1.3.1]{HTT}, we are reduced to showing that for any exact functor $X\colon\fromto{\mathscr{D}}{\mathscr{C}}$, the simplicial set
\begin{equation*}
\mathrm{H}(\mathscr{D},\mathscr{B}(\mathscr{C},w\mathscr{C}))\times_{w\Fun_{\Wald}(\mathscr{D},\mathscr{C})}w\Fun_{\Wald}(\mathscr{D},\mathscr{C})_{X/}
\end{equation*}
is weakly contractible. This simplicial set is the geometric realization of the simplicial space
\begin{equation*}
\goesto{\mathbf{n}}{\mathrm{H}_{1+n}(\mathscr{D},\mathscr{B}(\mathscr{C},w\mathscr{C}))\times_{w\Fun_{\Wald}(\mathscr{D},\mathscr{C})}\{X\}};
\end{equation*}
in particular, it may be identified with the path space of the fiber of the map
\begin{equation*}
\fromto{\mathrm{H}(\mathscr{D},\mathscr{B}(\mathscr{C},w\mathscr{C}))}{w\Fun_{\Wald}(\mathscr{D},\mathscr{C})}
\end{equation*}
over the vertex $X$.
\end{proof}
\begin{proof}[Proof B] Consider the ordinary category $\Delta_{w\Fun_{\Wald}^{\flat}(\mathscr{D},\mathscr{C})}$ of simplices of the simplicial set $w\Fun_{\Wald}(\mathscr{D},\mathscr{C})$. Corresponding to the natural map
\begin{equation*}
\fromto{N(\Delta_{w\Fun_{\Wald}(\mathscr{D},\mathscr{C})}^{\op}\times_{\Delta^{\op}}\mathrm{M}_{\dag})}{\Fun_{\Wald}(\mathscr{D},\mathscr{C})}
\end{equation*}
is a map
\begin{equation*}
\fromto{N\Delta_{w\Fun_{\Wald}(\mathscr{D},\mathscr{C})}^{\op}}{\mathrm{H}(\mathscr{D},\mathscr{B}(\mathscr{C},w\mathscr{C}))}.
\end{equation*}
This map identifies the nerve $N\Delta_{w\Fun_{\Wald}(\mathscr{D},\mathscr{C})}^{\op}$ with the simplicial subset of $\mathrm{H}(\mathscr{D},\mathscr{B}(\mathscr{C},w\mathscr{C}))$ whose simplices correspond to maps
\begin{equation*}
\fromto{\Delta^n\times_{\Delta^{\op}}\mathrm{M}_{\dag}}{\Fun_{\Wald}(\mathscr{D},\mathscr{C})}
\end{equation*}
that carry cocartesian edges (over $\Delta^n$) to degenerate edges. The composite
\begin{equation*}
N\Delta_{w\Fun_{\Wald}(\mathscr{D},\mathscr{C})}^{\op}\ \tikz[baseline]\draw[>=stealth,->](0,0.5ex)--(0.5,0.5ex);\ \mathrm{H}(\mathscr{D},\mathscr{B}(\mathscr{C},w\mathscr{C}))\ \tikz[baseline]\draw[>=stealth,->](0,0.5ex)--(0.5,0.5ex);\ w\Fun_{\Wald}(\mathscr{D},\mathscr{C})
\end{equation*}
is the ``initial vertex map,'' which is a well-known weak equivalence. A simple argument now shows that the map $\fromto{N\Delta_{w\Fun_{\Wald}(\mathscr{D},\mathscr{C})}^{\op}}{\mathrm{H}(\mathscr{D},\mathscr{B}(\mathscr{C},w\mathscr{C}))}$ is also a weak equivalence.
\end{proof}
\end{lem}

In other words, the virtual Waldhausen $\infty$-category $\mathscr{B}(\mathscr{C},w\mathscr{C})$ attached to a labeled Waldhausen $\infty$-category $(\mathscr{C},w\mathscr{C})$ is not itself representable, but it's close:
\begin{prp}\label{prp:BCwCiswFunblankC} The virtual Waldhausen $\infty$-category $\mathscr{B}(\mathscr{C},w\mathscr{C})$ attached to a labeled Waldhausen $\infty$-category $(\mathscr{C},w\mathscr{C})$ is equivalent to the functor
\begin{equation*}
\goesto{\mathscr{D}}{w\Fun_{\Wald}(\mathscr{D},\mathscr{C})}.
\end{equation*}
\end{prp}


\subsection*{Inverting labeled edges} Unfortunately, for a labeled Waldhausen $\infty$-category $(\mathscr{C},w\mathscr{C})$, the functor (Nt. \ref{ntn:sigmastar})
\begin{equation*}
\sigma^{\star}_{(\mathscr{C},w\mathscr{C})}\colon\fromto{\mathscr{B}(\mathscr{C},w\mathscr{C})}{(N\Delta^{\op})^{\flat}\times\mathscr{C}}
\end{equation*}
will typically fail to be a morphism of $\mathbf{Wald}_{\infty/N\Delta^{\op}}^{\cocart}$, because the cocartesian edges of $\mathscr{B}(\mathscr{C},w\mathscr{C})$ will be carried to labeled edges, but not necessarily to equivalences. Hence one may not regard $\sigma_{(\mathscr{C},w\mathscr{C})}^{\star}$ as a natural transformation of functors $\fromto{N\Delta^{\op}}{\Wald}$. To rectify this, we may formally invert the edges in $w\mathscr{C}$ in the $\infty$-categorical sense.

\begin{lem} The inclusion functor $\into{\Wald}{\lWald}$ defined by the assignment $\goesto{(\mathscr{C},\mathscr{C}_{\dag})}{(\mathscr{C},\mathscr{C}_{\dag},\iota\mathscr{C})}$ admits a left adjoint $\fromto{\lWald}{\Wald}$.
\begin{proof} The inclusion functor $\into{\Wald}{\lWald}$ preserves all limits and all filtered colimits. Now the result follows from the Adjoint Functor Theorem \cite[Cor. 5.5.2.9]{HTT} along with Pr. \ref{prp:lWaldpres}.
\end{proof}
\end{lem}

Let us denote by $w\mathscr{C}^{-1}\mathscr{C}$ the image of a labeled Waldhausen $\infty$-category $(\mathscr{C},w\mathscr{C})$ under the left adjoint above. The canonical exact functor $\fromto{\mathscr{C}}{w\mathscr{C}^{-1}\mathscr{C}}$ is initial with the property that it carries labeled edges to equivalences. As an example, let us consider the case of an ordinary category with cofibrations and weak equivalences in the sense of Waldhausen \cite[\S 1.2]{MR86m:18011}.

\begin{prp}\label{thm:classWaldareWald} If $(C,\cof C,wC)$ is a category with cofibrations and weak equivalences that is a \emph{partial model category} \cite{BarKan1102} in the sense that: \textup{(1)} the weak equivalences satisfy the two-out-of-six axiom \cite[9.1]{MR2102294}, and \textup{(2)} the weak equivalences and trivial cofibrations are part of a three-arrow calculus of fractions \cite[11.1]{MR2102294}, then the Waldhausen $\infty$-category $(NwC)^{-1}(NC)$ is equivalent to the relative nerve $N(C,wC)$, equipped with the smallest pair structure containing the images of $\cof C$ \textup{(Ex. \ref{exm:catwithcofibsisWaldinftycat})}.
\begin{proof} We first claim that $N(C,wC)$ is a Waldhausen $\infty$-category.

First, by \cite[38.3(iii)]{MR2102294}, the image of the zero object $0\in C$ is again a zero object of $N(C,wC)$. It is also an initial object of $N(C,wC)_{\dag}$, since for any object $X$, the mapping space $\Map_{N(C,wC)_{\dag}}(0,X)$ is a union of connected components of $\Map_{N(C,wC)}(0,X)$, whence it is either empty or contractible, but the image of the edge $\fromto{0}{X}$ is ingressive by definition.

Now let us see that pushouts along ingressives exist and are ingressives. The $\infty$-category $\Fun_{\Pair_{\infty}}(\Lambda_0\mathscr{Q}^2,N(C,wC))$ is the relative nerve of the full subcategory $C^{\ulcorner}$ of $\Fun(\mathbf{1}\cup^{\{0\}}\mathbf{1},C)$ spanned by those functors that carry the first arrow $\fromto{0}{1}$ to a cofibration, equipped with the objectwise weak equivalences. Similarly, $\Fun_{\Pair_{\infty}}(\mathscr{Q}^2,N(C,wC))$ is the relative nerve of the full subcategory $C^{\Box}$ of $\Fun(\mathbf{1}\times\mathbf{1},C)$ spanned by those functors that carry the arrows $\fromto{(0,0)}{(0,1)}$ and $\fromto{(1,0)}{(1,1)}$ each to cofibrations, equipped with the objectwise weak equivalences. The forgetful functor $U\colon\fromto{C^{\Box}}{C^{\ulcorner}}$ and its left adjoint $F\colon\fromto{C^{\ulcorner}}{C^{\Box}}$ are each relative functors, whence they descend to an adjunction
\begin{equation*}
\adjunct{F}{\Ho(C^{\ulcorner})}{\Ho(C^{\Box})}{U}
\end{equation*}
on the $\Ho s\Set$-enriched homotopy categories, using the description \cite[36.3]{MR2102294}. Furthermore, the unit is clearly an equivalence $\id\simeq UF$. Hence the forgetful functor
\begin{equation*}
\fromto{\Fun_{\Pair_{\infty}}(\mathscr{Q}^2,N(C,wC))}{\Fun_{\Pair_{\infty}}(\Lambda_0\mathscr{Q}^2,N(C,wC))}
\end{equation*}
admits a left adjoint, and the unit for this adjunction is an equivalence. This is precisely the condition that pushouts along ingressives exist and are ingressives. Thus $N(C,wC)$ is a Waldhausen $\infty$-category.

Moreover, if $\cofto{X}{Y}$ is a cofibration of $C$ and if $\fromto{X}{X'}$ is an arrow of $C$, a square
\begin{equation*}
\begin{tikzpicture} 
\matrix(m)[matrix of math nodes, 
row sep=4ex, column sep=4ex, 
text height=1.5ex, text depth=0.25ex] 
{X&Y\\ 
X'&Y'\\}; 
\path[>=stealth,->,font=\scriptsize] 
(m-1-1) edge[>->] (m-1-2) 
edge (m-2-1) 
(m-1-2) edge (m-2-2) 
(m-2-1) edge[>->] (m-2-2); 
\end{tikzpicture}
\end{equation*}
in $N(C,wC)$ is a pushout just in case it is is the essential image of the left adjoint above. This in turn holds just in case it is equivalent to the image of a pushout square in $C$.

Now suppose $\mathscr{D}$ a Waldhausen $\infty$-category. Since the canonical functor
\begin{equation*}
\fromto{NC}{N(C,wC)}
\end{equation*}
is exact, there is an induced functor
\begin{equation*}
R\colon\fromto{\Fun_{\Wald}(N(C,wC),\mathscr{D})}{\Fun'_{\Wald}(NC,\mathscr{D})},
\end{equation*}
where $\Fun'_{\Wald}(NC,\mathscr{D})\subset\Fun_{\Wald}(NC,\mathscr{D})$ is the full subcategory spanned by those exact functors that carry arrows in $wC$ to equivalences in $D$. The universal property of $N(C,wC)$, combined with the definition of its pair structure, guarantees an equivalence
\begin{equation*}
\equivto{\Fun_{\Pair_{\infty}}(N(C,wC),\mathscr{D})}{\Fun'_{\Pair_{\infty}}(NC,\mathscr{D})},
\end{equation*}
where $\Fun'_{\Pair_{\infty}}(NC,\mathscr{D})\subset\Fun_{\Pair_{\infty}}(NC,\mathscr{D})$ is the full subcategory spanned by those functors of pairs that carry arrows in $wC$ to equivalences in $\mathscr{D}$. Hence $R$ is fully faithful. Since an object (respectively, a morphism, a square) in $N(C,wC)$ is a zero object (resp., an ingressive morphism, a pushout square along an ingressive morphism) just in case it is equivalent to the image of one under the functor $\fromto{NC}{N(C,wC)}$, it follows that a functor of pairs $\fromto{N(C,wC)}{\mathscr{D}}$ that induces an exact functor $\fromto{C}{\mathscr{D}}$ is itself exact. Thus $R$ is essentially surjective.
\end{proof}
\end{prp}

Let us give another example of a situation in which we can identify the Waldhausen $\infty$-category $w\mathscr{C}^{-1}\mathscr{C}$, up to splitting certain idempotents. We thank an anonymous referee and Peter Scholze for identifying an error in the original formulation of this result.

\begin{dfn}\label{dfn:weakcofinal} We say that a full Waldhausen subcategory $\mathscr{C}'\subset\mathscr{C}$ of a Waldhausen $\infty$-category is \emph{strongly cofinal} if, for any object $X\in\mathscr{C}$, there exists an object $Y\in\mathscr{C}$ such that $X\vee Y\in\mathscr{C}'$.
\end{dfn}

We will show below in Th. \ref{thm:cofinality} that a strongly cofinal subcategory $\mathscr{C}'\subset\mathscr{C}$ of a Waldhausen $\infty$-category has the same algebraic $K$-theory as $\mathscr{C}$ in positive degrees.

\begin{prp}\label{prp:localizationofperfwaldcat} Suppose $C$ a compactly generated $\infty$-category containing a zero object, suppose $L\colon\fromto{C}{D}$ an accessible localization of $C$, and suppose the inclusion $\into{D}{C}$ preserves filtered colimits. Assume also that the class of all $L$-equivalences of $C$ is generated (as a strongly saturated class) by the $L$-equivalences between compact objects. Then if $wC^{\omega}\subset C^{\omega}$ is the subcategory consisting of $L$-equivalences between compact objects, then $D^{\omega}$ is the idempotent completion of $(wC^{\omega})^{-1}C^{\omega}$.

In particular, $C$ and $D$ are additive (Df. \ref{item:directsums}), then with their maximal pair structures, the inclusion $\into{(wC^{\omega})^{-1}C^{\omega}}{D^{\omega}}$ is strongly cofinal.
\begin{proof} Let us begin by giving, for any labeled Waldhausen $\infty$-category $A$ with a \emph{maximal} pair structure, a construction of $wA^{-1}A$. We begin by inverting the edges of $wA$ in $A$ as an $\infty$-category; the result is an $\infty$-category $A'$ and a functor $i\colon\fromto{A}{A'}$ that induces, for any $\infty$-category $B$, a fully faithful functor
\[\fromto{\Fun(A',B)}{\Fun(A,B)}\]
whose essential image is spanned by those functors that carry the edges in $wA$ to equivalences in $B$. Now we will use the ideas of \cite[\S 5.3.6]{HTT}. Consider the class $\mathscr{R}$ consisting of the following diagrams: the composite
\[\varnothing^{\rhd}\ \tikz[baseline]\draw[>=stealth,->,font=\scriptsize,inner sep=0.75pt](0,0.5ex)--node[above]{$z$}(0.5,0.5ex);\ A\ \tikz[baseline]\draw[>=stealth,->,font=\scriptsize,inner sep=0.75pt](0,0.5ex)--node[above]{$i$}(0.5,0.5ex);\ A',\]
in which $z$ is the inclusion of the zero object, and the composites
\[(\Lambda^2_0)^{\rhd}\ \tikz[baseline]\draw[>=stealth,->,font=\scriptsize,inner sep=0.75pt](0,0.5ex)--node[above]{$p$}(0.5,0.5ex);\ A\ \tikz[baseline]\draw[>=stealth,->,font=\scriptsize,inner sep=0.75pt](0,0.5ex)--node[above]{$i$}(0.5,0.5ex);\ A'\]
in which $p$ is a pushout square. Now let $\mathscr{F}$ denote the collection of all finite simplicial sets. In the notation of \cite[Pr. 5.3.6.2]{HTT}, we claim that $wA^{-1}A\simeq\mathscr{P}_{\mathscr{R}}^{\mathscr{F}}(A')$, where the latter $\infty$-category in endowed with its maximal pair structure.

To prove this claim, let us first note that the inclusion of the full subcategory $\Cat_{\infty}^{\Rex,z}\subset\Wald$ spanned by those Waldhausen $\infty$-categories equipped with the maximal pair structure admits a left adjoint. This much follows from the adjoint functor theorem, but in fact we can be more precise: it is the construction $\goesto{\mathscr{C}}{\mathscr{P}_{\mathscr{W}}^{\mathscr{F}}\mathscr{C}}$, where $\mathscr{W}$ consists of the initial object $\fromto{\varnothing^{\rhd}}{\mathscr{C}}$ and the pushouts $\fromto{(\Lambda^2_0)^{\rhd}}{\mathscr{C}}$ of cofibrations, and $\mathscr{F}$ consists of all finite simplicial sets. Note that since the diagrams of $\mathscr{W}$ are colimits in $\mathscr{C}$, it follows that the unit $j\colon\fromto{\mathscr{C}}{\mathscr{P}_{\mathscr{W}}^{\mathscr{F}}\mathscr{C}}$ is fully faithful.

Now for any Waldhausen $\infty$-category $\mathscr{C}$, let us consider the square
\begin{equation*}
\begin{tikzpicture}[baseline]
\matrix(m)[matrix of math nodes,
row sep=4ex, column sep=4ex,
text height=1.5ex, text depth=0.25ex]
{\Fun_{\Wald}(\mathscr{P}_{\mathscr{R}}^{\mathscr{F}}(A'),\mathscr{C}) & \Fun'_{\Wald}(A,\mathscr{C}) \\
\Fun_{\Wald}(\mathscr{P}_{\mathscr{R}}^{\mathscr{F}}(A'),\mathscr{P}_{\mathscr{W}}^{\mathscr{F}}(\mathscr{C})) & \Fun'_{\Wald}(A,\mathscr{P}_{\mathscr{W}}^{\mathscr{F}}(\mathscr{C})), \\};
\path[>=stealth,->,font=\scriptsize]
(m-1-1) edge node[above]{} (m-1-2)
edge node[left]{} (m-2-1)
(m-1-2) edge node[right]{} (m-2-2)
(m-2-1) edge node[below]{} (m-2-2);
\end{tikzpicture}
\end{equation*}
where $\Fun'$ denotes the full subcategory spanned by those exact functors that carry the edges of $wA$ to equivalences. Unwinding the universal properties, one sees immediately that the bottom horizontal functor is an equivalence; our claim is that the top horizontal functor is an equivalence. Hence we aim to show that the square above is homotopy cartesian; this amounts to the claim that in commutative diagram of exact functors
\begin{equation*}
\begin{tikzpicture}[baseline]
\matrix(m)[matrix of math nodes,
row sep=4ex, column sep=4ex,
text height=1.5ex, text depth=0.25ex]
{A & \mathscr{C} \\
\mathscr{P}_{\mathscr{R}}^{\mathscr{F}}(A') & \mathscr{P}_{\mathscr{W}}^{\mathscr{F}}(\mathscr{C}), \\ };
\path[>=stealth,->,font=\scriptsize]
(m-1-1) edge node[above]{} (m-1-2)
edge node[left]{$i$} (m-2-1)
(m-1-2) edge node[right]{$j$} (m-2-2)
(m-2-1) edge node[below]{$F$} (m-2-2);
\end{tikzpicture}
\end{equation*}
the functor $F$ factors through $j$. This now follows from the minimality of the construction of $\mathscr{P}_{\mathscr{R}}^{\mathscr{F}}(A')$, as in the proof of \cite[Pr. 5.3.6.2]{HTT}. This completes the proof that $wA^{-1}A\simeq\mathscr{P}_{\mathscr{R}}^{\mathscr{F}}(A')$.

Let us now note the inclusion of the full subcategory $\Cat_{\infty}^{\Rex,z,\vee}\subset\Cat_{\infty}^{\Rex,z}$ spanned by those Waldhausen $\infty$-categories equipped with the maximal pair structure also admits a left adjoint. This is given by the idempotent completion $\goesto{A}{A^{\vee}}$ of \cite[\S 5.1.4]{HTT}.

We turn to our localization. For any idempotent complete $\infty$-category $A$ that admits all finite colimits, the localization $\fromto{C^{\omega}}{D^{\omega}}$ induces an equivalence
\[\equivto{\Fun_{\Rex}(D^{\omega},A)}{\Fun_{\Rex}'(C^{\omega},A)},\]
where $\Fun_{\Rex}'(C^{\omega},A)\subset\Fun_{\Rex}(C^{\omega},A)$ is the full subcategory spanned by those finite colimit-preserving functors that carry $L$-equivalences to equivalences. (Here we are using the mutually inverse equivalences $\goesto{A}{\Ind(A)}$ and $\goesto{B}{B^{\omega}}$ of \cite[Pr. 5.5.7.10]{HTT}.) This target $\infty$-category is of course equivalent to the full subcategory of $\Fun_{\Wald}(C^{\omega},A)$ spanned by those exact functors that carries $L$-equivalences to equivalences. We therefore deduce that the natural functor $\fromto{(wC^{\omega})^{-1}C^{\omega}}{D^{\omega}}$ induces an equivalence
\[\equivto{\Fun_{\Wald}(D^{\omega},A)}{\Fun_{\Wald}((wC^{\omega})^{-1}C^{\omega},A)^{\vee}}.\]
Consequently, we deduce that
\[D^{\omega}\simeq\mathscr{P}_{\mathscr{W}}^{\mathscr{F}}((C^{\omega})')^{\vee}\simeq((wC^{\omega})^{-1}C^{\omega})^{\vee},\]
as desired.
\end{proof}
\end{prp}

\begin{ntn}\label{ntn:gammaC} Composing the canonical exact functor $\fromto{\mathscr{C}}{w\mathscr{C}^{-1}\mathscr{C}}$ with the functor
\begin{equation*}
\fromto{\mathscr{B}(\mathscr{C},w\mathscr{C})}{(N\Delta^{\op})^{\flat}\times\mathscr{C}},
\end{equation*}
we obtain a morphism of $\mathbf{Wald}_{\infty/N\Delta^{\op}}^{\cocart}$
\begin{equation*}
\fromto{\mathscr{B}(\mathscr{C},w\mathscr{C})}{(N\Delta^{\op})^{\flat}\times w\mathscr{C}^{-1}\mathscr{C}}
\end{equation*}
that carries cocartesian edges of $\mathscr{B}(\mathscr{C},w\mathscr{C})$ to equivalences. Applying the realization $|\cdot|_{N\Delta^{\op}}$ (Df. \ref{dfn:realization}), we obtain a morphism of $\VWald$
\begin{equation*}
\gamma_{(\mathscr{C},w\mathscr{C})}\colon\fromto{\mathscr{B}(\mathscr{C},w\mathscr{C})}{w\mathscr{C}^{-1}\mathscr{C}}.
\end{equation*}
\end{ntn}

We emphasize that for a general labeled Waldhausen $\infty$-category $(\mathscr{C},w\mathscr{C})$, the comparison morphism $\gamma_{(\mathscr{C},w\mathscr{C})}$ is not an equivalence of $\VWald$; nevertheless, we will find (Pr. \ref{prp:gammaonKisequiv}) that $\gamma_{(\mathscr{C},w\mathscr{C})}$ often induces an equivalence on $K$-theory. 


\subsection*{Waldhausen's Fibration Theorem, redux} We now aim to prove an analogue of Waldhausen's Generic Fibration Theorem \cite[Th. 1.6.4]{MR86m:18011}. For this we require a suitable analogue of Waldhausen's cylinder functor in the $\infty$-categorical context. This should reflect the idea that a labeled edge can, to some extent, be replaced by a labeled ingressive.

\begin{ntn}\label{ntn:wdagFunDC} To this end, for any labeled Waldhausen $\infty$-category $(\mathscr{A},\mathscr{A}_{\dag})$, write $w_{\dag}\mathscr{A}\coloneq w\mathscr{A}\cap\mathscr{A}_{\dag}$. The subcategory $w_{\dag}\mathscr{A}\subset\mathscr{A}$ defines a new pair structure, \emph{but not a new labeling}, of $\mathscr{A}$. Nevertheless, we may consider the full subcategory $\mathscr{B}(\mathscr{A},w_{\dag}\mathscr{A})\subset\mathscr{F}(\mathscr{A})$ spanned by those filtered objects
\begin{equation*}
X_0\ \tikz[baseline]\draw[>=stealth,>->,font=\scriptsize](0,0.5ex)--(0.75,0.5ex);\ X_1\ \tikz[baseline]\draw[>=stealth,>->,font=\scriptsize](0,0.5ex)--(0.75,0.5ex);\ \cdots\ \tikz[baseline]\draw[>=stealth,>->,font=\scriptsize](0,0.5ex)--(0.75,0.5ex);\ X_m
\end{equation*}
such that each ingressive $\cofto{X_{i}}{X_{i+1}}$ is labeled; we shall regard it as a subpair. One may verify that $\mathscr{B}_{m}(\mathscr{A},w_{\dag}\mathscr{A})\subset\mathscr{F}_m(\mathscr{A})$ is a Waldhausen subcategory, and $\fromto{\mathscr{B}(\mathscr{A},w_{\dag}\mathscr{A})}{N\Delta^{\op}}$ is a Waldhausen cocartesian fibration.

For any pair $\mathscr{D}$, write $w_{\dag}\Fun_{\Pair_{\infty}}(\mathscr{D},\mathscr{A})\subset\Fun_{\Pair_{\infty}}(\mathscr{D},\mathscr{A})$ for the following pair structure. A natural transformation
\begin{equation*}
\eta\colon\fromto{\mathscr{D}\times\Delta^1}{\mathscr{A}}
\end{equation*}
lies in $w_{\dag}\Fun_{\Pair_{\infty}}(\mathscr{D},\mathscr{A})$ if and only if the it satisfies the following two conditions.
\begin{enumerate}[(\ref{ntn:wdagFunDC}.1)]
\item For any object $X$ of $\mathscr{D}$, the edge $\fromto{\Delta^1\cong\Delta^1\times\{X\}\subset\Delta^1\times\mathscr{D}}{\mathscr{A}}$ is both ingressive and labeled.
\item For any ingressive $f\colon\cofto{X}{Y}$ of $\mathscr{D}$, the corresponding edge
\begin{equation*}
\fromto{\Delta^1}{\mathscr{F}_1(\mathscr{A})}
\end{equation*}
is ingressive in the sense of Df. \ref{dfn:injpairstruct}.
\end{enumerate}
If $\mathscr{D}$ is a Waldhausen $\infty$-category, write
\begin{equation*}
w_{\dag}\Fun_{\Wald}(\mathscr{D},\mathscr{A})\subset w_{\dag}\Fun_{\Pair_{\infty}}(\mathscr{D},\mathscr{A})
\end{equation*}
for the full subcategory spanned by the exact functors.
\end{ntn}

\begin{nul} Note that the proofs of Lm. \ref{lem:iotaBCwCiswC} apply also to the pair $(\mathscr{A},w_{\dag}\mathscr{A})$ to guarantee that for any compact Waldhausen $\infty$-category $\mathscr{D}$, the natural map
\begin{equation*}
\fromto{\mathrm{H}(\mathscr{D},(\mathscr{B}(\mathscr{A},w_{\dag}\mathscr{A})/N\Delta^{\op}))}{w_{\dag}\Fun_{\Wald}(\mathscr{D},\mathscr{A})}
\end{equation*}
induced by $\sigma$ is a weak homotopy equivalence.
\end{nul}

\begin{dfn}\label{dfn:enoughcofibrations} Suppose $(\mathscr{A},w\mathscr{A})$ a labeled Waldhausen $\infty$-category. We shall say that $(\mathscr{A},w\mathscr{A})$ \textbf{\emph{has enough cofibrations}} if for any small pair of $\infty$-categories $\mathscr{D}$, the inclusion
\begin{equation*}
\into{w_{\dag}\Fun_{\Pair_{\infty}}(\mathscr{D},\mathscr{A})}{w\Fun_{\Pair_{\infty}}(\mathscr{D},\mathscr{A})}
\end{equation*}
is a weak homotopy equivalence.
\end{dfn}
\noindent In particular, if every labeled edge of $(\mathscr{A},w\mathscr{A})$ is ingressive, then $(\mathscr{A},w\mathscr{A})$ has enough cofibrations. More generally, this may prove to be an extremely difficult condition to verify, but the following lemma simplifies matters somewhat.
\begin{lem}\label{lem:enoughcofs} Suppose $(\mathscr{A},\mathscr{A}_{\dag},w\mathscr{A})$ a labeled Waldhausen $\infty$-category. Suppose that there exists a functor
\[F\colon\fromto{\Fun(\Delta^1,\mathscr{A})}{\Fun(\Delta^1,\mathscr{A})}\]
along with a natural transformation $\eta\colon\fromto{\id}{F}$ such that:
\begin{enumerate}[(\ref{lem:enoughcofs}.1)]
\item The functor $F$ carries $\Fun(\Delta^1,w\mathscr{A})$ to $\Fun(\Delta^1,w_{\dag}\mathscr{A})$.
\item If $f$ is a labeled ingressive, then $\eta_f$ is an equivalence.
\item If $f$ is labeled, then $\eta_f$ is objectwise labeled.
\end{enumerate}
Then $(\mathscr{A},\mathscr{A}_{\dag},w\mathscr{A})$ has enough cofibrations.
\begin{proof} For any pair $\mathscr{D}$, the functor $F$ induces a functor
\begin{equation*}
\fromto{\Fun(\Delta^1,w\Fun_{\Pair_{\infty}}(\mathscr{D},\mathscr{A}))}{\Fun(\Delta^1,w_{\dag}\Fun_{\Pair_{\infty}}(\mathscr{D},\mathscr{A})),}
\end{equation*}
and $\eta$ induces natural transformations that exhibit this functor as a homotopy inverse to the inclusion
\begin{equation*}
\into{\Fun(\Delta^1,w_{\dag}\Fun_{\Pair_{\infty}}(\mathscr{D},\mathscr{A}))}{\Fun(\Delta^1,w\Fun_{\Pair_{\infty}}(\mathscr{D},\mathscr{A}))}.
\end{equation*}
The result now follow from the homotopy equivalence between a simplicial set and its (unbased) path space.
\end{proof}
\end{lem}

\begin{lem}\label{lem:wAhasenoughcofibs} If a labeled Waldhausen $\infty$-category $(\mathscr{A},w\mathscr{A})$ has enough cofibrations, then for any Waldhausen $\infty$-category $\mathscr{D}$, the inclusion
\begin{equation*}
\into{w_{\dag}\Fun_{\Wald}(\mathscr{D},\mathscr{A})}{w\Fun_{\Wald}(\mathscr{D},\mathscr{A})}
\end{equation*}
is a weak homotopy equivalence.
\begin{proof} For any Waldhausen $\infty$-category $\mathscr{B}$, the square
\begin{equation*}
\begin{tikzpicture} 
\matrix(m)[matrix of math nodes, 
row sep=4ex, column sep=4ex, 
text height=1.5ex, text depth=0.25ex] 
{w_{\dag}\Fun_{\Wald}(\mathscr{D},\mathscr{A})&w\Fun_{\Wald}(\mathscr{D},\mathscr{A})\\ 
w_{\dag}\Fun_{\Pair_{\infty}}(\mathscr{D},\mathscr{A})&w\Fun_{\Pair_{\infty}}(\mathscr{D},\mathscr{A})\\}; 
\path[>=stealth,->,font=\scriptsize] 
(m-1-1) edge (m-1-2) 
edge (m-2-1) 
(m-1-2) edge (m-2-2) 
(m-2-1) edge (m-2-2); 
\end{tikzpicture}
\end{equation*}
is a pullback, and the vertical maps are inclusions of connected components.
\end{proof}
\end{lem}

\begin{thm}[Generic Fibration Theorem II]\label{thm:fibration} Suppose $(\mathscr{A},w\mathscr{A})$ a labeled Waldhausen $\infty$-category that has enough cofibrations. Suppose $\phi\colon\fromto{\Wald}{\mathscr{E}}$ an additive theory with left derived functor $\Phi$. Write $\mathscr{A}^w\subset\mathscr{A}$ for the full subcategory spanned by those objects $X$ such that a map from a zero object to $X$ is labeled, with the pair structure inherited from $\mathscr{A}$. Then $\mathscr{A}^w$ is a Waldhausen $\infty$-category, the inclusion $i\colon\into{\mathscr{A}^w}{\mathscr{A}}$ is exact, and it along with the morphism of virtual Waldhausen $\infty$-categories $e\colon\fromto{\mathscr{A}}{\mathscr{B}(\mathscr{A},w\mathscr{A})}$ give rise to a fiber sequence
\begin{equation*}
\begin{tikzpicture} 
\matrix(m)[matrix of math nodes, 
row sep=4ex, column sep=4ex, 
text height=1.5ex, text depth=0.25ex] 
{\phi(\mathscr{A}^w)&\phi(\mathscr{A})\\ 
\ast&\Phi(\mathscr{B}(\mathscr{A},w\mathscr{A})).\\}; 
\path[>=stealth,->,font=\scriptsize] 
(m-1-1) edge (m-1-2) 
edge (m-2-1) 
(m-1-2) edge (m-2-2) 
(m-2-1) edge (m-2-2); 
\end{tikzpicture}
\end{equation*}
\begin{proof} It follows from Pr. \ref{prp:fibthmi} that it is enough to exhibit an equivalence between $\Phi(\mathscr{B}(\mathscr{A},w\mathscr{A}))$ and $\Phi(\mathscr{K}(i))$ as objects of $\mathscr{E}_{\phi(\mathscr{A})/}$.

The forgetful functor $\fromto{\mathscr{K}(i)}{\mathscr{F}\mathscr{A}}$ is fully faithful, and its essential image $\widetilde{\mathscr{F}}^w\mathscr{A}$ consists of those filtered objects
\begin{equation*}
X_0\ \tikz[baseline]\draw[>=stealth,>->,font=\scriptsize](0,0.5ex)--(0.75,0.5ex);\ X_1\ \tikz[baseline]\draw[>=stealth,>->,font=\scriptsize](0,0.5ex)--(0.75,0.5ex);\ \cdots\ \tikz[baseline]\draw[>=stealth,>->,font=\scriptsize](0,0.5ex)--(0.75,0.5ex);\ X_m
\end{equation*}
such that the induced ingressive $\cofto{X_{i}/X_0}{X_{i+1}/X_0}$ is labeled; this contains the subcategory $\mathscr{B}(\mathscr{A},w_{\dag}\mathscr{A})$. We claim that for any $m\geq 0$, the induced morphism $\fromto{\phi(\mathscr{B}_m(\mathscr{A},w_{\dag}\mathscr{A}))}{\phi(\widetilde{\mathscr{F}}_m^w(\mathscr{A}))}$ is an equivalence. Indeed, one may select an exact functor $p\colon\fromto{\mathscr{K}_m(i)}{\mathscr{B}_m(\mathscr{A},w_{\dag}\mathscr{A})}$ that carries an object
\begin{equation*}
\begin{tikzpicture} 
\matrix(m)[matrix of math nodes, 
row sep=4ex, column sep=4ex, 
text height=1.5ex, text depth=0.25ex] 
{X_0&X_1&X_2&\dots&X_m\\ 
0&U_1&U_2&\dots&U_m\\}; 
\path[>=stealth,->,font=\scriptsize] 
(m-1-1) edge[>->] (m-1-2) 
edge (m-2-1)
(m-1-2) edge[>->] (m-1-3)
edge (m-2-2)
(m-1-3) edge[>->] (m-1-4)
edge (m-2-3)
(m-1-4) edge[>->] (m-1-5)
(m-1-5) edge (m-2-5)
(m-2-1) edge[>->] (m-2-2) 
(m-2-2) edge[>->] (m-2-3)
(m-2-3) edge[>->] (m-2-4)
(m-2-4) edge[>->] (m-2-5); 
\end{tikzpicture}
\end{equation*}
to the filtered object
\begin{equation*}
X_0\ \tikz[baseline]\draw[>=stealth,>->,font=\scriptsize](0,0.5ex)--(0.75,0.5ex);\ X_0\vee U_1\ \tikz[baseline]\draw[>=stealth,>->,font=\scriptsize](0,0.5ex)--(0.75,0.5ex);\ X_0\vee U_2\ \tikz[baseline]\draw[>=stealth,>->,font=\scriptsize](0,0.5ex)--(0.75,0.5ex);\ \cdots\ \tikz[baseline]\draw[>=stealth,>->,font=\scriptsize](0,0.5ex)--(0.75,0.5ex);\ X_0\vee U_m.
\end{equation*}
When $m=0$, this functor is compatible with the canonical equivalences from $\mathscr{A}$. Additivity now guarantees that $p$ defines a (homotopy) inverse to the morphism $\fromto{\phi(\mathscr{B}_m(\mathscr{A},w_{\dag}\mathscr{A}))}{\phi(\widetilde{\mathscr{F}}_m^w\mathscr{A})}$.

One has an obvious forgetful functor $\fromto{\mathscr{B}(\mathscr{A},w_{\dag}\mathscr{A})}{\mathscr{B}(\mathscr{A},w\mathscr{A})}$ over $N\Delta^{\op}$. We claim that this induces an equivalence of virtual Waldhausen $\infty$-categories $\fromto{|\mathscr{B}(\mathscr{A},w_{\dag}\mathscr{A})|_{N\Delta^{\op}}}{|\mathscr{B}(\mathscr{A},w\mathscr{A})|_{N\Delta^{\op}}}$. So we wish to show that for any compact Waldhausen $\infty$-category $\mathscr{D}$, the morphism
\begin{equation*}
\fromto{\mathrm{H}(\mathscr{D},(\mathscr{B}(\mathscr{A},w_{\dag}\mathscr{A})/N\Delta^{\op}))}{\mathrm{H}(\mathscr{D},(\mathscr{B}(\mathscr{A},w\mathscr{A})/N\Delta^{\op}))}
\end{equation*}
of simplicial sets is a weak homotopy equivalence.

By Lm. \ref{lem:iotaBCwCiswC} and its extension to the pair $(\mathscr{A},w_{\dag}\mathscr{A})$, we have a square
\begin{equation*}
\begin{tikzpicture} 
\matrix(m)[matrix of math nodes, 
row sep=4ex, column sep=4ex, 
text height=1.5ex, text depth=0.25ex] 
{\mathrm{H}(\mathscr{D},(\mathscr{B}(\mathscr{A},w_{\dag}\mathscr{A})/N\Delta^{\op}))&\mathrm{H}(\mathscr{D},(\mathscr{B}(\mathscr{A},w\mathscr{A})/N\Delta^{\op}))\\ 
w_{\dag}\Fun_{\Wald}(\mathscr{D},\mathscr{C})&w\Fun_{\Wald}(\mathscr{D},\mathscr{C})\\}; 
\path[>=stealth,->,font=\scriptsize] 
(m-1-1) edge (m-1-2) 
edge (m-2-1) 
(m-1-2) edge (m-2-2) 
(m-2-1) edge (m-2-2); 
\end{tikzpicture}
\end{equation*}
in which the vertical maps are weak homotopy equivalences. Since $(\mathscr{A},w\mathscr{A})$ has enough cofibrations, the horizontal map along the bottom is a weak homotopy equivalence as well by Lm. \ref{lem:wAhasenoughcofibs}.
\end{proof}
\end{thm}


\part{Algebraic $K$-theory}

We are finally prepared to describe the Waldhausen $K$-theory of $\infty$-categories. We define (Df. \ref{dfn:Ktheory}) $K$-theory as the additivization of the the theory $\iota$ that assigns to any Waldhausen $\infty$-category the maximal $\infty$-groupoid (Nt. \ref{ntn:interior}) contained therein. Since the theory $\iota$ is representable by the particularly simple Waldhausen $\infty$-category $N\Gamma^{\op}$ of pointed finite sets (Pr. \ref{prp:iotaisrepresentable}), we obtain, for any additive theory $\phi$, a description of the space of natural transformations $\fromto{K}{\phi}$ as the value of $\phi$ on $N\Gamma^{\op}$.

Following this, we briefly describe two key examples that exploit certain features of the algebraic $K$-theory functor of which we are fond. The first of these (\S \ref{sect:example1}) lays the foundations for the algebraic $K$-theory of $E_1$-algebras in a variety of monoidal $\infty$-categories, and we prove a straightforward localization theorem. Second (\S \ref{sect:example2}), we extend algebraic $K$-theory to the context of spectral Deligne--Mumford stacks in the sense of Lurie, and we prove Thomason's ``proto-localization'' theorem in this context.


\section{The universal property of Waldhausen $K$-theory}\label{sect:univpropKthy} In this section, we \emph{define} algberaic $K$-theory as the additivization of the functor that assigns to any Waldhausen $\infty$-category its moduli space of objects. More precisely, the functor $\iota\colon\fromto{\Wald}{\Kan}$ that assigns to any Waldhausen $\infty$-category its interior $\infty$-groupoid (Nt. \ref{ntn:interior}) is a theory.
\begin{dfn}\label{dfn:Ktheory} The \textbf{\emph{algebraic $K$-theory functor}}
\begin{equation*}
K\colon\fromto{\Wald}{\Kan}
\end{equation*}
is defined as the additivization $K\coloneq D\iota$ of the interior functor $\iota\colon\fromto{\Wald}{\Kan}$. We denote by $\KK\colon\fromto{\Wald}{\Sp_{\geq 0}}$ its canonical connective delooping, as guaranteed by Cor. \ref{cor:additivevaluedinstab} and Pr. \ref{prp:deloopingofdfisconnective}.
\end{dfn}

Unpacking this definition, we obtain a global universal property of the natural morphism $\fromto{\iota}{K}$.
\begin{prp}\label{prp:univpropofk} For any additive theory $\phi$, the morphism $\fromto{\iota}{K}$ induces a natural homotopy equivalence
\begin{equation*}
\equivto{\Map(K,\phi)}{\Map(\iota,\phi)}.
\end{equation*}
\end{prp}

We will prove in Cor. \ref{cor:easycompare} and Cor. \ref{cor:hardcompare} that our definition extends Waldhausen's.

\begin{exm}\label{exm:Atheoryofinfintopoi} For any $\infty$-topos $\mathscr{E}$, one may define the $A$-theory space
\begin{equation*}
A(\mathscr{E})\coloneq K(\mathscr{E}_{\ast}^{\omega})
\end{equation*}
(Ex. \ref{exm:inftytopoigiveWaldcats}). In light of Ex. \ref{exm:stabinftytopoi}, we have
\begin{equation*}
A(\mathscr{E})\simeq K(\Sp(\mathscr{E}^{\omega})).
\end{equation*}
For any Kan simplicial set $X$, if
\begin{equation*}
\mathscr{E}=\Fun(X,\Kan)\simeq\Kan_{/X},
\end{equation*}
then it will follow from Cor. \ref{cor:fromWaldtoWald} that $A(\mathscr{E})$ agrees with Waldhausen's $A(X)$, where one defines the latter via the category $\mathscr{R}_{\mathrm{df}}(X)$ of finitely dominated retractive spaces over $X$ \cite[p. 389]{MR86m:18011}. Then of course one has $A(\mathscr{E})\simeq K(\Fun(X,\Sp^{\omega}))$.
\end{exm}


\subsection*{Representability of algebraic $K$-theory} Algebraic $K$-theory is controlled, as an additive theory, by the theory $\iota$. It is therefore valuable to study this functor as a theory. As a first step, we find that it is corepresentable.
\begin{ntn} For any finite set $I$, write $I_+$ for the finite set $I\sqcup\{\infty\}$. Denote by $\Gamma^{\op}$ the ordinary category of pointed finite sets. Denote by $\Gamma^{\op}_{\dag}\subset\Gamma^{\op}$ the subcategory comprised of monomorphisms $\fromto{J_+}{I_+}$.
\end{ntn}

\begin{prp}\label{prp:iotaisrepresentable} For any Waldhausen $\infty$-category $\mathscr{C}$, the inclusion
\begin{equation*}
\into{\{\ast\}}{N\Gamma^{\op}}
\end{equation*}
induces an equivalence of $\infty$-categories
\begin{equation*}
\equivto{\Fun_{\Wald}(N\Gamma^{\op},\mathscr{C})}{\mathscr{C}}.
\end{equation*}
In particular, the functor $\iota\colon\fromto{\Wald}{\Kan}$ is corepresented by the object $N\Gamma^{\op}$.
\begin{proof} Write $N\Gamma^{\op}_{\leq 1}$ for the full subcategory of $N\Gamma^{\op}$ spanned by the objects $\varnothing$ and $\ast$. Then it follows from Joyal's theorem \cite[Pr. 1.2.12.9]{HTT} that the inclusion $\into{\{\ast\}}{N\Gamma^{\op}}$ induces an equivalence between $\mathscr{C}$ and the full subcategory $\Fun^{\ast}(N\Gamma^{\op}_{\leq1},\mathscr{C})$ of $\Fun(N\Gamma^{\op}_{\leq1},\mathscr{C})$ spanned by functors $z\colon\fromto{N\Gamma^{\op}_{\leq1}}{\mathscr{C}}$ such that $z(\varnothing)$ is a zero object. Now the result follows from the observation that the $\infty$-category $\Fun_{\Wald}(N\Gamma^{\op},\mathscr{C})$ can be identified as the full subcategory of the $\infty$-category $\Fun(N\Gamma^{\op},\mathscr{C})$ spanned by those functors $Z\colon\fromto{N\Gamma^{\op}}{\mathscr{C}}$ such that (1) $Z(\varnothing)$ is a zero object, and (2) the identity exhibits $Z$ as a left Kan extension of $Z|_{(N\Gamma^{\op}_{\leq1})}$ along the inclusion $\into{N\Gamma^{\op}_{\leq1}}{N\Gamma^{\op}}$.
\end{proof}
\end{prp}
\noindent In the language of Cor. \ref{cor:Wconstruction}, we find that $W(\Delta^0)\simeq N\Gamma^{\op}$. Note also that it follows that the left derived functor $I\colon\fromto{\VWald}{\Kan}$ of $\iota$ is given by evaluation at $W(\Delta^0)\simeq N\Gamma^{\op}$. From this, the Yoneda lemma combines with Pr. \ref{prp:univpropofk} to imply the following.
\begin{cor}\label{cor:tracesarephiofgamma} For any additive theory $\phi\colon\fromto{\Wald}{\Kan_{\ast}}$, there is a homotopy equivalence
\begin{equation*}
\Map(K,\phi)\simeq\phi(N\Gamma^{\op}),
\end{equation*}
natural in $\phi$.
\end{cor}
\noindent In particular, the theorem of Barratt--Priddy--Quillen \cite{MR0358767} implies the following.
\begin{cor} The space of endomorphisms of the $K$-theory functor
\[K\colon\fromto{\Wald}{\Kan}\]
is given by
\begin{equation*}
\End(K)\simeq QS^0.
\end{equation*}
\end{cor}


\subsection*{The local universal property of algebraic $K$-theory} Though conceptually pleasant, the universal property of $K$-theory as an object of $\Add(\Kan)$ does not obviously provide an easy recognition principle for the $K$-theory of any \emph{particular} Waldhausen $\infty$-category. For that, we note that $\iota$ is pre-additive, and we appeal to Cor. \ref{for:preaddsareeasy} to obtain the following result.
\begin{prp} For any virtual Waldhausen $\infty$-category $\mathscr{X}$, the $K$-theory space $K(\mathscr{X})$ is homotopy equivalent to the loop space $\Omega I(\mathscr{S}(\mathscr{X}))$, where $I$ is the left derived functor of $\iota$.
\end{prp}
\noindent We observe that for any sifted $\infty$-category and any Waldhausen cocartesian fibration $\fromto{\mathscr{Y}}{S}$, the space $I(\mathscr{S}(|\mathscr{Y}|_S))$ may be computed as the underlying space of the subcategory $\iota_{N\Delta^{\op}\times S}\mathscr{S}(\mathscr{Y}/S)$ of the $\infty$-category $\mathscr{S}(\mathscr{Y})$ comprised of the cocartesian edges with respect to the cocartesian fibration $\fromto{\mathscr{S}(\mathscr{Y}/S)}{N\Delta^{\op}\times S}$ (Df. \ref{rec:leftfib}). This provides us with a (singly delooped) model of the algebraic $K$-theory space $K(|\mathscr{Y}|_S)$ as the underlying simplicial set of an $\infty$-category.
\begin{cor} For any sifted $\infty$-category $S$ and any Waldhausen cocartesian fibration $\fromto{\mathscr{Y}}{S}$, the $K$-theory space $K(|\mathscr{Y}|_S)$ is homotopy equivalent to the loop space $\Omega\iota_{(N\Delta^{\op}\times S)}\mathscr{S}(\mathscr{Y}/S)$.
\end{cor}
\noindent The total space of a left fibration is weakly equivalent to the homotopy colimit of the functor that classifies it. So the $K$-theory space $K(\mathscr{C})$ of a Waldhausen $\infty$-category is given by
\begin{equation*}
K(\mathscr{C})\simeq\Omega(\colim\iota\SSS_{\ast}(\mathscr{C})),
\end{equation*}
where
\begin{equation*}
\SSS_{\ast}(\mathscr{C})\colon\fromto{N\Delta^{\op}}{\Wald}
\end{equation*}
classifies the Waldhausen cocartesian fibration $\fromto{\mathscr{S}(\mathscr{C})}{N\Delta^{\op}}$. Since this is precisely how Waldhausen's $K$-theory is defined \cite[\S 1.3]{MR86m:18011}, we obtain a comparison between our $\infty$-categorical $K$-theory and Waldhausen $K$-theory.
\begin{cor}\label{cor:easycompare} If $(C,\cof C)$ is an ordinary category with cofibrations in the sense of Waldhausen \cite[\S 1.1]{MR86m:18011}, then the algebraic $K$-theory of the Waldhausen $\infty$-category $(NC,N(\cof C))$ is naturally equivalent to Waldhausen's algebraic $K$-theory of $(C,\cof C)$.
\end{cor}
\noindent The fact that the algebraic $K$-theory space $K(\mathscr{X})$ of a virtual Waldhausen $\infty$-category $\mathscr{X}$ can be exhibited as the loop space of the underlying simplicial set of an $\infty$-category permits us to find the following sufficient condition that a morphism of Waldhausen cocartesian fibrations induce an equivalence on $K$-theory.
\begin{cor}\label{cor:ThmAKtheory} For any sifted $\infty$-category $S$, a morphism $\fromto{(\mathscr{Y}'/S)}{(\mathscr{Y}'/S)}$ of Waldhausen cocartesian fibrations induces an equivalence
\begin{equation*}
\equivto{K(|\mathscr{Y}'|_S)}{K(|\mathscr{Y}|_S)}
\end{equation*}
if the following two conditions are satisfied.
\begin{enumerate}[(\ref{cor:ThmAKtheory}.1)]
\item For any object $X\in\iota_S\mathscr{Y}$, the simplicial set
\begin{equation*}
\iota_S\mathscr{Y}'\times_{\iota_S\mathscr{Y}}(\iota_S\mathscr{Y})_{/X}
\end{equation*}
is weakly contractible.
\item For any object $Y\in\iota_S\mathscr{F}_1(\mathscr{Y}/S)$, the simplicial set
\begin{equation*}
\iota_S\mathscr{F}_1(\mathscr{Y}'/S)\times_{\iota_S\mathscr{F}_1(\mathscr{Y}/S)}\iota_S\mathscr{F}_1(\mathscr{Y}/S)_{/Y}
\end{equation*}
is weakly contractible.
\end{enumerate}
\begin{proof} We aim to show that the map $\fromto{\iota_{N\Delta^{\op}\times S}\mathscr{S}(\mathscr{Y}'/S)}{\iota_{N\Delta^{\op}\times S}\mathscr{S}(\mathscr{Y}/S)}$ is a weak homotopy equivalence; it is enough to show that for any $\mathbf{n}\in\Delta$, the map $\fromto{\iota_S\mathscr{F}_n(\mathscr{Y}'/S)}{\iota_S\mathscr{F}_n(\mathscr{Y}/S)}$ is a weak homotopy equivalence. Since $\mathscr{F}(\mathscr{Y}'/S)$ and $\mathscr{F}(\mathscr{Y}/S)$ are each category objects (Pr. \ref{prp:Fisacategory}), it is enough to prove this claim for $n\in\{0,1\}$. The result now follows from Joyal's $\infty$-categorical version of Quillen's Theorem A \cite[Th. 4.1.3.1]{HTT}.
\end{proof}
\end{cor}

Using Pr. \ref{prp:localrecogofDf}, we further deduce the following recognition principle for the $K$-theory of a Waldhausen $\infty$-category.
\begin{prp}\label{prp:univpropKofC} For any Waldhausen $\infty$-category $\mathscr{C}$, and any functor
\begin{equation*}
\SSS_{\ast}(\mathscr{C})\colon\fromto{N\Delta^{\op}}{\Wald}
\end{equation*}
that classifies the Waldhausen cocartesian fibration $\fromto{\mathscr{S}(\mathscr{C})}{N\Delta^{\op}}$, the $K$-theory space $K(\mathscr{C})$ is the underlying space of the initial object of the $\infty$-category
\begin{equation*}
\Grp(\Kan)\times_{\Fun(N\Delta^{\op},\Kan)}\Fun(N\Delta^{\op},\Kan)_{\iota\SSS_{\ast}(\mathscr{C})/}.
\end{equation*}
\end{prp}


\subsection*{The algebraic $K$-theory of labeled Waldhausen $\infty$-category} We now study the $K$-theory of labeled Waldhausen $\infty$-categories.
\begin{dfn} Suppose $(\mathscr{C},w\mathscr{C})$ a labeled Waldhausen $\infty$-category (Df. \ref{dfn:labeling}). Then we define $K(\mathscr{C},w\mathscr{C})$ as the $K$-theory space $K(\mathscr{B}(\mathscr{C},w\mathscr{C}))$.
\end{dfn}

\begin{ntn} If $\mathscr{C}$ is a Waldhausen $\infty$-category, and if $w\mathscr{C}\subset\mathscr{C}$ is a labeling, then define $w_{N\Delta^{\op}}\mathscr{S}(\mathscr{C})\subset\mathscr{S}(\mathscr{C})$ as the smallest subcategory containing all cocartesian edges and all morphisms of the form $(\id,\psi)\colon\fromto{(\mathbf{m},Y)}{(\mathbf{m},X)}$, where for any integer $0\leq k\leq m$, the induced morphism $\fromto{Y_k}{X_k}$ is labeled.
\end{ntn}

In light of Lm. \ref{lem:iotaBCwCiswC}, we now immediately deduce the following.
\begin{prp} For any labeled Waldhausen $\infty$-category $(\mathscr{C},w\mathscr{C})$, the $K$-theory space $K(\mathscr{C},w\mathscr{C})$ is weakly homotopy equivalent to the loopspace
\begin{equation*}
\Omega(w_{N\Delta^{\op}}\mathscr{S}(\mathscr{C})).
\end{equation*}
\end{prp}
\noindent In other words, for any labeled Waldhausen $\infty$-category $(\mathscr{C},w\mathscr{C})$, the simplicial set $K(\mathscr{C},w\mathscr{C})$ is weakly homotopy equivalent to the loopspace
\begin{equation*}
\Omega\colim w\SSS_{\ast}(\mathscr{C}).
\end{equation*}
Since this again is precisely how Waldhausen's $K$-theory is defined \cite[\S 1.3]{MR86m:18011}, we obtain a further comparison between our $\infty$-categorical $K$-theory for labeled Waldhausen $\infty$-categories and Waldhausen $K$-theory, analogous to Cor. \ref{cor:easycompare}.
\begin{cor}\label{cor:hardcompare} If $(C,\cof C,wC)$ is an ordinary category with cofibrations and weak equivalences in the sense of Waldhausen \cite[\S 1.2]{MR86m:18011}, then the algebraic $K$-theory of the labeled Waldhausen $\infty$-category $(NC,N(\cof C),wC)$ is naturally equivalent to Waldhausen's algebraic $K$-theory of $(C,\cof C,wC)$.
\end{cor}
\noindent Using Cor. \ref{cor:ThmAKtheory}, we obtain the following.
\begin{cor}\label{prp:gammaonKisequiv} Suppose $(\mathscr{C},w\mathscr{C})$ a labeled Waldhausen $\infty$-category. Then the comparison morphism $\gamma_{(\mathscr{C},w\mathscr{C})}$ \textup{(Nt. \ref{ntn:gammaC})} induces an equivalence
\begin{equation*}
\fromto{K(\mathscr{C},w\mathscr{C})}{K(w\mathscr{C}^{-1}\mathscr{C})}
\end{equation*}
of $K$-theory spaces if the following conditions are satisfied.
\begin{enumerate}[(\ref{prp:gammaonKisequiv}.1)]
\item For any object $X$ of $w\mathscr{C}^{-1}\mathscr{C}$, the simplicial set
\begin{equation*}
w\mathscr{C}\times_{\iota(w\mathscr{C}^{-1}\mathscr{C})}\iota(w\mathscr{C}^{-1}\mathscr{C})_{/X}
\end{equation*}
is weakly contractible.
\item For any object $Y$ of $\mathscr{F}_1(w\mathscr{C}^{-1}\mathscr{C})$, the simplicial set
\begin{equation*}
w\mathscr{F}_1(\mathscr{C})\times_{\iota\mathscr{F}_1(w\mathscr{C}^{-1}\mathscr{C})}\iota\mathscr{F}_1(w\mathscr{C}^{-1}\mathscr{C})_{/Y}
\end{equation*}
is weakly contractible.
\end{enumerate}
\end{cor}

Pr. \ref{thm:classWaldareWald}, combined with Cor. \ref{prp:gammaonKisequiv}, yields a further corollary.
\begin{cor}\label{cor:fromWaldtoWald} Suppose $C$ a full subcategory of a model category $M$ that is stable under weak equivalences, then the Waldhausen $K$-theory of $(C,C\cap\cof M,C\cap wM)$ is canonically equivalent to the $K$-theory of a relative nerve $N(C,C\cap wM)$, equipped with the smallest pair structure containing the image of $\cof C$ \textup{(Ex. \ref{exm:catwithcofibsisWaldinftycat})}.
\begin{proof} The only nontrivial point is to check the conditions of Lm. \ref{lem:enoughcofs} for the labeled Waldhausen $\infty$-category $(NC,N(C\cap\cof M),N(C\cap wM))$. Fix a functorial factorization of any map of $C$ into a trivial cofibration followed by a fibration. The functor $F\colon\fromto{\Fun(\Delta^1,NC)}{\Fun(\Delta^1,NC)}$ that carries any map to the trivial cofibration in its factorization now does the job.
\end{proof}
\end{cor}


\subsection*{Cofinality and more fibration theorems} We may also use Cor. \ref{prp:gammaonKisequiv} in combination with Pr. \ref{prp:localizationofperfwaldcat} to specialize the second Generic Fibration Theorem (Th. \ref{thm:fibration}).  We first prove a cofinality result, which states that strongly cofinal inclusions (Df. \ref{dfn:weakcofinal}) of Waldhausen $\infty$-categories do not affect the $K$-theory in high degrees. We are thankful to Peter Scholze for noticing an error that necessitated the inclusion of this result. We follow closely the model of Staffeldt \cite[Th. 2.1]{MR990574}, which works in our setting with only superficial changes.
\begin{thm}[Cofinality]\label{thm:cofinality} The map on $K$-theory induced by the inclusion $i\colon\into{\mathscr{C}'}{\mathscr{C}}$ of a strongly cofinal subcategory fits into a fiber sequence
\begin{equation*}
K(\mathscr{C}')\to K(\mathscr{C})\to A,
\end{equation*}
where $A$ is the abelian group $K_0(\mathscr{C})/K_0(\mathscr{C}')$, regarded as a discrete simplicial set.
\begin{proof} It is convenient to describe the classifying space $BA$ in the following manner. Denote by $BA$ the nerve of the following ordinary category. An object $(m,(x_{i}))$ consists of an integer $m\geq 0$ and a tuple $(x_{i})_{i\in \{1,\dots,m\}}$, and a morphism
\begin{equation*}
\fromto{(m,(x_{i}))}{(n,(y_{j}))}
\end{equation*}
is a morphism $\phi\colon\fromto{\mathbf{n}}{\mathbf{m}}$ of $\Delta$ such that for any $j\in\{1,\dots,n\}$,
\begin{equation*}
y_{j}=\prod_{\phi(j-1)\leq i-1\leq\phi(j)-1}x_{i}.
\end{equation*}
The projection $\fromto{BA}{N\Delta^{\op}}$ clearly induces a left fibration, and the simplicial space $\fromto{N\Delta^{\op}}{\mathbf{Kan}}$ that classifies it visibly satisfies the Segal condition and thus exhibits $(BA)_1\cong A$ as the loop space $\Omega BA$.

We appeal to the Generic Fibration Theorem \ref{prp:fibthmi}. Consider the left fibration
\begin{equation*}
p\colon\fromto{\iota_{N(\Delta^{\op}\times\Delta^{\op})}\mathscr{S}\!\!\mathscr{K}\!(i)}{N(\Delta^{\op}\times\Delta^{\op})}
\end{equation*}
and more particularly its composite $q\coloneq\pr_2\circ p$ with the projection
\begin{equation*}
\pr_2\colon\fromto{N(\Delta^{\op}\times\Delta^{\op})}{N\Delta^{\op}}
\end{equation*}
(whose fiber over $\mathbf{n}\in\Delta$ is $\iota_{N(\Delta^{\op}\times\Delta^{\op})}\mathscr{S}_{\!n}\mathscr{K}\!(i)$). The Generic Fibration Theorem will imply the Cofinality Theorem once we have furnished an equivalence $\iota_{N(\Delta^{\op}\times\Delta^{\op})}\mathscr{S}\!\!\mathscr{K}\!(i)\simeq BA$ over $N\Delta^{\op}$.

Observe that an object $X$ of the $\infty$-category $\iota_{N(\Delta^{\op}\times\Delta^{\op})}\mathscr{S\!\!K}\!(i)$ consists of a diagram in $\mathscr{C}$ of the form
\begin{equation*}
\begin{tikzpicture} 
\matrix(m)[matrix of math nodes, 
row sep=4ex, column sep=4ex, 
text height=1.5ex, text depth=0.25ex] 
{0&0&\cdots&0\\ 
X_{01}&X_{11}&\cdots&X_{m1}\\
\vdots&\vdots&&\vdots\\
X_{0n}&X_{1n}&\cdots&X_{mn},\\}; 
\path[>=stealth,->,font=\scriptsize]
(m-1-1) edge[-,double distance=1.5pt] (m-1-2)
edge[>->] (m-2-1)
(m-1-2) edge[-,double distance=1.5pt] (m-1-3)
edge[>->] (m-2-2)
(m-1-3) edge[-,double distance=1.5pt] (m-1-4)
(m-1-4) edge[>->] (m-2-4)
(m-2-1) edge[>->] (m-2-2) 
edge[>->] (m-3-1)
(m-3-1) edge[>->] (m-4-1)
(m-2-2) edge[>->] (m-2-3)
edge[>->] (m-3-2)
(m-3-2) edge[>->] (m-4-2)
(m-2-4) edge[>->] (m-3-4)
(m-3-4) edge[>->] (m-4-4)
(m-4-1) edge[>->] (m-4-2)
(m-4-2) edge[>->] (m-4-3)
(m-2-3) edge[>->] (m-2-4)
(m-4-3) edge[>->] (m-4-4); 
\end{tikzpicture}
\end{equation*}
such that each $X_{k\ell}/X_{(k-1)\ell}\in\mathscr{C}'$ and the maps
\begin{equation*}
\fromto{X_{(k-1)\ell}\cup^{X_{(k-1)(\ell-1)}}X_{k(\ell-1)}}{X_{k\ell}}
\end{equation*}
are all ingressive. Consequently, we may define a map
\begin{equation*}
\Phi\colon\fromto{\iota_{N(\Delta^{\op}\times\Delta^{\op})}\mathscr{S\!\!K}\!(i)}{BA}
\end{equation*}
that carries an $n$-simplex
\begin{equation*}
X(0)\to\cdots\to X(n)
\end{equation*}
of $\iota_{N(\Delta^{\op}\times\Delta^{\op})}\mathscr{S\!\!K}\!(i)$ to the obvious $n$-simplex whose $i$-th vertex is
\begin{equation*}
\left(q(X(i)),([X(i)_{0\ell}/X(i)_{0(\ell-1)}])_{\ell\in\{1,\dots,q(X(i))\}}\right)
\end{equation*}
of $BA$, where $[Z]$ denotes the image of any object $Z\in\mathscr{C}$ in $K_0(\mathscr{C})/K_0(\mathscr{C}')$. This is easily seen to be a map of simplicial sets over $N\Delta^{\op}$.

Our aim is now to show that $\Phi$ is a fiberwise equivalence. Note that the target satisfies the Segal condition by construction, and the source satisfies it thanks to the Additivity Theorem. Consequently, we are reduced to checking that the induced map
\begin{equation*}
\Phi_1\colon\fromto{\iota_{N\Delta^{\op}}\mathscr{K}(i)}{(BA)_1\cong A}
\end{equation*}
is a weak equivalence. This is the unique map determined by the condition that it carry an object
\[X_0\ \tikz[baseline]\draw[>=stealth,>->](0,0.5ex)--(0.5,0.5ex);\ \cdots\ \tikz[baseline]\draw[>=stealth,>->](0,0.5ex)--(0.5,0.5ex);\ X_n\]
of $\iota\mathscr{K}(i)$ to the class $[X_0]=[X_1]=\cdots=[X_n]\in A$.

One may check that $\Phi_1$ induces a bijection $\equivto{\pi_0\iota_{N\Delta^{\op}}\mathscr{K}\!(i)}{A}$ exactly as in \cite[p. 517]{MR990574}.

Now fix an object $Z\in\mathscr{C}$, and write $\iota_{N\Delta^{\op}}\mathscr{K}\!(i)_{Z}\subset\iota_{N\Delta^{\op}}\mathscr{K}\!(i)$ for the connected component corresponding to the class $[Z]$. This is the full subcategory spanned by those objects
\begin{equation*}
X_0\ \tikz[baseline]\draw[>=stealth,>->](0,0.5ex)--(0.5,0.5ex);\ \cdots\ \tikz[baseline]\draw[>=stealth,>->](0,0.5ex)--(0.5,0.5ex);\ X_n
\end{equation*}
such that $[X_0]=[Z]$ in $A$. We may construct a functor
\begin{equation*}
T\colon\fromto{\iota_{N\Delta^{\op}}\mathscr{F}(\mathscr{C}')}{\iota_{N\Delta^{\op}}\mathscr{K}\!(i)_{Z}}
\end{equation*}
that carries an object
\begin{equation*}
Y_0\ \tikz[baseline]\draw[>=stealth,>->](0,0.5ex)--(0.5,0.5ex);\ \cdots\ \tikz[baseline]\draw[>=stealth,>->](0,0.5ex)--(0.5,0.5ex);\ Y_n
\end{equation*}
to an object
\begin{equation*}
Y_0\vee Z\ \tikz[baseline]\draw[>=stealth,>->](0,0.5ex)--(0.5,0.5ex);\ \cdots\ \tikz[baseline]\draw[>=stealth,>->](0,0.5ex)--(0.5,0.5ex);\ Y_n\vee Z.
\end{equation*}
In the other direction, choose an object $W\in\mathscr{C}$ such that $Z\vee W\in\mathscr{C}'$. Let $S\colon\fromto{\iota_{N\Delta^{\op}}\mathscr{K}\!(i)_{Z}}{\iota_{N\Delta^{\op}}\mathscr{F}(\mathscr{C}')}$ be the obvious functor that carries an object
\begin{equation*}
X_0\ \tikz[baseline]\draw[>=stealth,>->](0,0.5ex)--(0.5,0.5ex);\ \cdots\ \tikz[baseline]\draw[>=stealth,>->](0,0.5ex)--(0.5,0.5ex);\ X_n
\end{equation*}
to an object
\begin{equation*}
X_0\vee W\ \tikz[baseline]\draw[>=stealth,>->](0,0.5ex)--(0.5,0.5ex);\ \cdots\ \tikz[baseline]\draw[>=stealth,>->](0,0.5ex)--(0.5,0.5ex);\ X_n\vee W.
\end{equation*}

Now for any finite simplicial set $K$ and any map $g\colon\fromto{K}{\iota_{N\Delta^{\op}}\mathscr{F}(\mathscr{C}')}$, we construct a map
\begin{equation*}
G\colon\fromto{K\times\Delta^1}{\iota_{N\Delta^{\op}}\mathscr{F}(\mathscr{C}')}
\end{equation*}
such that
\begin{equation*}
G|(K\times\Delta^{\{0\}})\cong g\textrm{\quad and\quad}G|(K\times\Delta^{\{1\}})\cong S\circ T\circ g
\end{equation*}
in the following manner. We let the map $\fromto{K\times\Delta^1}{N\Delta^{\op}}$ induced by $G$ be the projection onto $K$ followed by the map $\fromto{K}{N\Delta^{\op}}$ induced by $g$. The natural transformation from the identity on $\mathscr{C}'$ to the functor $\goesto{X}{X\vee Z\vee W}$ now gives a map $\fromto{(K\times\Delta^1)\times_{N\Delta^{\op}}N\mathrm{M}}{\mathscr{C}}$, which by definition corresponds to the desired map $G$.

In almost exactly the same manner, for any map $f\colon\fromto{K}{\iota_{N\Delta^{\op}}\mathscr{K}\!(i)_{Z}}$, one may construct a map
\begin{equation*}
F\colon\fromto{K\times\Delta^1}{\iota_{N\Delta^{\op}}\mathscr{K}\!(i)_{Z}}
\end{equation*}
such that
\begin{equation*}
F|(K\times\Delta^{\{0\}})\cong f\textrm{\quad and\quad}F|(K\times\Delta^{\{1\}})\cong T\circ S\circ f.
\end{equation*}
We therefore conclude that for any simplicial set $K$, the functors $T$ and $S$ induce a bijection
\begin{equation*}
[K,\iota_{N\Delta^{\op}}\mathscr{F}(\mathscr{C}')]\cong[K,\iota_{N\Delta^{\op}}\mathscr{K}\!(i)_{Z}],
\end{equation*}
whence $S$ and $T$ are homotopy inverses. Now since $\iota_{N\Delta^{\op}}\mathscr{F}(\mathscr{C}')$ is contractible, it follows that $\iota_{N\Delta^{\op}}\mathscr{K}\!(i)_{Z}$ is as well. Thus $\iota_{N\Delta^{\op}}\mathscr{K}\!(i)$ is equivalent to the discrete simplicial set $A$, as desired.
\end{proof}
\end{thm}

In the situation of Pr. \ref{prp:localizationofperfwaldcat}, we find that the natural map
\begin{equation*}
\fromto{K((wC^{\omega})^{-1}C^{\omega})}{K(D^{\omega})}
\end{equation*}
is a homotopy monomorphism; that is, it induces an inclusion on $\pi_0$ and an isomorphism on $\pi_k$ for $k\geq 1$. We therefore obtain the following.
\begin{prp}[Special Fibration Theorem]\label{cor:specialfibration} Suppose $C$ a compactly generated $\infty$-category that is  additive (Df. \ref{item:directsums}). Suppose $L\colon\fromto{C}{D}$ an accessible localization, and suppose the inclusion $\into{D}{C}$ preserves filtered colimits. Assume also that the class of all $L$-equivalences of $C$ is generated (as a strongly saturated class) by the $L$-equivalences between compact objects. Then $L$ induces a pullback square of spaces
\begin{equation*}
\begin{tikzpicture} 
\matrix(m)[matrix of math nodes, 
row sep=4ex, column sep=4ex, 
text height=1.5ex, text depth=0.25ex] 
{K(E^{\omega})&K(C^{\omega})\\ 
\ast&K(D^{\omega}),\\}; 
\path[>=stealth,->,font=\scriptsize] 
(m-1-1) edge (m-1-2) 
edge (m-2-1) 
(m-1-2) edge (m-2-2) 
(m-2-1) edge (m-2-2); 
\end{tikzpicture}
\end{equation*}
where $C^{\omega}$ and $D^{\omega}$ are equipped with the maximal pair structure, and $E^{\omega}\subset C^{\omega}$ is the full subcategory spanned by those objects $X$ such that $LX\simeq 0$.
\end{prp}

A further specialization of this result is now possible. Suppose $C$ a compactly generated stable $\infty$-category. Then $C=\Ind(A)$ for some small $\infty$-category $A$, and so, since $\Ind(A)\subset\mathscr{P}(A)$ is closed under filtered colimits and finite limits, it follows that filtered colimits of $C$ are left exact \cite[Df. 7.3.4.2]{HTT}. Suppose also that $C$ is equipped with a t-structure such that $C_{\leq 0}\subset C$ is stable under filtered colimits. Then the localization $\tau_{\geq 1}\colon\fromto{C}{C}$, being the fiber of the natural transformation $\fromto{\id}{\tau_{\leq 0}}$, preserves filtered colimits as well. Now by \cite[Pr. 1.2.1.16]{HA}, the class $S$ of morphisms $f$ such that $\tau_{\leq 0}(f)$ is an equivalence is generated as a quasi-saturated class by the class $\{\fromto{0}{X}\ |\ X\in C_{\geq 1}\}$. But now writing $X$ as a filtered colimit of compact objects and applying $\tau_{\geq 1}$, we find that $S$ is generated under filtered colimits in $\Fun(\Delta^1,C)$ by the set $\{\fromto{0}{X}\ |\ X\in C^{\omega}\cap C_{\geq 1}\}$. Hence the $\tau_{\leq 0}$-equivalences are generated by $\tau_{\leq 0}$-equivalences between compact objects, and we have the following.
\begin{cor}\label{cor:tstruct} Suppose $C$ a compactly generated stable $\infty$-category. Suppose also that $C$ is equipped with a t-structure such that $C_{\leq 0}\subset C$ is stable under filtered colimits. Then the functor $\tau_{\leq 0}$ induces a pullback square
\begin{equation*}
\begin{tikzpicture} 
\matrix(m)[matrix of math nodes, 
row sep=4ex, column sep=4ex, 
text height=1.5ex, text depth=0.25ex] 
{K(C^{\omega}\cap C_{\geq 1})&K(C^{\omega})\\ 
\ast&K(C^{\omega}\cap C_{\leq 0}),\\}; 
\path[>=stealth,->,font=\scriptsize] 
(m-1-1) edge (m-1-2) 
edge (m-2-1) 
(m-1-2) edge (m-2-2) 
(m-2-1) edge (m-2-2); 
\end{tikzpicture}
\end{equation*}
where the $\infty$-categories that appear are equipped with their maximal pair structure.
\end{cor}
\noindent In particular, we can exploit the equivalence of \cite[Pr. 5.5.7.8]{HTT} to deduce the following.
\begin{cor} Suppose $A$ a small stable $\infty$-category that is equipped with a t-structure. Then the functor $\tau_{\leq 0}$ induces a pullback square
\begin{equation*}
\begin{tikzpicture} 
\matrix(m)[matrix of math nodes, 
row sep=4ex, column sep=4ex, 
text height=1.5ex, text depth=0.25ex] 
{K(A_{\geq 1})&K(A)\\ 
\ast&K(A_{\leq 0}),\\}; 
\path[>=stealth,->,font=\scriptsize] 
(m-1-1) edge (m-1-2) 
edge (m-2-1) 
(m-1-2) edge (m-2-2) 
(m-2-1) edge (m-2-2); 
\end{tikzpicture}
\end{equation*}
where the $\infty$-categories that appear are equipped with their maximal pair structure.
\begin{proof} If $A$ is idempotent-complete, then we can appeal to Cor. \ref{cor:tstruct} and \cite[Pr. 5.5.7.8]{HTT} directly. If not, then we may embed $A$ in its idempotent completion $A'$, and we extend the t-structure using the condition that any summand of an object $X\in A_{\leq 0}$ (respectively, $X\in A_{\geq 1}$) must lie in $A'_{\leq 0}$ (resp., $A'_{\geq 1}$). Now we appeal to the Cofinality Theorem \ref{thm:cofinality} to complete the proof.
\end{proof}
\end{cor}


\section{Example: Algebraic $K$-theory of $E_1$-algebras}\label{sect:example1} To any associative ring in any suitable monoidal $\infty$-category we can attach its $\infty$-category of modules. We may then impose suitable finiteness hypotheses on these modules and extract a $K$-theory spectrum. Here we identify some important examples of these $K$-theory spectra.

\begin{ntn} Suppose $\mathscr{A}$ a presentable, symmetric monoidal $\infty$-category \cite[Df. 2.0.0.7]{HA} with the property that the tensor product $\otimes\colon\fromto{\mathscr{A}\times\mathscr{A}}{\mathscr{A}}$ preserves (small) colimits separately in each variable; assume also that $\mathscr{A}$ is  additive (Df. \ref{item:directsums}). We denote by $\Alg(\mathscr{A})$ the $\infty$-category of $E_1$-algebras in $\mathscr{A}$, and we denote by $\Mod^{\ell}(\mathscr{A})$ the $\infty$-category $\mathrm{LMod}(\mathscr{A})$ defined in \cite[Df. 4.2.1.13]{HA}. We have the canonical presentable fibration
\begin{equation*}
\theta\colon\fromto{\Mod^{\ell}(\mathscr{A})}{\Alg(\mathscr{A})}
\end{equation*}
\cite[Cor. 4.2.3.7]{HA}, whose fiber over any $E_1$-algebra $\Lambda$ is the presentable $\infty$-category $\Mod^{\ell}_{\Lambda}$ of left $\Lambda$-modules. Informally, we describe the objects of $\Mod^{\ell}(\mathscr{A})$ as pairs $(\Lambda,E)$ consisting of an $E_1$-algebra $\Lambda$ in $\mathscr{A}$ and a left $\Lambda$-module $E$.
\end{ntn}

Our aim now is to impose hypotheses on the objects of $(\Lambda,E)$ and pair structures on the resulting full subcategories in order to ensure that the restriction of $\theta$ is a Waldhausen cocartesian fibration.
\begin{dfn}\label{ntn:perfectmodule} For any $E_1$-algebra $\Lambda$ in $\mathscr{A}$, a left $\Lambda$-module $E$ will be said to be \textbf{\emph{perfect}} if it satisfies the following two conditions.
\begin{enumerate}[(\ref{ntn:perfectmodule}.1)]
\item As an object of the $\infty$-category $\Mod^{\ell}_{\Lambda}$ of left $\Lambda$-modules, $E$ is compact. 
\item The functor $\fromto{\Mod^{\ell}_{\Lambda}}{\mathscr{A}}$ corepresented by $E$ is exact.
\end{enumerate}
Denote by $\Perf^{\ell}(\mathscr{A})\subset\Mod^{\ell}(\mathscr{A})$ the full subcategory spanned by those pairs $(\Lambda,E)$ in which $E$ is perfect.
\end{dfn}

These two conditions can be more efficiently expressed by saying that $E$ is perfect just in case the functor $\fromto{\Mod^{\ell}_{\Lambda}}{\mathscr{A}}$ corepresented by $E$ preserves all small colimits. Note that this is \emph{not} the same as \emph{complete compactness}, i.e., requiring that the functor $\fromto{\Mod^{\ell}_{\Lambda}}{\Kan}$ corepresented by $E$ preserves all small colimits. 

\begin{exm} When $\mathscr{A}$ is the nerve of the ordinary category of abelian groups, $\Alg(\mathscr{A})$ is the category of associative rings, and $\Mod^{\ell}(\mathscr{A})$ is the nerve of the ordinary category of pairs $(\Lambda,E)$ consisting of an associative ring $\Lambda$ and a left $\Lambda$-module $E$. An $\Lambda$-module $E$ is perfect just in case it is (1) finitely presented and (2) projective. Thus $\Perf^{\ell}_{\Lambda}$ is the nerve of the ordinary category of finitely generated projective $\Lambda$-modules.
\end{exm}

\begin{exm}\label{exm:summapyperfinconn} When $\mathscr{A}$ is the $\infty$-category of connective spectra, $\Alg(\mathscr{A})$ can be identified with the $\infty$-category of connective $E_1$-rings, and $\Mod^{\ell}(\mathscr{A})$ can be identified with the $\infty$-category of pairs $(\Lambda,E)$ consisting of a connective $E_1$-ring $\Lambda$ and a connective left $\Lambda$-module $E$. Since the functor $\Omega^{\infty}\colon\fromto{\Sp_{\geq 0}}{\Kan}$ is conservative \cite[Cor. 5.1.3.9]{HA} and preserves sifted colimits \cite[Pr. 1.4.3.9]{HA}, it follows using \cite[Lm. 1.3.3.10]{HA} the second condition of Df. \ref{ntn:perfectmodule} amounts to the requirement that $E$ be a projective object. Now \cite[Pr. 8.2.2.6 and Cor. 8.2.2.9]{HA} guarantees that the following are equivalent for a left $\Lambda$-module $E$.
\begin{enumerate}[(\ref{exm:summapyperfinconn}.1)]
\item The left $\Lambda$-module $E$ is perfect.
\item The left $\Lambda$-module $E$ is projective, and $\pi_0E$ is finitely generated as a $\pi_0\Lambda$-module.
\item The $\pi_0\Lambda$-module $\pi_0E$ is finitely generated, and for every $\pi_0A$-module $M$ and every integer $m\geq 1$, the abelian group $\Ext^m(E,M)$ vanishes.
\item There exists a finitely generated free $\Lambda$-module $F$ such that $E$ is a retract of $F$.
\end{enumerate}
\end{exm}

\begin{exm}\label{exm:summapyperfinsmod} The situation for modules over simplicial associative rings is nearly identical. When $\mathscr{A}$ is the $\infty$-category of simplicial abelian groups, $\Alg(\mathscr{A})$ can be identified with the $\infty$-category of simplicial associative rings, and $\Mod^{\ell}(\mathscr{A})$ can be identified with the $\infty$-category of pairs $(\Lambda,E)$ consisting of a simplicial associative ring $\Lambda$ and a left $\Lambda$-module $E$. Since the forgetful functor $\fromto{\mathscr{A}}{\Kan}$ is conservative and preserves sifted colimits, it follows that the second condition of Df. \ref{ntn:perfectmodule} amounts to the requirement that $E$ be a projective object. One may show that the following are equivalent for a left $\Lambda$-module $E$.
\begin{enumerate}[(\ref{exm:summapyperfinsmod}.1)]
\item The left $\Lambda$-module $E$ is perfect.
\item The left $\Lambda$-module $E$ is projective, and $\pi_0E$ is finitely generated as a $\pi_0A$-module.
\item The $\pi_0\Lambda$-module $\pi_0E$ is finitely generated, and for every $\pi_0\Lambda$-module $M$ and every integer $m\geq 1$, the abelian group $\Ext^m(E,M)$ vanishes.
\item There exists a finitely generated free $\Lambda$-module $F$ such that $E$ is a retract of $F$.
\end{enumerate}
\end{exm}

\begin{exm}\label{nul:summaryperfectmad} When $\mathscr{A}$ is the $\infty$-category of \emph{all} spectra, $\Alg(\mathscr{A})$ is the $\infty$-category of $E_1$-rings, and $\Mod^{\ell}(\mathscr{A})$ is the $\infty$-category of pairs $(\Lambda,E)$ consisting of an $E_1$-ring $\Lambda$ and a left $\Lambda$-module $E$. Suppose $\Lambda$ an $E_1$-ring. The second condition of Df. \ref{ntn:perfectmodule} is vacuous since $\mathscr{A}$ is stable. Hence by \cite[Pr. 8.2.5.4]{HA}, the following are equivalent for a left $\Lambda$-module $E$.
\begin{enumerate}[(\ref{nul:summaryperfectmad}.1)]
\item The left $\Lambda$-module $E$ is perfect.
\item The left $\Lambda$-module $E$ is contained in the smallest stable subcategory of the $\infty$-category $\Mod^{\ell}_{\Lambda}$ of left $\Lambda$-modules that contains $\Lambda$ itself and is closed under retracts.
\item The left $\Lambda$-module $E$ is compact as an object of $\Mod^{\ell}_{\Lambda}$.
\item There exists a right $\Lambda$-module $E^{\vee}$ such that the functor $\fromto{\Mod^{\ell}_{\Lambda}}{\Kan}$ informally written as $\Omega^{\infty}(E^{\vee}\otimes_{\Lambda}-)$ is corepresented by $E$.
\end{enumerate}
\end{exm}

Now we wish to endow $\Perf^{\ell}(\mathscr{A})$ with a suitable pair structure. In general, this may not be possible, but we can isolate those situations in which it is possible.
\begin{dfn}\label{cofsinperf} Denote by $S$ the class of morphisms $\fromto{(\Lambda',E')}{(\Lambda,E)}$ of the $\infty$-category $\Perf^{\ell}(\mathscr{A})$ with the following two properties.
\begin{enumerate}[(\ref{cofsinperf}.1)]
\item The morphism $\fromto{\Lambda'}{\Lambda}$ of $\Alg(\mathscr{A})$ is an equivalence.
\item Any pushout diagram
\begin{equation*}
\begin{tikzpicture}[baseline]
\matrix(m)[matrix of math nodes, 
row sep=4ex, column sep=4ex, 
text height=1.5ex, text depth=0.25ex] 
{(\Lambda',E')&(\Lambda,E)\\ 
(\Lambda',0)&(\Lambda,E'')\\}; 
\path[>=stealth,->,font=\scriptsize] 
(m-1-1) edge (m-1-2) 
edge (m-2-1) 
(m-1-2) edge (m-2-2) 
(m-2-1) edge (m-2-2); 
\end{tikzpicture}
\end{equation*}
in $\Mod^{\ell}(\mathscr{A})$ in which $0\in\Mod^{\ell}_{\Lambda'}$ is a zero object is also a pullback diagram, and the $\Lambda$-module $E''$ is perfect.
\end{enumerate}
We shall say that $\mathscr{A}$ is \emph{admissible} if the class $S$ is stable under pushout in $\Perf^{\ell}(\mathscr{A})$ and composition. (Note that pushouts in $\Perf^{\ell}(\mathscr{A})$ are )
\end{dfn}

\begin{exm} When $\mathscr{A}$ is the nerve of the category of abelian groups, $S$ is the class of morphisms $\fromto{(\Lambda',E')}{(\Lambda,E)}$ such that $\fromto{\Lambda'}{\Lambda}$ is an isomorphism, and the induced map of $\Lambda'$-modules $\fromto{E'}{E}$ is an admissible monomorphism. It is a familiar fact that these are closed under pushout and composition, so that the nerve of the category of abelian groups is admissible.
\end{exm}

\begin{exm} When $\mathscr{A}$ is the $\infty$-category of connective spectra or the $\infty$-category of simplicial abelian groups, $S$ is the class of morphisms $\fromto{(\Lambda',E')}{(\Lambda,E)}$ such that $\fromto{\Lambda'}{\Lambda}$ is an equivalence, and the induced homomorphism
\begin{equation*}
\fromto{\Ext^0(E,M)}{\Ext^0(E',M)}
\end{equation*}
is a surjection for every $\pi_0\Lambda'$-module $M$. This is visibly closed under composition. To see that these are closed under pushouts, let us proceed in two steps. First, for any morphism $\fromto{\Lambda}{\Lambda'}$ of $\Alg(\mathscr{A})$, the functor informally described as $\goesto{E}{E\otimes_{\Lambda}\Lambda'}$ clearly carries morphisms of $\Perf^{\ell}_{\Lambda}$ that lie in $S$ to morphisms of $\Perf^{\ell}_{\Lambda'}$ that lie in $S$. Now, for a fixed $E_1$-algebra $\Lambda$ in $\mathscr{A}$, suppose
\begin{equation*}
\begin{tikzpicture} 
\matrix(m)[matrix of math nodes, 
row sep=4ex, column sep=4ex, 
text height=1.5ex, text depth=0.25ex] 
{E'&E\\ 
F'&F\\}; 
\path[>=stealth,->,font=\scriptsize] 
(m-1-1) edge (m-1-2) 
edge (m-2-1) 
(m-1-2) edge (m-2-2) 
(m-2-1) edge (m-2-2); 
\end{tikzpicture}
\end{equation*}
a pushout square in $\Perf^{\ell}_{\Lambda}$ in which $\fromto{E'}{E}$ lies in the class $S$, and suppose $M$ a $\pi_0\Lambda$-module. For any morphism $\fromto{F'}{M}$, one may precompose to obtain a morphism $\fromto{E'}{M}$. Our criterion on the morphism $\fromto{E'}{E}$ now guarantees that there is a commutative square
\begin{equation*}
\begin{tikzpicture} 
\matrix(m)[matrix of math nodes, 
row sep=4ex, column sep=4ex, 
text height=1.5ex, text depth=0.25ex] 
{E'&E\\ 
F'&M\\}; 
\path[>=stealth,->,font=\scriptsize] 
(m-1-1) edge (m-1-2) 
edge (m-2-1) 
(m-1-2) edge (m-2-2) 
(m-2-1) edge (m-2-2); 
\end{tikzpicture}
\end{equation*}
up to homotopy. Now the universal property of the pushout yields a morphism $\fromto{F}{M}$ that extends the morphism $\fromto{F'}{M}$, up to homotopy. Thus both connective spectra and simplicial abelian groups are admissible $\infty$-categories.
\end{exm}

\begin{exm} When $\mathscr{A}$ is the $\infty$-category of all spectra, every morphism is contained in the class $S$. Hence the $\infty$-category of all spectra is an admissible $\infty$-category.
\end{exm}

\begin{ntn} If $\mathscr{A}$ is admissible, denote by $\Perf^{\ell}_{\dag}(\mathscr{A})$ the subcategory of $\Perf^{\ell}(\mathscr{A})$ whose morphisms are those that lie in the class $S$. With this pair structure, the $\infty$-category $\Perf^{\ell}(\mathscr{A})$ is a Waldhausen $\infty$-category.
\end{ntn}

\begin{lem} If $\mathscr{A}$ is admissible, then the functor $\fromto{\Perf^{\ell}(\mathscr{A})}{\Alg(\mathscr{A})}$ is a Waldhausen cocartesian fibration.
\begin{proof} It is clear that the fibers of this cocartesian fibration are Waldhausen $\infty$-categories. We claim that for any morphism $\fromto{\Lambda'}{\Lambda}$ of $E_1$-algebras, the corresponding functor
\begin{equation*}
\fromto{\Mod^{\ell}_{\Lambda'}}{\Mod^{\ell}_{\Lambda}}
\end{equation*}
given informally by the assignment $\goesto{E'}{\Lambda\otimes_{\Lambda'}E'}$ carries perfect modules to perfect modules. Indeed, it is enough to show that the right adjoint functor
\begin{equation*}
\fromto{\Mod^{\ell}_{\Lambda}}{\Mod^{\ell}_{\Lambda'}}
\end{equation*}
preserves small colimits. This is immediate, since colimits are computed in the underlying $\infty$-category $\mathscr{A}$ \cite[Pr. 3.2.3.1]{HA}.

The induced functor $\fromto{\Perf^{\ell}_{\Lambda'}}{\Perf^{\ell}_{\Lambda}}$ carries an ingressive morphism $\cofto{F'}{E'}$ to the morphism of left $\Lambda$-modules $\fromto{F'\otimes_{\Lambda'}\Lambda}{E'\otimes_{\Lambda'}\Lambda}$, which fits into a pushout square
\begin{equation*}
\begin{tikzpicture} 
\matrix(m)[matrix of math nodes, 
row sep=4ex, column sep=4ex, 
text height=1.5ex, text depth=0.25ex] 
{(\Lambda',F')&(\Lambda',E')\\ 
(\Lambda,F'\otimes_{\Lambda'}\Lambda)&(\Lambda,E'\otimes_{\Lambda'}\Lambda)\\}; 
\path[>=stealth,->,font=\scriptsize] 
(m-1-1) edge (m-1-2) 
edge (m-2-1) 
(m-1-2) edge (m-2-2) 
(m-2-1) edge (m-2-2); 
\end{tikzpicture}
\end{equation*}
in $\Perf^{\ell}(\mathscr{A})$; hence $\fromto{F'\otimes_{\Lambda'}\Lambda}{E'\otimes_{\Lambda'}\Lambda}$ is ingressive.
\end{proof}
\end{lem}

\begin{dfn} The \textbf{\emph{algebraic $K$-theory of $E_1$-rings}}, which we will abusively denote
\begin{equation*}
\KK\colon\fromto{\Alg(\mathscr{A})}{\Sp_{\geq0}},
\end{equation*}
is the composite functor $\KK\circ P$, where $P\colon\fromto{\Alg(\mathscr{A})}{\Wald}$ is the functor classified by the Waldhausen cocartesian fibration $\fromto{\Perf^{\ell}(\mathscr{A})}{\Alg(\mathscr{A})}$.
\end{dfn}

\begin{cnstr} The preceding definition ensures that $K$ is well-defined up to a contractible ambiguity. To obtain an explicit model of $K$, we proceed in the following manner. Apply $\mathscr{S}$ to $\fromto{\Perf^{\ell}(\mathscr{A})}{\Alg(\mathscr{A})}$ in order to obtain a Waldhausen cocartesian fibration $\fromto{\mathscr{S}(\Perf^{\ell}(\mathscr{A}))}{N\Delta^{\op}\times\Alg(\mathscr{A})}$. Now consider the subcategory
\begin{equation*}
\iota_{(N\Delta^{\op}\times\Alg(\mathscr{A}))}\mathscr{S}(\Perf^{\ell}(\mathscr{A}))\subset\mathscr{S}(\Perf^{\ell}(\mathscr{A}))
\end{equation*}
consisting of cocartesian edges. The composite
\begin{equation*}
\iota_{(N\Delta^{\op}\times\Alg(\mathscr{A}))}\mathscr{S}\Perf^{\ell}(\mathscr{A})\ \tikz[baseline]\draw[>=stealth,->,font=\scriptsize](0,0.5ex)--(0.75,0.5ex);\ N\Delta^{\op}\times\Alg(\mathscr{A})\ \tikz[baseline]\draw[>=stealth,->,font=\scriptsize](0,0.5ex)--(0.75,0.5ex);\ \Alg(\mathscr{A})
\end{equation*}
is now a left fibration with a contractible space of sections given by
\begin{equation*}
\Alg(\mathscr{A})\cong\{0\}\times\Alg(\mathscr{A})\ \tikz[baseline]\draw[>=stealth,<-,font=\scriptsize,inner sep=0.5pt](0,0.5ex)--node[above]{$\sim$}(0.5,0.5ex);\ \iota\mathscr{S}_0\Perf^{\ell}(\mathscr{A})\ \tikz[baseline]\draw[>=stealth,right hook->,font=\scriptsize](0,0.5ex)--(0.75,0.5ex);\ \iota_{(N\Delta^{\op}\times\Alg(\mathscr{A}))}\mathscr{S}\Perf^{\ell}(\mathscr{A}).
\end{equation*}
It is clear by construction that this left fibration classifies a functor
\begin{equation*}
L\colon\fromto{\Alg(\mathscr{A})}{\Kan}
\end{equation*}
such that $K\simeq\Omega\circ L$. 
\end{cnstr}

Let us now concentrate on the case in which $\mathscr{A}$ is the $\infty$-category of spectra.
\begin{exm} Combining Ex. \ref{exm:stabinftytopoi}, Ex. \ref{exm:Atheoryofinfintopoi}, and the identification of $\Fun(X,\Sp)$ with $\Mod^{\ell}(\Sigma^{\infty}_+X)$, we obtain the well-known equivalence
\[A(X)\simeq K(\Sigma^{\infty}_+X).\]
\end{exm}

\begin{prp}\label{prp:protolocforE1rings} Suppose $\Lambda$ an $E_1$ ring spectrum, and suppose $S\subset\pi_{\ast}\Lambda$ a collection of homogeneous elements satisfying the left Ore condition \cite[Df. 8.2.4.1]{HA}. Then the morphism $\fromto{\Lambda}{\Lambda[S^{-1}]}$ of $\Alg(\Sp)$ induces a fiber sequence of connective spectra
\begin{equation*}
\KK(\mathbf{Nil}_{(\Lambda,S)}^{\ell,\omega})\ \tikz[baseline]\draw[>=stealth,->](0,0.5ex)--(0.5,0.5ex);\ \KK(\Lambda)\ \tikz[baseline]\draw[>=stealth,->](0,0.5ex)--(0.5,0.5ex);\ \KK(\Lambda[S^{-1}]),
\end{equation*}
where $\mathbf{Nil}_{(\Lambda,S)}^{\ell,\omega}\subset\Perf^{\ell}_{\Lambda}$ is the full subcategory spanned by those perfect left $\Lambda$-modules that are $S$-nilpotent.
\begin{proof} Consider the t-structure
\begin{equation*}
(\mathbf{Nil}_{(\Lambda,S)}^{\ell},\mathbf{Loc}_{(\Lambda,S)}^{\ell}),
\end{equation*}
where $\mathbf{Nil}_{(\Lambda,S)}^{\ell}\subset\Mod^{\ell}_{\Lambda}$ is the full subcategory spanned by the $S$-nilpotent left $\Lambda$-modules, and $\mathbf{Loc}_{(\Lambda,S)}^{\ell}\subset\Mod^{\ell}_A$ is the full subcategory spanned by the $S$-local left $\Lambda$-modules. We claim that this t-structure restricts to one on $\Perf_{\Lambda}^{\ell}$. To this end, we note that $\Mod^{\ell}_{\Lambda}$ is compactly generated, and $\mathbf{Loc}_{(\Lambda,S)}^{\ell}\subset\Mod^{\ell}_{\Lambda}$ is in fact stable under all colimits \cite[Rk. 8.2.4.16]{HA}. Now we apply Cor. \ref{cor:tstruct}, and our description of the cofiber term now follows from the discussion preceding \cite[Rk. 8.2.4.26]{HA}.
\end{proof}
\end{prp}
\noindent Such a result is surely well-known among experts; see for example \cite[Pr. 1.4 and Pr. 1.5]{BM}.

\begin{exm}\label{exm:ARconj} For a prime $p$ (suppressed from the notation) and an integer $n\geq 0$, the truncated Brown--Peterson spectra $\mathrm{BP}\langle n\rangle$, with coefficient ring
\begin{equation*}
\pi_{\ast}\mathrm{BP}\langle n\rangle\cong\ZZ_{(p)}[v_1,v_2,\dots,v_n]
\end{equation*}
admit compatible $E_1$ structures \cite[p. 506]{MR1990937}. We may consider the multiplicative system $S\subset\pi_{\ast}\mathrm{BP}\langle n\rangle$ of homogeneous elements generated by $v_n$. Then $\mathrm{BP}\langle n\rangle[v_n^{-1}]$ is an $E_1$-algebra equivalent to the Johnson--Wilson spectrum $E(n)$. The exact sequence above yields a fiber sequence of connective spectra
\begin{equation*}
\KK(\mathbf{Nil}_{(\mathrm{BP}\langle n\rangle,S)}^{\ell,\omega})\ \tikz[baseline]\draw[>=stealth,->](0,0.5ex)--(0.5,0.5ex);\ \KK(\mathrm{BP}\langle n\rangle)\ \tikz[baseline]\draw[>=stealth,->](0,0.5ex)--(0.5,0.5ex);\ \KK(E(n)).
\end{equation*}
The content of a well-known conjecture of Ausoni--Rognes \cite[(0.2)]{MR1947457} identifies the fiber term (possibly after p-adic completion) as $\KK(\mathrm{BP}\langle n-1\rangle)$. In light of results such as \cite[Lm. 8.4.2.13]{HA}, such a result will follow from a suitable form of a \emph{D\'evissage Theorem} \cite[Th. 4]{MR0338129}; we hope to return to such a result in later work (cf. \cite[1.11.1]{MR92f:19001}).

Of course, when $n=1$, such a  D\'evissage Theorem has already been provided thanks to beautiful work of Andrew Blumberg and Mike Mandell \cite{BM}. They prove that the $K$-theory of the $\infty$-category of perfect, $\beta$-nilpotent modules over the $p$-local Adams summand can be identified with the $K$-theory of $\ZZ_{(p)}$. Consequently, they provide a fiber sequence
\begin{equation*}
K(\ZZ_{(p)})\ \tikz[baseline]\draw[>=stealth,->](0,0.5ex)--(0.5,0.5ex);\ K(\ell)\ \tikz[baseline]\draw[>=stealth,->](0,0.5ex)--(0.5,0.5ex);\ K(L).
\end{equation*}
\end{exm}


\section{Example: Algebraic $K$-theory of derived stacks}\label{sect:example2} Here we introduce the algebraic $K$-theory of spectral Deligne--Mumford stacks in the sense of Lurie, and we prove an easy localization theorem (analogous to what Thomason called the ``Proto-localization Theorem'') in this context.

We appeal here to the theory of nonconnective spectral Deligne--Mumford stacks and their module theory as exposed in \cite{DAGV,DAGVIII}. Much of what we will say can probably be done in other contexts of derived algebraic geometry as well, such as \cite{MR2137288,MR2394633}; we have opted to use Lurie's approach only because that is the one with which we are least unfamiliar. We begin by summarizing some general facts about quasicoherent modules over nonconnective spectral Deligne--Mumford stacks. Since Lurie at times concentrates on connective Deligne--Mumford stacks, we will at some points comment on how to extend the relevant definitions and results to the nonconnective case.

\begin{ntn}\label{nul:summary} Recall from \cite[\S 2.3, Pr. 2.5.1]{DAGVIII} that the functor
\begin{equation*}
\fromto{\Sch(\mathscr{G}_{\textrm{\'et}}^{n\mathrm{M}})^{\op}}{\Stk^{\mathrm{nc}}}
\end{equation*}
is a cocartesian fibration, and its fiber over a nonconnective spectral Deligne--Mumford stack $(\mathscr{E},\mathscr{O})$ is the stable, presentable $\infty$-category $\QCoh(\mathscr{E},\mathscr{O})$ of \emph{quasicoherent} $\mathscr{O}$-modules.

For any nonconnective Deligne--Mumford stack $(\mathscr{E},\mathscr{O})$, the following are equivalent for an $\mathscr{O}$-module $\mathscr{M}$.
\begin{enumerate}[(\ref{nul:summary}.1)]
\item The $\mathscr{O}$-module $\mathscr{M}$ is quasicoherent.
\item For any morphism $\fromto{U}{V}$ of $\mathscr{E}$ such that $(\mathscr{X}_{/U},\mathscr{O}_{|U})$ and $(\mathscr{X}_{/V},\mathscr{O}_{|V})$ are affine, the natural morphism $\fromto{\mathscr{M}(V)\otimes_{\mathscr{O}(V)}\mathscr{O}(U)}{\mathscr{M}(U)}$ is an equivalence.
\item\label{item:simplecharqcoh} The following conditions obtain.
\begin{enumerate}[(\ref{nul:summary}.\ref{item:simplecharqcoh}.1)]
\item For every integer $n$, the homotopy sheaf $\pi_n\mathscr{M}$ is a quasicoherent module on the underlying ordinary Deligne--Mumford stack of $(\mathscr{E},\mathscr{O})$
\item The object $\Omega^{\infty}\mathscr{M}$ is hypercomplete in the $\infty$-topos $\mathscr{E}$.
\end{enumerate}
\end{enumerate}
\end{ntn}

Using ideas from \cite[\S 2.7]{DAGVIII}, we shall now make sense of the notion of quasicoherent module over any functor $\fromto{\mathbf{CAlg}}{\Kan(\kappa_1)}$. As suggested in \cite[Rk. 2.7.9]{DAGVIII}, write
\begin{equation*}
\QCoh\colon\fromto{\Fun(\mathbf{CAlg},\Kan(\kappa_1))^{\op}}{\Cat_{\infty}(\kappa_1)}
\end{equation*}
for the right Kan extension of the functor $\fromto{\mathbf{CAlg}}{\Cat_{\infty}(\kappa_1)}$ that classifies the cocartesian fibration $\fromto{\Mod}{\mathbf{CAlg}}$. Then for any functor $X\colon\fromto{\mathbf{CAlg}}{\Kan(\kappa_1)}$, we obtain the $\infty$-category of \emph{quasicoherent modules} $\QCoh(X)$ on the functor $X$. Many of the results of \S 2.7 of loc. cit. hold in this context with precisely the same proofs, including the following brace of results.
\begin{prp}[cf. \protect{\cite[Rk. 2.7.17]{DAGVIII}}] For any functor $X\colon\fromto{\mathbf{CAlg}}{\Kan(\kappa_1)}$, the $\infty$-category $\QCoh(X)$ is stable.
\end{prp}

\begin{prp}[cf. \protect{\cite[Rk. 2.7.18]{DAGVIII}}]\label{prp:nonconnversionof2718} Suppose $(\mathscr{E},\mathscr{O})$ a nonconnective Deligne--Mumford stack representing a functor $X\colon\fromto{\mathbf{CAlg}}{\Kan(\kappa_1)}$. Then there is a canonical equivalence of $\infty$-categories
\begin{equation*}
\QCoh(\mathscr{E},\mathscr{O})\simeq\QCoh(X).
\end{equation*}
\end{prp}

\begin{dfn} Suppose $X\colon\fromto{\mathbf{CAlg}}{\Kan(\kappa_1)}$ a functor. We say that a quasicoherent module $\mathscr{M}$ on $X$ is \textbf{\emph{perfect}} if for any $E_{\infty}$ ring $A$ and any point $x\in X(A)$, the $A$-module $\mathscr{M}(x)$ is perfect (Df. \ref{ntn:perfectmodule}). Write $\Perf(X)\subset\QCoh(X)$ for the full subcategory spanned by the perfect modules.
\end{dfn}
\noindent In particular, we can now use Pr. \ref{prp:nonconnversionof2718} to specialize the notion of perfect module to the setting of nonconnective Deligne--Mumford stacks.
\begin{ntn} Denote by $\Perf\subset\Sch(\mathscr{G}_{\textrm{\'et}}^{\mathrm{M}})^{\op}$ the full subcategory of those objects $(\mathscr{E},\mathscr{O},\mathscr{M})$ such that $\mathscr{M}$ is perfect.
\end{ntn}

\begin{nul} For any functor $X\colon\fromto{\mathbf{CAlg}}{\Kan(\kappa_1)}$, the $\infty$-category $\QCoh(X)$ admits a symmetric monoidal structure \cite[Nt. 2.7.27]{DAGVIII}. Moreover, this is functorial, yielding a functor
\begin{equation*}
\QCoh^{\otimes}\colon\fromto{\Fun(\mathbf{CAlg},\Kan(\kappa_1))^{\op}}{\mathbf{CAlg}(\Cat_{\infty}(\kappa_1))}.
\end{equation*}
\end{nul}

\begin{prp}[cf. \protect{\cite[Pr. 2.7.28]{DAGVIII}}] For any functor $X\colon\fromto{\mathbf{CAlg}}{\Kan(\kappa_1)}$ a quasicoherent module $\mathscr{M}$ on $X$ is perfect if and only if it is a dualizable object of $\QCoh(X)$.
\end{prp}
\noindent Since the pullback functors are symmetric monoidal, they preserve dualizable objects. This proves the following.
\begin{cor} The functor $\fromto{\Perf}{\Stk^{\mathrm{nc}}}$ is a cocartesian fibration.
\end{cor}

We endow $\Perf$ with a pair structure by $\Perf_{\dag}\coloneq\Perf\times_{\Stk^{\mathrm{nc}}}\iota\Stk^{\mathrm{nc}}$, so that the fibers are equipped with the maximal pair structure.
\begin{prp} The functor $\fromto{\Perf}{\Stk^{\mathrm{nc}}}$ is a Waldhausen cocartesian fibration.
\end{prp}
In fact, the fiber over a nonconnective Deligne--Mumford stack $(\mathscr{E},\mathscr{O})$ is a stable $\infty$-category $\Perf(\mathscr{E},\mathscr{O})$.

\begin{dfn}\label{dfn:algKthyDMstacks} We now define the \textbf{\emph{algebraic $K$-theory of nonconnective Deligne--Mumford stacks}} as a functor that we abusively denote
\begin{equation*}
\KK\colon\fromto{\Stk^{\mathrm{nc}}}{\Sp_{\geq 0}}
\end{equation*}
given by the composite $\KK\circ P$, where $P$ is the functor $\fromto{\Stk^{\mathrm{nc},\op}}{\Wald}$ classified by the Waldhausen cocartesian fibration $\fromto{\Perf}{\Stk^{\mathrm{nc}}}$.
\end{dfn}

\begin{lem} For any open immersion of quasicompact nonconnective spectral Deligne--Mumford stacks $j\colon\fromto{\mathscr{U}}{\mathscr{X}}$, the induced functor
\begin{equation*}
j_{\star}\colon\fromto{\QCoh(\mathscr{U})}{\QCoh(\mathscr{X})}
\end{equation*}
is fully faithful.
\begin{proof} When $\mathscr{X}$ is of the form $\Spec^{\textrm{\'et}}A$, this is proved in \cite[Cor. 2.4.6]{DAGVIII}. For any map $x\colon\fromto{\Spec^{\textrm{\'et}}A}{\mathscr{X}}$, we have the open immersion
\begin{equation*}
\fromto{\mathscr{U}\times_{\mathscr{X}}\Spec^{\textrm{\'et}}A}{\Spec^{\textrm{\'et}}A},
\end{equation*}
which induces a fully faithful functor
\begin{equation*}
\fromto{\QCoh(\mathscr{U}\times_{\mathscr{X}}\Spec^{\textrm{\'et}}A)}{\QCoh(\Spec^{\textrm{\'et}}A)}.
\end{equation*}
Now letting $A$ vary and applying \cite[Pr. 2.4.5(3)]{DAGVIII}, we obtain a functor
\begin{equation*}
\fromto{\mathbf{CAlg}_{\mathscr{X}/}}{\mathscr{O}(\Cat_{\infty}(\kappa_1))}
\end{equation*}
whose values are all fully faithful functors. Thanks to Pr. \ref{prp:nonconnversionof2718}, the limit of this functor is then equivalent to a functor
\begin{equation*}
\alpha\colon\fromto{\lim_{A\in\mathbf{CAlg}_{\mathscr{X}/}}\QCoh(\mathscr{U}\times_{\mathscr{X}}\Spec^{\textrm{\'et}}A)}{\QCoh(\mathscr{X})},
\end{equation*}
which is thus fully faithful. We aim to identify this functor with $j_{\star}$.

Since each of the $\infty$-categories $\QCoh(\mathscr{U}\times_{\mathscr{X}}\Spec^{\textrm{\'et}}A)$ can itself be described as the limit of the $\infty$-categories $\Mod_B$ for $B\in\mathbf{CAlg}_{\mathscr{U}\times_{\mathscr{X}}\Spec^{\textrm{\'et}}A/}$, it follows that the source of $\alpha$ can be expressed as the limit of the $\infty$-categories $\Mod_B$ over the $\infty$-category $C$ of squares of nonconnective Deligne--Mumford stacks of the form
\begin{equation*}
\begin{tikzpicture} 
\matrix(m)[matrix of math nodes, 
row sep=4ex, column sep=4ex, 
text height=1.5ex, text depth=0.25ex] 
{\Spec^{\textrm{\'et}}B&\Spec^{\textrm{\'et}}A\\ 
\mathscr{U}&\mathscr{X}.\\}; 
\path[>=stealth,->,font=\scriptsize] 
(m-1-1) edge (m-1-2) 
edge (m-2-1) 
(m-1-2) edge (m-2-2) 
(m-2-1) edge node[below]{$j$} (m-2-2); 
\end{tikzpicture}
\end{equation*}
Now there is a forgetful functor $g\colon\fromto{C}{\mathbf{CAlg}_{\mathscr{U}/}}$ that carries an object as above to the morphism $\fromto{\Spec^{\textrm{\'et}}B}{\mathscr{U}}$. This is the functor that induces the canonical functor
\begin{equation*}
\fromto{\lim_{A\in\mathbf{CAlg}_{\mathscr{X}/}}\QCoh(\mathscr{U}\times_{\mathscr{X}}\Spec^{\textrm{\'et}}A)}{\QCoh(\mathscr{U})};
\end{equation*}
hence it suffices to show that $g$ is right cofinal. This now follows from the fact that the functor $g$ admits a right adjoint $\fromto{\mathbf{CAlg}_{\mathscr{U}/}}{C}$, which carries a morphism $x\colon\fromto{\Spec^{\textrm{\'et}}C}{\mathscr{U}}$ to the object
\begin{equation*}
\begin{tikzpicture} 
\matrix(m)[matrix of math nodes, 
row sep=4ex, column sep=4ex, 
text height=1.5ex, text depth=0.25ex] 
{\Spec^{\textrm{\'et}}C&\Spec^{\textrm{\'et}}C\\ 
\mathscr{U}&\mathscr{X}.\\}; 
\path[>=stealth,->,font=\scriptsize] 
(m-1-1) edge[-,double distance=1.5pt] (m-1-2) 
edge node[left]{$x$} (m-2-1) 
(m-1-2) edge node[right]{$j\circ x$} (m-2-2) 
(m-2-1) edge node[below]{$j$} (m-2-2); 
\end{tikzpicture}
\end{equation*}
The proof is complete.
\end{proof}
\end{lem}

\begin{ntn} For any open immersion $j\colon\fromto{\mathscr{U}}{\mathscr{X}}$ of quasicompact nonconnective spectral Deligne--Mumford stacks, let us write $\Perf(\mathscr{X},\mathscr{X}\setminus\mathscr{U})$ for the full subcategory of $\Perf(\mathscr{X})$ spanned by those perfect modules $\mathscr{M}$ on $\mathscr{X}$ such that $j^{\star}\mathscr{M}\simeq0$. Write
\begin{equation*}
\KK(\mathscr{X},\mathscr{X}\setminus\mathscr{U})\coloneq\KK(\Perf(\mathscr{X},\mathscr{X}\setminus\mathscr{U})).
\end{equation*}
\end{ntn}

\begin{prp}[``Proto-localization,'' cf. \protect{\cite[Th. 5.1]{MR92f:19001}}]\label{prp:protolocforDMstacks} For any quasicompact open immersion $j\colon\fromto{\mathscr{U}}{\mathscr{X}}$ of quasicompact, quasiseparated spectral algebraic spaces \cite[Df. 1.3.1]{DAGVIII} and \cite[Df. 1.3.1]{DAGXII}, the functor $j^{\star}\colon\fromto{\Perf(\mathscr{X})}{\Perf(\mathscr{U})}$ induces a fiber sequence of connective spectra
\begin{equation*}
\KK(\mathscr{X},\mathscr{X}\setminus\mathscr{U})\ \tikz[baseline]\draw[>=stealth,->,font=\scriptsize](0,0.5ex)--(0.75,0.5ex);\ \KK(\mathscr{X})\ \tikz[baseline]\draw[>=stealth,->,font=\scriptsize](0,0.5ex)--(0.75,0.5ex);\ \KK(\mathscr{U}).
\end{equation*}
\begin{proof} We wish to employ the Special Fibration Theorem \ref{cor:specialfibration}. We note that by \cite[Cor. 1.5.12]{DAGXII}, the $\infty$-category $\QCoh(\mathscr{X})$ is compactly generated, and one has $\Perf(\mathscr{X})=\QCoh(\mathscr{X})^{\omega}$; the analogous claim holds for $\mathscr{U}$. It thus remains to show that $j^{\star}$-equivalences of $\QCoh(\mathscr{X})$ --- i.e., the class of morphisms of $\QCoh(\mathscr{X})$ whose restriction to $\mathscr{U}$ is an equivalence --- is generated (as a strongly saturated class) by $j^{\star}$-equivalences between compact objects. Since $\QCoh(\mathscr{X})$ is stable, we find that it suffices to show that the full subcategory $\QCoh(\mathscr{X},\mathscr{X}-\mathscr{U})$ of $\QCoh(\mathscr{X})$ spanned by the $j^{\star}$-acyclics --- i.e., those quasicoherent modules $M$ such that $j^{\star}M\simeq 0$ --- is generated by compact objects of $\QCoh(\mathscr{X})$. This will follow from \cite[Th. 1.5.10]{DAGXII} once we know that the quasicoherent stack $\Phi_{\mathscr{X}}(\QCoh(\mathscr{X},\mathscr{X}-\mathscr{U}))$ of \cite[Constr. 8.5]{DAGXI} is locally compactly generated.

So suppose $R$ a connective $E_{\infty}$ ring spectrum, and suppose $\eta\in\mathscr{X}(R)$. We wish to show that the $\infty$-category
\begin{equation*}
\Phi_{\mathscr{X}}(\QCoh(\mathscr{X},\mathscr{X}-\mathscr{U}))(\eta)\simeq\Mod_R\otimes_{\QCoh(\mathscr{X})}\QCoh(\mathscr{X},\mathscr{X}\setminus\mathscr{U})
\end{equation*}
is compactly generated. It easy to see that this $\infty$-category can be identified with the full subcategory of $\Mod_R$ spanned by those modules $M$ that are carried to zero by the functor
\begin{equation*}
\fromto{\Mod_R\simeq\Mod_R\otimes_{\QCoh(\mathscr{X})}\QCoh(\mathscr{X})}{\Mod_R\otimes_{\QCoh(\mathscr{X})}\QCoh(\mathscr{U})}.
\end{equation*}
By a theorem of Ben Zvi, Francis, and Nadler \cite[Cor. 8.22]{DAGXI}, this functor may be identified with the restriction functor along the open embedding
\begin{equation*}
j'\colon\into{\mathscr{U}'\coloneq\Spec^{\textrm{\'et}}R\times_{\mathscr{X}}\mathscr{U}}{\Spec^{\textrm{\'et}}R}.
\end{equation*}
The open immersion $j'$ is determined by a quasicompact open $U\subset\Spec^{\mathrm{Z}}A$, which consists of those prime ideals of $\pi_0A$ that do not contain a finitely generated ideal $I$. The proof is now completed by \cite[Pr. 4.1.15 and Pr. 5.1.3]{DAGVIII}.
\end{proof}
\end{prp}

When $j$ is the open complement of a closed immersion $i\colon\fromto{\mathscr{Z}}{\mathscr{X}}$, one may ask whether $\KK(\mathscr{X}\setminus\mathscr{U})$ can be identified with $\KK(\mathscr{Z})$. In general, the answer is no, but in special situations, such an identification is possible. Classically, this is the result of a \emph{D\'evissage Theorem} \cite[Th. 4]{MR0338129}; we hope to return to a higher categorical analogue of such a result in later work (cf. \cite[1.11.1]{MR92f:19001})


\bibliographystyle{amsplain}
\bibliography{kthy}

\end{document}